%% file: eigenWave.tex
\documentclass[preprint]{elsarticle}
\usepackage[bookmarks=true,colorlinks=true,linkcolor=blue]{hyperref}

\usepackage{amsmath,amssymb,mathabx,mathtools}
\usepackage{graphicx,esint,ulem}

\biboptions{sort&compress}
\usepackage[usenames,dvipsnames,svgnames,table]{xcolor}

\usepackage{color}

\newcommand{\fb}{{\mathbf f}}

\usepackage{calc}

\usepackage{fancyvrb}

\usepackage{empheq}
\usepackage[most]{tcolorbox}

\newtcbox{\mymath}[1][]{%
    nobeforeafter, math upper, tcbox raise base,
    enhanced, colframe=blue!60!black,
    colback=blue!20, boxrule=1pt,
    #1}

\usepackage{algorithm,algpseudocode}
\usepackage{algorithmicx}
\algrenewcommand\alglinenumber[1]{\footnotesize #1:} 
\newcommand{\algFontSize}{\footnotesize}

\usepackage{tcolorbox}
\usepackage[english]{babel}
\usepackage{amsmath}
\usepackage{amsfonts}
\usepackage{fullpage}
\usepackage{graphicx}
\usepackage{setspace}
\usepackage{float}

\usepackage{tabu}

\usepackage{multirow}

\usepackage{mdframed}

\allowdisplaybreaks

\usepackage{longtable}

\usepackage{tikz}

\newtheorem{theorem}{Theorem}

\newenvironment{proof}[1][Proof]{\begin{trivlist}
\item[\hskip \labelsep {\bfseries #1.}]}{\end{trivlist}}

\input tex/defs.tex
\input tex/trimFig.tex
\usepackage{xargs}
\input tex/plotFigureMacros.tex

\usepackage{todonotes}

\definecolor{jwbGreen}{rgb}{0, .6, 0}
\definecolor{purple}{rgb}{.7, 0., .8}
\definecolor{pinegreen}{rgb}{0.0, 0.47, 0.44}
\definecolor{darkyellow}{rgb}{0.7, 0.5, 0.0}

\newcommand{\blue}{\color{blue}}

\begin{document}


\begin{frontmatter}
 \title{EigenWave: An Optimal O(N) Method for Computing Eigenvalues\\ and Eigenvectors by Time-Filtering the Wave Equation}

\author[vtu]{Daniel Appel\"o\fnref{DanielThanks}}
\ead{appelo@vt.edu}

\author[rpi]{Jeffrey W.~Banks\fnref{NSFgrants}}
\ead{banksj3@rpi.edu}

\author[rpi]{William D.~Henshaw\fnref{NSFgrants}}
\ead{henshw@rpi.edu}

\author[rpi]{Ngan~Le\fnref{NSFgrants}}
\ead{leb2@rpi.edu}


\author[rpi]{Donald~W.~Schwendeman\corref{cor}\fnref{NSFgrants}}
\ead{schwed@rpi.edu}

\address[vtu]{Department of Mathematics, Virginia Tech, Blacksburg, VA 24061 U.S.A.}
\address[rpi]{Department of Mathematical Sciences, Rensselaer Polytechnic Institute, Troy, NY 12180, USA}

\cortext[cor]{Corresponding author}

\fntext[DanielThanks]{Research supported by National Science Foundation under grant DMS-2345225, and Virginia Tech.}
\fntext[NSFgrants]{Research supported by the National Science Foundation under grants DMS-1519934 and DMS-1818926.}

\begin{abstract}
   An algorithm named EigenWave is described to compute eigenvalues and eigenvectors of elliptic boundary value problems.
The algorithm, based on the recently developed WaveHoltz scheme, solves a related time-dependent wave equation
as part of an iteration. 
At each iteration, the solution to the wave equation is filtered in time.
As the iteration progresses, the filtered solution generally contains relatively larger and larger proportions of eigenmodes whose eigenvalues are near a chosen target frequency (target eigenvalue).
The ability to choose an arbitrary target frequency enables the computation of eigenvalues
anywhere in the spectrum, without the need to invert an indefinite matrix, as is common with other approaches. Furthermore, the iteration can be embedded within a matrix-free Arnoldi algorithm, which enables the efficient computation
of multiple eigenpairs near the target frequency. 
For efficiency, the time-dependent wave equation can be solved with implicit time-stepping and only about $10$ time-steps
per-period are needed, independent of the mesh spacing. When the (definite) implicit time-stepping equations are solved with a multigrid algorithm, the cost of the resulting EigenWave scheme
scales linearly with the number of grid points $N$ as the mesh is refined, giving an optimal $O(N)$ algorithm.
The approach is demonstrated by finding eigenpairs of the Laplacian in complex geometry
using overset grids. Results in two and three space dimensions are presented using 
second-order and fourth-order accurate approximations.
\end{abstract}


\end{frontmatter}

\tableofcontents

\clearpage
\input tex/intro


\input tex/algorithm

\input tex/analysis

\input tex/discreteApproximations

\input tex/usingExistingEigenvalueRoutines

\input tex/numericalResults

\input tex/optimalAlgorithm

\input tex/conclusions

\appendix


\input tex/discreteAnalysis

\input tex/arnoldi

\input tex/properties

\input tex/eigenpairTables
\clearpage
\bibliographystyle{elsart-num}
\bibliography{journal-ISI,henshaw,henshawPapers,eigen,WaveHoltz,helmholtz,eigenV2}
 
\end{document}

%% file: tex/defs.tex
\newcommand{\citeCount}[1]{}

\newcommand{\bogus}[1]{{}}

\newcommand{\mni}{\medskip\noindent}

\DeclareMathOperator{\sign}{sign}

\newcommand{\CR}{{\rm CR}}

\newlength{\ycbTop}
\newlength{\ycbMid}%

\newcommand{\p}{\partial}

\newcommand{\Np}{{N_p}}

\newcommand{\tableFontSize}{\tableFont}

\newcommand{\f}[2]{\frac{#1}{#2}}

\def\ba#1\ea{\begin{align}#1\end{align}}

\def\bas#1\eas{\begin{align*}#1\end{align*}}

\def\bat#1\eat{\begin{alignat}{3}#1\end{alignat}}

\def\bats#1\eats{\begin{alignat*}{3}#1\end{alignat*}}

\newcommand{\bse}{\begin{subequations}}
\newcommand{\ese}{\end{subequations}}

\newcommand{\Nt}{N_t}

\newcommand{\Dzt}{D_{0t}}
\newcommand{\Dpt}{D_{+t}}
\newcommand{\Dmt}{D_{-t}}
\newcommand{\dt}{\Delta t}

\newcommand{\eqdef}{\overset{{\rm def}}{=}}

\newcommand{\ev}{\mathbf{ e}}
\newcommand{\fv}{\mathbf{ f}}

\newcommand{\iv}{\mathbf{ i}}

\newcommand{\rv}{\mathbf{ r}}

\newcommand{\vv}{\mathbf{ v}}
\newcommand{\wv}{\mathbf{ w}}
\newcommand{\xv}{\mathbf{ x}}
\newcommand{\yv}{\mathbf{ y}}
\newcommand{\zv}{\mathbf{ z}}
\newcommand{\Av}{\mathbf{ A}}

\newcommand{\Gv}{\mathbf{ G}}
\newcommand{\Hv}{\mathbf{ H}}
\newcommand{\Iv}{\mathbf{ I}}

\newcommand{\Qv}{\mathbf{ Q}}
\newcommand{\Rv}{\mathbf{ R}}

\newcommand{\Vv}{\mathbf{ V}}
\newcommand{\Wv}{\mathbf{ W}}

\newcommand{\half}{{1\over2}}

\newcommand{\Real}{{\mathbb R}}

\newcommand{\Bc}{{\mathcal B}}

\newcommand{\Ec}{{\mathcal E}}

\newcommand{\Gc}{{\mathcal G}}

\newcommand{\Kc}{{\mathcal K}}

\newcommand{\Lc}{{\mathcal L}}

\newcommand{\nd}{n_d}

\newcommand{\ds}{{\Delta s}}

\newcommand{\ssf}{\scriptscriptstyle}

\newcommand{\sinc}{{\rm sinc}}
\newcommand{\Tbar}{T_f}
\newcommand{\Tf}{T_f}

\newcommand{\omegaTilde}{\tilde{\omega}}
\newcommand{\vHat}{\hat{v}}

\newcommand{\wHat}{\hat{w}}

\newcommand{\Gcd}{\Gc_d}
\newcommand{\NITS}{N_{\ssf ITS}}
\newcommand{\Nits}{\NITS}

\newcommand{\betad}{\beta_d}

\newcommand{\betadTilde}{{\tilde\beta_d}}

\newcommand{\EigenWaveb}{EigenWave\ }

\newcommand{\Aw}{S}

\newcommand{\Active}{\Omega_h^a}

\newcommand{\Ntotal}{N_a}
\newcommand{\Nreq}{N_r}
\newcommand{\Nextra}{N_e}

\newcommand{\mTotal}{{N_a}}
\newcommand{\kReq}{{N_r}}
\newcommand{\pExtra}{{N_e}}

\newcommand{\Gcde}{\Gc_{\rm de}}

\newcommand{\lambdaExplicit}{\lambda^{(e)}_{h,j}}
\newcommand{\lambdaImplicit}{\lambda^{(i)}_{h,j}}

\newcommand{\NT}{\Nt}

\newcommand{\Phiv}{\boldsymbol{\Phi}}
\newcommand{\mvCol}{\blue}

\newcommand{\mr}{{m_{r}}}
\newcommand{\mTheta}{{m_{\theta}}}

\newcommand{\br}{\beta_r}
\newcommand{\dlam}{\Delta\lambda}
\newcommand{\dw}{\delta}

\newcommand{\What}{\hat{W}}
\newcommand{\Vhat}{\hat{V}}
\newcommand{\NGd}{N_h}
\newcommand{\Lam}{\tilde\Lambda}
\newcommand{\dtTilde}{\widetilde{\dt}}

%% file: tex/trimFig.tex
\newlength{\tfwidth}
\newlength{\tfheight}
\newlength{\tfxa}
\newlength{\tfxb}
\newlength{\tfya}
\newlength{\tfyb}
\newcommand{\trimFigWithBox}[6]{%
\setlength\fboxsep{0pt}%
\setlength\fboxrule{1.0pt}
\fbox{\includegraphics[width=#2, clip, trim=#3 #4 #5 #6]{#1}}%
}
\newcommand{\trimFigNoBox}[6]{%
\setlength\fboxsep{1pt}
\setlength\fboxrule{0.0pt}
\fbox{\includegraphics[width=#2, clip, trim=#3 #4 #5 #6]{#1}}%
}
\newcommand{\trimFigHeightWithBox}[6]{%
\setlength\fboxsep{0pt}%
\setlength\fboxrule{1.0pt}
\fbox{\includegraphics[height=#2, clip, trim=#3 #4 #5 #6]{#1}}%
}
\newcommand{\trimFigHeightNoBox}[6]{%
\setlength\fboxsep{1pt}
\setlength\fboxrule{0.0pt}
\fbox{\includegraphics[height=#2, clip, trim=#3 #4 #5 #6]{#1}}%
}

\newsavebox\figBox

\newcommand{\trimw}[6]{%
\sbox\figBox{\includegraphics{#1}}
\setlength{\tfwidth}{\the\wd\figBox}
\setlength{\tfheight}{\the\ht\figBox}
\setlength{\tfxa}{\tfwidth*\real{#3}}%
\setlength{\tfxb}{\tfwidth*\real{#4}}%
\setlength{\tfya}{\tfheight*\real{#5}}%
\setlength{\tfyb}{\tfheight*\real{#6}}%
\trimFigNoBox{#1}{#2}{\tfxa}{\tfya}{\tfxb}{\tfyb}%
}

\newcommand{\trimwb}[6]{%

\sbox\figBox{\includegraphics{#1}}
\setlength{\tfwidth}{\the\wd\figBox}
\setlength{\tfheight}{\the\ht\figBox}
\setlength{\tfxa}{\tfwidth*\real{#3}}%
\setlength{\tfxb}{\tfwidth*\real{#4}}%
\setlength{\tfya}{\tfheight*\real{#5}}%
\setlength{\tfyb}{\tfheight*\real{#6}}%
\trimFigWithBox{#1}{#2}{\tfxa}{\tfya}{\tfxb}{\tfyb}%
}

\newcommand{\trimh}[6]{%
\sbox\figBox{\includegraphics{#1}}
\setlength{\tfwidth}{\the\wd\figBox}
\setlength{\tfheight}{\the\ht\figBox}
\setlength{\tfxa}{\tfwidth*\real{#3}}%
\setlength{\tfxb}{\tfwidth*\real{#4}}%
\setlength{\tfya}{\tfheight*\real{#5}}%
\setlength{\tfyb}{\tfheight*\real{#6}}%
\trimFigHeightNoBox{#1}{#2}{\tfxa}{\tfya}{\tfxb}{\tfyb}%
}

\newcommand{\trimhb}[6]{%

\sbox\figBox{\includegraphics{#1}}
\setlength{\tfwidth}{\the\wd\figBox}
\setlength{\tfheight}{\the\ht\figBox}
\setlength{\tfxa}{\tfwidth*\real{#3}}%
\setlength{\tfxb}{\tfwidth*\real{#4}}%
\setlength{\tfya}{\tfheight*\real{#5}}%
\setlength{\tfyb}{\tfheight*\real{#6}}%
\trimFigHeightWithBox{#1}{#2}{\tfxa}{\tfya}{\tfxb}{\tfyb}%
}

%% file: tex/plotFigureMacros.tex
\newcommandx{\figByHeight}[9][5=0, 6=0, 7=0, 8=0,9=]{
\draw (#1,#2) node[anchor=south west,xshift=-16pt,yshift=-4pt] {\trimh{#3}{#4}{#5}{#6}{#7}{#8}};}
\newcommandx{\figByHeightb}[9][5=0, 6=0, 7=0, 8=0,9=]{
\draw (#1,#2) node[anchor=south west,xshift=-16pt,yshift=-4pt] {\trimhb{#3}{#4}{#5}{#6}{#7}{#8}};}

\newcommandx{\figByHeightWithLabel}[9][5=0, 6=0, 7=0, 8=0,9=]{
\draw (#1,#2) node[anchor=south west,xshift=-16pt,yshift=-4pt] {\trimh{#3}{#4}{#5}{#6}{#7}{#8}} node[draw=white,fill=white,inner sep=1pt,anchor=south west] {#9};}
\newcommandx{\figByHeightWithLabelb}[9][5=0, 6=0, 7=0, 8=0,9=]{
\draw (#1,#2) node[anchor=south west,xshift=-16pt,yshift=-4pt] {\trimhb{#3}{#4}{#5}{#6}{#7}{#8}} node[draw=white,fill=white,inner sep=1pt,anchor=south west] {#9};}

\newcommandx{\figByWidth}[9][5=0, 6=0, 7=0, 8=0,9=]{
\draw (#1,#2) node[anchor=south west,xshift=-16pt,yshift=-4pt] {\trimw{#3}{#4}{#5}{#6}{#7}{#8}};}
\newcommandx{\figByWidthb}[9][5=0, 6=0, 7=0, 8=0,9=]{
\draw (#1,#2) node[anchor=south west,xshift=-16pt,yshift=-4pt] {\trimwb{#3}{#4}{#5}{#6}{#7}{#8}};}

\newcommandx{\figByWidthWithLabel}[9][5=0, 6=0, 7=0, 8=0,9=]{
\draw (#1,#2) node[anchor=south west,xshift=-16pt,yshift=-4pt] {\trimw{#3}{#4}{#5}{#6}{#7}{#8}} node[draw=white,fill=white,inner sep=1pt,anchor=south west] {#9};}
\newcommandx{\figByWidthWithLabelb}[9][5=0, 6=0, 7=0, 8=0,9=]{
\draw (#1,#2) node[anchor=south west,xshift=-16pt,yshift=-4pt] {\trimwb{#3}{#4}{#5}{#6}{#7}{#8}} node[draw=white,fill=white,inner sep=1pt,anchor=south west] {#9};}

%% file: tex/intro.tex
\section{Introduction}  \label{sec:intro}

We describe an algorithm, called EigenWave, designed to compute accurate approximations of eigenvalues and
eigenfunctions associated with elliptic boundary-value problems. The algorithm, based on the recently developed
WaveHoltz scheme~\cite
{appelo2020waveholtz,EMWaveHoltzPengAppeloIEEE2022,appeloElWaveHoltz2022,overHoltzArXiv2025,overHoltzPartOne,overHoltzPartTwo},
solves a related initial-boundary-value problem (IBVP) for a time-dependent wave equation as part of an iteration that successively updates an initial condition of the IBVP. At each step of the iteration, the solution of the IBVP
is filtered in time over a period (or possibly multiple periods) associated with a chosen target frequency resulting
in an update of the initial data.  The filter function is chosen so that the relative contribution of the eigenmodes  of the updated
initial data whose eigenvalues are near the target frequency is enhanced. The
iteration is embedded within a matrix-free Arnoldi algorithm which enables the efficient computation of multiple
eigenpairs near the target frequency. The ability to select an arbitrary target frequency enables the computation of
eigenvalues anywhere in the spectrum, without the need to invert an indefinite matrix as is common with other
approaches. While EigenWave can accommodate a variety of approaches to approximate the spatial elliptic operator, we
employ finite differences at various orders of accuracy on overset grids for complex geometry. For example, Figure~\ref
{fig:rpiEigenvectors} displays some sample eigenvectors of the Laplacian operator computed using EigenWave and Figure~\ref
{fig:rpiGrid} shows an overset grid for this geometry. Importantly, the solution of the time-dependent wave equation of the IBVP
at each iteration can be obtained efficiently using implicit time-stepping which leads to an $O(N)$ algorithm, where $N$ is the
number of spatial grid points.  The approach thus provides a powerful new
tool for solving large-scale eigenvalue problems arising in continuum mechanics.

\input tex/rpiGridAndContoursFig

Several excellent algorithms exist for finding all eigenvalues of a not-too-large symmetric matrix~$A$, such as
the QR algorithm, divide-and-conquer method, or Jacobi algorithm~\cite{BaiDemmel2000}.  In addition,
many excellent schemes have been developed for computing a few eigenpairs of large-scale problems, see the survey articles~\cite{GolubVanDerVorst2000,Sorensen2002} and the book by Saad~\cite{NumericalMethodsForLargeEigenvalueProblemsSaadBook}, and
high quality software packages also exist, such as ARPACK~\cite{ARPACK}, SLEPSc~\cite{SLEPc2005}, Ansaszi~\cite{Anasazi2009}, PRIMME~\cite{PRIMME2010}, and EVSL~\cite{EVSL2019}.
Popular schemes for computing a few eigenpairs include those based on Arnoldi, such as the 
explicitly and implicitly restarted Arnoldi method~\cite{LehoucqSoresen1996,Morgan2000,EmadPetitonEdjlali2005,MorganZeng2006,FreitagSpence2009,Evstigneev2017,XueFeiElman2021}, 
the Krylov-Schur algorithm~\cite{Stewart2002,Lehoucq2002,BaglamaCalvettiReichel2003,ZhouSaad2008}, schemes based on the 
Davidson and Jacobi-Davidson methods~\cite{StathopoulosSaadWu1998,SleijpenVanDerVorst2000,Notay2005,Zhou2006,StathopoulosMcCombs2007,Zhou2010,RomeroRoman2014,Miao2018,Miao2020,Miao2023},
and subspace iteration schemes, such as FEAST which uses contour integration of the resolvent to select eigenvalues of interest~\cite{Polizzi2009,TangPolizzi2014,KestynPolizziTang2016,YinChanYeung2017,YeXiaChanCauleyBalakrishnan2017,GavinPolizzi2018,HorningTownsend2020,GopalakrishnanGrubisicOvall2020}.
For interior eigenvalues, many of these algorithms rely on inverting a shifted-matrix $A-\sigma I$. 
However, this indefinite matrix can be very difficult to invert other than by direct sparse factorization~\cite{BaiDemmel2000,Galgon2015,LiYang2021}, and the resulting cost can become prohibitive. 
Another approach, more closely related to EigenWave is to use polynomial preconditioners~\cite{IterativeMethodsBookSaad2003,Evstigneev2017,AustinTrefethen2015,EVSL2019,LiYang2021} whereby a Krylov-based method is applied to some polynomial of the matrix, where 
the polynomial is chosen to transform the spectrum into a more suitable form while keeping the eigenvectors the same.
EigenWave with explicit time-stepping of the IBVP followed by time-filtering leads to a high-degree polynomial preconditioner that importantly damps unwanted high-frequency modes which 
are sometimes an issue with polynomial preconditioners~\cite{EmbreeLoeMorgan2021}. 
This damping property is inherited from the stability of the explicit time-stepping, i.e.~the time-step is chosen
to keep the scheme stable. 
With implicit time-stepping, EigenWave can be interpreted as a rational polynomial preconditioner 
where importantly the rational polynomial can be applied with an $O(N)$ computational cost.
For implicit time-stepping the unwanted high-frequency modes are also damped provided at least $5$ time-steps per period, 
although in practice using $10$ time-steps per period seems to be a good choice (this is discussed later in this article).


In recent work, performed independently from our work, Nannen and Wess~\cite{NannenWess2024} have also used ideas from the WaveHoltz scheme~\cite{appelo2020waveholtz}
to develop a Krylov subspace iteration based on filtering solutions to the wave equation. This work is a nice complement to our work since
it uses explicit time-stepping, a different filter, and finite element approximations.
However, a key difference is that we also consider the use of implicit time-stepping in addition to explicit time-stepping.
Even though explicit time-stepping may be faster in some cases, 
it is the implicit time-stepping with a large time-step that leads to an $O(N)$ algorithm for large $N$. 
In addition, we show how the approach can be used in a \textit{matrix-free} manner 
with existing high-quality eigenvalue software, such as that found in SLEPSc~\cite{SLEPc2005} and ARPACK~\cite{ARPACK}. 
This can be important in practice since 
these existing software packages are accurate and robust due to a very sophisticated implementation\footnote{In their description of the implicitly restarted 
Arnoldi method~\cite{LehoucqSorensen2000}, Lehoucq and Sorensen comment that ``While we have presented the IRAM as much as possible in template form, we do not recommend implementation from this description. High quality software is freely available in the form of ARPACK. The fine detail of implementation is quite important to the robustness and ultimate success of the method.
While a working method could be obtained from the description given here, \textbf{it would almost certainly be deficient in some respect}'' (emphasis added).}.
We will show that when combined with Krylov-Schur or IRAM algorithms in a \textsl{matrix-free} manner, 
EigenWave provides an efficient algorithm to compute multiple eigenpairs to 
high accuracy using
just a few (e.g.~$3$--$5$) wave-solves\footnote{A wave-solve is defined as a solution of the IBVP followed by an application of the time filter.} per eigenpair. 

The EigenWave algorithm is analyzed, and its behavior demonstrated, by 
computing eigenpairs of the Laplacian operator in various geometries
using overset grids~\cite{CGNS}. Overset grids are used to efficiently discretize the wave equation to high-order accuracy 
in space using finite-differences~\cite{adegdmi2020}. Second-order accuracy in time is sufficient for EigenWave since time accuracy does not effect
the spatial accuracy of the computed eigenvectors.
Interestingly, EigenWave first computes the desired eigenvectors but with different eigenvalues that lie
in the interval $[-\half,1]$; the WaveHoltz time-filter has shifted all eigenvalues of the Laplacian to this new interval while the eigenvectors are unaffected and computed accurately.
The desired eigenvalues of the Laplacian are then computed in a post-processing step using a Rayleigh quotient involving the computed eigenvectors.
The EigenWave filter can thus be viewed as a \textsl{spectral transform} in the parlance of eigenvalue algorithms. The most common spectral transform is the shift-and-invert transform, $(A-\sigma I)^{-1}$, which has a one-to-one and onto mapping between eigenvalues of the transformed and un-transformed problems. There is no such 
mapping between eigenvalues when using EigenWave, and instead the desired eigenvalue is determined after computing the corresponding eigenvector by using the Rayleigh quotient.
EigenWave is particularly efficient when used in conjunction with implicit time-stepping since, amazingly, only a few
time-steps are required (e.g.~10 total time-steps) when integrating the IBVP over one period of the target frequency.
Further, the matrix arising from implicit time-stepping is well suited for solution with fast methods
such as multigrid~\cite{automg,multigridWithNonstandardCoarsening2023}.
It is also interesting to note that although the solution of the wave equation on overset grids 
normally requires upwind dissipation for stability~\cite{ssmx2023,mxsosup2017},
it has been found in practice that no upwind dissipation is generally needed with EigenWave.  It is also worth noting that EigenWave is useful to obtain a set of selected eigenvectors used for deflation in the newly developed OverHoltz algorithm~\cite{overHoltzArXiv2025,overHoltzPartOne,overHoltzPartTwo}.

For reference, Table~\ref{tab:notation} provides a summary of some of the symbols and notation that will be introduced in subsequent sections.
\begin{table}\footnotesize
\begin{center}
  \begin{tabular}{lcl} 
     \hline 
    $\Lc$      &:& elliptic operator in the eigenvalue problem $\Lc \phi=-\lambda^2 \phi$ and wave equation $\p_t^2 = \Lc w$. \\
    $\Bc$      &:& boundary condition operator. \\
    $(\lambda_j,\phi_j(\xv))$ &:& continuous eigenvalues and eigenfunctions of $(\Lc,\Bc)$, $j=1,2,\ldots$. \\
    $L_{ph}$   &:& p-th order accurate approximation to $\Lc$. \\
    $B_{ph}$   &:& p-th order accurate boundary conditions. \\
    $N$        &:& total number of grid points. \\
    $\NGd$     &:& total number of eigenvectors in the discrete problem. \\
    $(\lambda_{h,j},\Phiv_j)$ &:& discrete eigenvalues and eigenvectors of $( L_{ph},B_{ph})$, $j=1,2,\ldots,\NGd$.  \\
    $w(\xv,t)$   &:& solution to the continuous wave equation. \\
    $W_\iv^n$  &:& discrete solution $W_\iv^n\approx w(\xv_\iv,t^n)$, for time-level $t^n=n\dt$ and grid index $\iv=[i_1,i_2,i_3]$. \\
    $\omega$   &:& target frequency (target eigenvalue). \\
    $T=(2\pi)/\omega$   &:& period corresponding to $\omega$. \\
    $v^{(k)}(\xv)$  &:& initial condition for $w(\xv,0)$ and approximation to an eigenvector. \\
    $\beta(\lambda,\omega)$ &:& time-continuous filter function. \\
    $N_p$      &:& number of time periods over which the filter is integrated. \\
    $\Nits$    &:& number of implicit time-steps per-period. \\
    $\Tf = \Np T$ &:& final time for each wave solve. \\
    $\Aw$      &:& continuous EigenWave linear operator with eigenvalues $\beta_j=\beta(\lambda_j;\omega)$\\
               & & and eigenfunctions $\phi_j$. \\
    $\Aw_{ph}$ &:& p-th order accurate approximation to $\Aw$ with eigenvalues $\beta_{h,j} \approx \beta(\lambda_{h,j};\omega)$  \\
               & & and eigenvectors $\Phiv_j$ .\\
    IRAM       &:& Implicitly Restarted Arnoldi Method, a sophisticated algorithm for \\
               & &  finding several eigenpairs of a (large) matrix. Implemented in ARPACK.\\
    Krylov-Schur (KS) &:& Variation of the IRAM algorithm. Implemented in SLEPSc. \\
    $N_r$      &:& number of requested eigenvalues, input to IRAM or KS. \\
    $N_c$      &:& number of computed eigenvalues returned from IRAM or KS. \\
    $N_a$      &:& dimension of the Krylov subspace for IRAM and KS, typically $N_a=2 N_r +1$\\
               & & (keeping additional vectors improves convergence of the $N_r$ requested). \\
    \hline 
  \end{tabular}
  \end{center}
  \caption{Nomenclature} \label{tab:notation}
\end{table}

%% file: tex/rpiGridAndContoursFig.tex
{
\newcommand{\drawContour}[7]{%
\begin{scope}[#1]
\draw(0.0,0) node[anchor=south west,xshift=-4pt,yshift=+10pt] {\trimfiga{fig/#2}{\figWidtha}};
\end{scope}
}
\newcommand{\cbWidth}{.2cm}
\newcommand{\cbHeight}{4cm}
\newcommand{\xcb}{.5cm}
\newcommand{\ycb}{-.2cm}
\setlength{\ycbTop}{\ycb+\cbHeight}
\setlength{\ycbMid}{\ycb+\cbHeight*\real{.5}}
\newcommand{\trimfigcb}[3]{\includegraphics[width=#2, height=#3, clip, trim=17cm 2.35cm 1.65cm 2.35cm]{#1}}
\newcommand{\figWidth}{5.25cm}
\newcommand{\figHeight}{2.65cm}
\newcommand{\figWidtha}{5.25cm}
\newcommand{\trimfiga}[2]{\trimw{#1}{#2}{.12}{.117}{.32}{.32}}
\begin{figure}[htb]
\begin{center}
\begin{tikzpicture}
   \useasboundingbox (0,.60) rectangle (3.06*\figWidth,2.1*\figHeight);  

   \begin{scope}[xshift=-.5cm,yshift=\figHeight]
     \drawContour{xshift=0.00*\figWidth}{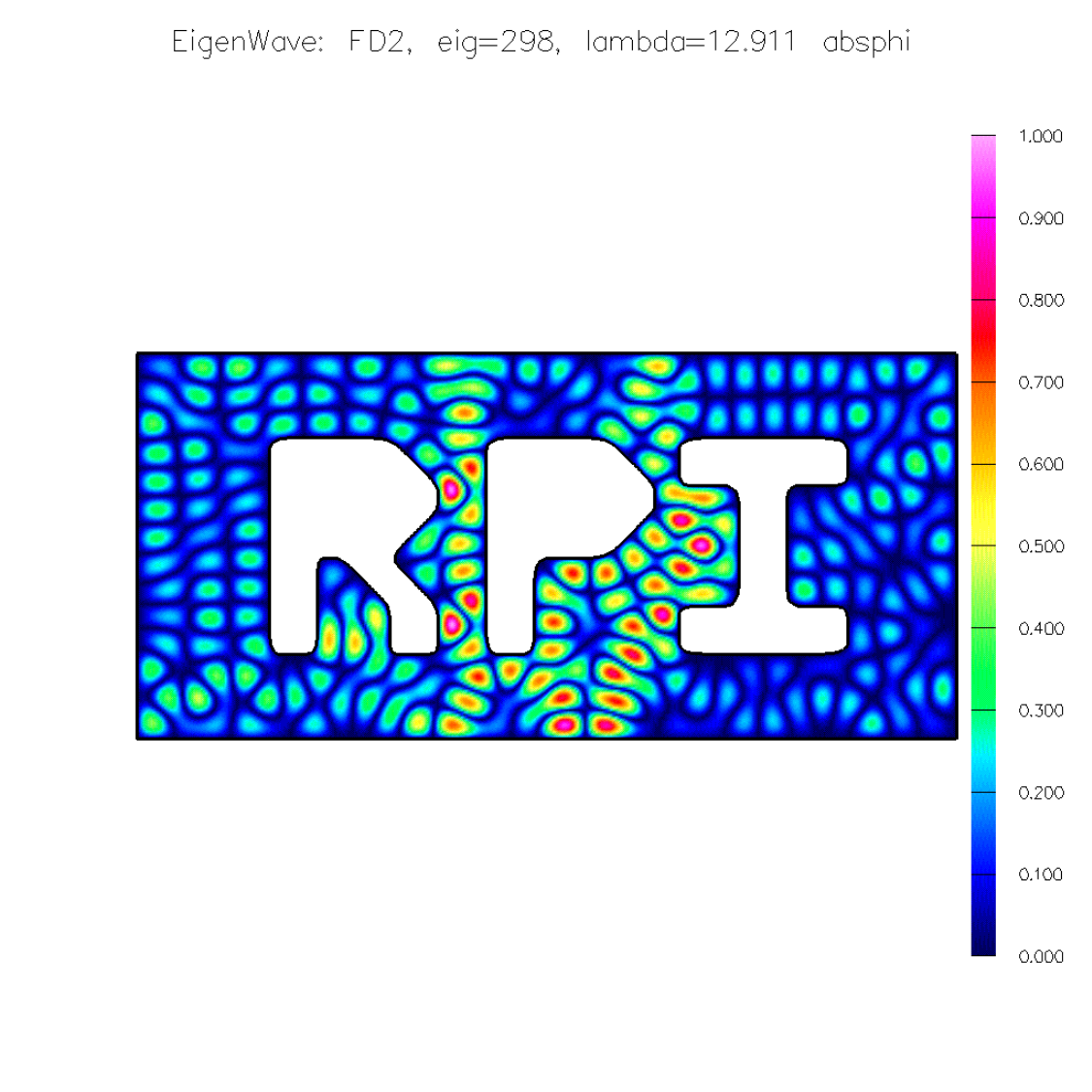}{$|\phi|$}{$v$}{$t=0.3$}{$0$}{$1.0$}     
     \drawContour{xshift=1.02*\figWidth}{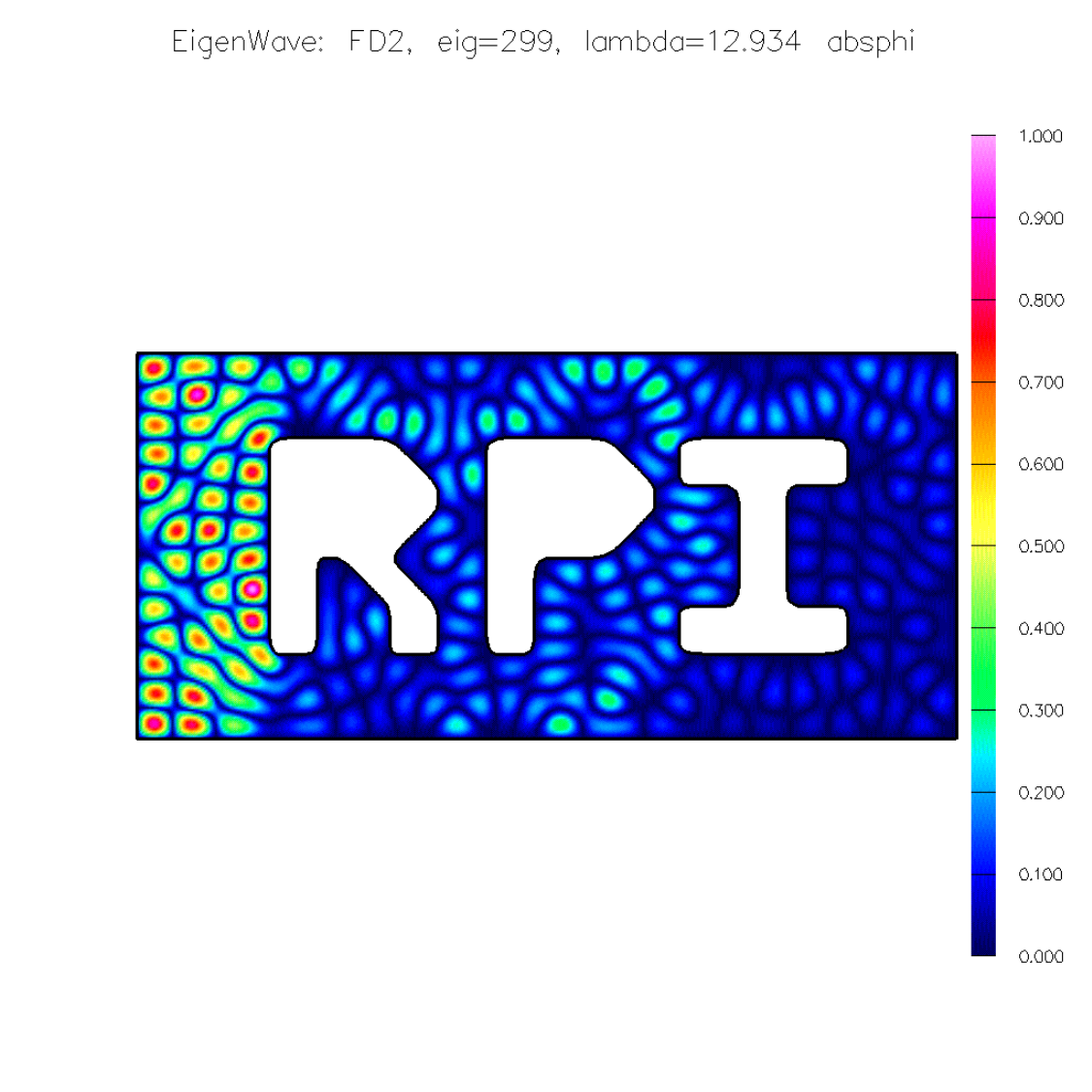}{$|\phi|$}{$v$}{$t=0.3$}{$0$}{$1.0$}    
     \drawContour{xshift=2.04*\figWidth}{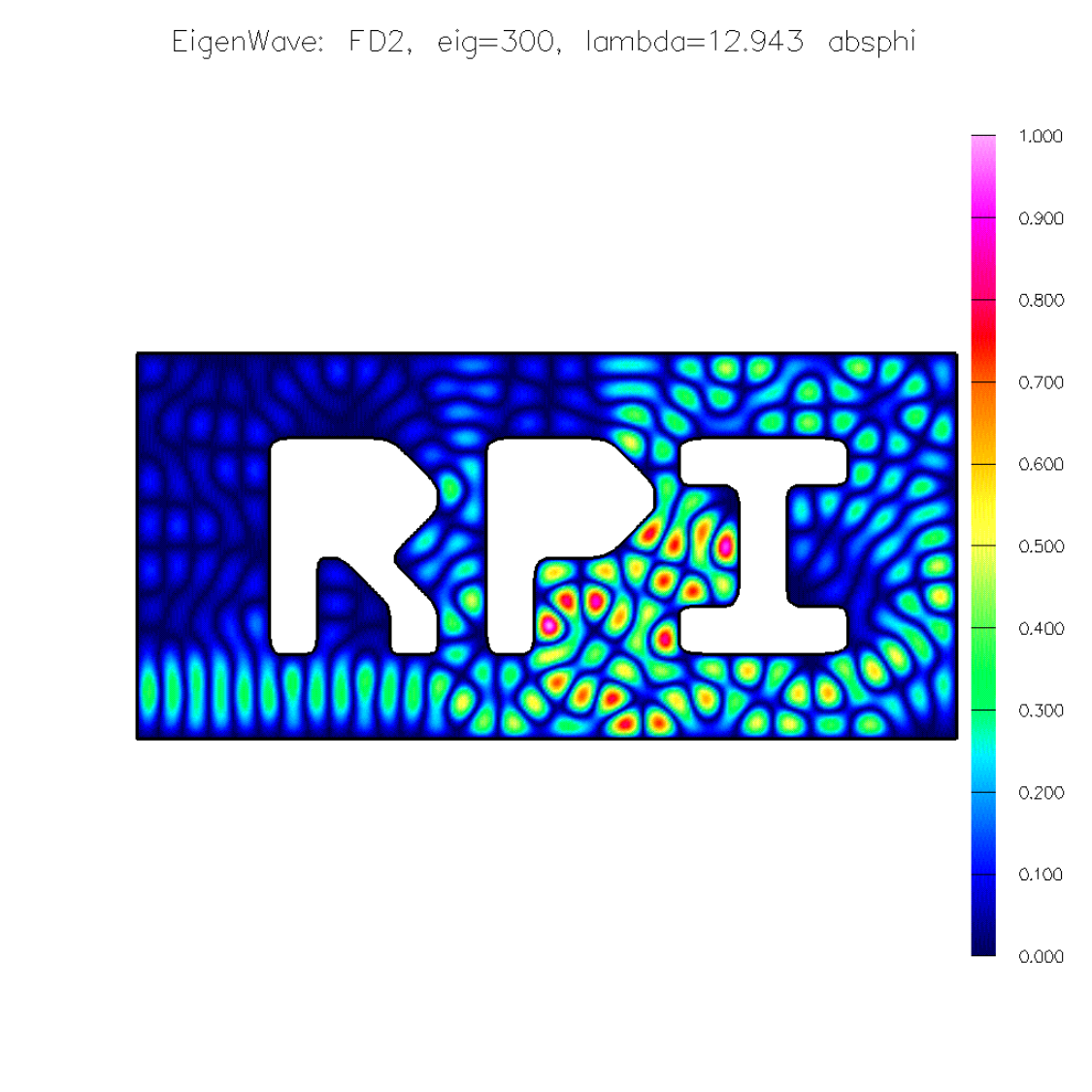}{$|\phi|$}{$v$}{$t=0.3$}{$0$}{$1.0$}    
   \end{scope}
   
   \begin{scope}[xshift=-.5cm,yshift=0cm]
     \drawContour{xshift=0.00*\figWidth}{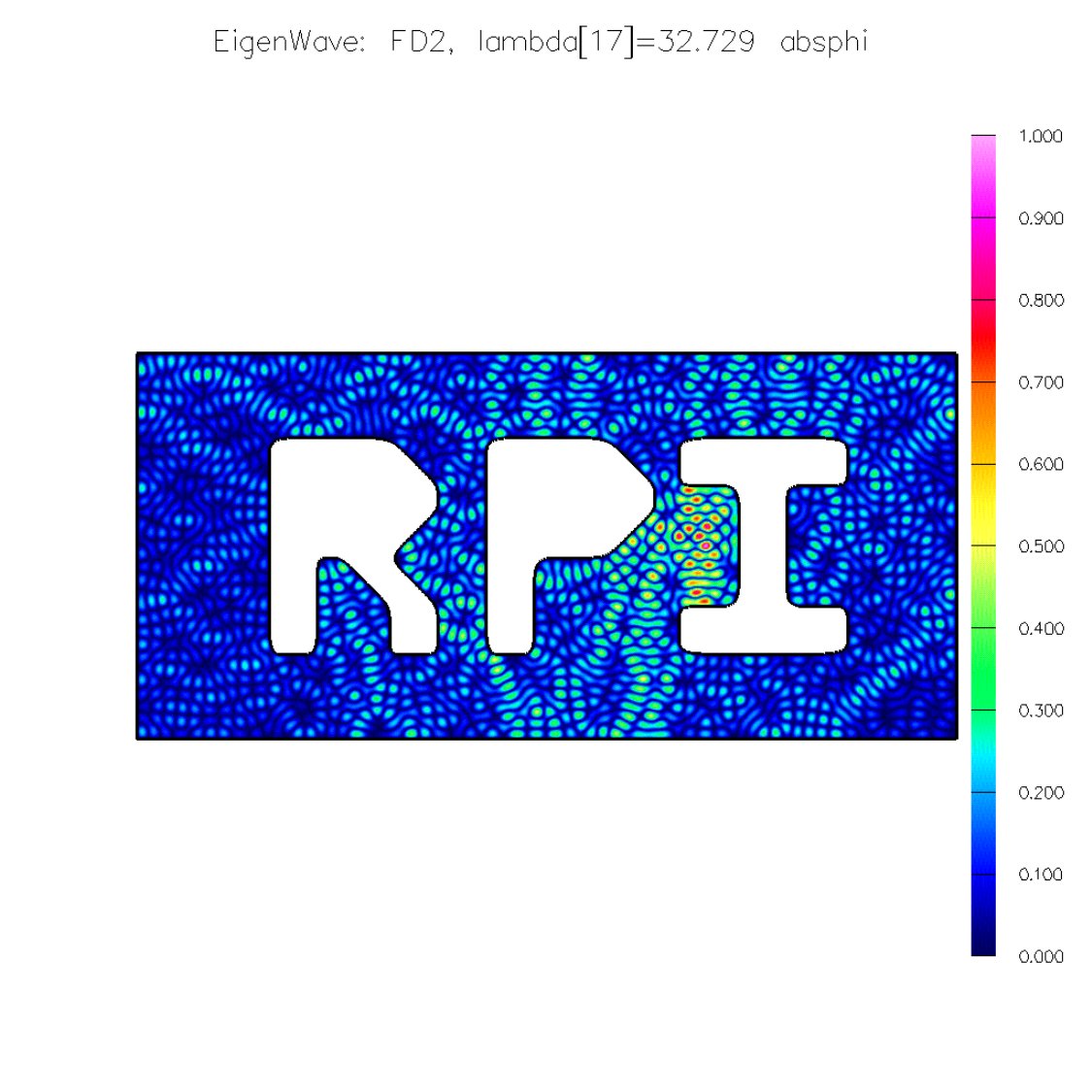}{$|\phi|$}{$v$}{$t=0.3$}{$0$}{$1.0$}     
     \drawContour{xshift=1.02*\figWidth}{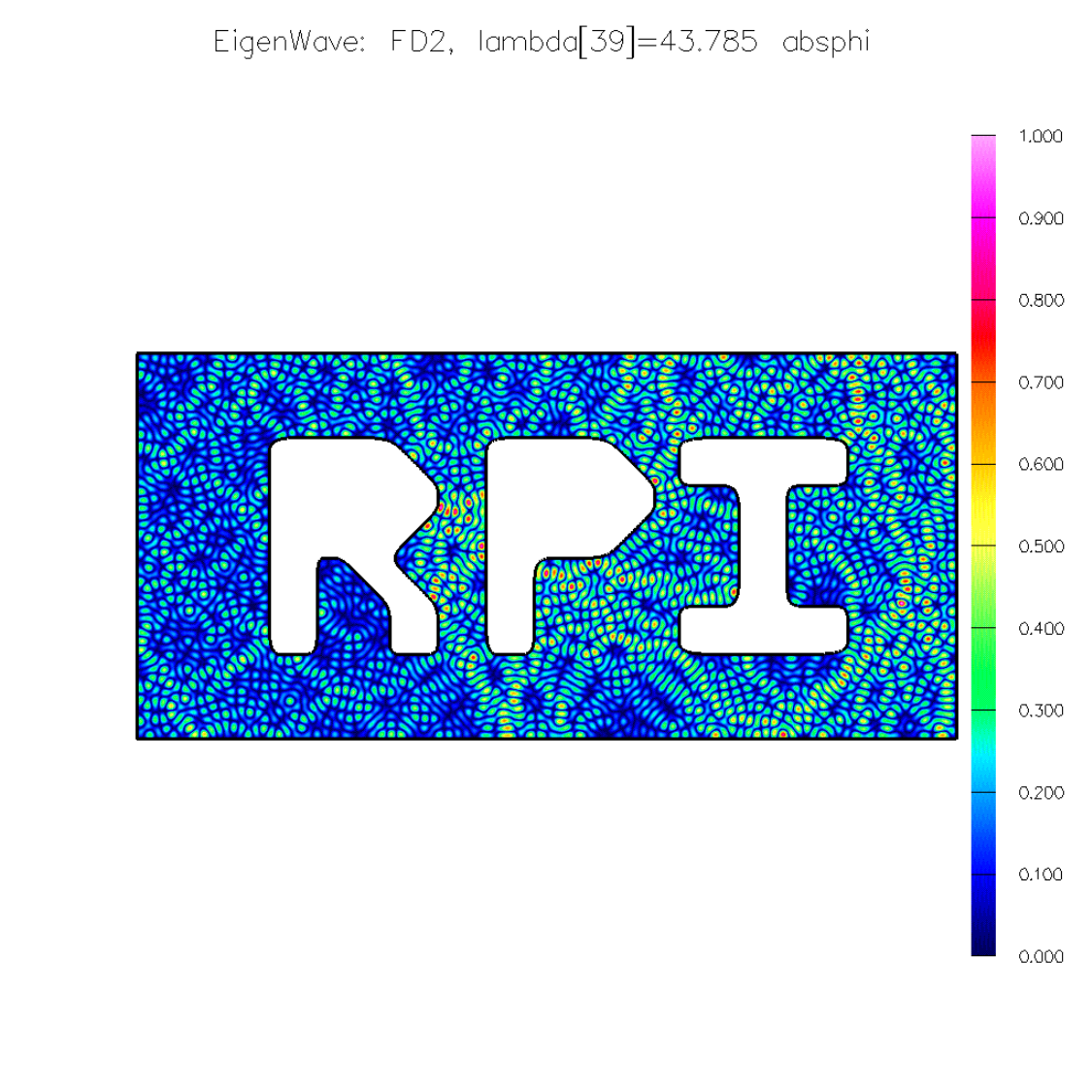}{$|\phi|$}{$v$}{$t=0.3$}{$0$}{$1.0$}     
     \drawContour{xshift=2.04*\figWidth}{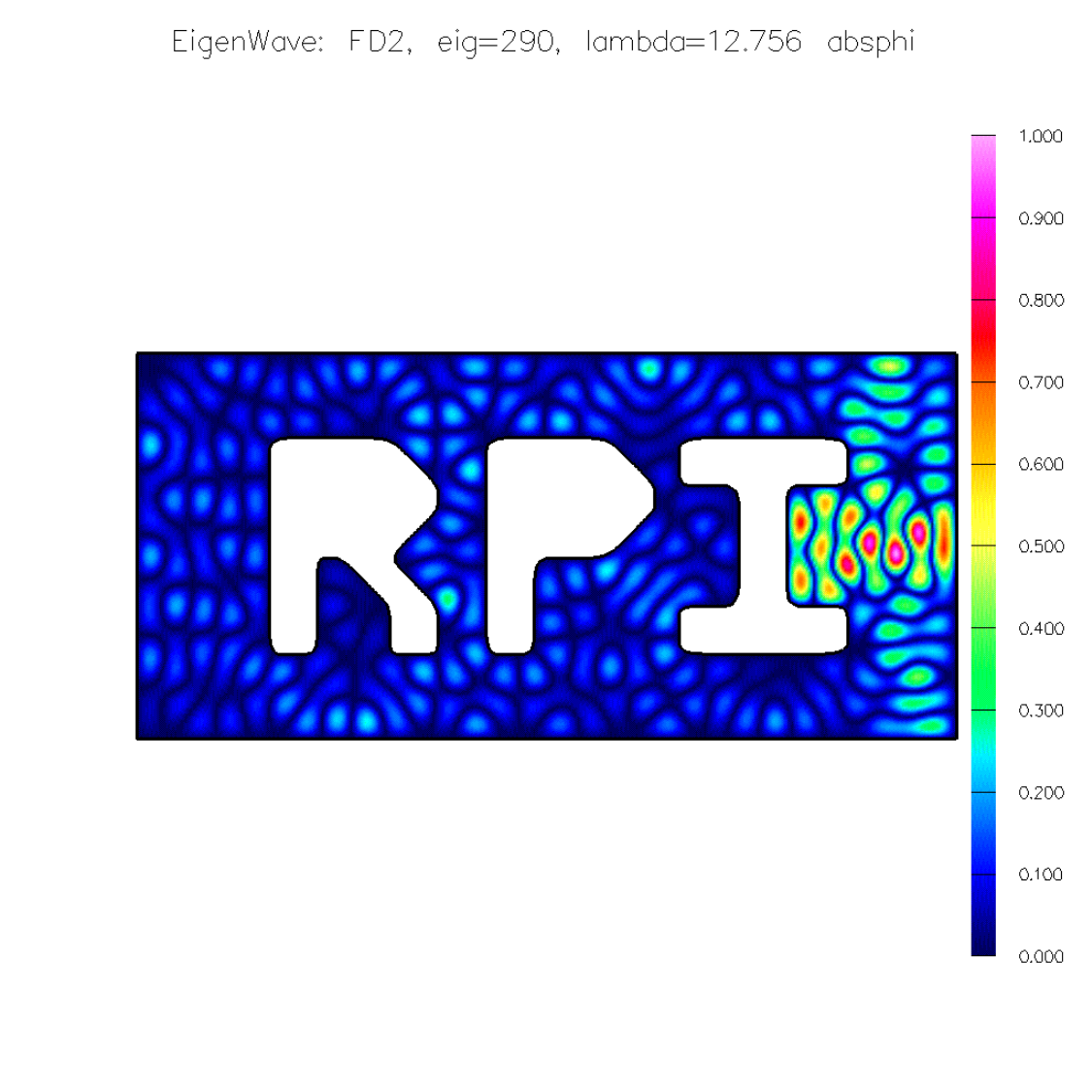}{$|\phi|$}{$v$}{$t=0.3$}{$0$}{$1.0$}     
   \end{scope}

  \begin{scope}[xshift=5cm,yshift=.2cm]
    \draw (\xcb,\ycb) node[anchor=south west,xshift=0.25cm,yshift=.5cm,rotate=-90] {\trimfigcb{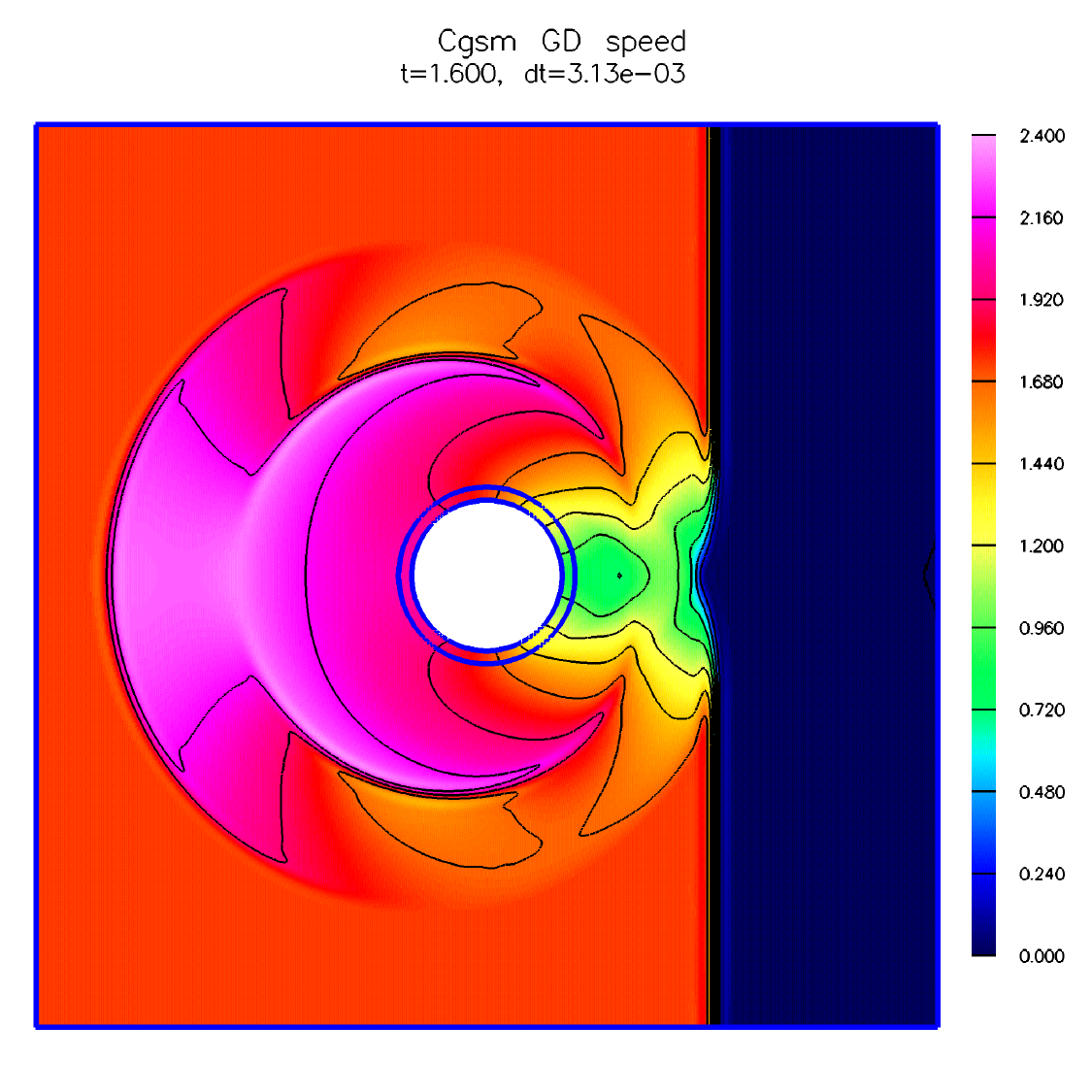}{\cbWidth}{\cbHeight}};
    \draw (.8,0) node[anchor=north,xshift=+3pt,yshift=+2pt] {\scriptsize $0.0$};
    \draw (4.8,0) node[anchor=north,xshift=+0pt,yshift=+2pt] {\scriptsize $1.0$};
  \end{scope}

\end{tikzpicture}
\end{center}
\caption{Absolute value of selected eigenvectors of the Laplacian with Dirichlet boundary conditions computed with the EigenWave algorithm.
    }
\label{fig:rpiEigenvectors}
\end{figure}
}

{
\newcommand{\figWidth}{6cm}
\newcommand{\figHeight}{2.8cm}
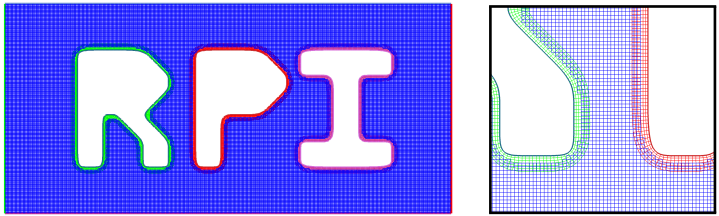
\begin{figure}[htb]
\begin{center}
\begin{tikzpicture}
   \useasboundingbox (0,.65) rectangle (1.65*\figWidth,1*\figHeight);  

   \begin{scope}[xshift=0cm,yshift=0cm]
     \figByWidth{   0}{0.}{fig/rpiGridG4}{6.cm}[0.025][0.025][0.275][0.275]
     \figByWidthb{ 6.5}{.05}{fig/rpiGridG4Zoom}{2.95cm}[0.1][0.1][0.165][0.1]
   \end{scope}  

\end{tikzpicture}
\end{center}
\caption{Overset grid (and magnification) for the letters in RPI. 
    }
\label{fig:rpiGrid}
\end{figure}
}


\bogus{
\newcommand{\drawContour}[7]{%
\begin{scope}[#1]
\draw(0.0,0) node[anchor=south west,xshift=-4pt,yshift=+10pt] {\trimfiga{fig/#2}{\figWidtha}};
\begin{scope}[xshift=.9cm,yshift=+7pt]
  \draw (\xcb,\ycb) node[anchor=south west,xshift=0.25cm,yshift=.5cm,rotate=-90] {\trimfigcb{fig/colourBarLines}{\cbWidth}{\cbHeight}};
  \draw (.8,0) node[anchor=north,xshift=+3pt,yshift=+2pt] {\scriptsize $#6$};
  \draw (4.8,0) node[anchor=north,xshift=+0pt,yshift=+2pt] {\scriptsize $#7$};
\end{scope}
\end{scope}
}
\newcommand{\cbWidth}{.2cm}
\newcommand{\cbHeight}{4cm}
\newcommand{\xcb}{.5cm}
\newcommand{\ycb}{-.2cm}
\setlength{\ycbTop}{\ycb+\cbHeight}
\setlength{\ycbMid}{\ycb+\cbHeight*\real{.5}}
\newcommand{\trimfigcb}[3]{\includegraphics[width=#2, height=#3, clip, trim=17cm 2.35cm 1.65cm 2.35cm]{#1}}
\newcommand{\figWidtha}{6.5cm}
\newcommand{\trimfiga}[2]{\trimw{#1}{#2}{.12}{.117}{.32}{.32}}
\begin{figure}[htb]
\begin{center}
\begin{tikzpicture}
   \useasboundingbox (0,.75) rectangle (13.5,6.5);  

   \begin{scope}[xshift=-.5cm,yshift=0cm]
     \drawContour{xshift=0cm,yshift=0.0cm}{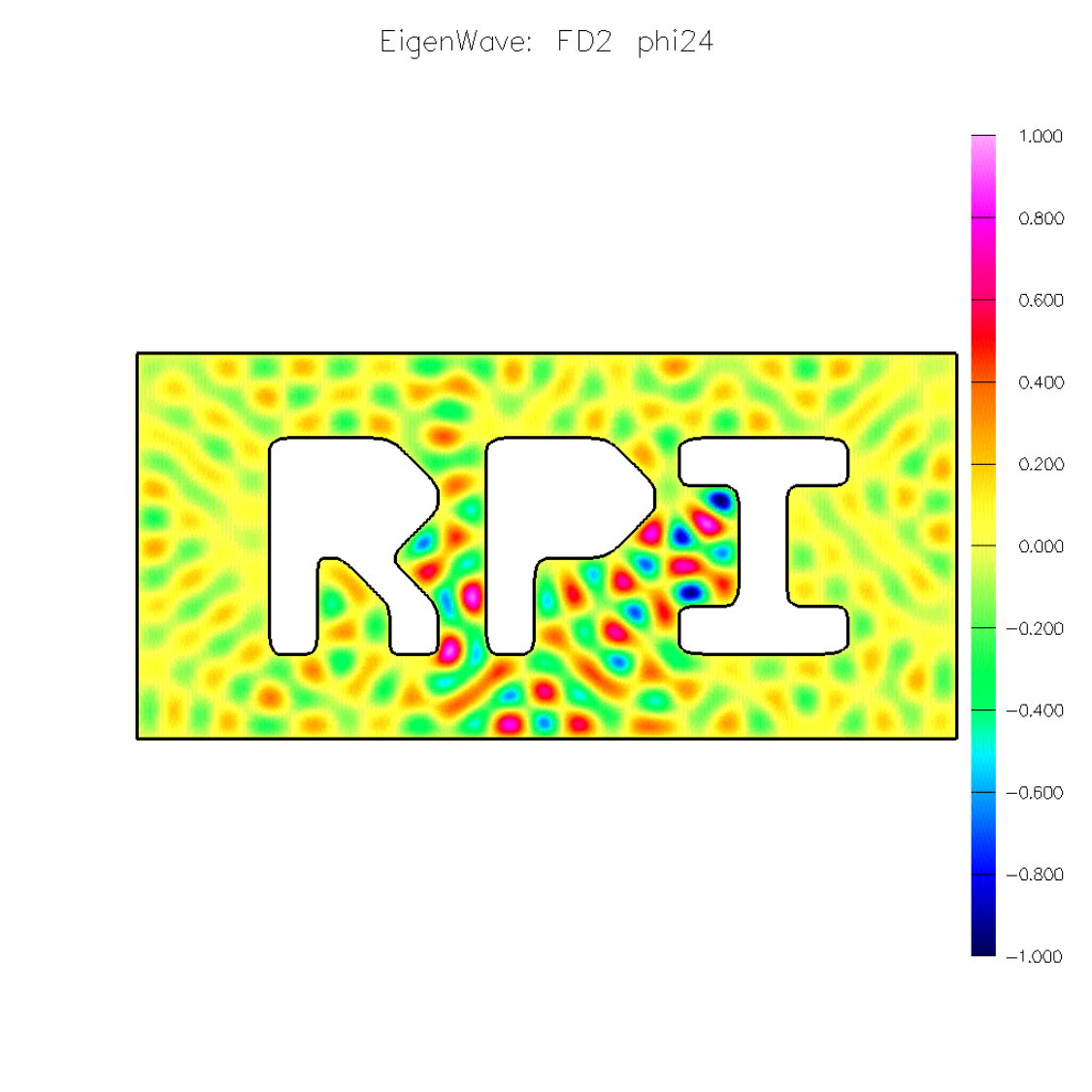}{$\phi^{(308)}$}{$v$}{$t=0.3$}{$-1.0$}{$1.0$};     
     \drawContour{xshift=7cm,yshift=0.0cm}{rpiG16O2phi13}{$\phi^{(334)}$}{$v$}{$t=0.3$}{$-1.0$}{$1.0$};     
   \end{scope}

   \begin{scope}[xshift=1.75cm,yshift=3.15cm]
     \figByWidth{   0}{0.5}{fig/rpiGridG4}{6.cm}[0.025][0.025][0.275][0.275]
     \figByWidthb{ 6.5}{.55}{fig/rpiGridG4Zoom}{3cm}[0.1][0.1][0.165][0.1]
   \end{scope}   

\end{tikzpicture}
\end{center}
\caption{At top are overset grid (and zoom) for the letters in RPI. At bottom are selected eigenvectors for the Laplacian computed with the \EigenWaveb algorithm.
    }
\label{fig:rpiGridAndContoursII}
\end{figure}



}

%% file: tex/algorithm.tex
\section{Problem specification and the EigenWave algorithm} \label{sec:eigenWaveAlgorithm}

Consider the problem of computing numerical approximations of 
selected eigenvalues and eigenfunctions associated with a boundary-value problem (BVP) involving an elliptic PDE. 
Let~$\Omega \subset \Real^{\nd}$ be a bounded domain in~$\nd$ space dimensions with boundary~$\partial\Omega$.
The BVP is defined in terms of an elliptic operator~$\Lc$,  
along with boundary conditions given by a boundary operator~$\Bc$ representing Dirichlet, Neumann or Robin boundary conditions.
If $c>0$ denotes a wave-speed then $\Lc$ is taken here to be $\Lc=c^2 \Delta$ where $\Delta$ is the Laplacian. 
The eigenvalue problem for~$\Lc$, with homogeneous boundary conditions specified~by $\Bc$, is given by
\bse
\label{eq:eigBVP}
\bat
  &  \Lc \phi = - \lambda^2 \, \phi, \qquad&&  \xv\in\Omega, \qquad \\
  &  \Bc \phi = 0,                   \qquad && \xv\in\partial\Omega ,
\eat
\ese
where $\phi=\phi(\xv)$ is an eigenfunction corresponding to the eigenvalue\footnote{Here $\mu=-\lambda^2$ is used to be consistent with the discussion in previous WaveHoltz articles.}  $\mu=-\lambda^2$.  For the problems of interest here, the eigenvalues $\mu$ 
are real and non-positive\footnote{This places some restrictions on the coefficients in the Robin boundary condition~\cite{StraussBook2007}.} 
so that we can take $\lambda\ge 0$ without loss of generality.  Also, while the set of eigenfunctions of~\eqref{eq:eigBVP} are linearly independent, the eigenvalues need not be distinct.\footnote{Note that the EigenWave algorithm can be extended to more general $\Lc$, variable coefficients, and more general boundary conditions that may lead to complex-valued eigenpairs.}

The EigenWave algorithm is based on the solution of a related time-dependent initial-boundary-value problem (IBVP).  Let
$w(\xv,t)$ solve the IBVP for the wave equation over a time interval $[0,T_f]$ given by 
\bse
\label{eq:waveHoltz}
\bat
  &  \p_t^2 w =  \Lc w ,                         \quad && \xv\in \Omega,\quad \;\;0<t<T_f,        \\
  &  \Bc w(\xv,t) = 0  ,                         \quad && \xv\in \partial\Omega,\quad 0<t<T_f, \label{eq:waveBC}\\
  &  w(\xv,0) = v^{(k)}(\xv),                          \quad && \xv\in \Omega, \\
  &  \p_t w(\xv,0)=0,                            \quad && \xv\in \Omega,
\eat
\ese
for some initial function $v^{(k)}(\xv)$, where $k$ denotes an iteration number. For a given solution of \eqref{eq:waveHoltz}, define the time filter
\ba
  v^{(k+1)}(\xv) 
    = \frac{2}{T_f}\int_0^{T_f}\left(\cos(\omega t)-\frac14\right) \, w(\xv,t; v^{(k)} ) \,dt ,
     \label{eq:waveHoltzFilterStep}
\ea
where we have indicated that $w$ depends on the choice of the initial function $v^{(k)}(\xv)$ and $\omega$ is a chosen target frequency (target eigenvalue).  The final time $T_f$ is related to $\omega$ by 
\ba
  T_f \eqdef \Np \, \f{2\pi}{\omega} ,
\ea
where $\Np$ is an integer defining the number of periods over which the IBVP is integrated.
The procedure of solving the IBVP in~\eqref{eq:waveHoltz} with initial condition $v^{(k)}(\xv)$
and then applying the time filter in~\eqref{eq:waveHoltzFilterStep} corresponds to one step in the WaveHoltz iteration, and this
defines a linear operator $\Aw=\Aw(\omega)$ that maps $v^{(k)}$ to $v^{(k+1)}$, 
\ba
    v^{(k+1)} = \Aw\, v^{(k)}.
    \label{eq:Amap}
\ea
As shown in Section~\ref{sec:analysis}, the operator $\Aw$ has exactly the same eigenfunctions as $\Lc$ in~\eqref{eq:eigBVP}, while the eigenvalues of $\Aw$ are different. The 
eigenvalues of $\Aw$ all lie in the interval $[-\half,\,1]$, and importantly, 
the largest eigenvalues of $\Aw$  generally correspond to eigenvalues $\lambda$ near the chosen target frequency~$\omega$.
The WaveHoltz step (called a \textsl{wave-solve}) consisting of the solution of the wave equation in~\eqref{eq:waveHoltz}
 and an application of the time filter in~\eqref{eq:waveHoltzFilterStep} has thus transformed
the eigenvalue problem for $\Lc$ into a new eigenvalue problem for $\Aw$ whose
spectrum is more easily computed by standard eigenvalue algorithms. 
Note that the target frequency $\omega$ can be adjusted
to select eigenvalues in different intervals.

\input tex/eigenWaveAlgorithm

The basic EigenWave algorithm, given in Algorithm~\ref{alg:EigenWave}, employs a power iteration on the operator $\Aw$ in~\eqref{eq:Amap} 
to determine an eigenpair $(\lambda,\phi(x))$ for $\lambda$ near the target frequency $\omega$. 
Note that $(u,v)$ denotes the usual $L_2$ inner product on $\Omega$ and $\Vert u\Vert\sp2=(u,u)$.
The algorithm starts from some initial guess $v^{(0)}$ and computes a sequence of approximations $v^{(k)}$, $k=1,2,\ldots$.
At each iteration the current guess $v^{(k)}$ is used as the initial condition to the wave equation solver.
The new iterate $v^{(k+1)} = \Aw v^{(k)}$ is computed as the weighted time-integral of the wave equation solution.
Note that the power iteration does not directly compute estimates to an eigenvalue $\lambda_j$ of $\Lc$, 
but rather estimates to an eigenvalue $\beta_j$ of $\Aw$ given by $\beta\sp{(k+1)}$, $k=0,1,2,\ldots$, in Step~8 of the algorithm.
However, upon convergence, an approximation to the eigenvalue $\lambda_j$ 
can be computed from the approximate eigenfunction $\phi(\xv)$ using a Rayleigh quotient involving~$\Lc$, see Steps~14 and~15.
Note that in practice the time-integral on line~\ref{alg:filter} can be accumulated inside the time-stepping loop to avoid storing the solution to the wave equation over time.

%% file: tex/eigenWaveAlgorithm.tex
\renewcommand{\algFontSize}{\small}
\begin{algorithm}
\algFontSize 
\caption{EigenWave algorithm - power iteration on $\Aw$ to compute one eigenpair $(\lambda,\phi)$.}
\begin{algorithmic}[1]

  \Function{$[\lambda,\phi]$ = EigenWave}{$\omega$,$v^{(0)}$,$\Np$}  
    \State // Input: target frequency~$\omega$, initial guess~$v^{(0)}$ with norm one, number of periods $\Np$
    \State $T=2\pi/\omega$,~~$T_f=\Np T$ \Comment Period and final time.  \label{alg:init}
    \For{k=0,1,\ldots} \Comment Start EigenWave~iterations.

      \State $w^{(k)}(\xv,0) = v^{(k)}(\xv)$ \Comment Initial condition for wave equation solve.
      \State $w^{(k)}(\xv,t)$ = \Call{solveWaveEquation}{$w^{(k)}(\xv,0)$,\,$T_f$} \Comment Solve for $\wv(\xv,t)$ for $t\in[0,T_f]$.\label{alg:solveWave}
      \State $\displaystyle v^{(k+1)}(\xv) = \f{2}{T_f} \int_{0}^{T_f} \left( \cos(\omega t) - \f{1}{4} \right) \, w^{(k)}(\xv,t; v^{(k)}) \, dt$
           \Comment Time filter $w(\xv,t)$ to give $v^{(k+1)} = \Aw v^{(k)}$. \label{alg:filter}
      \State $\beta^{(k+1)} = (v^{(k+1)},v^{(k)})$  \Comment Rayleigh quotient estimate for eigenvalue of $\Aw$
      \State  $v^{(k+1)} =  v^{(k+1)}/\|  v^{(k+1)} \|$ \Comment Normalize
      \If{ $\| v^{(k+1)} - \sign(\beta^{(k+1)}) \,  v^{(k)} \| < {\rm tolerance}$ } \Comment $\sign(\beta^{(k+1)})=\pm 1$
        \State break from loop
      \EndIf
    \EndFor    \Comment End EigenWave iterations.
    \State $\displaystyle \phi(\xv) = v^{(k+1)}(\xv) $ \Comment Approximate eigenfunction.
     \State $\displaystyle \lambda = \sqrt{ (\phi, -\Lc \phi) }$   \Comment Approximate eigenvalue of $\Lc$ from a Rayleigh quotient.
 \EndFunction
\end{algorithmic} 
\label{alg:EigenWave}
\end{algorithm}

%% file: tex/analysis.tex
\newcommand{\eHat}{\hat{e}}
\section{Analysis of the continuous EigenWave algorithm}  \label{sec:analysis}

In this section some useful properties of the EigenWave~algorithm are established.
The analysis here closely follows the analysis of the WaveHoltz algorithm as given in~\cite{appelo2020waveholtz}.

\begin{theorem}[Eigenvalues of the EigenWave operator]
The EigenWave operator $\Aw$, defined in~\eqref{eq:Amap}, has the same eigenfunctions $\phi_j(\xv)$, $j=0,1,2,\ldots$,
as the operator $\Lc$ (and boundary conditions) in~\eqref{eq:eigBVP} but with different eigenvalues 
\ba 
   \beta_j \eqdef \beta(\lambda_j;\omega)\in\hbox{$[-\half,\,1]$} ,
\ea
where $\beta=\beta(\lambda;\omega)$ is the WaveHoltz filter function defined by 
\ba
  \beta(\lambda;\omega) \eqdef \frac{2}{T_f}\int_0^{T_f}\left(\cos(\omega t)-\frac14\right) \, \cos(\lambda t) \, dt. 
  \label{eq:filterFunction}
\ea
\end{theorem}
\begin{proof}

The eigenvalue problem in~\eqref{eq:eigBVP}
has a complete set of eigenfunctions, $\phi_j(\xv)$, $j=0,1,2,\ldots$, for eigenvalues $\lambda_j \ge 0$.
The solution to the wave equation, $w(\xv,t)$, initial condition function, $v^{(k)}$, and next iterate $v^{(k+1)} =\Aw \, v^{(k)}$ in~\eqref{eq:Amap} 
can be expanded in terms of the eigenfunctions,
\ba
  & w(\xv,t) = \sum_{j=0}^\infty \wHat_{j}(t)\, \phi_j(\xv), \quad v^{(k)}(\xv) = \sum_{j=0}^\infty \vHat_{j}^{(k)} \, \phi_j(\xv),
   \quad v^{(k+1)}(\xv) = \sum_{j=0}^\infty \vHat_j^{(k+1)} \, \phi_j(\xv), \label{eq:eigenExpansion}
\ea
where $\wHat_{j}(t)$, $\vHat_{j}^{(k)}$, and $\vHat^{(k+1)}$ denote generalized Fourier coefficients in the expansions for $w$, $v^{(k)}$, and $v^{(k+1)}$, respectively.
Substituting~\eqref{eq:eigenExpansion} into the wave equation IBVP~\eqref{eq:waveHoltz}, and using $\Lc \phi_j = - \lambda_j^2 \, \phi_j$, leads to an initial-value problem for each time-dependent Fourier coefficient $\wHat_j (t)$, $j=0,1,2,\ldots$, given by
\bse
\label{eq:wHatEquation}
\ba
  & \p_t^2 \wHat_j = -\lambda_j^2 \, \wHat_j , \\
  & \wHat_j (0) = \vHat_{j}^{(k)}, \\
  & \p_t\wHat_j (0) =0,
\ea
\ese
whose solution is easily found to be 
\ba
   \wHat_j(t) = \vHat_{j}^{(k)} \, \cos(\lambda_j t) \, . \label{eq:wHatSolution}
\ea
Substituting the expansions~\eqref{eq:eigenExpansion} into the time filter~\eqref{eq:waveHoltzFilterStep}
implies that, in terms of the generalized Fourier coefficients, the time filtering step takes the form
\ba
  \vHat_j^{(k+1)} =  \frac{2}{T_f}\int_0^{T_f}\left(\cos(\omega t)-\frac{1}{4}\right) \, \wHat_j(t)  \,dt . \label{eq:timeFilterFourier}
\ea
Using~\eqref{eq:wHatSolution} in~\eqref{eq:timeFilterFourier} together with the formula for the filter function~\eqref{eq:filterFunction} gives 
\ba
 &  \vHat_j^{(k+1)}  = \beta(\lambda_j;\omega)  \, \vHat_j^{(k)},  \qquad j=0,1,2,\ldots \label{eq:AvFourier}
\ea
Recall that $v^{(k+1)} = \Aw \, v^{(k)}$ and so~\eqref{eq:AvFourier} shows 
\ba
    \Aw  \sum_{j=0}^\infty \vHat_{j}^{(k)}\, \phi_j(\xv)  = \sum_{j=0}^\infty \beta(\lambda_j;\omega) \, \vHat_{j}^{(k)} \, \phi_j(\xv) ,
\ea
for any coefficients $\vHat_{j}^{(k)}$. 
Therefore by setting $\vHat_i^{(k)}=1$ and $\vHat_j^{(k)}=0$ for $i\ne j$, it follows that 
\ba
   \Aw \, \phi_i = \beta(\lambda_i,\omega) \,\phi_i, \qquad i=0,1,2\ldots 
\ea
Thus $\Aw$ has eigenfunctions $\phi_j(\xv)$, $j=0,1,2,\ldots$, with corresponding eigenvalues $\beta_j \eqdef \beta(\lambda_j;\omega).$
As shown in Figure~\ref{fig:waveHoltzBeta} and discussed further below, $\beta(\lambda;\omega)\in[-\half,1]$.
Thus, all eigenvalues $\beta_j$ of $\Aw$ are real and lie in the interval~$[-\half,\,1]$. This completes the proof.
\qed
\end{proof}

\input tex/filterFunctionFig

The WaveHoltz filter function~\eqref{eq:filterFunction} 
can be written as the sum of sinc functions,
\ba
   \beta(\lambda; \omega)    &= \sinc\bigl((\omega-\lambda)T_f\bigr) + \sinc\bigl((\omega+\lambda)T_f\bigr) - \half \, \sinc\bigl(\lambda T_f\bigr),
   \label{eq:filterThreeSincs} 
\ea 
with one centered at $\lambda=0$ and the others centered at $\lambda=\pm\omega$.
As shown in Figure~\ref{fig:waveHoltzBeta}, $\beta$ has a maximum at $\lambda=\omega$ (note that $\beta$ is a function of $\lambda/\omega$ for a given integer $\Np$).
The widths of the peaks and valleys of this oscillatory function
can be decreased by increasing the number of periods, $\Np$.  
Note that the largest eigenvalues $\beta_j$ generally correspond to
values of $\lambda_j$ closest to the target frequency $\omega$.

A simple power iteration on $\Aw$ as given in Algorithm~\ref{alg:EigenWave} 
 can be used to compute an eigenpair $(\beta_j,\phi_j)$ 
for the $\beta_j$ with the largest magnitude (if that eigenvalue is isolated).
The convergence rate is then given by the ratio of the largest magnitude $\beta_j$ to the next largest in magnitude according to the standard analysis of the power method for computing eigenvalues.  
Since the eigenfunctions are the same, an application of a Rayleigh quotient involving the elliptic operator $\Lc$ can  be used to compute the corresponding value for~$\lambda_j$.  More sophisticated algorithms can be used to compute one or more eigenvalue-eigenfunction pairs
as discussed in Section~\ref{sec:Arnoldi}.

%% file: tex/filterFunctionFig.tex
{
\newcommand{\figw}{7cm}
\newcommand{\figh}{6cm}
\begin{figure}[htb]
\begin{center}
\begin{tikzpicture}
  \useasboundingbox (0,.65) rectangle (\figw,.95*\figh);  
  \begin{scope}[yshift=0*\figh]
    \figByWidth{0}{0}{fig/waveHoltzBetaFunction}{\figw}[0.][0.][0.][0.]
  \end{scope}  
\end{tikzpicture}
\end{center}
\caption{
  WaveHoltz filter function $\beta$ for $\Np=1$, $\Np=2$, and $\Np=3$ periods per time-interval.
   }
\label{fig:waveHoltzBeta}
\end{figure}
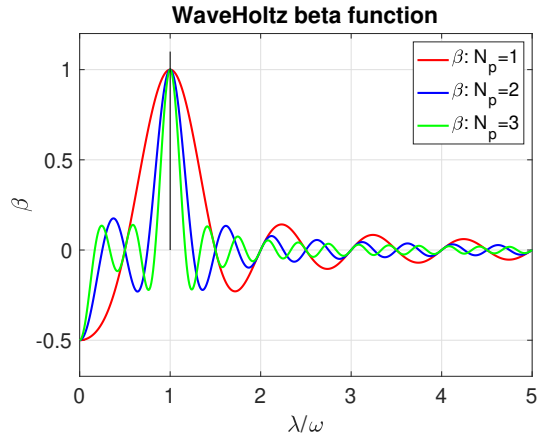
}

%% file: tex/discreteApproximations.tex
\section{Discrete approximations with explicit and implicit time-stepping}  \label{sec:discreteApproximations}

The EigenWave algorithm can be implemented with any number of numerical schemes such as those
based on finite differences, finite volumes, or finite elements.
The algorithm is developed in this article using finite difference methods on overset grids. 
An overset grid is a collection of curvilinear grids that cover a chosen problem domain $\Omega$ and overlap
where they meet, recall Figure~\ref{fig:rpiGrid} for example.  Thus, we may consider a discretization of a problem on a single curvilinear grid for ease of discussion, and note that the extension to a full overset grid is straightforward~\cite{CGNS,max2006b}.

Consider first the discrete approximation to the eigenvalue problem~\eqref{eq:eigBVP}.
We assume there is a smooth and invertible mapping, $\xv=\Gv(\rv)$,  from the unit square coordinates $\rv\in\Real^\nd$ in $\nd$-dimensions to the physical 
domain $\xv\in\Real^d$.  
The \textsl{mapping method} uses the chain rule to transform the governing equations from derivatives in $\xv$ to derivatives in $\rv$.
The transformed equations in the unit square coordinates can be discretized using standard 
centred or conservative finite difference approximations (see for example~\cite{max2006b}). 
Let $\xv_\iv$ denote the grid points where the subscript $\iv=[i_1,i_2,i_3]$ denotes a multi-index.
Let $\Phi_\iv \approx \phi(\xv_\iv)$ denote the grid function approximation to an eigenfunction $\phi$ at point $\xv=\xv_\iv$.
The discretized form of the eigenvalue BVP~\eqref{eq:eigBVP} is given by 
\bse
\label{eq:eigBVPdiscrete}
\bat
   & L_{ph} \Phi_\iv  = - \lambda_h^2 \, \Phi_\iv,  \qquad && \iv\in \Omega_h^a ,  \\
   & B_{ph} \Phi_\iv  = 0 ,                         \qquad && \iv\in \partial\Omega_h, \label{eq:eigBVPdiscreteBC}
\eat
\ese
where $L_{ph}$ denotes a $p^{\rm th}$-order accurate approximation to $\Lc$ 
and $B_{ph}$ denotes a a $p^{\rm th}$-order accurate approximation to the boundary conditions.
Here $\Omega_h^a$ denotes the set of \textsl{active} grid points where the interior equations are applied and $\partial\Omega_h$
denotes the grid points where the boundary conditions are applied. 
Note that in general we use additional compatibility boundary conditions
as numerical boundary conditions to treat the wide stencils associated with high-order finite difference approximations to $L_{ph}$, but these details are suppressed.  (A more detailed discussion of compatibility conditions is given in~\cite{lcbc2022} for example.)

Now consider discretizing the IBVP for the wave equation~\eqref{eq:waveHoltz} using the same spatial curvilinear grid.
Let $W_\iv^n \approx w(\xv_\iv,t^n)$ with $t^n = n\dt$, where $\dt$ is the time-step and the integer $n$ (and superscript $n$) denotes the time-level.
The explicit schemes employed in this article use a three-level approximation in time 
to discretize the wave equation IBVP~\eqref{eq:waveHoltz} as  
\bse
\label{eq:explicitScheme}
\bat
  & \Dpt\Dmt W_\iv^{n} = L_{ph} W_\iv^n , \quad&&\iv\in \Omega_h^a ,   \quad&& n=0,1,2,\ldots \,, \label{eq:explicitSchemeInterior} \\
  & B_{ph} W_\iv^{n}  = 0 ,                                            \quad&& \iv\in \partial\Omega_h, \quad&& n=1,2,\ldots \, ,   \label{eq:explicitSchemeBC} \\
  & W_\iv^0 = V_\iv^{(k)},  \quad&&\iv\in \Omega_h ,    \label{eq:explicitSchemeIC1} \\
  & \Dzt W_\iv^0 = 0  \quad&&\iv\in \Omega_h ,    \label{eq:explicitSchemeIC2}
\eat
\ese
where $\Omega_h$ is the set of all grid points,
and where $\Dpt$, $\Dmt$, and $\Dzt$ are the forward, backward, and centered divided difference operators in time, 
$\Dpt W_\iv^n\eqdef(W_\iv^{n+1}-W_\iv^{n})/\dt$, 
$\Dmt W_\iv^n\eqdef(W_\iv^{n}-W_\iv^{n-1})/\dt$,
$\Dzt W_\iv^n\eqdef(W_\iv^{n+1}-W_\iv^{n-1})/(2\dt)$.
Combining the initial conditions~\eqref{eq:explicitSchemeIC1} and~\eqref{eq:explicitSchemeIC2} with the interior scheme~\eqref{eq:explicitSchemeInterior} for $n=0$ 
gives the formula for the first time-step, 
\bat
  & W_\iv^1 = V_\iv^{(k)} + \half \dt^2 L_{ph} W_\iv^0, \quad&&\iv\in \Omega_h^a . \label{eq:explicitStep1}
\eat
The implicit schemes used here are also three-level schemes and they employ a trapezoidal-type approximation in time as \bse
\label{eq:implicitScheme}
\bat
  &  \Dpt\Dmt W_\iv^{n} = \half L_{ph} \bigl( W_\iv^{n+1} + W_\iv^{n-1}  \bigr) , \quad&&\iv\in \Omega_h^a ,   \quad&& n=0,1,2,\ldots\, , \label{eq:implicitSchemeInterior}  \\
  & B_{ph} W_\iv^{n}  = 0 ,    \quad&& \iv\in \partial\Omega_h, \quad&& n=1,2,\ldots \, ,            \label{eq:implicitSchemeBC} \\
  & W_\iv^0 = V_\iv^{(k)},           \quad&&\iv\in \Omega_h ,                                                       \label{eq:implicitSchemeIC1} \\
  & \Dzt W_\iv^0 = 0  ,         \quad&&\iv\in \Omega_h .                                                  \label{eq:implicitSchemeIC2}
\eat
The implicit scheme is unconditionally stable as discussed in~\cite{wimp2025}. 
The implicit matrix that needs to be inverted at each time step is (apart from boundary conditions) has the form 
\ba
   M_{ph} = I -{\dt^2\over2} L_{ph} , \label{eq:impMatrix}
\ea
where $I$ is the identity matrix.
Note that $M_{ph}$ with boundary conditions is a definite matrix (having eigenvalues $1+ \lambda_h^2 \dt^2 /2$ 
that is well suited to solution by fast iterative methods such as multigrid.
Combining the initial conditions~\eqref{eq:implicitSchemeIC1} and~\eqref{eq:implicitSchemeIC2} 
with the interior scheme~\eqref{eq:implicitSchemeInterior} for $n=0$ 
gives the formula for the first implicit time-step,
\ba
  &  M_{ph} W_\iv^1 = W_\iv^0, \qquad\iv\in \Omega_h^a . \label{eq:implicitStep1}
\ea
\ese
Note that the implicit system in~\eqref{eq:implicitStep1} involves the same coefficient matrix $M_{ph}$ used for time-stepping; 
this avoids the need to form and invert a different matrix for the first step.

Discrete solutions are computed to a time $T_f=\Np(2\pi/\omega)$ using $\NT$ time-steps with either the explicit or implicit time-stepping schemes.
The time filter in~\eqref{eq:waveHoltzFilterStep} is then applied using a trapezoidal rule quadrature
\bse
\ba
    V_\iv^{(k+1)} = \f{2}{T_f} \sum_{n=0}^{\NT} \sigma_{n} \left( \cos(\omega t^n) - \f{\alpha_d}{2} \right)  W_\iv^n, \label{eq:filterQuadrature}
\ea
where $\sigma_n$ are the weights in the quadrature and $\alpha_d$ is
\ba
  \alpha_d = \alpha_d(\omega\dt)  \eqdef \f{\tan(\omega\dt/2)}{\tan(\omega\dt)}.   \label{eq:alphad}
\ea
\ese
The value of $\alpha_d$ in~\eqref{eq:alphad} is chosen so that the discrete $\beta$ function using trapezoidal quadrature reaches a maximum of one at $\lambda=\omega$,
see~\cite{overHoltzArXiv2025,overHoltzPartOne,overHoltzPartTwo} for details.

Note that the explicit and implicit schemes and the trapezoidal quadrature rule 
are only second-order accurate in time, while the order of accuracy of the spatial operator may be second order or higher.
As shown in~\ref{sec:discreteAnalysis},
errors associated with the discretization in time do not affect the spatial accuracy of the discrete eigenvectors.
The time approximations can, however, affect the convergence of the EigenWave algorithm for large $\dt$ as discussed
in~\ref{sec:discreteAnalysis}.

\medskip
After discretization, the iteration in~\eqref{eq:Amap} takes the discrete form
\ba
    \Vv^{(k+1)} = \Aw_{ph} \Vv^{(k)} ,
    \label{eq:AmapDiscrete}
\ea
where $\Vv^{(k)}$ denotes the vector of unknowns $V_\iv^{(k)}$ on the grid, and $\Aw_{ph}$ denotes the matrix corresponding to the approximation of~$\Aw$.  We note that while $\Aw_{ph}$ exists, it is never formed explicitly in the EigenWave algorithm.

%% file: tex/usingExistingEigenvalueRoutines.tex
\section{Using EigenWave with existing high quality eigenvalue solvers} \label{sec:eigSoftware}

EigenWave can be combined with existing eigenvalue solvers 
that support a matrix free option. 
It has been found that the Krylov-Schur algorithm from SLEPc~\cite{SLEPc2005} and the 
IRAM (Implicitly Restarted Arnoldi Method) algorithm from ARPACK~\cite{ARPACK}, or the \texttt{eigs} function in Matlab give similarly
good results. 
These solvers are highly sophisticated and generally perform remarkable well, even for eigenvalues of high multiplicity. 
For more details on the IRAM and the Krylov-Schur algorithms, 
see the discussion in~\cite{BaiDemmel2000} and~\cite{Stewart2002}. 
See also~\ref{sec:Arnoldi} for further discussion of Arnoldi-based solvers.

\renewcommand{\algFontSize}{\small}
\begin{algorithm}
\algFontSize 
\caption{Matrix free Krylov-Schur or IRAM Eigenvalue Solver (e.g. from SLEPSc, ARPACK, Matlab).}
\begin{algorithmic}[1]
%
  \Function{[$N_c$,Y,E] = KyrlovSchur}{ \textsc{{\mvCol matVec}}, $N_{\Aw}$, $N_r$, $N_a$, whichEigs, tolerance }  
    \State Parameters:
    \State\quad Input: \textsc{{\mvCol matVec}} : matrix-vector multiply function, $\tilde\Vv = \Aw_{ph} \Vv$. 
    \State\quad Input: $N_{\Aw}$ : dimension of the matrix $\Aw_{ph} \in \Real^{N_{\Aw} \times N_{\Aw}}$ and vectors $\Vv,\tilde\Vv\in\Real^{N_{\Aw}}$.
    \State\quad Input: $N_r$ : request this many eigenpairs be computed.
    \State\quad Input: $N_a$ : size of Krylov subspace (total Arnoldi vectors kept), typically $N_a=2 N_r +1$.
    \State\quad Input: whichEigs : specify which eigenvalues to find, e.g. `largest', `smallest', `largestAbs'.
    \State\quad Input: tolerance : convergence tolerance for the eigenpairs.
    \State\quad Output: $N_c$ : number of converged eigenpairs (can be larger or smaller than $N_r$).
    \State\quad Output: Y(1:$N_{\Aw}$,1:$N_c$) : eigenvectors as columns of a matrix.
    \State\quad Output : E(1:$N_c$) : vector of eigenvalues. 
 \EndFunction
\end{algorithmic} 
\label{alg:KrylovSchurAlgorithm}
\end{algorithm}

In order to help the reader understand the results discussed in Section~\ref{sec:numericalResults} we present,
in Algorithm~\ref{alg:KrylovSchurAlgorithm}, the typical input and output for a Krylov-Schur or IRAM eigenvalue solver (arbitrarily called \textsc{KyrlovSchur} to be concrete).
A matrix free solver requires a black-box function, here called \textsc{{\mvCol matVec}},
to compute the product of the matrix times a vector. For EigenWave this function performs a wave-solve for a generic initial condition $\Vv$ and returns a corresponding next iterate $\tilde\Vv$ such that $\tilde\Vv = \Aw_{ph} \Vv$. 
An example description of \textsc{{\mvCol matVec}} is given in Algorithm~\ref{alg:matVec} which shows how to exclude constraint equations such as boundary conditions.
The input also includes the requested number, $N_r$, of eigenvalues, the number of vectors, $N_a$, to keep in the Krylov subspace, 
and a description of which eigenvalues to find, in our case we choose the largest in absolute value.
The output consists of the number, $N_c$, of converged eigenpairs\footnote{$N_c$ is sometimes less than, and sometimes more than $N_r$, depending the convergence behaviour of the KrylovSchur algorithm and the user convergence tolerances.}
as well as the eigenvalues and eigenvectors (these may be requested individually instead of all being returned).

For a large number of requested eigenvalues, the storage requirements of the EigenWave algorithm is
normally dominated by the storage requirements of the IRAM or Krylov-Schur algorithms which involves at least $(2 N_r+1) N$ floating point numbers, where $N$ is the total number of grid points.
Exceptions to this rule would be if a direct sparse solver is used to solve the implicit time-stepping equations; 
these require significant storage for fill-in, especially in three dimensions.
On the other hand, the matrix-free multigrid solver we use is quite memory efficient especially 
when most grid points belong to Cartesian grids (see Section~\ref{sec:implicitTimeSteppingWithMultigrid} for more details on the multigrid solver we use).

%% file: tex/numericalResults.tex
\section{Numerical results}  \label{sec:numericalResults}

Numerical results are presented showing the basic behavior of the EigenWave algorithm for 
eigenvalue problems 
in various two and three-dimensional domains with Dirichlet boundary conditions. 
Results using Neumann boundary are similar.
The schemes use second and fourth-order accurate spatial discretizations on overset grids.
In each case, the eigenvalues and eigenvectors obtained using the EigenWave algorithm are compared to the true discrete eigen-pairs, which are computed directly from the discrete problem \eqref{eq:eigBVPdiscrete} (see \ref{sec:accuracyOfEigenpairs} for details on how they are computed).
The results demonstrate that the EigenWave algorithm can compute eigen-pairs to near full machine precision using just a few
wave-solves per eigen-pair.  Moreover, with implicit time-stepping, 
only $10$ time-steps per wave-solve are needed.
This is true for simple problem domains, such as a square or annulus, as well as for more complex domains that use overset grids.
Unless otherwise stated, the implicit scheme uses a direct sparse solver since this is often the fastest approach for smaller problem sizes. 
Further properties of EigenWave are covered in~\ref{sec:properties}.

\medskip
The accuracy of the eigen-pairs computed using the EigenWave algorithm is measured in three ways: the relative
error in the eigenvalue, the relative error in the eigenvector, and the relative residual, each defined, respectively, by 
\ba
  &  \texttt{eig-err} =  \f{|\lambda_{h,j} - \lambda_{h,j}^{\rm true}|}{\lambda_{h,j}^{\rm true}} , 
  \quad
   \texttt{evect-err} =  \f{\| V_{\iv,j} - V_{\iv,j}^{\rm true} \|_\infty }{\| V_{\iv,j}^{\rm true} \|_\infty}, 
   \quad
   \texttt{eig-res} =  \f{ \, \| L_{ph} V_{\iv,j} + \lambda_{h,j}^2 V_{\iv,j} \|_\infty }{ \lambda_{h,j}^2 },
\ea
where $(\lambda_{h,j},V_{\iv,j})$ denotes the $j^{\rm th}$ discrete eigen-pair computed using EigenWave, while corresponding eigen-pairs with the ``${\rm true}$'' superscript denote the true discrete values computed separately to a very small error tolerance. 
For eigenvalues with multiplicity greater than one, the eigenvectors are not unique, although the eigenspace is.
  For multiple eigenvalues the error in the eigenvector is computed as the distance to the closest vector in the true discrete eigenspace as follows.
Given an approximate eigenvector $\vv$ and an eigenspace $\Ec ={\rm span}\{ \wv_1, \ldots,\wv_m \}$, 
we find the orthogonal projection of $\vv $ onto $\Ec$, denoted as $\tilde{\vv}$, by solving a least squares problem.
The max-norm distance is then computed as $\| \vv - \tilde{\vv} \|_{\infty}$. 



\input tex/squareEigenpairs

\input tex/diskEigenpairs

\input tex/circleInAChannel


\input tex/threeShapesEigenpairs

\input tex/rpiEigenpairs

\input tex/doubleEllipseEigenpairs

\input tex/boxEigenpairs

\input tex/pipeEigenpairs

\input tex/sphereEigenpairs

\input tex/doubleEllipsoidEigenpairs.tex

%% file: tex/squareEigenpairs.tex
\subsection{Eigenmodes of a square} \label{sec:squareEigenpairs}

\input tex/squareEigsFig.tex

Here, eigenpairs on the unit square $\Omega=[0,1]\times[0,1]$ are computed on a Cartesian grid with $128$ cells in each direction (subsequently called \texttt{square128}). The target frequency is chosen as $\omega=12$. Table~\ref{tab:square128Summary} summarizes the results using the second and fourth-order accurate discretizations. In both cases, twenty-four eigenpairs are requested ($N_r=24$ in Algorithm~\ref{alg:KrylovSchurAlgorithm}), and 
twenty-seven converged eigenpairs  are found ($N_c=27$ in Algorithm~\ref{alg:KrylovSchurAlgorithm}) by the KrylovSchur
algorithm\footnote{Recall that the Krylov Schur algorithm actually just computes the eigenvectors while the eigenvalues are computed subsequently by a Rayleigh quotient.}
in a total of $89$ wave-solves. This corresponds to approximately 
$3.3$ wave-solves per eigenpair found.
Since there are ten implicit time-steps per wave-solve, i.e.~$\Nits=10$, this implies 
about $33$ implicit solves per eigenpair found.
Figure~\ref{fig:squareEigs} shows the filter-functions $\beta(\lambda;\omega)$ for both orders of accuracy with the computed
eigenvalues marked. 
Note that here, and in subsequent graphs like those in Figure~\ref{fig:squareEigs}, the curves and marked eigenvalues are adjusted 
to account for time-discretization errors 
as discussed in~\ref{sec:discreteAnalysis}\footnote{Note that the adjusted eigenvalues appearing in the graphs are only used for convergence theory, 
the eigenvalues computed by EigenWave have no errors due to time discretizations as evidenced, for example, in Table~\ref{tab:square128Summary}.}.
Table~\ref{tab:square128Order2} in~\ref{sec:eigenpairTables} gives further details of the 27 eigenvalues computed for this example, including their multiplicity and accuracy. Importantly, we note that this problem has many multiple eigenvalues, and EigenWave coupled with KrylovSchur still performs very well. Also note that these results indicate that the behavior of the algorithm 
at fourth-order accuracy is almost identical with that for second-order accuracy. This is not unexpected
since the convergence of EigenWave should depend primarily on the distribution of the eigenvalues and not on
the order of accuracy or grid spacing.

\bogus{\begin{table}[hbt]
\begin{center}\tableFontSize
\begin{tabular}{|c|c|c|c|c|c|c|c|} \hline
 \multicolumn{8}{|c|}{EigenWave: grid=square128, ts=implicit, order=2, $\omega=  12$, $N_p=1$, KrylovSchur } \\ \hline 
   num   &  wave       & time-steps  & wave-solves &  time-steps &   max      &  max      &  max       \\ 
   eigs  &  solves     & per period  &  per eig    &  per-eig    &   eig-err  & evect-err & eig-res    \\ 
 \hline
    27    &       89     &      10     &    3.3    &     32     &  7.99e-15 &  4.89e-13 &  2.60e-12 \\
 \hline 
\end{tabular}
\end{center}
\vspace*{-1\baselineskip}
\caption{Summary of EigenWave performance for \texttt{square128} grid using the KrylovSchur algorithm and implicit time-stepping.  The spatial order of accuracy is 2 and the wave-solves use $N_p=1$ to determine the final time.}
\label{tab:square128Order2Summary}
\end{table}



\begin{table}[hbt]
\begin{center}\tableFontSize
\begin{tabular}{|c|c|c|c|c|c|c|c|} \hline
 \multicolumn{8}{|c|}{EigenWave: grid=square128, ts=implicit, order=4,  $\omega=  12$, $N_p=1$, KrylovSchur } \\ \hline 
   num   &  wave       & time-steps  & wave-solves &  time-steps &   max      &  max      &  max       \\ 
   eigs  &  solves     & per period  &  per eig    &  per-eig    &   eig-err  & evect-err & eig-res    \\ 
 \hline
    27    &       89     &      10     &    3.3    &     32     &  7.38e-14 &  1.46e-12 &  7.42e-12 \\
 \hline 
\end{tabular}
\end{center}
\vspace*{-1\baselineskip}
\caption{Summary of EigenWave performance for \texttt{square128} grid using the KrylovSchur algorithm and implicit time-stepping.  The spatial order of accuracy is 4 and the wave-solves use $N_p=1$ to determine the final time.}
\label{tab:square128Order4Summary}
\end{table}}

\begin{table}[hbt]
\begin{center}\tableFontSize
\begin{tabular}{|c|c|c|c|c|c|c|c|c|} \hline
 \multicolumn{9}{|c|}{EigenWave: grid=square128, ts=implicit, $\omega=  12$, $N_p=1$, KrylovSchur } \\ \hline 
   order & num   &  wave       & time-steps  & wave-solves &  time-steps &   max      &  max      &  max       \\ 
   & eigs  &  solves     & per period  &  per eig    &  per-eig    &   eig-err  & evect-err & eig-res    \\ 
 \hline
    2 & 27    &       89     &      10     &    3.3    &     32     &  7.99e-15 &  4.89e-13 &  2.60e-12 \\
    4 & 27    &       89     &      10     &    3.3    &     32     &  7.38e-14 &  1.46e-12 &  7.42e-12 \\
 \hline 
\end{tabular}
\end{center}
\vspace*{-1\baselineskip}
\caption{Summary of EigenWave performance for \texttt{square128} grid using the KrylovSchur algorithm and implicit time-stepping.  The spatial order of accuracy is 2 for the top row and 4 for the bottom row, and the wave-solves use $N_p=1$ to determine the final time.}
\label{tab:square128Summary}
\end{table}

%% file: tex/squareEigsFig.tex
{
\newcommand{\figw}{6.5cm}
\newcommand{\figh}{5.5cm}

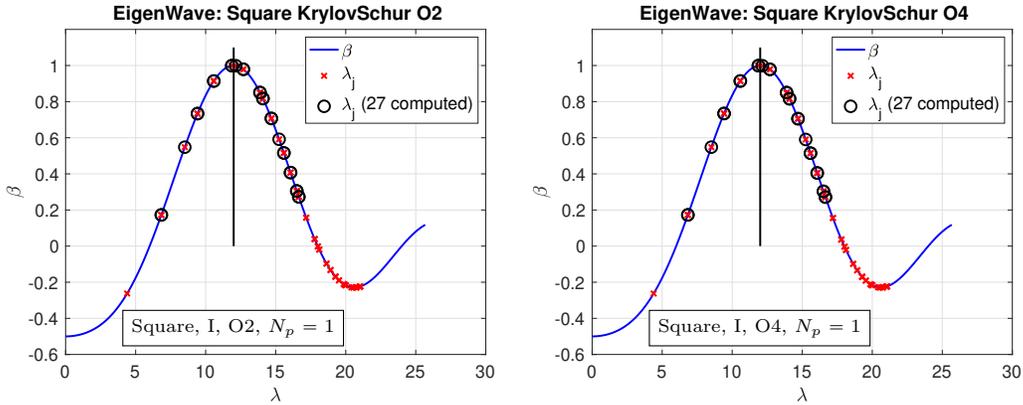
\begin{figure}[htb]
\begin{center}
\begin{tikzpicture}[scale=1]
  \useasboundingbox (0,.7) rectangle (14,5.25);  

  \begin{scope}[xshift=0cm]
    \figByWidth{0.0}{0}{fig/eigSquareO2ImpEig12Krylov}{\figw}[0][0][0][0]
    \draw (1.5,1) node[draw,fill=white,anchor=west,xshift=0pt,yshift=0pt] {\scriptsize Square, I, O2, $\Np=1$};
  \end{scope}
  \begin{scope}[xshift=7cm]
    \figByWidth{0.0}{0}{fig/eigSquareO4ImpEig12Krylov}{\figw}[0][0][0][0]
    \draw (1.5,1) node[draw,fill=white,anchor=west,xshift=0pt,yshift=0pt] {\scriptsize Square, I, O4, $\Np=1$};
  \end{scope}  

\end{tikzpicture}
\end{center}
\caption{Eigenpairs of a square with target frequency $\omega=12$. The discrete time corrected  filter function $\beta=\beta(\lambda;\omega)$ is plotted in blue, the true discrete eigenvalues $\lambda_j$ are marked with red x's, and the computed eigenvalues are marked with black circles.
At left are second-order accurate results, and at right are fourth-order accurate results.
}
\label{fig:squareEigs}
\end{figure}
}

%% file: tex/diskEigenpairs.tex
\renewcommand{\Gcd}{\Gc_{\rm disk}}
\subsection{Eigenmodes of a circular disk} \label{sec:diskEigenpairs}

\input tex/diskGridAndFiltersFig

\input tex/diskContoursFig.tex

In this section, eigenpairs of the Laplacian operator for a circular disk in two dimensions with Dirichlet boundary conditions are computed using the EigenWave algorithm.
The exact continuous eigenvalues and eigenfunctions for this problem are given by
\bse
\label{eq:diskTrueEigenpairs}
\ba
   \lambda_{0,\mr}={q_{0,\mr}\over a},\qquad \phi(r,\theta)=J_0\left(\lambda_{0,\mr} r\right),\qquad \mr=1,2,\ldots
\ea
and
\ba
  \lambda_{\mTheta,\mr} ={q_{\mTheta,\mr}\over a},\qquad \phi(r,\theta)
     =\begin{cases}
     J_\mTheta\left(\lambda_{\mTheta,\mr} r\right)\cos(\mTheta \theta) \\
     J_\mTheta\left(\lambda_{\mTheta,\mr} r\right)\sin(\mTheta \theta)
  \end{cases} , 
  \quad \mTheta,\mr=1,2,\ldots
\ea
\ese
where $(r,\theta)$ are the usual polar coordinates, $a$ is the radius of the disk (here taken as $a=1$) and $q_{\mTheta,\mr}$ is the $\mr\sp{{\rm th}}$ zero of $J_\mTheta$, the first kind Bessel function of integer order~$\mTheta$.
Note that there 
are many double eigenvalues due to the rotational symmetry of the geometry, and so this problem is potentially challenging.  This problem also serves to test the EigenWave algorithm for a domain using an overset grid and for which the exact eigenvalues and eigenfunctions are known.

Solutions are computed using an overset grid for the disk of radius one, denoted by $\Gcd^{(j)}$, consisting of a Cartesian grid covering the central portion of the domain and
an annular boundary-fitted grid as shown in Figure~\ref{fig:diskFig} (left). 
This grid is constructed to have typical grid spacing $\ds^{(j)}=1/(10 j)$, where $j$ is a positive integer specifying the grid resolution.
Table~\ref{tab:sice4Summary} summarizes results of computing eigenpairs using second and fourth-order accurate spatial discretizations on grid $\Gcd^{(4)}$. The target frequency is $\omega=10$, and implicit time-stepping is used with $\Nits=10$ time-steps per period.
Forty-two eigenpairs are requested, and EigenWave using the Krylov-Schur algorithm returns~$43$ and~$44$ eigenpairs for second and fourth-order accurate discretizations, respectively. In both cases the algorithm used a total of $169$ wave-solves, corresponding to approximately 
$3.9$ wave-solves per computed eigenpair at second order, and $3.8$ wave-solves per eigenpair  at fourth order. The middle and right panels in Figure~\ref{fig:diskFig} shows the filter function $\beta(\lambda;\omega)$ for orders two and four, respectively, both with the computed
eigenvalues marked. Overall the convergence behavior of the EigenWave algorithm is seen to be almost identical between second and fourth-order accuracy. Furthermore, the Eigenwave algorithm is seen to provide very good approximations to the discrete eigenpairs, even for this difficult problem with many duplicate eigenvalues.
Figure~\ref{fig:diskEigenfunctionSolution} shows contours of selected eigenvectors (e.g. $\phi^{(46)}$ denotes the $46$-th eigenvector when \textsl{all} the eigenvalues are ordered from smallest to largest, including those not computed by EigenWave).

\input tex/diskShortTables

Table~\ref{tab:sice4Order4} in~\ref{sec:eigenpairTables} gives more details on the accuracy of the computed eigenvalues and eigenvectors.
There are many eigenpairs with multiplicity equal to 2 as indicated in the table.  The multiplicity is estimated numerically by
checking the distance between nearby eigenvalues against a chosen tolerance.  We note that the largest errors in the eigenvectors tend to occur
for eigenpairs associated with certain exact eigenvalues with multiplicity 2, but where the difference between the computed eigenvalue pairs,
both for $\lambda_{h,j}$ and $\lambda_{h,k}^{\rm true}$, is large enough that the multiplicity
estimator reports the multiplicity as only one (e.g.~$\lambda_{h,j}$ with $j=19$ and~$j=20$ in the table).  Since
the error in $\lambda_{h,k}^{\rm true}$ with respect to the exact eigenvalues for these `near-multiple' eigenvalues is larger,
it is not surprising that the error in $\lambda_{h,j}$ with respect to $\lambda_{h,k}^{\rm true}$ is larger as well.  

%% file: tex/diskGridAndFiltersFig.tex
{
\newcommand{\figw}{5.5cm}
\newcommand{\figh}{5.5cm}
\begin{figure}[htb]
\begin{center}
\begin{tikzpicture}
  \useasboundingbox (0,.5) rectangle (16,4.5);  

  \figByWidth{0.0}{0}{fig/sicGridG2}{4.5cm}[0.1][0.1][0.1][0.1]

  \begin{scope}[xshift=4.75cm]
     \figByWidth{0}{0}{fig/diskG4O2Freq10Krylov}{\figw}[0][0][0][0]
     \draw (1.5,1) node[draw,fill=white,anchor=west,xshift=0pt,yshift=0pt] {\scriptsize Disk, I, O2, $\Np=1$};
  \end{scope} 
  \begin{scope}[xshift=10.5cm]
     \figByWidth{0}{0}{fig/diskG4O4Freq10Krylov}{\figw}[0][0][0][0]
     \draw (1.5,1) node[draw,fill=white,anchor=west,xshift=0pt,yshift=0pt] {\scriptsize Disk, I, O4, $\Np=1$};
  \end{scope}    
  
\end{tikzpicture}
\end{center}
\caption{At left is overset grid $\Gcd^{(2)}$ for a disk. The middle and right show graphs of the filter function $\beta$ with the  computed eigenvalues marked with black circles for $2nd$ and $4th$-order accurate discretizations respectively, both using grid $\Gcd^{(4)}$. The target frequency $\omega=10$ is marked as a vertical black line.
 }
\label{fig:diskFig}
\end{figure}
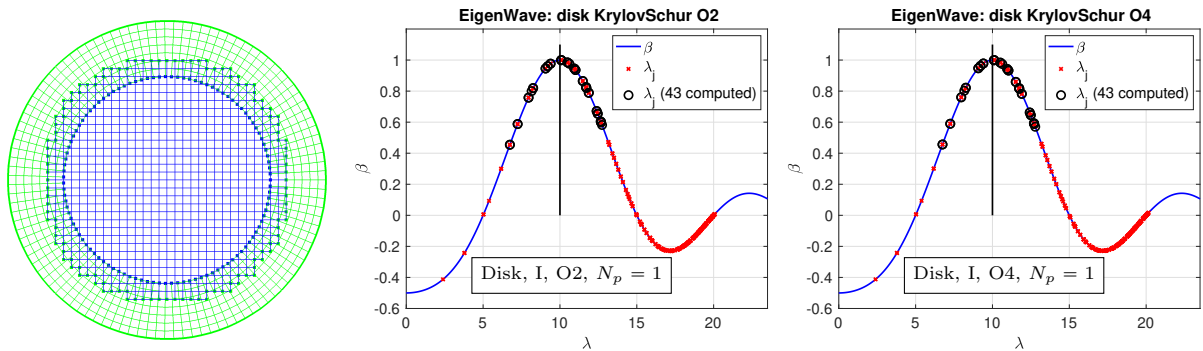
}

%% file: tex/diskContoursFig.tex
{
\newcommand{\drawContour}[7]{%
\begin{scope}[#1]
\draw(0.0,0) node[anchor=south west,xshift=-4pt,yshift=+0pt] {\trimfiga{fig/#2}{\figWidtha}};
  \draw(.5,.5) node[draw,fill=white,anchor=west,xshift=2pt,yshift=1pt] {\scriptsize #3};
\begin{scope}[xshift=-.3cm,yshift=-2pt]
  \draw (\xcb,\ycb) node[anchor=south west,xshift=0.25cm,yshift=.5cm,rotate=-90] {\trimfigcb{fig/colourBarLines}{\cbWidth}{\cbHeight}};
  \draw (.8,0) node[anchor=north,xshift=+3pt,yshift=+2pt] {\scriptsize $#6$};
  \draw (4.8,0) node[anchor=north,xshift=+0pt,yshift=+2pt] {\scriptsize $#7$};
\end{scope}
\end{scope}
}
\newcommand{\cbWidth}{.2cm}
\newcommand{\cbHeight}{4cm}
\newcommand{\xcb}{.5cm}
\newcommand{\ycb}{-.2cm}
\setlength{\ycbTop}{\ycb+\cbHeight}
\setlength{\ycbMid}{\ycb+\cbHeight*\real{.5}}
\newcommand{\trimfigcb}[3]{\includegraphics[width=#2, height=#3, clip, trim=17cm 2.35cm 1.65cm 2.35cm]{#1}}
\newcommand{\figWidtha}{4.25cm}
\newcommand{\trimfiga}[2]{\trimw{#1}{#2}{.11}{.117
}{.11}{.11}}
\begin{figure}[htb]
\begin{center}
\begin{tikzpicture}
   \useasboundingbox (0,.3) rectangle (15,4.5);  

   \begin{scope}[yshift=0cm]
     \drawContour{xshift=0.cm,yshift=0.00cm}{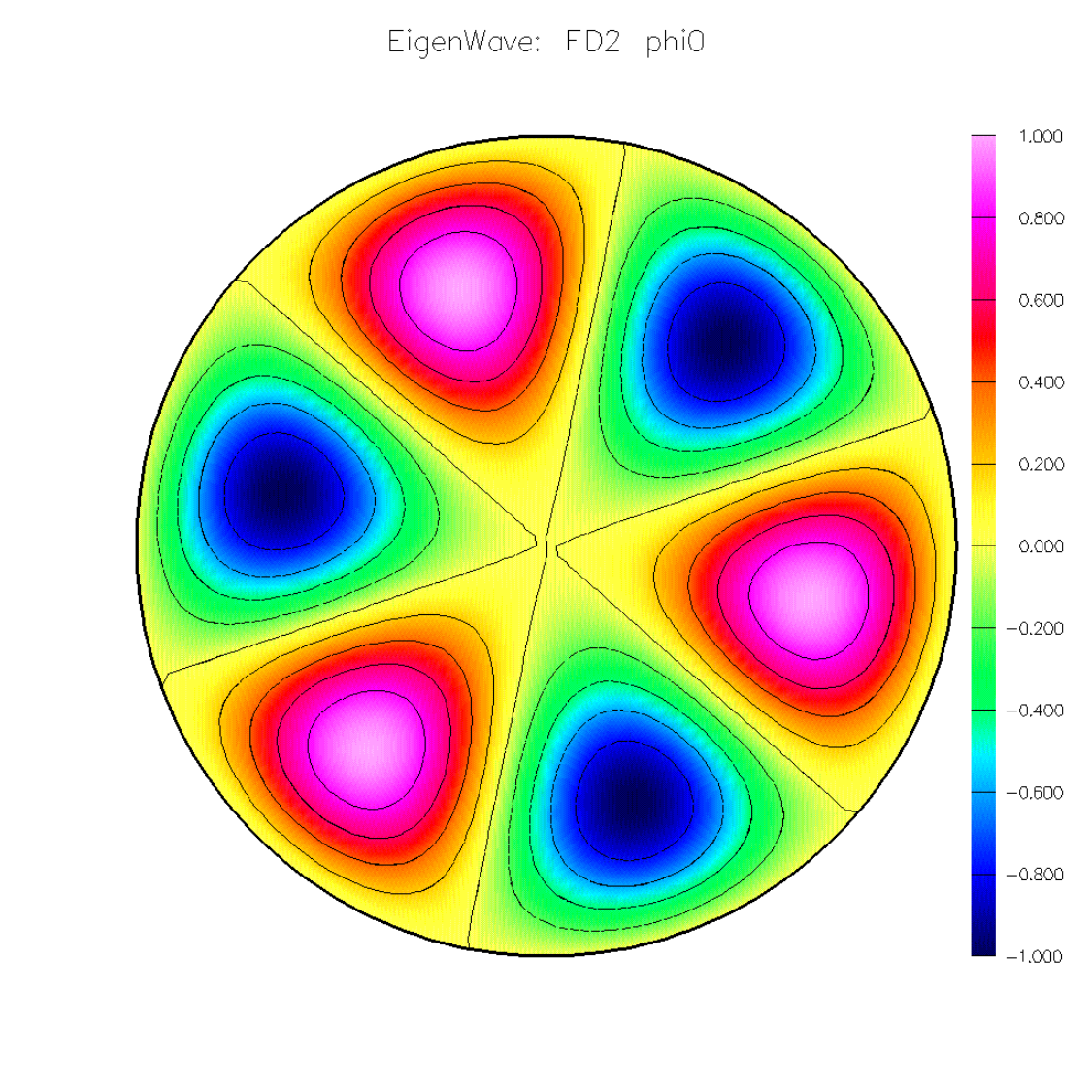}{$\phi^{(6)}$}{$v$}{$t=0.3$}{$-1.0$}{$1.0$}     
     \drawContour{xshift=5cm,yshift=0.00cm}{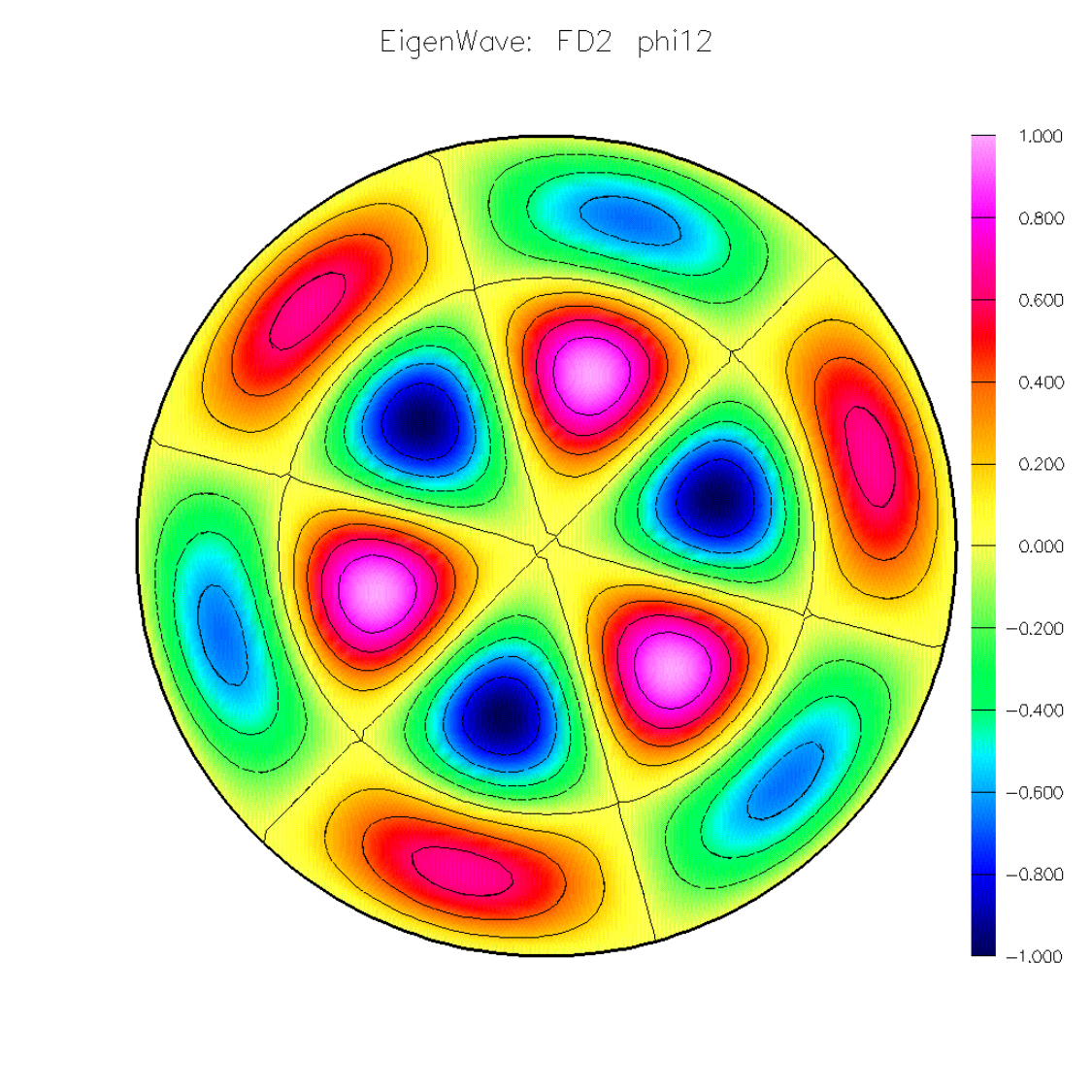}{$\phi^{(18)}$}{$v$}{$t=0.3$}{$-1.0$}{$1.0$}    
     \drawContour{xshift=10cm,yshift=0.00cm}{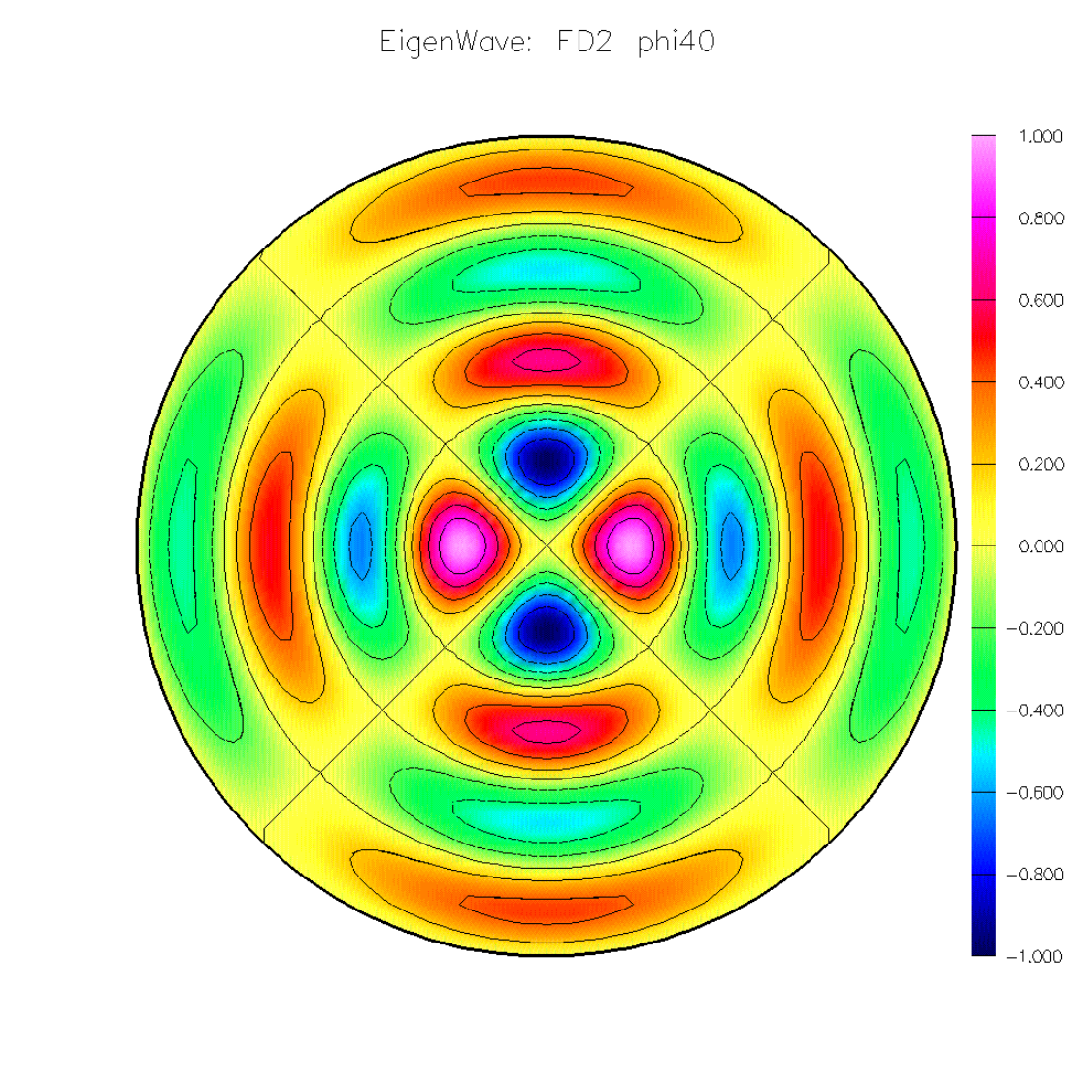}{$\phi^{(46)}$}{$v$}{$t=0.3$}{$-1.0$}{$1.0$}    
   \end{scope}

\end{tikzpicture}
\end{center}
\caption{Some computed eigenfunctions of a disk where $\phi^{(m)}$ denotes the eigenvector corresponding to the m-th eigenvalue, sorted from smallest to largest.
    }
\label{fig:diskEigenfunctionSolution}
\end{figure}
}

%% file: tex/diskShortTables.tex
\begin{table}[hbt]
\begin{center}\tableFontSize
\begin{tabular}{|c|c|c|c|c|c|c|c|c|} \hline
 \multicolumn{9}{|c|}{EigenWave: disk, ts=implicit, $\omega=  10.0$, $N_p=1$, KrylovSchur } \\ \hline 
   order & num   &  wave       & time-steps  & wave-solves &  time-steps &   max      &  max      &  max       \\ 
   & eigs  &  solves     & per period  &  per eig    &  per-eig    &   eig-err  & evect-err & eig-res    \\ 
 \hline
    2 & 43    &       169     &      10     &    3.9    &     39     &  5.83e-15 &  1.10e-11 &  1.07e-12 \\
    4 & 44    &       169     &      10     &    3.8    &     38     &  6.78e-15 &  1.42e-11 &  9.13e-13 \\
 \hline 
\end{tabular}
\end{center}
\vspace*{-1\baselineskip}
\caption{Summary of EigenWave performance for disk grid $\Gcd^{(4)}$ using the KrylovSchur algorithm and implicit time-stepping.  The top row uses second-order accurate discretization while bottom row uses fourth-order accurate discretization. In all cases the wave-solves use $N_p=1$ to determine the final time.}
\label{tab:sice4Summary}
\end{table}

%% file: tex/circleInAChannel.tex
\newcommand{\Gcic}{\Gc_{\rm cic}}
\subsection{Eigenmodes for a circle-in-a-channel domain} \label{sec:circleInAChannel}

In this section, eigenpairs of a square domain with a circular hole are computed using the \EigenWaveb algorithm.
The overset grid for the square domain $[-2,2]^2$ with circular hole of
radius one-half, consists of a background Cartesian grid and
an annular grid as shown in Figure~\ref{fig:cicFig} (left). 
The grid with target grid spacing $\ds^{(j)}=1/(10 j)$ is denoted by $\Gcic^{(j)}$.

{
\newcommand{\figw}{6.5cm}
\newcommand{\figh}{5.5cm}
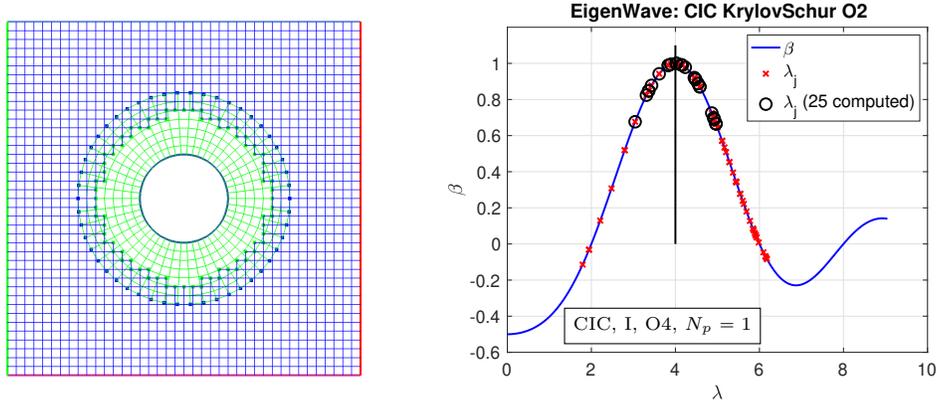
\begin{figure}[htb]
\begin{center}
\begin{tikzpicture}
  \useasboundingbox (0,.5) rectangle (12.5,5.25);  
  \figByWidth{0.0}{0}{fig/cicGrid}{5cm}[0.1][0.1][0.1][0.1]
   \begin{scope}[xshift=6cm,yshift=-5pt]
     \figByWidth{0}{0}{fig/cicG4O2ImpEig4Krylov}{\figw}[0][0][0][0]
     \draw (1.5,1) node[draw,fill=white,anchor=west,xshift=0pt,yshift=0pt] {\scriptsize CIC, I, O4, $\Np=1$};
  \end{scope}  
  
\end{tikzpicture}
\end{center}
\caption{Left: overset grid $\Gcic^{(2)}$ for a circle in a channel.
Right: graph of the filter function $\beta$ with the  
computed eigenvalues marked with black circles for a fourth-order accurate computation on grid $\Gcic^{(4)}$.
The target frequency $\omega=4$ is marked as a vertical black line.
 }
\label{fig:cicFig}
\end{figure}
}

{

\input tables/cicG4O2ImpEig4Krylov.tex
  \eigenWaveSummaryTable  
  \eigenWaveLongTable
}
Table~\ref{tab:cice4Order2Summary} summarizes results of computing eigenpairs to second-order accuracy on grid $\Gcd^{(4)}$.
The target frequency was $\omega=4$. Twenty-four eigenpairs were requested and twenty-five eigenpairs were accurately found by the KrylovSchur
algorithm in a total of $136$ WaveHoltz solves. This corresponds to approximately 
$5.4$ wave-solves per eigenpair found.
The right graph in Figure~\ref{fig:cicFig} shows the filter-function $\beta(\lambda;\omega)$ with the computed
eigenvalues marked. 
Table~\ref{tab:cice4Order2} shows more detailed results.
Contours of selected eigenvectors are shown in Figure~\ref{fig:cicEigenVectors}.


\input tex/cicContoursFig.tex

%% file: tables/cicG4O2ImpEig4Krylov.tex
\newcommand{\eigenWaveLongTable}{
\begin{table}[hbt]
\begin{center}\tableFontSize
\begin{tabular}{|c|r|r|c|c|c|c|c|} \hline
 \multicolumn{8}{|c|}{EigenWave: grid=cice4, ts=implicit, order=2, $N_p=1$, KrylovSchur } \\ \hline 
$j$  & \multicolumn{1}{c|}{$\lambda_{h,j}$} &  \multicolumn{1}{c|}{$\lambda_{h,k}^{\rm true}$}  & $k$ &  mult &  eig-err & evect-err & eig-res  \\ \hline   
   0   &   3.196154 &   3.196154 &     7 &    1  &  3.89e-15 &  1.48e-13 &  9.96e-13 \\ 
   1   &   3.514270 &   3.514270 &     8 &    1  &  7.08e-15 &  2.23e-13 &  1.05e-12 \\ 
   2   &   3.576084 &   3.576084 &     9 &    2  &  1.99e-15 &  1.89e-13 &  8.17e-13 \\ 
   3   &   3.576084 &   3.576084 &    10 &    2  &  1.24e-15 &  2.15e-13 &  7.72e-13 \\ 
   4   &   3.660054 &   3.660054 &    11 &    1  &  2.79e-15 &  3.27e-13 &  9.24e-13 \\ 
   5   &   3.882355 &   3.882355 &    12 &    1  &  3.43e-15 &  1.36e-13 &  6.93e-13 \\ 
   6   &   4.157318 &   4.157318 &    13 &    2  &  2.14e-15 &  5.91e-13 &  5.46e-13 \\ 
   7   &   4.157319 &   4.157319 &    14 &    2  &  1.92e-15 &  8.79e-13 &  7.44e-13 \\ 
   8   &   4.232627 &   4.232627 &    15 &    1  &  2.10e-15 &  6.80e-13 &  6.17e-13 \\ 
   9   &   4.405174 &   4.405174 &    16 &    2  &  3.43e-15 &  1.19e-12 &  4.75e-11 \\ 
  10   &   4.405176 &   4.405176 &    17 &    2  &  1.81e-15 &  3.84e-13 &  1.15e-11 \\ 
  11   &   4.554997 &   4.554997 &    18 &    1  &  1.56e-15 &  1.03e-12 &  6.82e-13 \\ 
  12   &   4.591306 &   4.591306 &    19 &    1  &  5.03e-15 &  8.14e-13 &  7.19e-13 \\ 
  13   &   4.686505 &   4.686505 &    20 &    1  &  3.60e-15 &  2.20e-13 &  6.89e-13 \\ 
  14   &   4.983295 &   4.983295 &    21 &    1  &  4.81e-15 &  5.62e-13 &  4.72e-13 \\ 
  15   &   5.020320 &   5.020320 &    22 &    2  &  2.65e-15 &  4.87e-13 &  4.27e-13 \\ 
  16   &   5.020323 &   5.020323 &    23 &    2  &  4.07e-15 &  1.42e-13 &  4.17e-13 \\ 
  17   &   5.110229 &   5.110229 &    24 &    1  &  2.43e-15 &  2.90e-13 &  5.05e-13 \\ 
  18   &   5.169540 &   5.169540 &    25 &    2  &  1.72e-15 &  2.86e-13 &  4.09e-13 \\ 
  19   &   5.169541 &   5.169541 &    26 &    2  &  8.59e-16 &  2.80e-13 &  2.81e-13 \\ 
  20   &   5.589282 &   5.589282 &    27 &    1  &  4.77e-16 &  1.76e-13 &  3.60e-13 \\ 
  21   &   5.656494 &   5.656494 &    28 &    1  &  7.85e-16 &  1.57e-13 &  3.38e-13 \\ 
  22   &   5.680284 &   5.680284 &    29 &    1  &  3.60e-15 &  1.72e-13 &  4.58e-13 \\ 
  23   &   5.740765 &   5.740765 &    30 &    2  &  2.32e-15 &  7.22e-14 &  4.54e-13 \\ 
  24   &   5.740766 &   5.740766 &    31 &    2  &  1.08e-15 &  1.59e-13 &  3.02e-13 \\ 
 \hline 
\end{tabular}
\end{center}
\caption{grid=cice4, method=KrylovSchur, ts=implicit, order=2, $N_p=1$.}
\label{tab:cice4Order2}
\end{table}
}

\newcommand{\eigenWaveSummaryTable}{
\begin{table}[hbt]
\begin{center}\tableFontSize
\begin{tabular}{|c|c|c|c|c|c|c|c|} \hline
 \multicolumn{8}{|c|}{EigenWave: grid=cice4, ts=implicit, order=2, $N_p=1$, KrylovSchur } \\ \hline 
   num   &  wave       & time-steps  & wave-solves &  time-steps &   max      &  max      &  max       \\ 
   eigs  &  solves     & per period  &  per eig    &  per-eig    &   eig-err  & evect-err & eig-res    \\ 
 \hline
    25    &       136     &      10     &    5.4    &     54     &  7.08e-15 &  1.19e-12 &  4.75e-11 \\
 \hline 
\end{tabular}
\end{center}
\caption{EigenWave: grid=cice4, method=KrylovSchur, ts=implicit, order=2, $N_p=1$.
 The index $k$ denotes the closest true discrete eigenvalue $\lambda_{h,k}^{\rm true}$ to the EigenWave value $\lambda_{h,j}$.
 }
\label{tab:cice4Order2Summary}
\end{table}
}

%% file: tex/cicContoursFig.tex
{
\newcommand{\drawContour}[7]{%
\begin{scope}[#1]
\draw(0.0,0) node[anchor=south west,xshift=-4pt,yshift=+0pt] {\trimfiga{fig/#2}{\figWidtha}};
  \draw(.5,.5) node[draw,fill=white,anchor=west,xshift=2pt,yshift=1pt] {\scriptsize #3};
\begin{scope}[xshift=-.2cm,yshift=-2pt]
  \draw (\xcb,\ycb) node[anchor=south west,xshift=0.25cm,yshift=.5cm,rotate=-90] {\trimfigcb{fig/colourBarLines}{\cbWidth}{\cbHeight}};
  \draw (.8,0) node[anchor=north,xshift=+3pt,yshift=+2pt] {\scriptsize $#6$};
  \draw (4.8,0) node[anchor=north,xshift=+0pt,yshift=+2pt] {\scriptsize $#7$};
\end{scope}
\end{scope}
}
\newcommand{\cbWidth}{.2cm}
\newcommand{\cbHeight}{4cm}
\newcommand{\xcb}{.5cm}
\newcommand{\ycb}{-.2cm}
\setlength{\ycbTop}{\ycb+\cbHeight}
\setlength{\ycbMid}{\ycb+\cbHeight*\real{.5}}
\newcommand{\trimfigcb}[3]{\includegraphics[width=#2, height=#3, clip, trim=17cm 2.35cm 1.65cm 2.35cm]{#1}}
\newcommand{\figWidtha}{4.25cm}
\newcommand{\trimfiga}[2]{\trimw{#1}{#2}{.11}{.117}{.11}{.11}}
\begin{figure}[htb]
\begin{center}
\begin{tikzpicture}
   \useasboundingbox (0,.3) rectangle (15,5);  

   \begin{scope}[yshift=0cm]
    \drawContour{xshift= 0.cm,yshift=0.00cm}{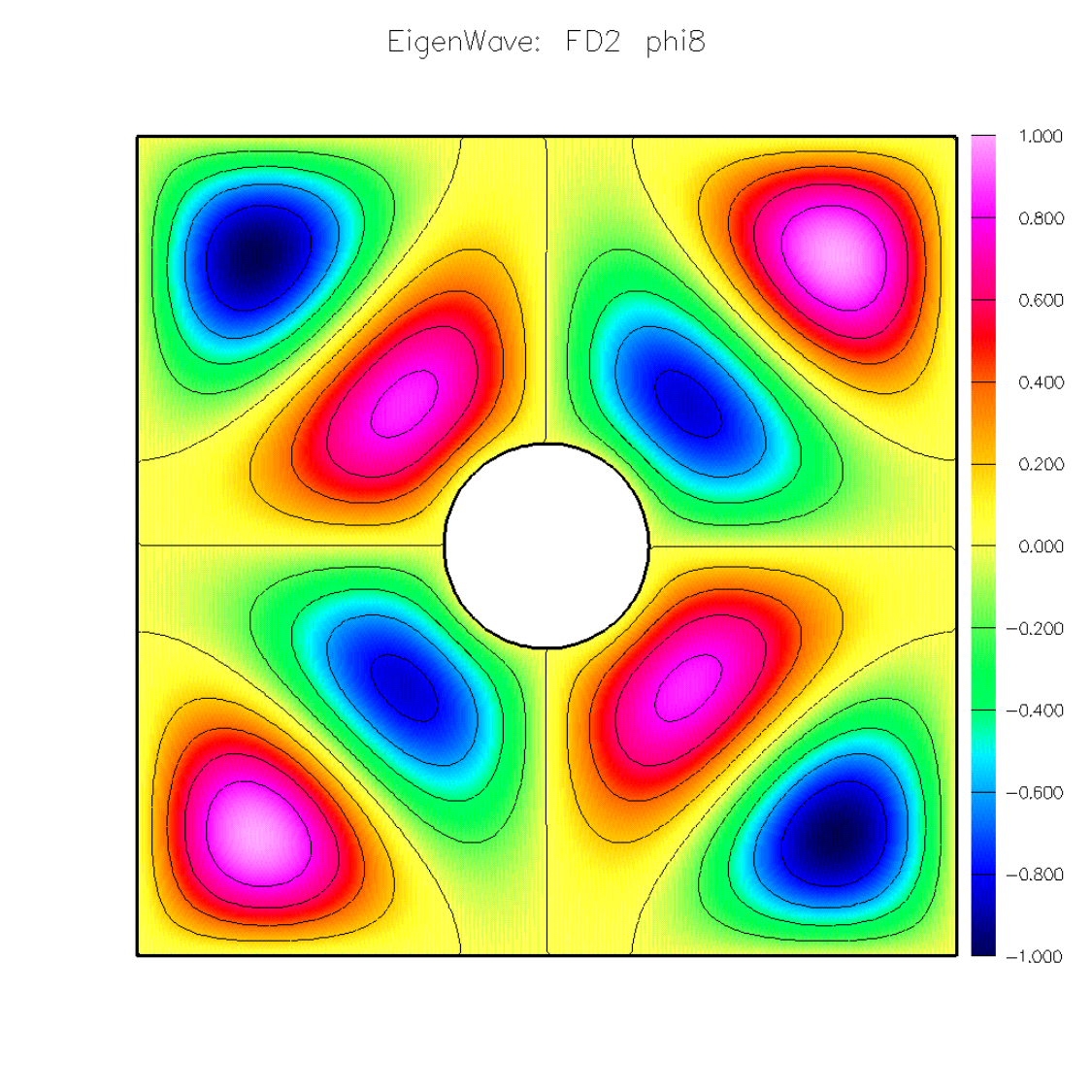}{$\phi^{(11)}$}{$v$}{$t=0.3$}{$-1.0$}{$1.0$}    
    \drawContour{xshift= 5.cm,yshift=0.00cm}{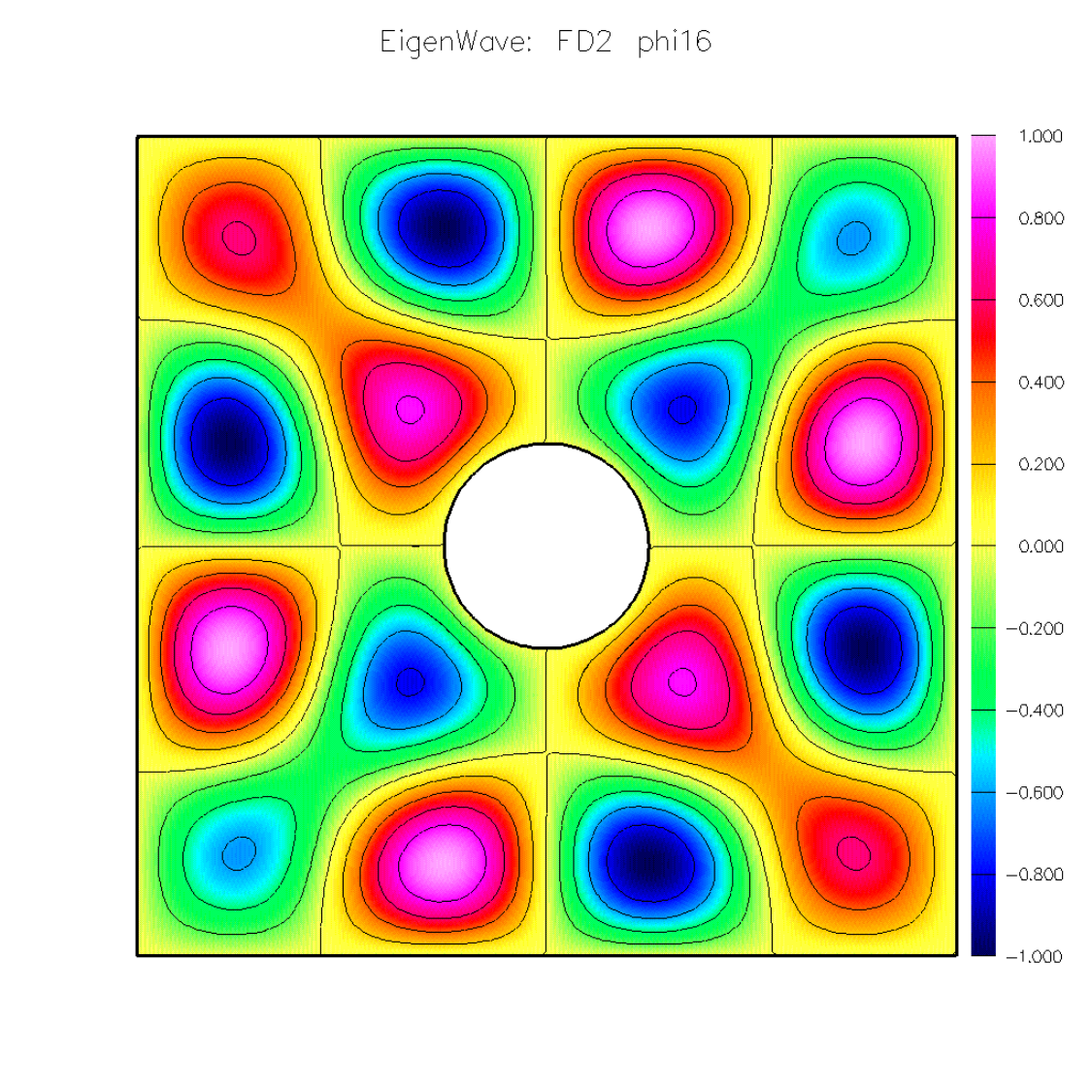}{$\phi^{(19)}$}{$v$}{$t=0.3$}{$-1.0$}{$1.0$}    
    \drawContour{xshift=10.cm,yshift=0.00cm}{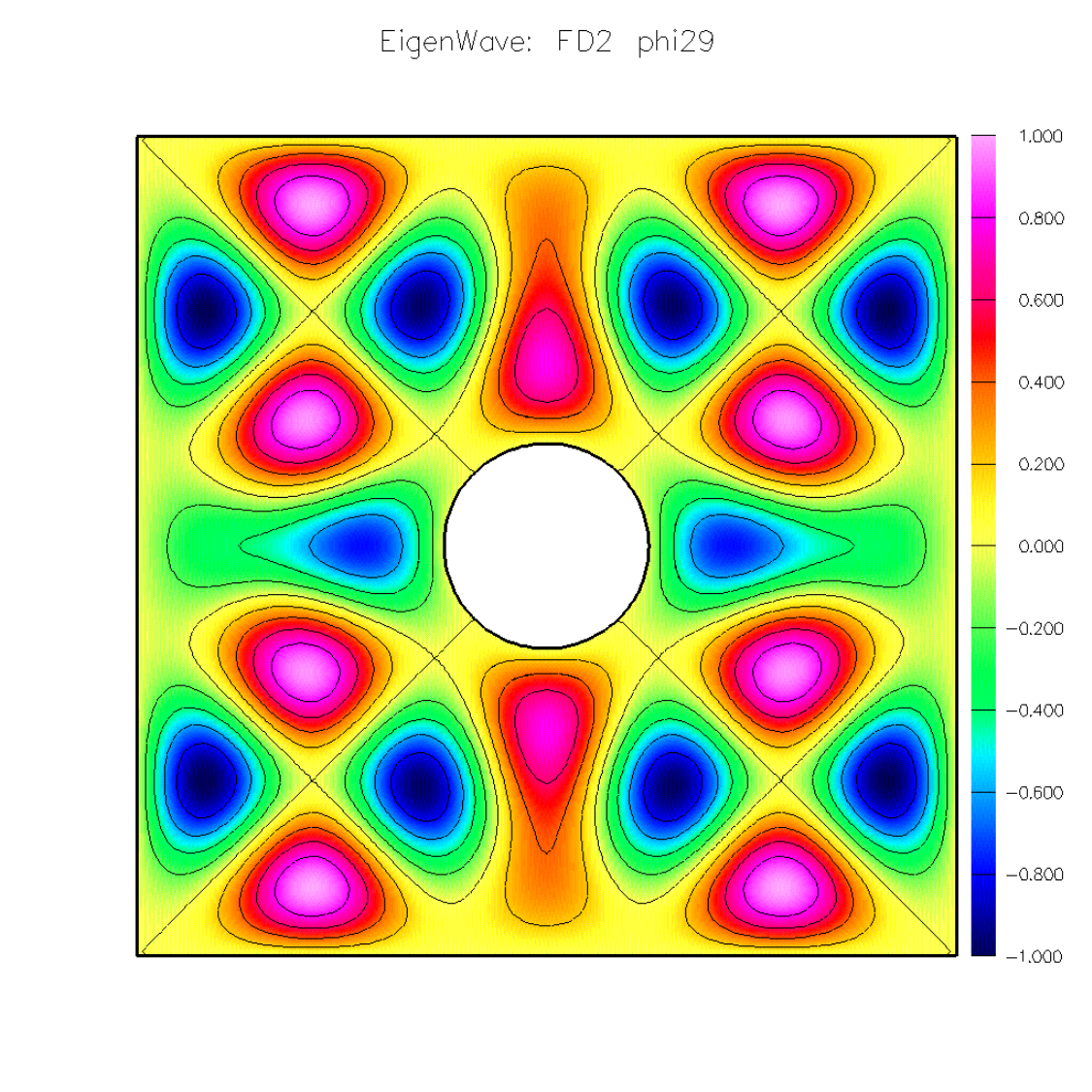}{$\phi^{(32)}$}{$v$}{$t=0.3$}{$-1.0$}{$1.0$}    
   \end{scope}

\end{tikzpicture}
\end{center}
\caption{Circle in a square: some computed eigenvectors. 
    }
\label{fig:cicEigenVectors}
\end{figure}
}

%% file: tex/threeShapesEigenpairs.tex
\newcommand{\Gcts}{\Gc_{ts}}
\subsection{Eigenmodes of a domain with three shapes}

\input tex/threeShapesGridFig.tex

In this study we compute eigenpairs on a domain with cutout regions in the form of three shapes, a circle, rectangle and triangle.
The \textsl{three-shapes} overset grid, denoted by $\Gcts^{(j)}$, consists of four component grids as shown
in Figure~\ref{fig:shapesGrid}. A blue background Cartesian grid covers the domain $[-1.25,1.25]\times[-1,1]$.
A red annular grid lies adjacent to the circle of radius $0.25$ and center $(-.5,.35)$. 
A green boundary fitted grid lies adjacent to a triangle 
with rounded corners; the vertices (before rounding) are 
$(-.866025,0)$, $(.866025,.5)$ and $(.866025,0)$.
A pink boundary fitted
grid lies adjacent to a rectangle with rounded corners; the vertices (before rounding) are located at 
$(-1,-.25)$, $(-1,-.75)$, $(-.25,-.75)$, and $(-.25,-.25)$.
The component grids in $\Gcts^{(j)}$ have a target grid spacing of $\ds^{(j)}=1/(10 j)$. 
Each boundary fitted curvilinear grid has $8$ grid lines in the normal direction.

\input tex/threeShapesContoursFig

Table~\ref{tab:shapese4Order2Summary} summarizes results using the KrylovSchur algorithm 
and indicates that $54$ eigenpairs were accurately computed with a total of $151$ wave solves for an
average of approximately $2.8$ wave-solves per eigenpair.
The right graph in Figure~\ref{fig:shapesGrid} shows the filter-function $\beta(\lambda;\omega)$ with the computed
eigenvalues marked. 
{
\renewcommand{\tableFontSize}{\scriptsize}

\input tables/shapesG4O2ImpEig10Krylov.tex
  \eigenWaveSummaryTable
}

%% file: tex/threeShapesGridFig.tex
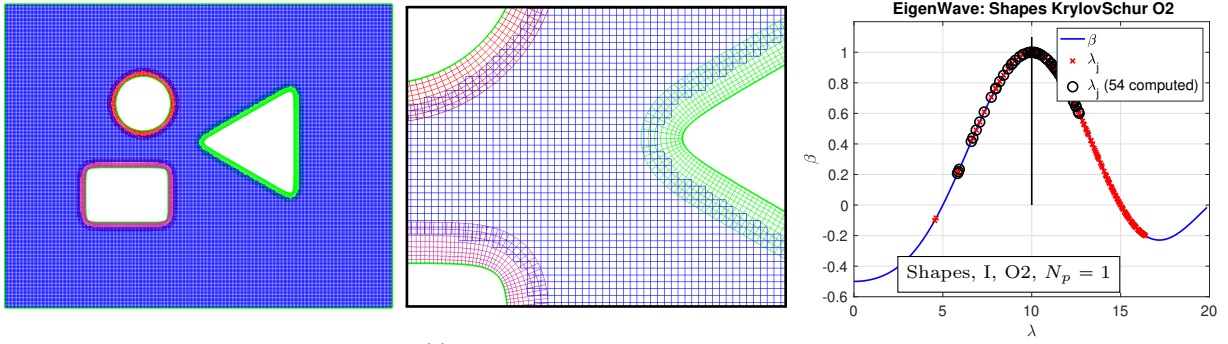
\begin{figure}[htb]
\begin{center}
\begin{tikzpicture}
  \useasboundingbox (0,.25) rectangle (16,4);  
  \begin{scope}[xshift=-4pt,yshift=-2pt]
    \figByWidth{0.0}{0}{fig/shapesGridG4}{5.2cm}[0.025][0.025][0.125][0.125]
  \end{scope}
  \figByWidthb{5.25}{0}{fig/shapesGridG4Zoom}{5cm}[0.025][0.025][0.125][0.125]

  \begin{scope}[xshift=10.5cm,yshift=-12pt]
     \figByWidth{0.0}{0}{fig/shapesG4O2ImpEig10Krylov}{5.5cm}[0][0][0][0];
     \draw (1.2,.85) node[draw,fill=white,anchor=west,xshift=0pt,yshift=0pt] {\scriptsize Shapes, I, O2, $\Np=1$};
  \end{scope}  

\end{tikzpicture}
\end{center}
\caption{Left: \textsl{three-shapes} overset grid $\Gcts^{(4)}$. Middle: magnified view of the grid. 
Right: graph of the filter function $\beta$ with the  
computed eigenvalues marked with black circles.
 }
\label{fig:shapesGrid}
\end{figure}

%% file: tex/threeShapesContoursFig.tex
{
\newcommand{\drawContour}[7]{%
\begin{scope}[#1]
\draw(0.0,0) node[anchor=south west,xshift=-4pt,yshift=+0pt] {\trimfiga{fig/#2}{\figWidtha}};
  \draw(.5,.5) node[draw,fill=white,anchor=west,xshift=2pt,yshift=2pt] {\scriptsize #3};
\begin{scope}[xshift=-.15cm,yshift=0pt]
  \draw (\xcb,\ycb) node[anchor=south west,xshift=0.25cm,yshift=.5cm,rotate=-90] {\trimfigcb{fig/colourBarLines}{\cbWidth}{\cbHeight}};
  \draw (.8,0) node[anchor=north,xshift=+3pt,yshift=+2pt] {\scriptsize $#6$};
  \draw (4.8,0) node[anchor=north,xshift=+0pt,yshift=+2pt] {\scriptsize $#7$};
\end{scope}
\end{scope}
}
\newcommand{\cbWidth}{.2cm}
\newcommand{\cbHeight}{4cm}
\newcommand{\xcb}{.5cm}
\newcommand{\ycb}{-.2cm}
\setlength{\ycbTop}{\ycb+\cbHeight}
\setlength{\ycbMid}{\ycb+\cbHeight*\real{.5}}
\newcommand{\trimfigcb}[3]{\includegraphics[width=#2, height=#3, clip, trim=17cm 2.35cm 1.65cm 2.35cm]{#1}}
\newcommand{\figWidtha}{4.75cm}
\newcommand{\trimfiga}[2]{\trimw{#1}{#2}{.11}{.117}{.18}{.18}}
\begin{figure}[htb]
\begin{center}
\begin{tikzpicture}
   \useasboundingbox (0,.25) rectangle (15,4.5);  

   \begin{scope}[yshift=0cm]
    \drawContour{xshift= 0.cm,yshift=0.00cm}{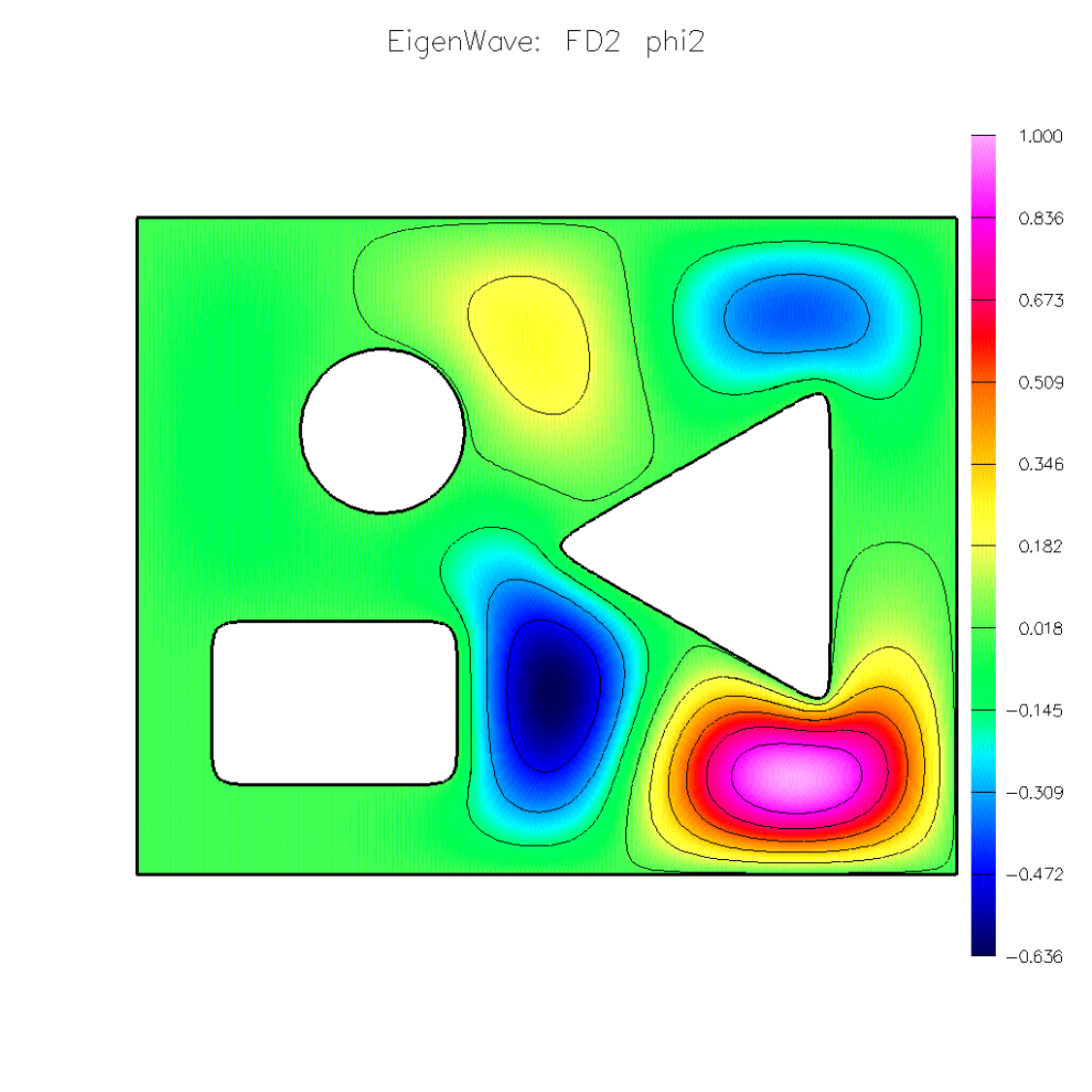}{$\phi^{(4)}$}{$v$}{$t=0.3$}{$-1.0$}{$1.0$}    
    \drawContour{xshift= 5.cm,yshift=0.00cm}{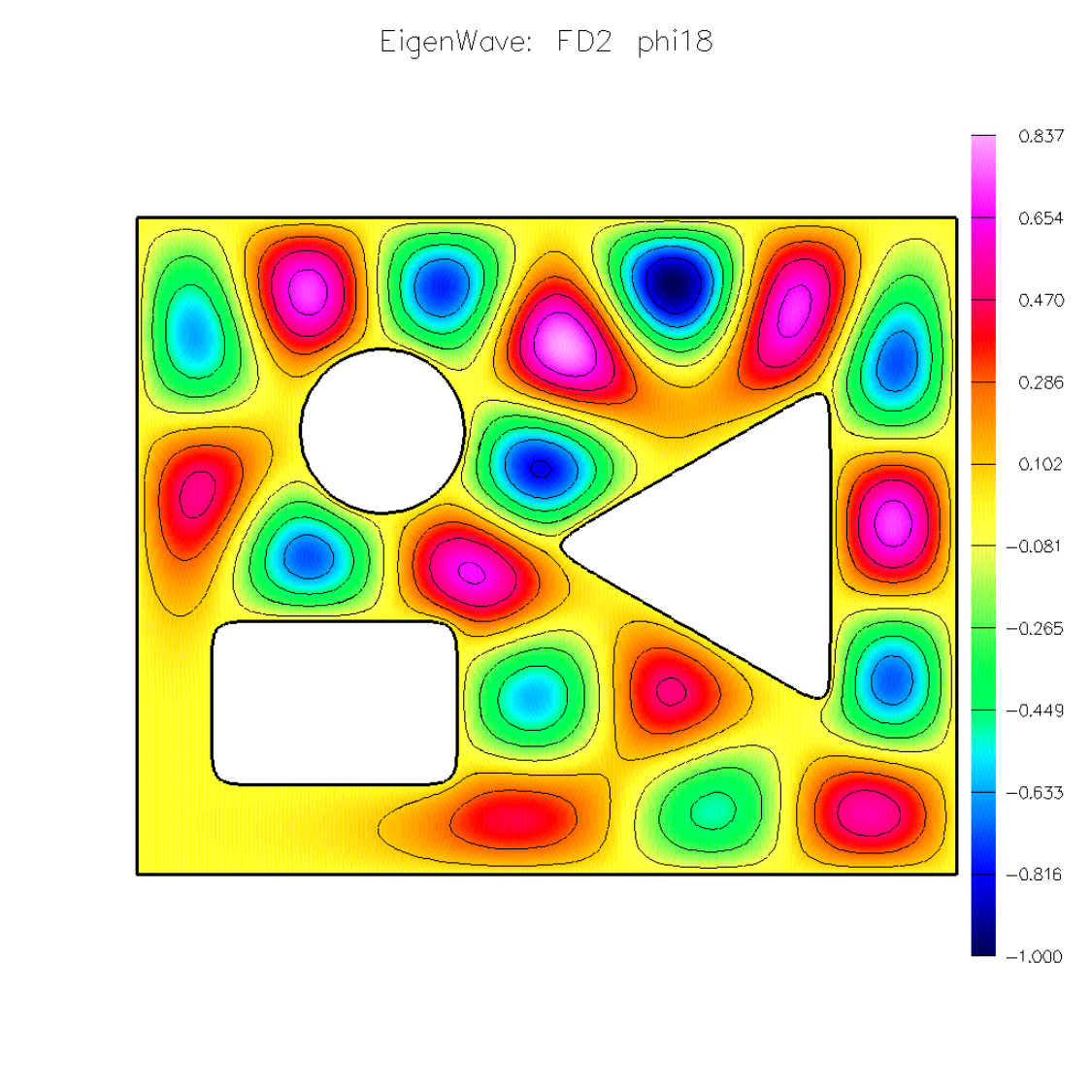}{$\phi^{(20)}$}{$v$}{$t=0.3$}{$-1.0$}{$1.0$}    
    \drawContour{xshift=10.cm,yshift=0.00cm}{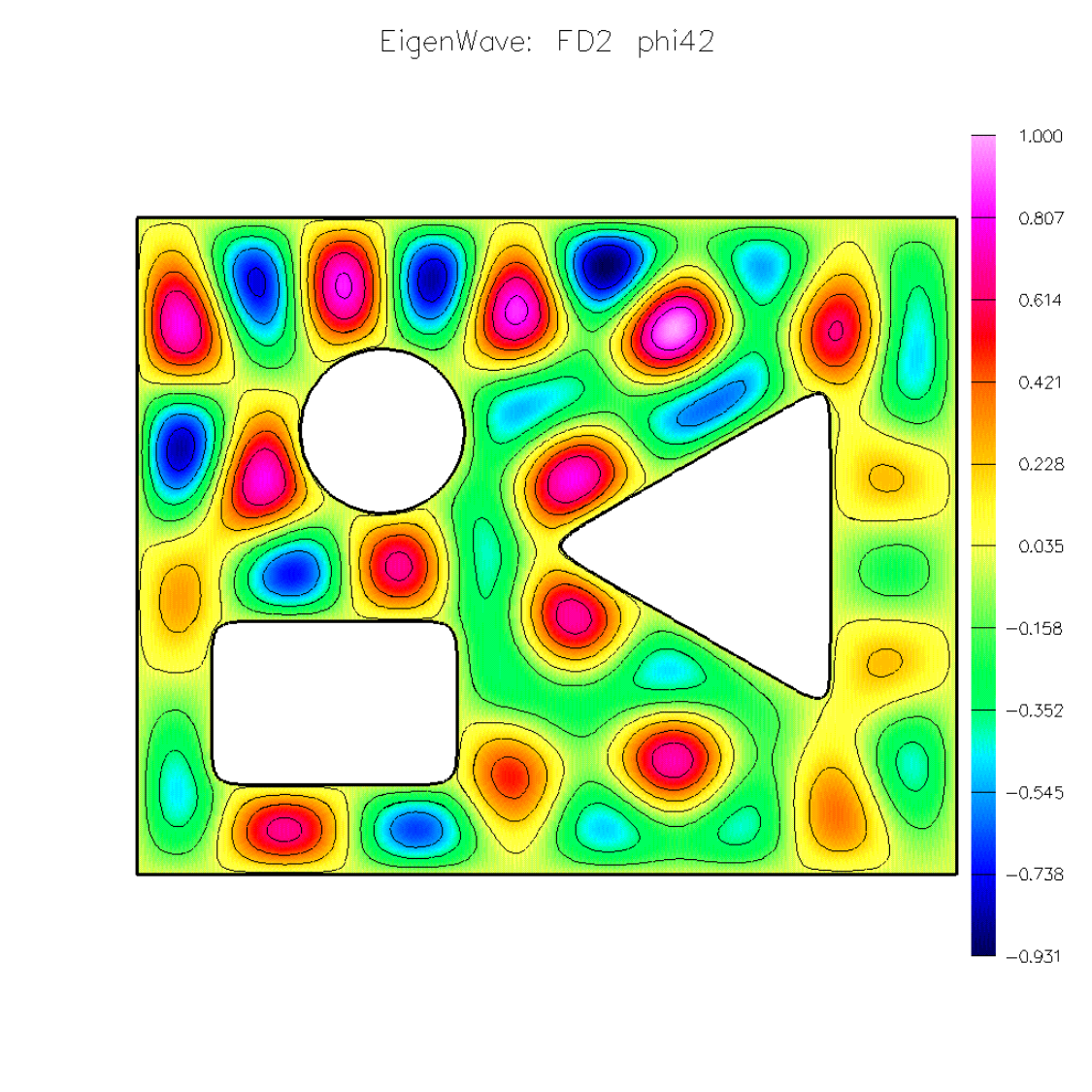}{$\phi^{(44)}$}{$v$}{$t=0.3$}{$-1.0$}{$1.0$}    
   \end{scope}

\end{tikzpicture}
\end{center}
\caption{Three shapes. Selected computed eigenvectors.
    }
\label{fig:threeShapesEigenVectors}
\end{figure}
}

%% file: tables/shapesG4O2ImpEig10Krylov.tex
\newcommand{\eigenWaveLongTable}{
\begin{table}[hbt]
\begin{center}\tableFontSize
\begin{tabular}{|c|r|r|c|c|c|c|c|} \hline
 \multicolumn{8}{|c|}{EigenWave: grid=shapese4, ts=implicit, order=2, $N_p=1$, KrylovSchur } \\ \hline 
   $j$  & \multicolumn{1}{c|}{$\lambda_j$} &  \multicolumn{1}{c|}{$\lambda_{\rm true}$}  & eig &  mult &  eig-err & evect-err & eig-res  \\ \hline
   0   &   6.029340 &   6.029340 &     2 &    1  &  1.92e-15 &  7.47e-13 &  1.58e-12 \\ 
   1   &   6.080277 &   6.080277 &     3 &    1  &  4.38e-15 &  8.39e-13 &  1.39e-12 \\ 
   2   &   6.143754 &   6.143754 &     4 &    1  &  9.40e-15 &  2.69e-13 &  1.46e-12 \\ 
   3   &   6.902001 &   6.902001 &     5 &    1  &  5.53e-15 &  3.01e-13 &  1.23e-12 \\ 
   4   &   7.000978 &   7.000978 &     6 &    1  &  1.90e-15 &  2.84e-12 &  7.31e-12 \\ 
   5   &   7.218703 &   7.218703 &     7 &    1  &  1.13e-14 &  1.08e-13 &  7.44e-13 \\ 
   6   &   7.450822 &   7.450822 &     8 &    1  &  7.15e-16 &  1.77e-13 &  9.77e-13 \\ 
   7   &   7.748758 &   7.748758 &     9 &    1  &  8.02e-15 &  8.09e-13 &  1.38e-12 \\ 
   8   &   8.220809 &   8.220809 &    10 &    1  &  6.27e-15 &  2.88e-13 &  7.73e-13 \\ 
   9   &   8.517618 &   8.517618 &    11 &    1  &  4.38e-15 &  3.91e-12 &  3.23e-12 \\ 
  10   &   8.543484 &   8.543484 &    12 &    1  &  8.32e-16 &  3.24e-13 &  4.40e-13 \\ 
  11   &   8.798151 &   8.798151 &    13 &    1  &  1.21e-15 &  3.70e-13 &  7.42e-13 \\ 
  12   &   9.069582 &   9.069582 &    14 &    1  &  2.35e-15 &  6.41e-13 &  5.68e-13 \\ 
  13   &   9.379993 &   9.379993 &    15 &    1  &  1.89e-15 &  2.47e-12 &  9.62e-13 \\ 
  14   &   9.632984 &   9.632984 &    16 &    1  &  1.29e-14 &  2.97e-11 &  3.75e-12 \\ 
  15   &   9.704143 &   9.704143 &    17 &    1  &  1.83e-15 &  2.16e-12 &  7.46e-13 \\ 
  16   &   9.903703 &   9.903703 &    18 &    1  &  1.79e-15 &  3.52e-13 &  5.19e-13 \\ 
  17   &  10.235274 &  10.235274 &    19 &    1  &  8.33e-15 &  2.61e-11 &  3.06e-12 \\ 
  18   &  10.457170 &  10.457170 &    20 &    1  &  9.51e-15 &  3.10e-13 &  4.48e-13 \\ 
  19   &  10.532813 &  10.532813 &    21 &    1  &  2.53e-15 &  1.74e-12 &  6.14e-13 \\ 
  20   &  10.852120 &  10.852120 &    22 &    1  &  1.64e-16 &  2.83e-13 &  4.39e-13 \\ 
  21   &  11.035549 &  11.035549 &    23 &    1  &  2.09e-15 &  8.79e-13 &  4.04e-13 \\ 
  22   &  11.053457 &  11.053457 &    24 &    1  &  1.12e-15 &  1.05e-12 &  2.86e-13 \\ 
  23   &  11.273879 &  11.273879 &    25 &    1  &  5.36e-15 &  1.81e-13 &  2.90e-13 \\ 
  24   &  11.393003 &  11.393003 &    26 &    1  &  1.25e-15 &  3.54e-12 &  1.64e-12 \\ 
  25   &  11.405776 &  11.405776 &    27 &    1  &  1.71e-15 &  3.20e-13 &  2.28e-13 \\ 
  26   &  11.497074 &  11.497074 &    28 &    1  &  4.94e-15 &  1.16e-13 &  2.19e-13 \\ 
  27   &  11.605039 &  11.605039 &    29 &    1  &  7.04e-15 &  1.69e-13 &  3.58e-13 \\ 
  28   &  11.774829 &  11.774829 &    30 &    1  &  1.66e-15 &  1.53e-13 &  2.82e-13 \\ 
  29   &  11.955016 &  11.955016 &    31 &    1  &  2.23e-15 &  7.37e-14 &  3.75e-13 \\ 
  30   &  12.271069 &  12.271069 &    32 &    1  &  5.65e-15 &  3.48e-13 &  3.19e-13 \\ 
  31   &  12.441207 &  12.441207 &    33 &    1  &  2.57e-15 &  2.09e-13 &  2.36e-13 \\ 
  32   &  12.697583 &  12.697583 &    34 &    1  &  7.13e-15 &  1.87e-13 &  3.73e-13 \\ 
  33   &  12.784269 &  12.784269 &    35 &    1  &  1.11e-15 &  1.20e-13 &  3.06e-13 \\ 
  34   &  12.913410 &  12.913410 &    36 &    1  &  5.50e-15 &  1.83e-13 &  1.76e-13 \\ 
  35   &  12.998371 &  12.998371 &    37 &    1  &  4.10e-15 &  2.39e-13 &  2.55e-13 \\ 
  36   &  13.095528 &  13.095528 &    38 &    1  &  7.87e-15 &  1.92e-12 &  1.52e-12 \\ 
  37   &  13.200508 &  13.200508 &    39 &    1  &  4.44e-15 &  1.30e-13 &  1.43e-13 \\ 
  38   &  13.418508 &  13.418508 &    40 &    1  &  1.85e-15 &  4.82e-14 &  2.92e-13 \\ 
  39   &  13.627207 &  13.627207 &    41 &    1  &  3.65e-15 &  3.37e-13 &  2.22e-13 \\ 
  40   &  13.653144 &  13.653144 &    42 &    1  &  4.29e-15 &  1.70e-13 &  2.38e-13 \\ 
  41   &  13.792644 &  13.792644 &    43 &    1  &  1.55e-15 &  1.39e-13 &  1.99e-13 \\ 
  42   &  13.922781 &  13.922781 &    44 &    1  &  5.61e-15 &  1.33e-13 &  5.60e-13 \\ 
  43   &  14.084873 &  14.084873 &    45 &    1  &  5.04e-16 &  4.58e-13 &  2.64e-13 \\ 
  44   &  14.128528 &  14.128528 &    46 &    1  &  5.28e-15 &  4.82e-13 &  3.02e-13 \\ 
  45   &  14.387199 &  14.387199 &    47 &    1  &  7.41e-16 &  1.43e-13 &  1.77e-13 \\ 
  46   &  14.499515 &  14.499515 &    48 &    1  &  2.45e-16 &  9.69e-14 &  1.78e-13 \\ 
  47   &  14.705445 &  14.705445 &    49 &    1  &  4.71e-15 &  1.85e-13 &  2.34e-13 \\ 
  48   &  14.757477 &  14.757477 &    50 &    1  &  3.01e-15 &  3.47e-13 &  2.39e-13 \\ 
  49   &  14.806384 &  14.806384 &    51 &    1  &  6.00e-16 &  3.51e-13 &  2.40e-13 \\ 
  50   &  14.962685 &  14.962685 &    52 &    1  &  4.51e-15 &  6.24e-13 &  1.94e-13 \\ 
  51   &  14.981501 &  14.981501 &    53 &    1  &  2.13e-15 &  5.55e-13 &  1.90e-13 \\ 
  52   &  15.205037 &  15.205037 &    54 &    1  &  2.57e-15 &  2.29e-13 &  1.50e-13 \\ 
  53   &  15.286595 &  15.286595 &    55 &    1  &  6.97e-15 &  5.18e-13 &  1.55e-12 \\ 
 \hline 
\end{tabular}
\end{center}
\caption{grid=shapese4, method=KrylovSchur, ts=implicit, order=2, $N_p=1$.}
\label{tab:shapese4Order2}
\end{table}
}

\newcommand{\eigenWaveSummaryTable}{
\begin{table}[hbt]
\begin{center}\tableFontSize
\begin{tabular}{|c|c|c|c|c|c|c|c|} \hline
 \multicolumn{8}{|c|}{EigenWave: grid=shapese4, ts=implicit, order=2, $N_p=1$, KrylovSchur } \\ \hline 
   num   &  wave       & time-steps  & wave-solves &  time-steps &   max      &  max      &  max       \\ 
   eigs  &  solves     & per period  &  per eig    &  per-eig    &   eig-err  & evect-err & eig-res    \\ 
 \hline
    54    &       151     &      10     &    2.8    &     27     &  1.29e-14 &  2.97e-11 &  7.31e-12 \\
 \hline 
\end{tabular}
\end{center}
\caption{EigenWave: grid=shapese4, method=KrylovSchur, ts=implicit, order=2, $N_p=1$.}
\label{tab:shapese4Order2Summary}
\end{table}
}

%% file: tex/rpiEigenpairs.tex
\newcommand{\Gcrpi}{\Gc_{\rm rpi}}
\subsection{Eigenmodes of the RPI grid} \label{sec:rpiEigenpairs}

Eigenpairs are computed for a domain with the letters R, P and I (also mentioned previously in the Section~\ref{sec:intro})
The letters are enclosed within the rectangular region $[0,8.5]\times[-2,2]$.
The overset grid for the domain, denoted by $\Gcrpi^{(j)}$ consists of four component grids, each with
a target grid spacing of $\ds^{(j)}=1/(10 j)$.
As shown in Figure~\ref{fig:rpiGridAndFilter} narrow body-fitted grids are used to form the outline of each letter, and these are
embedded in a background Cartesian grid.

\input tex/rpiGridFig.tex

{

\input tables/rpiG4O2ImpEigKrylov.tex
  \eigenWaveSummaryTable
}

Table~\ref{tab:rpiGride4Order2Summary} summarizes results of computing eigenpairs to second-order accuracy 
on grid $\Gcrpi^{(4)}$.
The target frequency was $\omega=4$. Forty-eight eigenpairs were requested and sixty eigenpairs were accurately found by the KrylovSchur
algorithm in a total of $177$ wave-solves. This corresponds to approximately 
$3$ wave-solves per eigenpair found.
The right graph in Figure~\ref{fig:rpiGrid} shows the filter-function $\beta(\lambda;\omega)$ with the computed
eigenvalues marked. 
Figure~\ref{fig:rpiEigenfunctionSolution} plots contours of some of the computed eigenvectors.

\input tex/rpiContoursFig.tex

\input tex/rpiHighFreqContoursFig.tex

Eigenpairs are also computed using a higher target frequency and a larger value for $\Np$, the number of periods
over which the wave equation is integrated.
This makes the main peak of the filter function narrower with fewer eigenvalues 
lying near the peak. As a result, a small number of eigenpairs can be found in an efficient way.
Reducing the number of computed eigenpairs may be desired, for example,
to avoid the excessive storage of keeping a large number of eigenvectors.
The disadvantage of increasing $\Np$ is that the cost per wave-solve increases by a factor of $\Np$.

{

\input tables/rpiG4O2ImpEig12Krylov.tex

\eigenWaveSummaryTable
\renewcommand{\tableFontSize}{\scriptsize}
}

Table~\ref{tab:rpiGride4Order2Eig12Summary} summarizes results of computing eigenpairs to second-order accuracy
 on grid $\Gcrpi^{(4)}$.
The target frequency was $\omega=12$. 
Fourty-eight eigenpairs were requested and fifty-one eigenpairs were found by EigenWave with the KrylovSchur
algorithm in a total of $192$ wave-solves. This corresponds to approximately 
$3.8$ wave-solves per eigenpair found.  
The graph in Figure~\ref{fig:rpiGridEigsNp} shows the filter-function $\beta(\lambda;\omega)$ with the computed
eigenpairs marked. 
Figure~\ref{fig:rpiHighFreqEigenvectors} shows contours of some eigenvectors.

{
\newcommand{\figw}{7.5cm}
\newcommand{\figh}{6.25cm}

\begin{figure}[htb]
\begin{center}
\begin{tikzpicture}[scale=1]
  \useasboundingbox (0,.7) rectangle (\figw,\figh);  

  \begin{scope}[xshift=0cm]
    \figByWidth{0.0}{0}{fig/rpiG4O2ImpEig12Krylov}{\figw}[0][0][0][0]
    \draw (1.5,1) node[draw,fill=white,anchor=west,xshift=0pt,yshift=0pt] {\scriptsize RPI, I, O2, $\Np=8$};
  \end{scope}

\end{tikzpicture}
\end{center}
\caption{Computing eigenpairs of the RPI domain with more periods (larger $\Np$). 
Graphs of the filter function $\beta=\beta(\lambda;\omega)$ 
with the true discrete eigenvalues $\lambda_j$ marked
with red x's and the computed eigenvalues marked with black circles.
The graph shows second-order accurate results with target frequency $\omega=12$.
}
\label{fig:rpiGridEigsNp}
\end{figure}
}

%% file: tex/rpiGridFig.tex
{
\newcommand{\figw}{6.5cm}
\newcommand{\figh}{5.5cm}
\begin{figure}[htb]
\begin{center}
\begin{tikzpicture}
   \useasboundingbox (0,.) rectangle (15,4.5);  
  \figByWidth{  0}{0.5}{fig/rpiGridG4Small}{7.65cm}[0.035][0.035][0.275][0.275]
  \figByWidthb{2.}{-.7}{fig/rpiGridZoom}{2.5cm}[0.0][0.0][0.][0.]
  \begin{scope}[xshift=8cm,yshift=-0.75cm]
     \figByWidth{0.0}{0}{fig/rpiG4O2ImpEigKrylov}{\figw}[0][0][0][0]
     \draw (1.5,1) node[draw,fill=white,anchor=west,xshift=0pt,yshift=0pt] {\scriptsize Disk, I, O2, $\Np=1$};
  \end{scope}  
\end{tikzpicture}
\end{center}
\caption{Left: overset grid for the letters R, P, and I, and a magnified portion. 
Right: graph of the filter function $\beta$ with the  
computed eigenvalues marked with black circles for a second-order accurate computation on grid $\Gcrpi^{(4)}$.
}
\label{fig:rpiGridAndFilter}
\end{figure}
}

%% file: tables/rpiG4O2ImpEigKrylov.tex
\newcommand{\eigenWaveLongTable}{
\begin{table}[hbt]
\begin{center}\tableFontSize
\begin{tabular}{|c|r|r|c|c|c|c|c|} \hline
 \multicolumn{8}{|c|}{EigenWave: grid=rpiGride4, ts=implicit, order=2, $N_p=1$, KrylovSchur } \\ \hline 
   $j$  & \multicolumn{1}{c|}{$\lambda_j$} &  \multicolumn{1}{c|}{$\lambda_{\rm true}$}  & eig &  mult &  eig-err & evect-err & eig-res  \\ \hline
   0   &   2.441115 &   2.441115 &     1 &    1  &  1.46e-15 &  1.07e-13 &  2.02e-12 \\ 
   1   &   2.456254 &   2.456254 &     2 &    1  &  7.41e-15 &  1.98e-13 &  1.77e-12 \\ 
   2   &   2.655458 &   2.655458 &     3 &    1  &  6.69e-16 &  1.95e-13 &  2.16e-12 \\ 
   3   &   2.962038 &   2.962038 &     4 &    1  &  1.05e-14 &  1.01e-13 &  1.82e-12 \\ 
   4   &   3.107722 &   3.107722 &     5 &    1  &  2.00e-15 &  4.91e-13 &  1.41e-12 \\ 
   5   &   3.182482 &   3.182482 &     6 &    1  &  7.12e-15 &  3.54e-13 &  1.19e-12 \\ 
   6   &   3.205546 &   3.205546 &     7 &    1  &  1.11e-15 &  3.14e-13 &  1.14e-12 \\ 
   7   &   3.299822 &   3.299822 &     8 &    1  &  1.87e-14 &  9.98e-14 &  1.04e-12 \\ 
   8   &   3.374167 &   3.374167 &     9 &    1  &  5.66e-15 &  1.47e-13 &  1.09e-12 \\ 
   9   &   3.406837 &   3.406837 &    10 &    1  &  2.09e-15 &  4.64e-13 &  8.17e-13 \\ 
  10   &   3.570304 &   3.570304 &    11 &    1  &  3.86e-15 &  9.80e-14 &  1.03e-12 \\ 
  11   &   3.663343 &   3.663343 &    12 &    1  &  7.03e-15 &  6.79e-13 &  1.20e-12 \\ 
  12   &   3.713957 &   3.713957 &    13 &    1  &  5.26e-15 &  2.07e-13 &  8.88e-13 \\ 
  13   &   3.769488 &   3.769488 &    14 &    1  &  1.18e-16 &  4.97e-13 &  8.15e-13 \\ 
  14   &   3.851635 &   3.851635 &    15 &    1  &  4.84e-15 &  1.19e-12 &  1.13e-12 \\ 
  15   &   3.906854 &   3.906854 &    16 &    1  &  1.14e-16 &  4.10e-13 &  1.10e-12 \\ 
  16   &   3.984613 &   3.984613 &    17 &    1  &  1.10e-14 &  9.89e-13 &  5.58e-13 \\ 
  17   &   4.003862 &   4.003862 &    18 &    1  &  5.32e-15 &  1.43e-12 &  6.22e-13 \\ 
  18   &   4.103577 &   4.103577 &    19 &    1  &  4.11e-15 &  4.42e-13 &  8.47e-13 \\ 
  19   &   4.137022 &   4.137022 &    20 &    1  &  6.23e-15 &  3.80e-12 &  1.01e-12 \\ 
  20   &   4.243583 &   4.243583 &    21 &    1  &  7.53e-15 &  1.43e-12 &  5.96e-13 \\ 
  21   &   4.368943 &   4.368943 &    22 &    1  &  2.03e-15 &  1.25e-12 &  8.92e-13 \\ 
  22   &   4.451446 &   4.451446 &    23 &    1  &  1.14e-14 &  3.94e-12 &  6.89e-13 \\ 
  23   &   4.536284 &   4.536284 &    24 &    1  &  0.00e+00 &  6.43e-13 &  7.70e-13 \\ 
  24   &   4.562011 &   4.562011 &    25 &    1  &  5.06e-15 &  2.73e-12 &  6.77e-13 \\ 
  25   &   4.637474 &   4.637474 &    26 &    1  &  5.36e-15 &  1.63e-13 &  8.43e-13 \\ 
  26   &   4.724604 &   4.724604 &    27 &    1  &  1.13e-15 &  3.34e-13 &  4.74e-13 \\ 
  27   &   4.773788 &   4.773788 &    28 &    1  &  6.33e-15 &  6.66e-13 &  5.50e-13 \\ 
  28   &   4.844252 &   4.844252 &    29 &    1  &  3.85e-15 &  4.60e-12 &  1.13e-12 \\ 
  29   &   4.859014 &   4.859014 &    30 &    1  &  1.02e-14 &  4.45e-13 &  2.24e-13 \\ 
  30   &   4.950589 &   4.950589 &    31 &    1  &  7.36e-15 &  1.34e-12 &  5.11e-13 \\ 
  31   &   4.955494 &   4.955494 &    32 &    1  &  6.09e-15 &  1.06e-12 &  3.98e-13 \\ 
  32   &   5.060037 &   5.060037 &    33 &    1  &  1.76e-15 &  3.79e-13 &  3.05e-13 \\ 
  33   &   5.099385 &   5.099385 &    34 &    1  &  3.48e-15 &  2.20e-13 &  4.28e-13 \\ 
  34   &   5.187012 &   5.187012 &    35 &    1  &  1.37e-15 &  2.10e-13 &  3.77e-13 \\ 
  35   &   5.222875 &   5.222875 &    36 &    1  &  1.70e-16 &  2.50e-13 &  5.04e-13 \\ 
  36   &   5.266958 &   5.266958 &    37 &    1  &  2.70e-15 &  1.06e-13 &  3.03e-13 \\ 
  37   &   5.333580 &   5.333580 &    38 &    1  &  5.50e-15 &  1.04e-12 &  7.40e-13 \\ 
  38   &   5.341157 &   5.341157 &    39 &    1  &  3.16e-15 &  3.95e-13 &  3.05e-13 \\ 
  39   &   5.430771 &   5.430771 &    40 &    1  &  1.39e-14 &  1.25e-13 &  3.84e-13 \\ 
  40   &   5.532926 &   5.532926 &    41 &    1  &  9.63e-16 &  2.72e-13 &  3.54e-13 \\ 
  41   &   5.575845 &   5.575845 &    42 &    1  &  3.82e-15 &  3.98e-13 &  6.02e-13 \\ 
  42   &   5.593307 &   5.593307 &    43 &    1  &  3.97e-15 &  5.20e-13 &  4.55e-13 \\ 
  43   &   5.685843 &   5.685843 &    44 &    1  &  2.81e-15 &  5.64e-13 &  5.43e-13 \\ 
  44   &   5.733256 &   5.733256 &    45 &    1  &  6.35e-15 &  1.41e-13 &  3.32e-13 \\ 
  45   &   5.781582 &   5.781582 &    46 &    1  &  1.54e-16 &  2.27e-13 &  4.00e-13 \\ 
  46   &   5.895501 &   5.895501 &    47 &    1  &  1.96e-15 &  1.54e-13 &  3.45e-13 \\ 
  47   &   5.936179 &   5.936179 &    48 &    1  &  1.32e-14 &  1.33e-13 &  3.74e-13 \\ 
  48   &   5.983383 &   5.983383 &    49 &    1  &  4.45e-16 &  3.41e-13 &  2.29e-13 \\ 
  49   &   6.012838 &   6.012838 &    50 &    1  &  7.98e-15 &  3.08e-13 &  3.47e-13 \\ 
  50   &   6.054506 &   6.054506 &    51 &    1  &  7.33e-16 &  2.86e-13 &  3.30e-13 \\ 
  51   &   6.119700 &   6.119700 &    52 &    1  &  2.90e-15 &  5.51e-14 &  2.27e-13 \\ 
  52   &   6.181501 &   6.181501 &    53 &    1  &  2.16e-15 &  1.58e-13 &  2.53e-13 \\ 
  53   &   6.210201 &   6.210201 &    54 &    1  &  2.29e-15 &  2.57e-13 &  2.64e-13 \\ 
  54   &   6.275184 &   6.275184 &    55 &    1  &  4.25e-16 &  2.29e-13 &  2.91e-13 \\ 
  55   &   6.338780 &   6.338780 &    56 &    1  &  3.64e-15 &  1.90e-13 &  2.57e-13 \\ 
  56   &   6.399247 &   6.399247 &    57 &    1  &  3.89e-15 &  5.35e-14 &  2.84e-13 \\ 
  57   &   6.447315 &   6.447315 &    58 &    1  &  2.62e-15 &  2.79e-13 &  2.62e-13 \\ 
  58   &   6.471012 &   6.471012 &    59 &    1  &  8.10e-15 &  2.33e-13 &  2.17e-13 \\ 
  59   &   6.515043 &   6.515043 &    60 &    1  &  4.77e-15 &  5.86e-13 &  2.66e-13 \\ 
 \hline 
\end{tabular}
\end{center}
\caption{grid=rpiGride4, method=KrylovSchur, ts=implicit, order=2, $N_p=1$.}
\label{tab:rpiGride4Order2}
\end{table}
}

\newcommand{\eigenWaveSummaryTable}{
\begin{table}[hbt]
\begin{center}\tableFontSize
\begin{tabular}{|c|c|c|c|c|c|c|c|} \hline
 \multicolumn{8}{|c|}{EigenWave: grid=rpiGride4, ts=implicit, order=2, $N_p=1$, KrylovSchur } \\ \hline 
   num   &  wave       & time-steps  & wave-solves &  time-steps &   max      &  max      &  max       \\ 
   eigs  &  solves     & per period  &  per eig    &  per-eig    &   eig-err  & evect-err & eig-res    \\ 
 \hline
    60    &       177     &      10     &    3.0    &     29     &  1.87e-14 &  4.60e-12 &  2.16e-12 \\
 \hline 
\end{tabular}
\end{center}
\caption{EigenWave: grid=rpiGride4, method=KrylovSchur, ts=implicit, order=2, $N_p=1$.}
\label{tab:rpiGride4Order2Summary}
\end{table}
}

%% file: tex/rpiContoursFig.tex
{
\newcommand{\drawContour}[7]{%
\begin{scope}[#1]
\draw(0.0,0) node[anchor=south west,xshift=-4pt,yshift=+10pt] {\trimfiga{fig/#2}{\figWidtha}};
  \draw(.5,.5) node[draw,fill=white,anchor=west,xshift=2pt,yshift=1pt] {\scriptsize #3};
\begin{scope}[xshift=.3cm,yshift=-2pt]
  \draw (\xcb,\ycb) node[anchor=south west,xshift=0.25cm,yshift=.5cm,rotate=-90] {\trimfigcb{fig/colourBarLines}{\cbWidth}{\cbHeight}};
  \draw (.8,0) node[anchor=north,xshift=+3pt,yshift=+2pt] {\scriptsize $#6$};
  \draw (4.8,0) node[anchor=north,xshift=+0pt,yshift=+2pt] {\scriptsize $#7$};
\end{scope}
\end{scope}
}
\newcommand{\cbWidth}{.2cm}
\newcommand{\cbHeight}{4cm}
\newcommand{\xcb}{.5cm}
\newcommand{\ycb}{-.2cm}
\setlength{\ycbTop}{\ycb+\cbHeight}
\setlength{\ycbMid}{\ycb+\cbHeight*\real{.5}}
\newcommand{\trimfigcb}[3]{\includegraphics[width=#2, height=#3, clip, trim=17cm 2.35cm 1.65cm 2.35cm]{#1}}
\newcommand{\figWidtha}{5.25cm}
\newcommand{\trimfiga}[2]{\trimw{#1}{#2}{.12}{.117}{.32}{.32}}
\begin{figure}[htb]
\begin{center}
\begin{tikzpicture}
   \useasboundingbox (0,.3) rectangle (16,3);  

   \begin{scope}[xshift=-.5cm,yshift=0cm]
     \drawContour{xshift=0.cm,yshift=0.00cm}{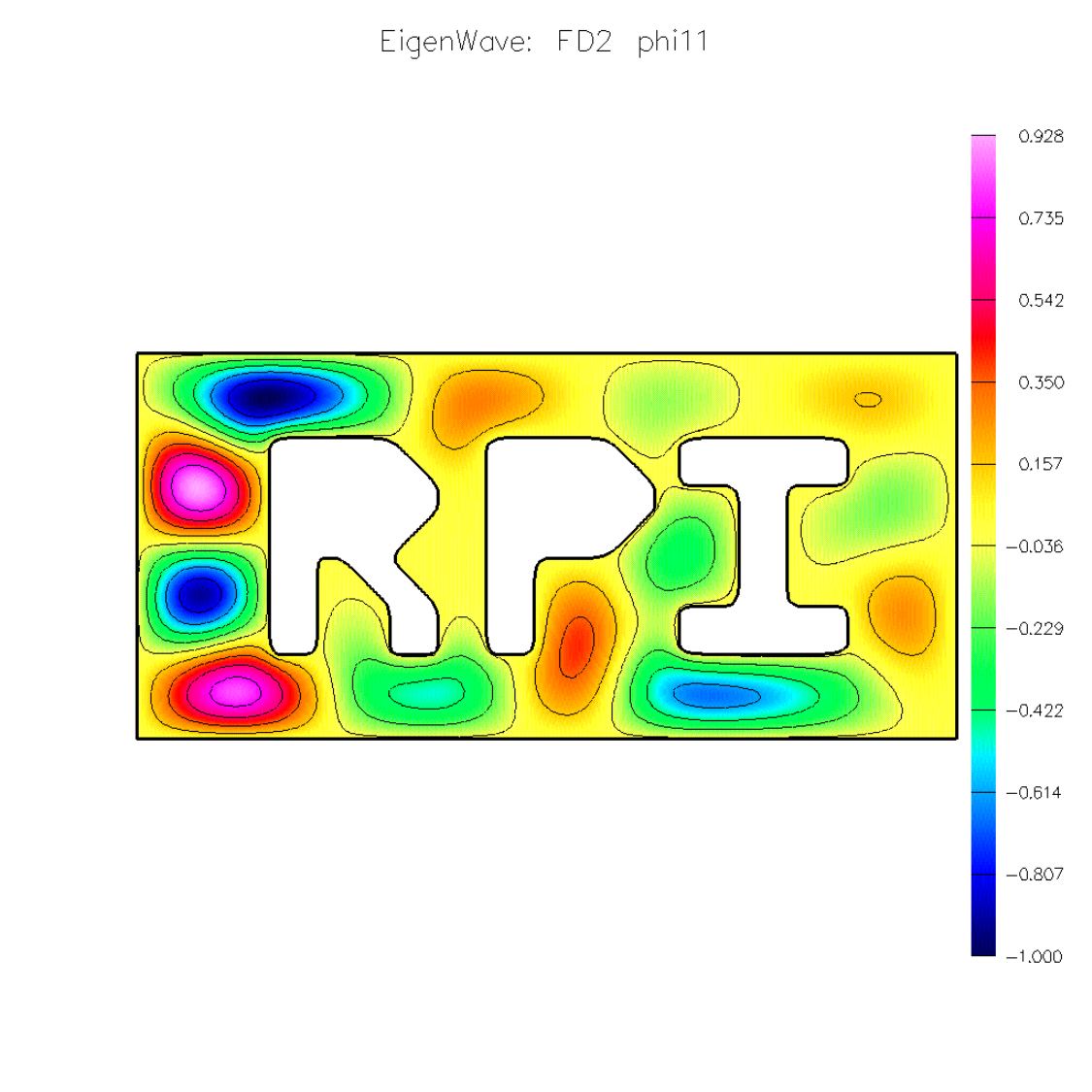}{$\phi^{(12)}$}{$v$}{$t=0.3$}{$-1.0$}{$1.0$}    
     \drawContour{xshift=5.5cm,yshift=0.00cm}{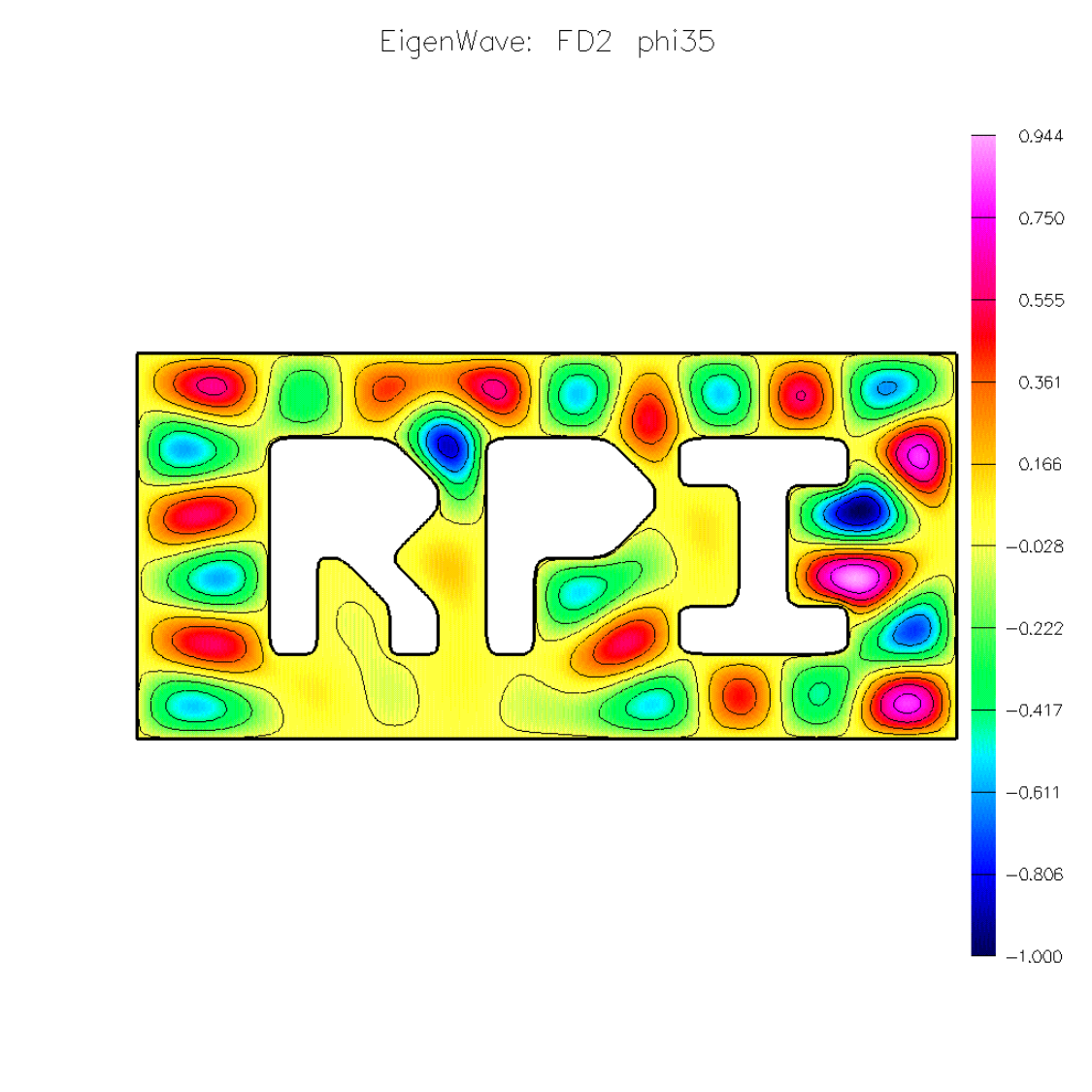}{$\phi^{(36)}$}{$v$}{$t=0.3$}{$-1.0$}{$1.0$}    
     \drawContour{xshift=11cm,yshift=0.00cm}{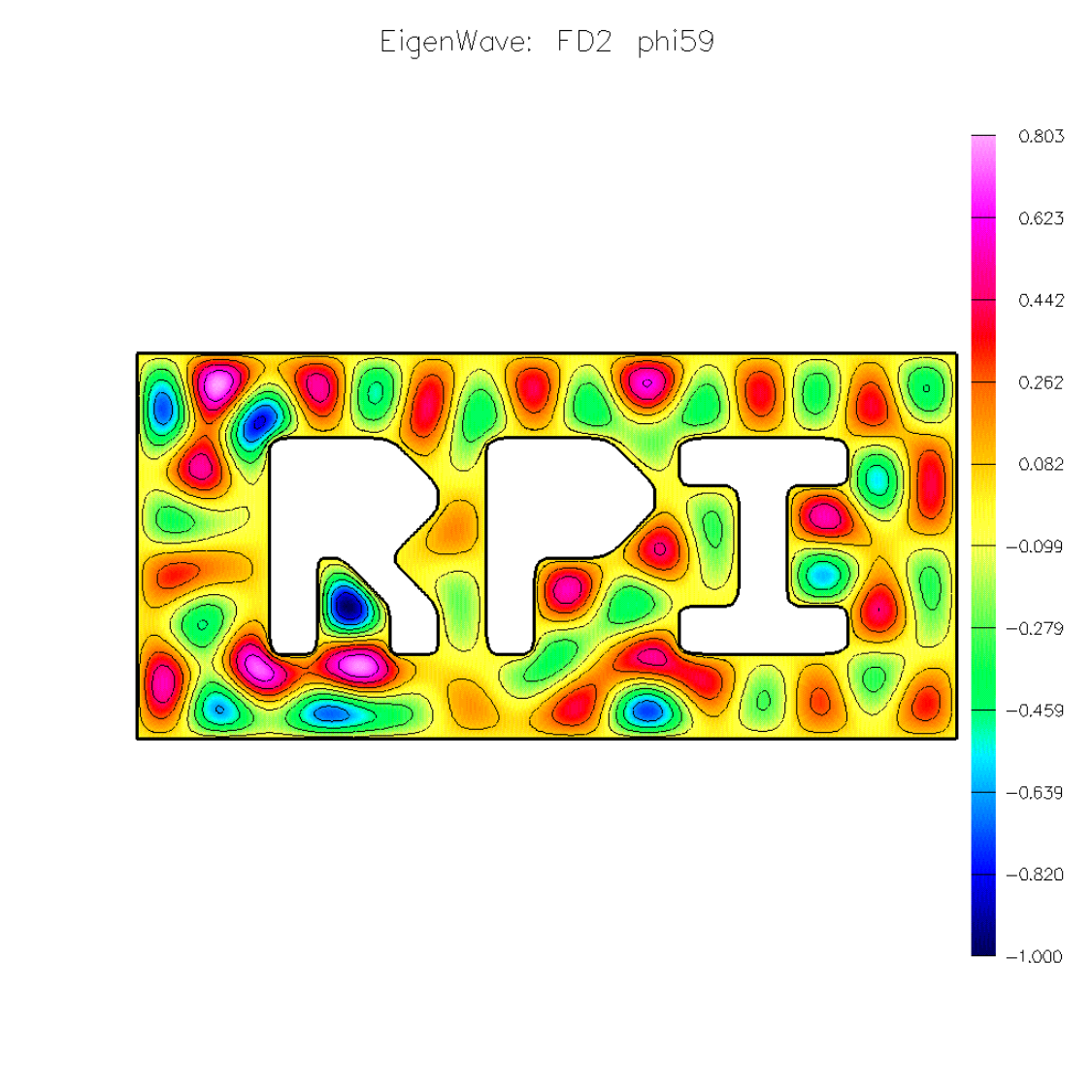}{$\phi^{(60)}$}{$v$}{$t=0.3$}{$-1.0$}{$1.0$}     
   \end{scope}

\end{tikzpicture}
\end{center}
\caption{RPI grid: computed eigenvectors.
    }
\label{fig:rpiEigenfunctionSolution}
\end{figure}
}

%% file: tex/rpiHighFreqContoursFig.tex
{
\newcommand{\drawContour}[7]{%
\begin{scope}[#1]
\draw(0.0,0) node[anchor=south west,xshift=-4pt,yshift=+10pt] {\trimfiga{fig/#2}{\figWidtha}};
  \draw(.5,.5) node[draw,fill=white,anchor=west,xshift=2pt,yshift=1pt] {\scriptsize #3};
\begin{scope}[xshift=.3cm,yshift=-2pt]
  \draw (\xcb,\ycb) node[anchor=south west,xshift=0.25cm,yshift=.5cm,rotate=-90] {\trimfigcb{fig/colourBarLines}{\cbWidth}{\cbHeight}};
  \draw (.8,0) node[anchor=north,xshift=+3pt,yshift=+2pt] {\scriptsize $#6$};
  \draw (4.8,0) node[anchor=north,xshift=+0pt,yshift=+2pt] {\scriptsize $#7$};
\end{scope}
\end{scope}
}
\newcommand{\cbWidth}{.2cm}
\newcommand{\cbHeight}{4cm}
\newcommand{\xcb}{.5cm}
\newcommand{\ycb}{-.2cm}
\setlength{\ycbTop}{\ycb+\cbHeight}
\setlength{\ycbMid}{\ycb+\cbHeight*\real{.5}}
\newcommand{\trimfigcb}[3]{\includegraphics[width=#2, height=#3, clip, trim=17cm 2.35cm 1.65cm 2.35cm]{#1}}
\newcommand{\figWidtha}{5.25cm}
\newcommand{\trimfiga}[2]{\trimw{#1}{#2}{.12}{.117}{.32}{.32}}
\begin{figure}[htb]
\begin{center}
\begin{tikzpicture}
   \useasboundingbox (0,.3) rectangle (16,3);  

   \begin{scope}[xshift=-.5cm,yshift=0cm]
     \drawContour{xshift=0.cm,yshift=0.00cm}{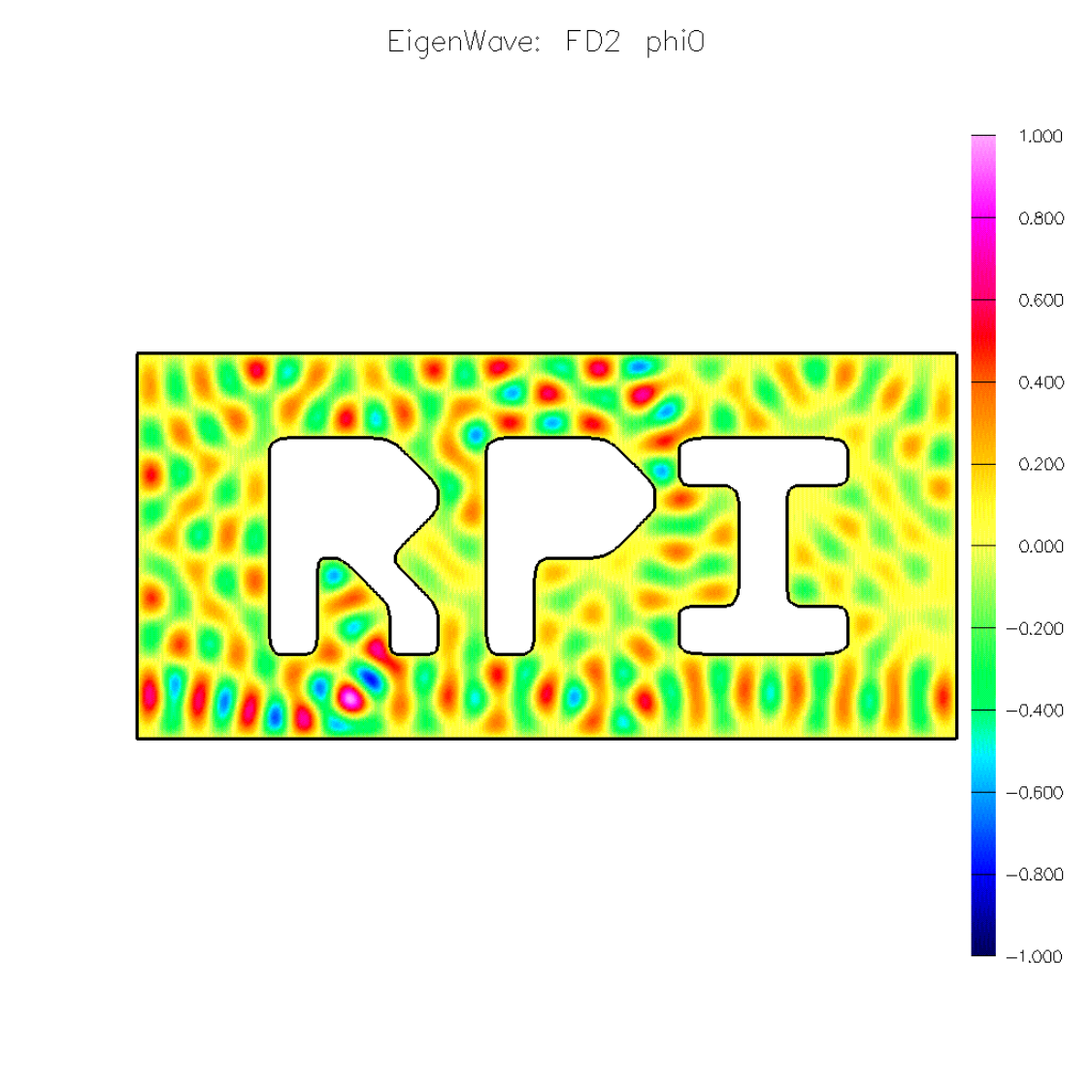}{$\phi^{(284)}$}{$v$}{$t=0.3$}{$-1.0$}{$1.0$}     
     \drawContour{xshift=5.5cm,yshift=0.00cm}{rpiG4O2Eig12phi24}{$\phi^{(308)}$}{$v$}{$t=0.3$}{$-1.0$}{$1.0$}     
     \drawContour{xshift=11cm,yshift=0.00cm}{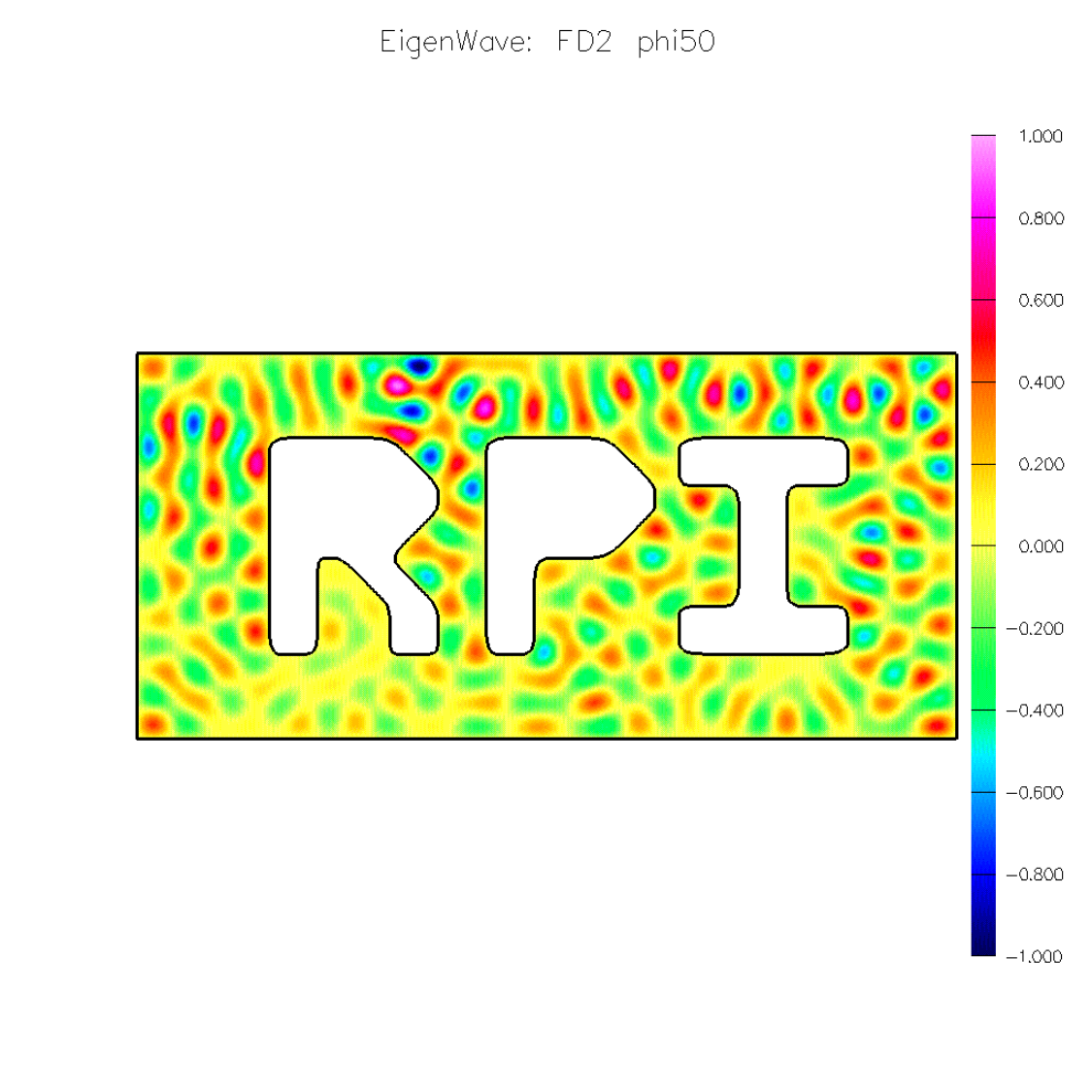}{$\phi^{(334)}$}{$v$}{$t=0.3$}{$-1.0$}{$1.0$}    
   \end{scope}

\end{tikzpicture}
\end{center}
\caption{RPI grid: Selected eigenvectors computed with the \EigenWaveb algoritm.
    }
\label{fig:rpiHighFreqEigenvectors}
\end{figure}
}

%% file: tables/rpiG4O2ImpEig12Krylov.tex
\newcommand{\eigenWaveLongTable}{
\begin{table}[hbt]
\begin{center}\tableFontSize
\begin{tabular}{|c|r|r|c|c|c|c|c|} \hline
 \multicolumn{8}{|c|}{EigenWave: grid=rpiGride4, ts=implicit, order=2, $N_p=8$, KrylovSchur } \\ \hline 
   $j$  & \multicolumn{1}{c|}{$\lambda_j$} &  \multicolumn{1}{c|}{$\lambda_{\rm true}$}  & eig &  mult &  eig-err & evect-err & eig-res  \\ \hline
   0   &  12.648077 &  12.648077 &   284 &    1  &  1.40e-16 &  2.57e-12 &  5.87e-14 \\ 
   1   &  12.664390 &  12.664390 &   285 &    1  &  1.68e-15 &  1.89e-12 &  6.74e-14 \\ 
   2   &  12.679320 &  12.679320 &   286 &    1  &  1.68e-15 &  5.78e-13 &  7.99e-14 \\ 
   3   &  12.692089 &  12.692089 &   287 &    1  &  5.88e-15 &  7.11e-13 &  9.39e-14 \\ 
   4   &  12.710202 &  12.710202 &   288 &    1  &  5.03e-15 &  1.49e-12 &  2.12e-13 \\ 
   5   &  12.730258 &  12.730258 &   289 &    1  &  4.19e-16 &  7.05e-13 &  7.17e-14 \\ 
   6   &  12.756340 &  12.756340 &   290 &    1  &  1.25e-15 &  7.22e-13 &  9.16e-14 \\ 
   7   &  12.764142 &  12.764142 &   291 &    1  &  8.07e-15 &  9.12e-13 &  5.93e-14 \\ 
   8   &  12.782401 &  12.782401 &   292 &    1  &  7.92e-15 &  7.64e-13 &  6.89e-14 \\ 
   9   &  12.801372 &  12.801372 &   293 &    1  &  9.71e-16 &  3.97e-13 &  6.17e-14 \\ 
  10   &  12.805995 &  12.805995 &   294 &    1  &  5.55e-16 &  8.97e-13 &  6.34e-14 \\ 
  11   &  12.822388 &  12.822388 &   295 &    1  &  4.43e-15 &  7.36e-13 &  7.46e-14 \\ 
  12   &  12.857674 &  12.857674 &   296 &    1  &  4.28e-15 &  8.84e-13 &  6.09e-14 \\ 
  13   &  12.862770 &  12.862770 &   297 &    1  &  3.45e-15 &  1.96e-12 &  7.88e-14 \\ 
  14   &  12.911115 &  12.911115 &   298 &    1  &  4.95e-15 &  7.96e-13 &  7.66e-14 \\ 
  15   &  12.934101 &  12.934101 &   299 &    1  &  3.85e-15 &  1.79e-12 &  7.64e-14 \\ 
  16   &  12.942594 &  12.942594 &   300 &    1  &  4.94e-15 &  1.45e-12 &  6.75e-14 \\ 
  17   &  12.953265 &  12.953265 &   301 &    1  &  2.74e-15 &  3.61e-12 &  7.64e-14 \\ 
  18   &  12.955215 &  12.955215 &   302 &    1  &  9.60e-16 &  3.19e-12 &  9.20e-14 \\ 
  19   &  12.990125 &  12.990125 &   303 &    1  &  5.47e-16 &  8.06e-13 &  7.39e-14 \\ 
  20   &  13.016261 &  13.016261 &   304 &    1  &  2.87e-15 &  7.82e-13 &  6.11e-14 \\ 
  21   &  13.031452 &  13.031452 &   305 &    1  &  7.09e-15 &  1.30e-12 &  8.85e-14 \\ 
  22   &  13.099484 &  13.099484 &   306 &    1  &  0.00e+00 &  1.36e-11 &  1.12e-13 \\ 
  23   &  13.122889 &  13.122889 &   307 &    1  &  6.50e-15 &  7.17e-12 &  6.10e-14 \\ 
  24   &  13.131260 &  13.131260 &   308 &    1  &  3.11e-15 &  5.21e-12 &  1.05e-13 \\ 
  25   &  13.147313 &  13.147313 &   309 &    1  &  8.24e-15 &  9.94e-12 &  1.02e-13 \\ 
  26   &  13.165265 &  13.165265 &   310 &    1  &  8.10e-16 &  2.35e-12 &  7.68e-14 \\ 
  27   &  13.187200 &  13.187200 &   311 &    1  &  2.29e-15 &  1.17e-12 &  6.51e-14 \\ 
  28   &  13.214572 &  13.214572 &   312 &    1  &  2.29e-15 &  1.73e-12 &  7.55e-14 \\ 
  29   &  13.230422 &  13.230422 &   313 &    1  &  6.85e-15 &  2.54e-12 &  8.10e-14 \\ 
  30   &  13.279307 &  13.279307 &   314 &    1  &  8.03e-16 &  7.44e-13 &  6.06e-14 \\ 
  31   &  13.287231 &  13.287231 &   315 &    1  &  6.15e-15 &  9.94e-13 &  7.53e-14 \\ 
  32   &  13.297599 &  13.297599 &   316 &    1  &  2.67e-15 &  1.30e-12 &  6.56e-14 \\ 
  33   &  13.313440 &  13.313440 &   317 &    1  &  1.33e-15 &  2.27e-12 &  8.21e-14 \\ 
  34   &  13.325555 &  13.325555 &   318 &    1  &  2.27e-15 &  2.49e-12 &  8.45e-14 \\ 
  35   &  13.371306 &  13.371306 &   319 &    1  &  2.13e-15 &  8.76e-13 &  5.91e-14 \\ 
  36   &  13.408881 &  13.408881 &   320 &    1  &  3.84e-15 &  3.74e-13 &  7.82e-14 \\ 
  37   &  13.428053 &  13.428053 &   321 &    1  &  2.65e-16 &  9.55e-13 &  6.43e-14 \\ 
  38   &  13.439632 &  13.439632 &   322 &    1  &  6.08e-15 &  6.46e-13 &  4.67e-14 \\ 
  39   &  13.457707 &  13.457707 &   323 &    1  &  2.11e-15 &  3.61e-13 &  8.02e-14 \\ 
  40   &  13.461542 &  13.461542 &   324 &    1  &  5.28e-16 &  9.63e-13 &  6.87e-14 \\ 
  41   &  13.475752 &  13.475752 &   325 &    1  &  1.58e-15 &  3.84e-13 &  4.08e-14 \\ 
  42   &  13.510229 &  13.510229 &   326 &    1  &  3.81e-15 &  1.95e-12 &  8.98e-14 \\ 
  43   &  13.517628 &  13.517628 &   327 &    1  &  4.07e-15 &  2.19e-12 &  3.97e-14 \\ 
  44   &  13.525088 &  13.525088 &   328 &    1  &  2.89e-15 &  3.54e-12 &  4.46e-14 \\ 
  45   &  13.539547 &  13.539547 &   329 &    1  &  3.28e-15 &  4.24e-12 &  2.48e-13 \\ 
  46   &  13.554331 &  13.554331 &   330 &    1  &  5.11e-15 &  1.19e-12 &  1.00e-13 \\ 
  47   &  13.570583 &  13.570583 &   331 &    1  &  2.36e-15 &  8.13e-13 &  7.49e-14 \\ 
  48   &  13.591812 &  13.591812 &   332 &    1  &  4.31e-15 &  1.48e-12 &  8.63e-14 \\ 
  49   &  13.601874 &  13.601874 &   333 &    1  &  8.88e-15 &  6.28e-13 &  9.17e-14 \\ 
  50   &  13.612912 &  13.612912 &   334 &    1  &  1.30e-14 &  1.23e-12 &  6.33e-14 \\ 
 \hline 
\end{tabular}
\end{center}
\caption{grid=rpiGride4, method=KrylovSchur, ts=implicit, order=2, $N_p=8$.}
\label{tab:rpiGride4Order2}
\end{table}
}

\newcommand{\eigenWaveSummaryTable}{
\begin{table}[hbt]
\begin{center}\tableFontSize
\begin{tabular}{|c|c|c|c|c|c|c|c|} \hline
 \multicolumn{8}{|c|}{EigenWave: grid=rpiGride4, ts=implicit, order=2, $N_p=8$, KrylovSchur } \\ \hline 
   num   &  wave       & time-steps  & wave-solves &  time-steps &   max      &  max      &  max       \\ 
   eigs  &  solves     & per period  &  per eig    &  per-eig    &   eig-err  & evect-err & eig-res    \\ 
 \hline
    51    &       192     &      10     &    3.8    &     301     &  1.30e-14 &  1.36e-11 &  2.48e-13 \\
 \hline 
\end{tabular}
\end{center}
\caption{EigenWave: grid=rpiGride4, method=KrylovSchur, ts=implicit, order=2, $N_p=8$.}
\label{tab:rpiGride4Order2Summary}
\end{table}
}

%% file: tex/doubleEllipseEigenpairs.tex
\input tex/doubleEllipseGridFig.tex

\subsection{Eigenmodes of the Penrose unilluminable room} \label{sec:doubleEllipseEigenpairs}

As a next example, we compute eigenpairs for the Penrose unilluminable room~\cite{Fukushima2015LightPI}.
The geometry, shown in Figure~\ref{fig:doubleEllipseGridFig}, is designed so that 
some of the alcoves, two at the top and two at the bottom, remain dark (or quiet) when there is a light source (or sound source) in the interior.
The design is based on two ellipses of different sizes. 
Two smaller half-ellipses, with semi-axes $(a_1,b_1)=(2,1)$, are located at the top and bottom. 
Two larger half-ellipses, with semi-axes $(a_2,b_2)=(3,6)$, are placed on the left and right.
The left and right ends of the smaller ellipses are located at the foci of the larger ellipses.
The overset grid for the domain is shown in Figure~\ref{fig:doubleEllipseGridFig} (left and middle).
The grid, denoted by $\Gcde^{(j)}$ and built with typical grid spacing $\ds^{(j)}=1/(10 j)$, consists of a total of nine component grids.  Four component grids are placed to fit the curved elliptical boundaries with four small Cartesian grids used to fit the straight portions of the boundaries in the alcoves (see middle image).  The ninth component grid is a large background Cartesian grid covering the bulk of the domain.

\input tex/doubleEllipseContoursFigII.tex

\begin{table}[hbt]
\begin{center}\tableFontSize
\begin{tabular}{|c|c|c|c|c|c|c|c|c|} \hline
 \multicolumn{9}{|c|}{EigenWave: double ellipse, ts=implicit, $\omega=11$, $N_p=6$, KrylovSchur } \\ \hline 
   order & num   &  wave       & time-steps  & wave-solves &  time-steps &   max      &  max      &  max       \\ 
   & eigs  &  solves     & per period  &  per eig    &  per-eig    &   eig-err  & evect-err & eig-res    \\ 
 \hline
    2 & 331    &       768     &      10     &    2.3    &     139     &  8.72e-14 &  7.29e-10 &  6.01e-10 \\
    4 & 324    &       768     &      10     &    2.4    &     142     &  6.95e-13 &  2.46e-10 &  1.71e-09 \\
 \hline 
\end{tabular}
\end{center}
\vspace*{-1\baselineskip}
\caption{Summary of EigenWave performance for double ellipse grid $\Gcde^{(4)}$ using the KrylovSchur algorithm and implicit time-stepping.  The spatial order of accuracy is 2 for the top row and 4 for the bottom row, and the wave-solves use $N_p=6$ to determine the final time.}
\label{tab:darkCornerRoomGride4Summary}
\end{table}

In this example, a large number of eigenpairs are desired for the purpose of investigating different normal modes of the room.
Table~\ref{tab:darkCornerRoomGride4Summary} summarizes results of computing eigenpairs to second and fourth-order accuracy on grid~$\Gcde^{(4)}$.
Implicit time-stepping is used with $\Nits=10$ time-steps per period and $\Np=6$ periods per wave-solve.
The target frequency is $\omega=11$ and $256$ eigenpairs are requested. For the second-order accurate discretization $331$ eigenpairs are found accurately, while for the fourth-order accurate discretization $324$ eigenpairs are found accurately. This corresponds to approximately 
$2.3$  and $2.4$ wave-solves per computed eigenpair, respectively.
Figure~\ref{fig:doubleEllipseEigenVectorsII} shows contours of selected eigenvectors.
The eigenvectors shown in the top left and middle correspond to the smallest and second smallest eigenvalues.
The eigenvector shown on the bottom left of the figure is a surface mode that is primarily restricted 
to the neighborhood of the left boundary.
The eigenvector shown on the top right is a mode that is restricted to the central portion of the domain while the remaining two eigenvectors show interesting modes
in the central region of the geometry.
For this problem with approximately $210,000$ grid points the 
implicit scheme with a direct sparse solver was about $1.8$ (or $1.7$) times faster then the explicit scheme at second (or fourth) order accuuracy.

%% file: tex/doubleEllipseGridFig.tex
{
\newcommand{\figw}{6.25cm}
\newcommand{\figh}{6cm}
\begin{figure}[htb]
\begin{center}
\begin{tikzpicture}
  \useasboundingbox (0,.5) rectangle (16,\figh);  

  \figByWidth{0.0}{-6pt}{fig/darkCornerRoomGridG2}{4.7cm}[0.1][0.1][0.][0.]
  \draw[thick,black,->,yshift=0pt] (0,0) node[anchor=north] {\scriptsize $-5$} -- (4.65,0.00) node[anchor=north] {\scriptsize $+5$};  
  \draw[thick,black,->,xshift=1pt] (0,0) node[anchor=east,yshift=4pt]  {\scriptsize $-6$} -- (0.00,5.47) node[anchor=east ] {\scriptsize $+6$};  

  \figByWidthb{4.95}{.75}{fig/darkCornerRoomGridG2Zoom}{4.5cm}[0.][0.][0.][0.15]

  \begin{scope}[scale=0.46,xshift=7.02cm,yshift=5.97cm]
    \draw (0,0) ellipse (3cm and 6cm);    
    \draw[black,fill=black] (0,5.196cm) circle (3pt); 
  \end{scope}
  \begin{scope}[scale=0.46,xshift=5cm,yshift=11.2cm]
    \draw (0,0) ellipse (2cm and 1cm);    
  \end{scope}

  \begin{scope}[xshift=9.5cm,yshift=.25cm]
     \figByWidth{0}{0}{fig/darkCornerRoomG4O4EigsKrylov}{\figw}[0][0][0][0]
     \draw (1.5,1) node[draw,fill=white,anchor=west,xshift=0pt,yshift=0pt] {\scriptsize Double-Ellipse, I, O4, $\Np=6$};
  \end{scope} 

\end{tikzpicture}
\end{center}
\caption{At left is double ellipse geometry and overset grid $\Gcde^{(2)}$, and at center is a closeup of a portion of the grid. At right is the filter function and computed eigenvalues.
 }
\label{fig:doubleEllipseGridFig}
\end{figure}
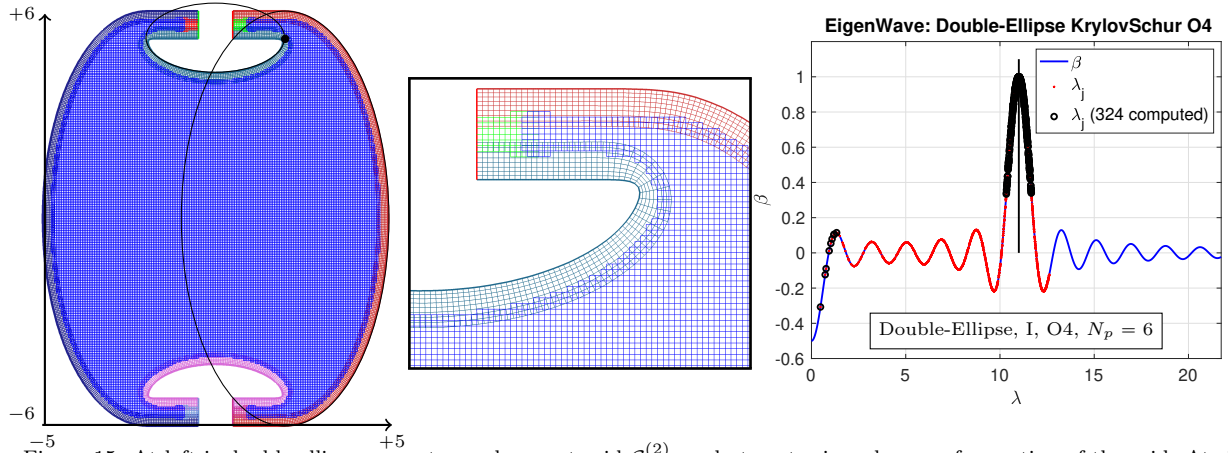
}

%% file: tex/doubleEllipseContoursFigII.tex
{
\newcommand{\drawContour}[7]{%
\begin{scope}[#1]
\draw(0.0,0) node[anchor=south west,xshift=-4pt,yshift=+0pt] {\trimfiga{fig/#2}{\figWidtha}};
  \draw(.35,.35) node[draw,fill=white,anchor=west,xshift=2pt,yshift=2pt] {\scriptsize #3};
\end{scope}
}
\newcommand{\drawContourWithColourBar}[7]{%
\begin{scope}[#1]
   \draw(0.0,0) node[anchor=south west,xshift=-4pt,yshift=+0pt] {\trimfiga{fig/#2}{\figWidtha}};
  \draw(.35,.35) node[draw,fill=white,anchor=west,xshift=2pt,yshift=2pt] {\scriptsize #3};
\begin{scope}[xshift=-.5cm,yshift=-5pt]
  \draw (\xcb,\ycb) node[anchor=south west,xshift=0.25cm,yshift=.5cm,rotate=-90] {\trimfigcb{fig/colourBarLines}{\cbWidth}{\cbHeight}};
  \draw (.8,0) node[anchor=north,xshift=+3pt,yshift=+2pt] {\scriptsize $#6$};
  \draw (4.8,0) node[anchor=north,xshift=+0pt,yshift=+2pt] {\scriptsize $#7$};
\end{scope}
\end{scope}
}
\newcommand{\cbWidth}{.2cm}
\newcommand{\cbHeight}{4cm}
\newcommand{\xcb}{.5cm}
\newcommand{\ycb}{-.2cm}
\setlength{\ycbTop}{\ycb+\cbHeight}
\setlength{\ycbMid}{\ycb+\cbHeight*\real{.5}}
\newcommand{\trimfigcb}[3]{\includegraphics[width=#2, height=#3, clip, trim=17cm 2.35cm 1.65cm 2.35cm]{#1}}
\newcommand{\figWidtha}{4cm}
\newcommand{\trimfiga}[2]{\trimw{#1}{#2}{.175}{.175}{.1}{.12}}
\begin{figure}[htb]
\begin{center}
\begin{tikzpicture}
   \useasboundingbox (0,.25) rectangle (13.5,10);  

   \begin{scope}[yshift=5cm]
    \drawContour{xshift= 0.cm,yshift=0.00cm}{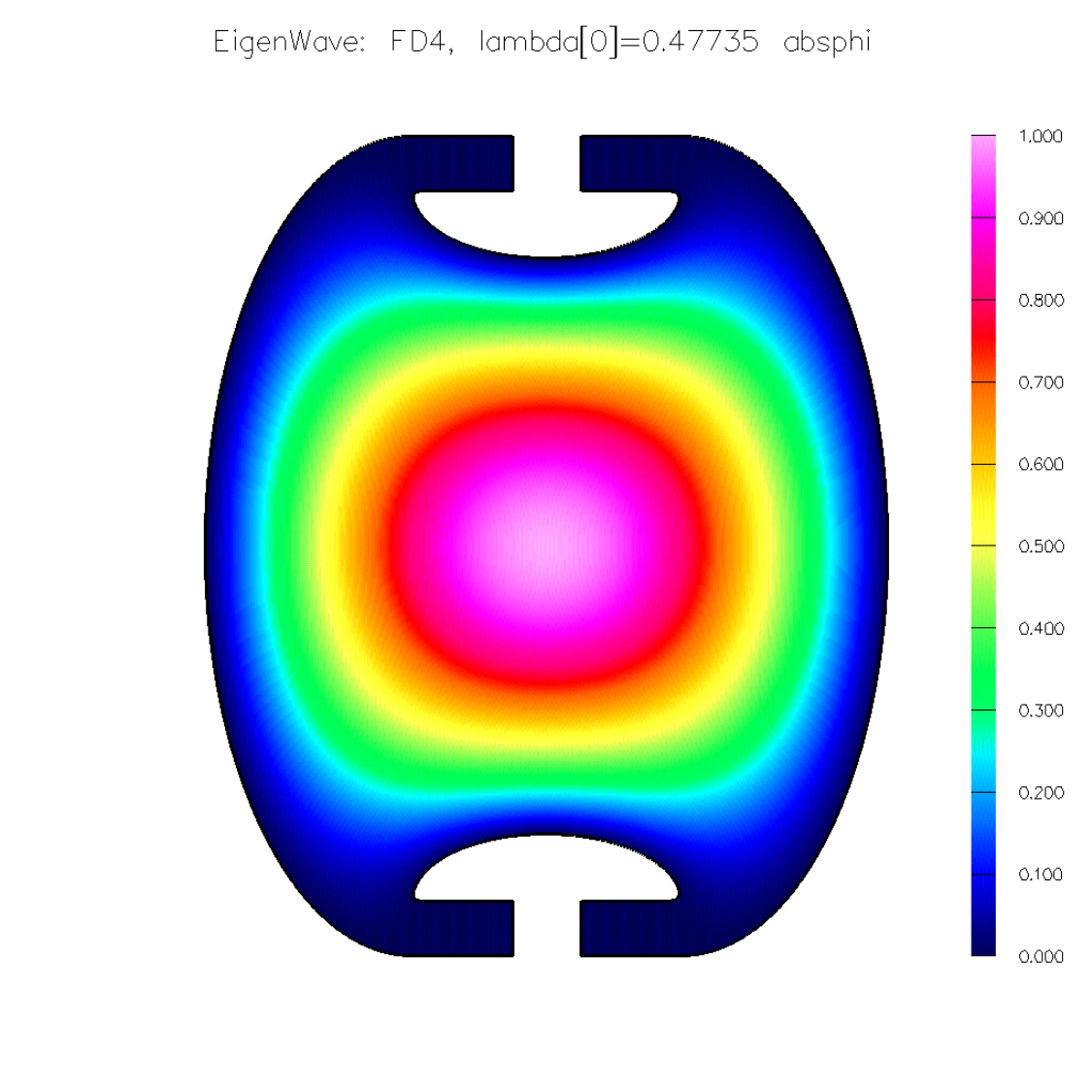}{$|\phi|$}{$v$}{$t=0.3$}{$-1.0$}{$1.0$}     
    \drawContour{xshift=4.5cm,yshift=0.00cm}{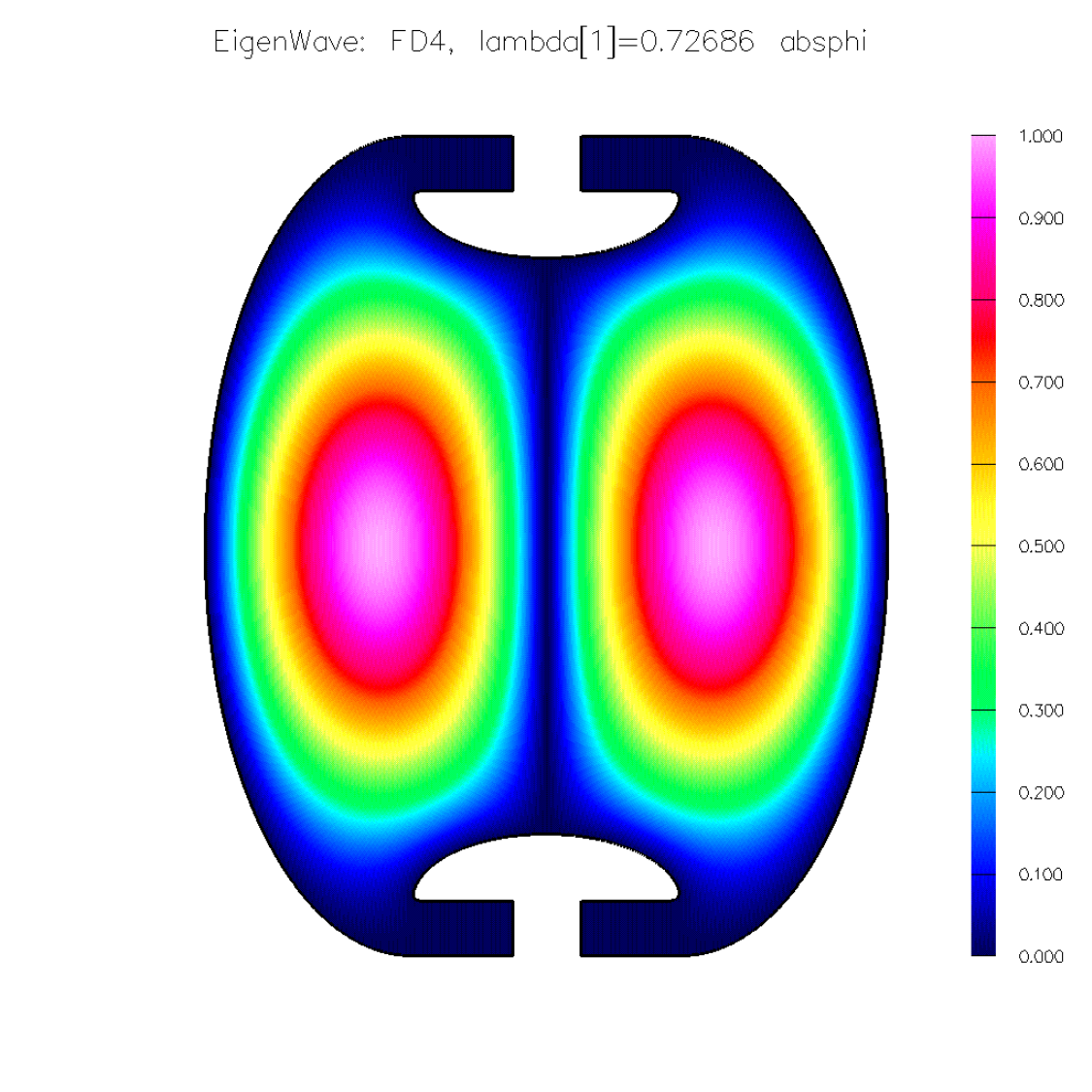}{$|\phi|$}{$v$}{$t=0.3$}{$-1.0$}{$1.0$}    
    \drawContour{xshift=9.0cm,yshift=0.00cm}{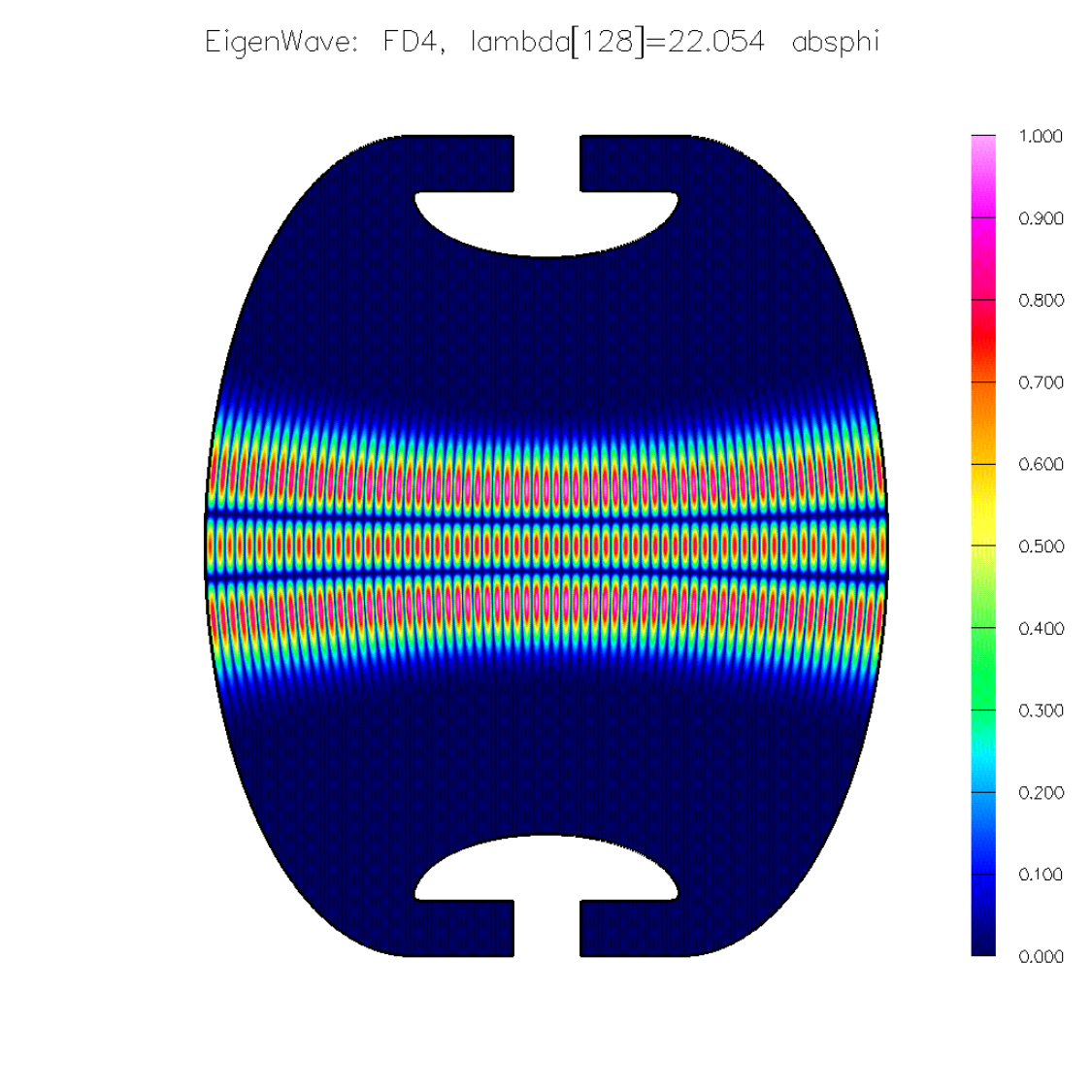}{$|\phi|$}{$v$}{$t=0.3$}{$-1.0$}{$1.0$}     
   \end{scope}
   \begin{scope}[yshift=0cm]
    \drawContour{xshift= 0.cm,yshift=0.00cm}{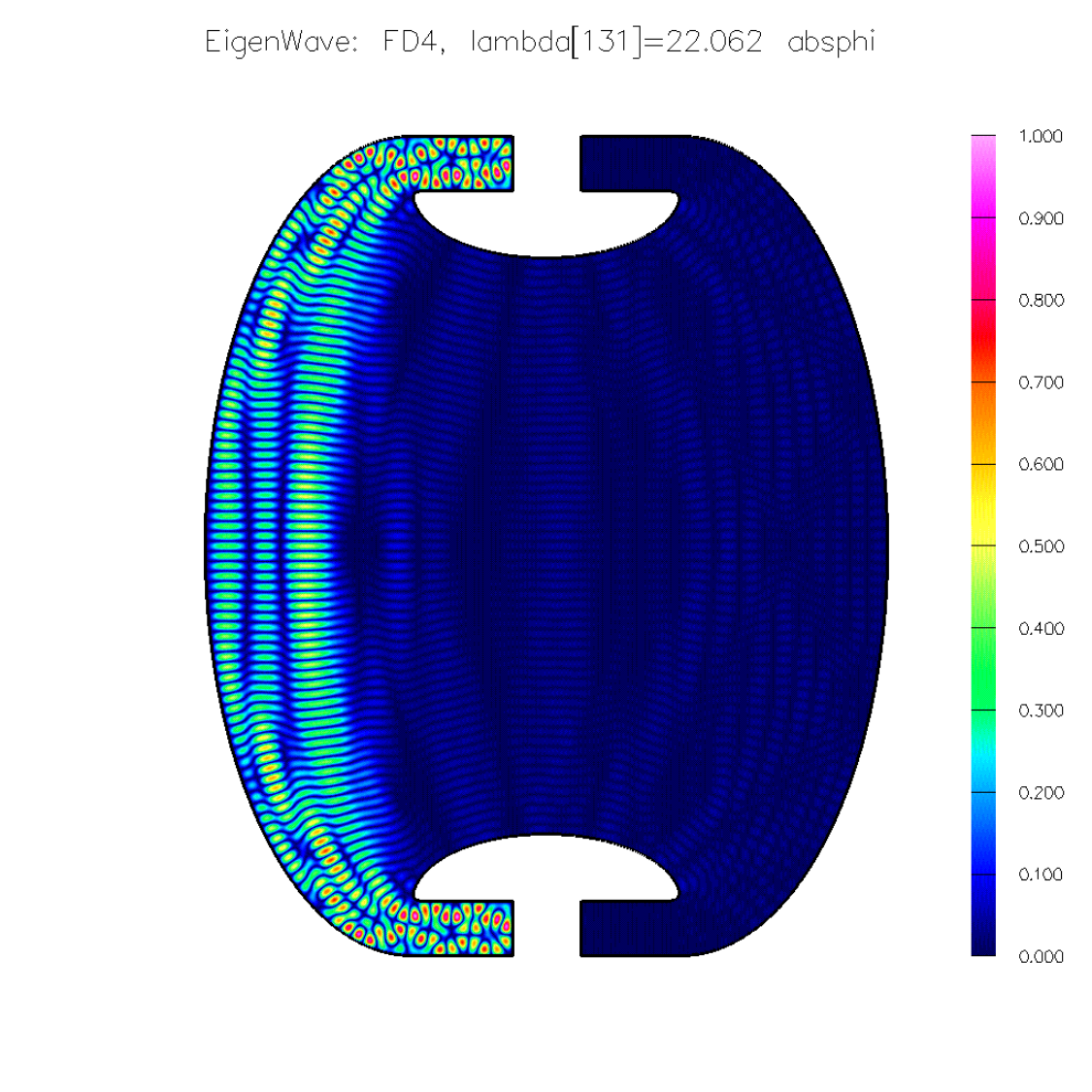}{$|\phi|$}{$v$}{$t=0.3$}{$-1.0$}{$1.0$}    
    \drawContourWithColourBar{xshift=4.5cm,yshift=0.00cm}{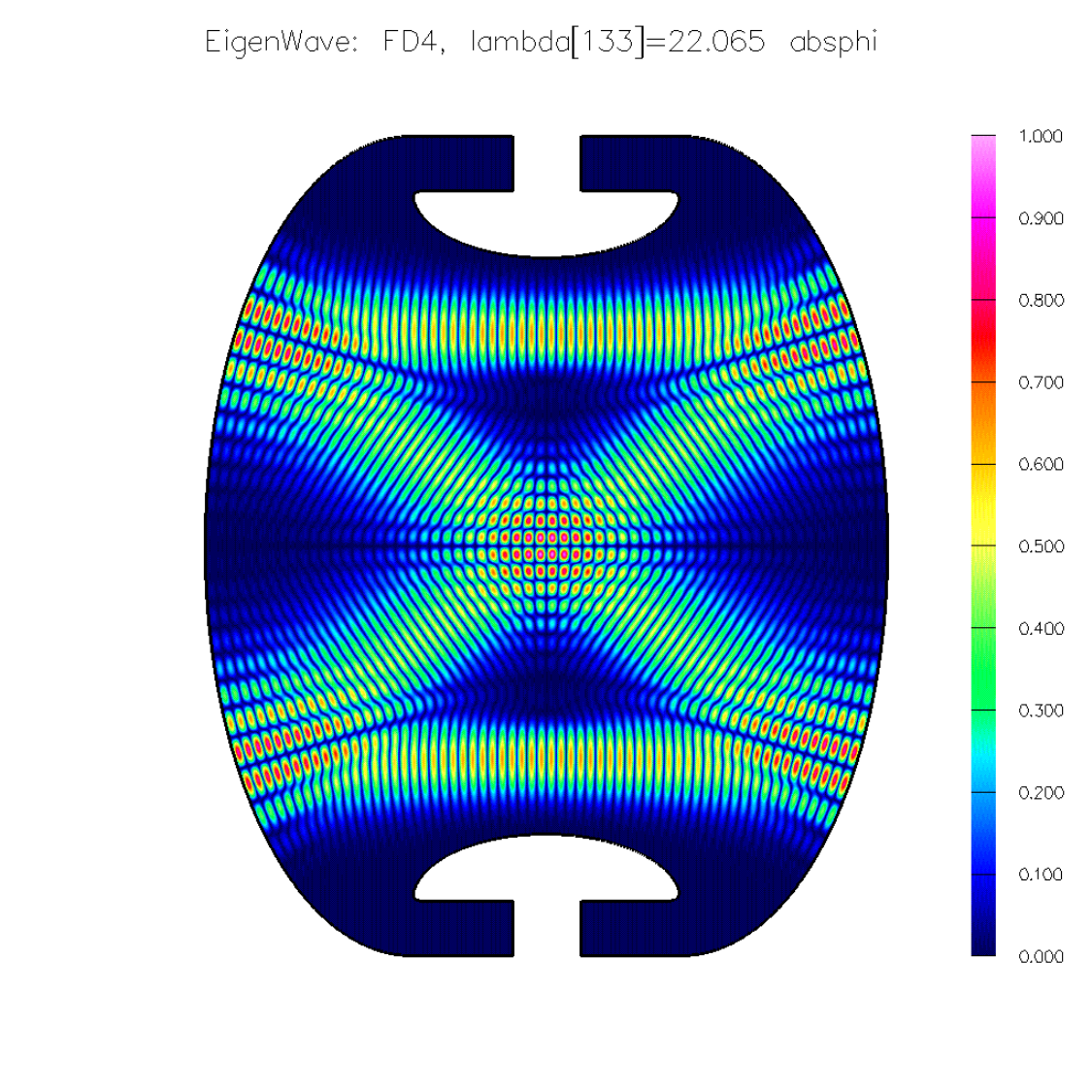}{$|\phi|$}{$v$}{$t=0.3$}{$0.0$}{$1.0$}    
    \drawContour{xshift=9.0cm,yshift=0.00cm}{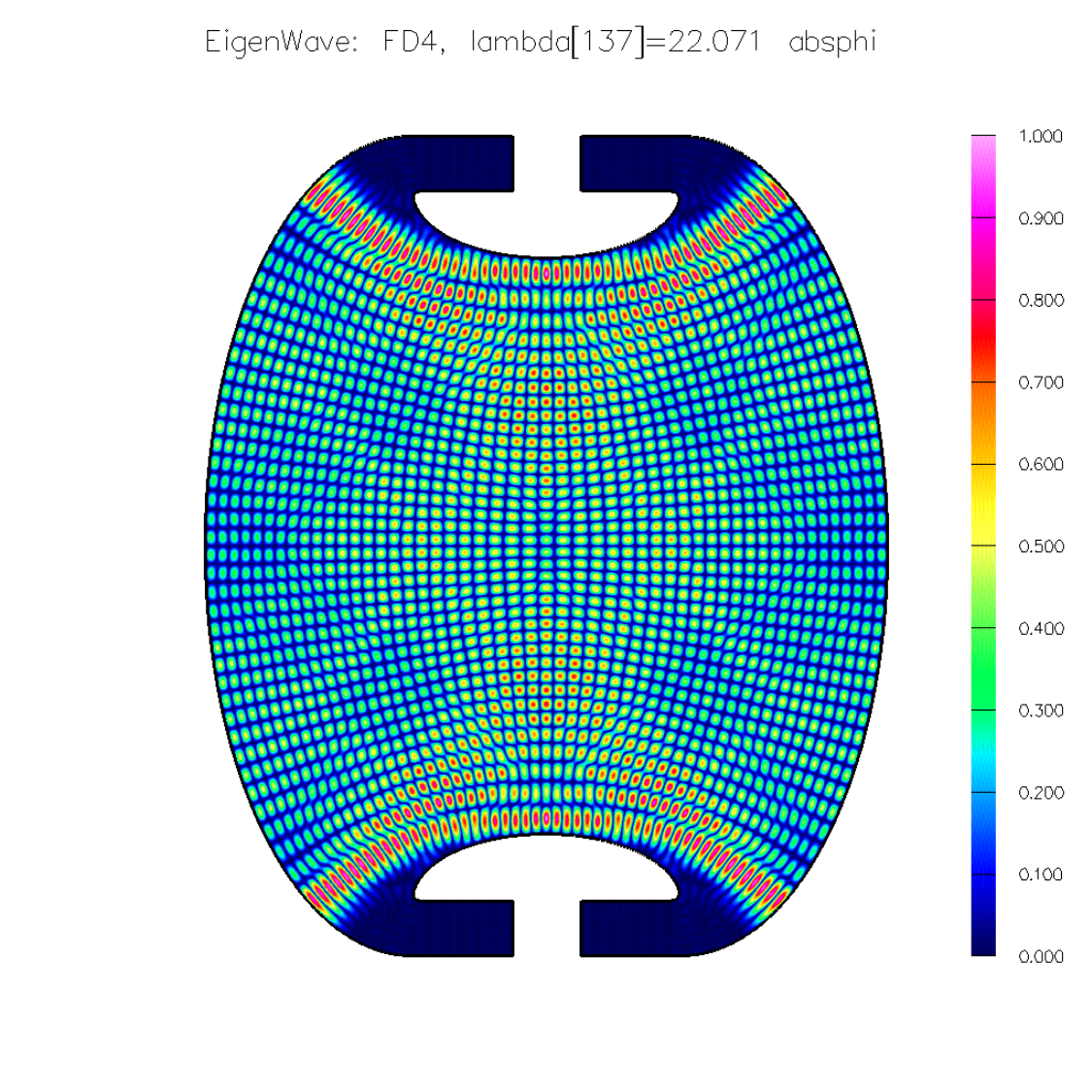}{$|\phi|$}{$v$}{$t=0.3$}{$-1.0$}{$1.0$}    
   \end{scope}  

\end{tikzpicture}
\end{center}
\caption{Selected eigenvectors of the double ellipse geometry on $\Gcde^{(8)}$ using the $4th$-order accurate discretization.
    }
\label{fig:doubleEllipseEigenVectorsII}
\end{figure}
}

%% file: tex/boxEigenpairs.tex
\newcommand{\Gcb}{\Gc_{\rm box}}
\subsection{Eigenmodes of a box} \label{sec:boxEigenpairs}

In this section eigenpairs of a three-dimensional unit box, $\Omega=[0,1]^3$, are computed.
The grid for the box, denoted by $\Gcb^{(j)}$,
consists of a single Cartesian grid with grid spacing of $\ds^{(j)}=1/(10 j)$.
Note that a box with all equal length sides has many eigenvalues of multiplicity 3 and 6.
The EigenWave algorithm is still able to accurately compute these high-multiplicity eigenpairs.

{
\newcommand{\figw}{6.5cm}
\newcommand{\figh}{5.5cm}

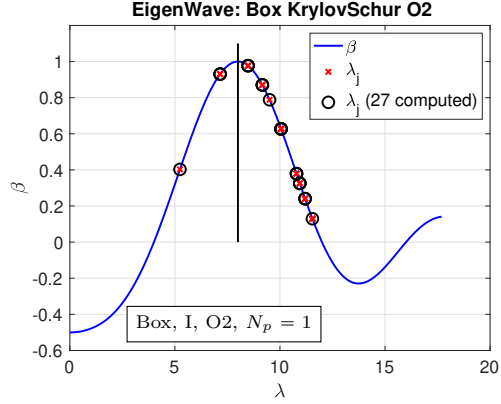
\begin{figure}[htb]
\begin{center}
\begin{tikzpicture}[scale=1]
  \useasboundingbox (0,.7) rectangle (7,5.5);  

  \begin{scope}[yshift=0cm]
     \figByWidth{0.0}{0}{fig/box2O2ImpEigKrylov}{\figw}[0][0][0][0]
     \draw (1.5,1) node[draw,fill=white,anchor=west,xshift=0pt,yshift=0pt] {\scriptsize Box, I, O2, $\Np=1$};
  \end{scope}

\end{tikzpicture}
\end{center}
\caption{Box: computing multiple eigenpairs, order=2. The computed eigenvalues are marked with black circles.
    } 
\label{fig:boxEigs}
\end{figure}
}

{

\input tables/box2O2ImpEigKrylov.tex
  \eigenWaveSummaryTable
}

Table~\ref{tab:box2Order2Summary} summarizes results of computing eigenpairs to second-order accuracy on grid $\Gcb^{(2)}$.
The target frequency was $\omega=8$. Twenty eigenpairs were requested and twenty-seven 
eigenpairs were accurately found by the KrylovSchur
algorithm in a total of $139$ wave-solves. This corresponds to approximately 
$5.1$ wave-solves per eigenpair found.
The graph in Figure~\ref{fig:boxEigs} shows the filter-function $\beta(\lambda;\omega)$ with the computed
eigenvalues marked.

%% file: tables/box2O2ImpEigKrylov.tex
\newcommand{\eigenWaveLongTable}{
\begin{table}[hbt]
\begin{center}\tableFontSize
\begin{tabular}{|c|r|r|c|c|c|c|c|} \hline
 \multicolumn{8}{|c|}{EigenWave: grid=box2, ts=implicit, order=2, $N_p=1$, KrylovSchur } \\ \hline 
   $j$  & \multicolumn{1}{c|}{$\lambda_j$} &  \multicolumn{1}{c|}{$\lambda_{\rm true}$}  & eig &  mult &  eig-err & evect-err & eig-res  \\ \hline
   0   &   5.435806 &   5.435806 &     0 &    1  &  3.92e-15 &  1.79e-14 &  2.38e-13 \\ 
   1   &   7.671600 &   7.671600 &     1 &    3  &  6.95e-16 &  3.74e-14 &  8.58e-14 \\ 
   2   &   7.671600 &   7.671600 &     1 &    3  &  1.16e-16 &  4.29e-14 &  8.28e-14 \\ 
   3   &   7.671600 &   7.671600 &     3 &    3  &  1.16e-16 &  8.86e-14 &  8.81e-14 \\ 
   4   &   9.389297 &   9.389297 &     4 &    3  &  2.65e-15 &  4.95e-14 &  7.59e-14 \\ 
   5   &   9.389297 &   9.389297 &     4 &    3  &  2.08e-15 &  3.96e-14 &  5.61e-14 \\ 
   6   &   9.389297 &   9.389297 &     4 &    3  &  9.46e-16 &  6.20e-14 &  6.77e-14 \\ 
   7   &  10.338928 &  10.338928 &     9 &    3  &  1.55e-15 &  2.75e-14 &  3.76e-14 \\ 
   8   &  10.338928 &  10.338928 &     9 &    3  &  3.09e-15 &  1.01e-14 &  3.52e-14 \\ 
   9   &  10.338928 &  10.338928 &     9 &    3  &  3.44e-15 &  3.24e-14 &  3.79e-14 \\ 
  10   &  10.838098 &  10.838098 &    10 &    1  &  5.24e-15 &  1.91e-14 &  3.88e-14 \\ 
  11   &  11.670428 &  11.670428 &    11 &    6  &  1.22e-15 &  1.47e-14 &  3.93e-14 \\ 
  12   &  11.670428 &  11.670428 &    11 &    6  &  1.83e-15 &  9.10e-15 &  2.68e-14 \\ 
  13   &  11.670428 &  11.670428 &    13 &    6  &  9.13e-16 &  1.54e-14 &  3.67e-14 \\ 
  14   &  11.670428 &  11.670428 &    13 &    6  &  3.04e-16 &  1.06e-14 &  4.12e-14 \\ 
  15   &  11.670428 &  11.670428 &    15 &    6  &  6.09e-16 &  1.33e-14 &  2.71e-14 \\ 
  16   &  11.670428 &  11.670428 &    15 &    6  &  9.13e-16 &  1.45e-14 &  4.34e-14 \\ 
  17   &  12.864850 &  12.864850 &    19 &    3  &  1.38e-15 &  2.28e-14 &  3.27e-14 \\ 
  18   &  12.864850 &  12.864850 &    19 &    3  &  1.66e-15 &  1.79e-14 &  2.91e-14 \\ 
  19   &  12.864850 &  12.864850 &    19 &    3  &  5.94e-15 &  2.38e-14 &  2.52e-14 \\ 
  20   &  13.133357 &  13.133357 &    20 &    3  &  1.89e-15 &  3.86e-14 &  5.77e-14 \\ 
  21   &  13.133357 &  13.133357 &    21 &    3  &  1.49e-15 &  3.40e-14 &  4.87e-14 \\ 
  22   &  13.133357 &  13.133357 &    21 &    3  &  5.41e-16 &  6.30e-14 &  3.38e-14 \\ 
  23   &  13.573463 &  13.573463 &    23 &    3  &  3.01e-15 &  4.31e-14 &  3.96e-14 \\ 
  24   &  13.573463 &  13.573463 &    25 &    3  &  2.62e-16 &  2.96e-14 &  5.17e-14 \\ 
  25   &  13.573463 &  13.573463 &    25 &    3  &  2.62e-15 &  2.04e-14 &  3.42e-14 \\ 
  26   &  14.205299 &  14.205299 &    27 &    6  &  3.75e-16 &  1.34e-13 &  1.65e-13 \\ 
 \hline 
\end{tabular}
\end{center}
\caption{grid=box2, method=KrylovSchur, ts=implicit, order=2, $N_p=1$.}
\label{tab:box2Order2}
\end{table}
}

\newcommand{\eigenWaveSummaryTable}{
\begin{table}[hbt]
\begin{center}\tableFontSize
\begin{tabular}{|c|c|c|c|c|c|c|c|} \hline
 \multicolumn{8}{|c|}{EigenWave: grid=box2, ts=implicit, order=2, $N_p=1$, KrylovSchur } \\ \hline 
   num   &  wave       & time-steps  & wave-solves &  time-steps &   max      &  max      &  max       \\ 
   eigs  &  solves     & per period  &  per eig    &  per-eig    &   eig-err  & evect-err & eig-res    \\ 
 \hline
    27    &       139     &      10     &    5.1    &     51     &  5.94e-15 &  1.34e-13 &  2.38e-13 \\
 \hline 
\end{tabular}
\end{center}
\caption{EigenWave: grid=box2, method=KrylovSchur, ts=implicit, order=2, $N_p=1$.}
\label{tab:box2Order2Summary}
\end{table}
}

%% file: tex/pipeEigenpairs.tex
\newcommand{\Gcp}{\Gc_{\rm pipe}}
\subsection{Eigenmodes of a pipe in three dimensions} \label{sec:pipeEigenpairs}

In this section eigenpairs of a three-dimensional cylindrical solid pipe are computed.
The pipe has a radius of $0.5$ and an axial length of $1$.
The composite grid for the solid cylinder, denoted by $\Gcp^{(j)}$, consists of two component grids,
each with grid spacings approximately equal to $\ds^{(j)}=1/(10 j)$ in all directions.
One component grid is a narrow boundary-fitted cylindrical shell, while the other component grid is a background Cartesian grid covering the interior of the cylindrical domain (see Figure~\ref{fig:pipeFig}).

\input tex/pipeGridFig.tex

\input tex/pipeContoursFig

{

\input tables/pipeG2O2ImpEigKrylov.tex
\eigenWaveSummaryTable
}
Table~\ref{tab:pipeze2Order2Summary} summarizes results of computing eigenpairs to second-order accuracy on grid $\Gcp^{(2)}$.
The target frequency was $\omega=10$. Thirty-two eigenpairs were requested and thirty-five 
eigenpairs were accurately found by the KrylovSchur
algorithm in a total of $124$ wave-solves. This corresponds to approximately 
$3.5$ wave-solves per eigenpair found.
The middle graph in Figure~\ref{fig:pipeFig} shows the filter-function $\beta(\lambda;\omega)$ with the computed
eigenvalues marked. 

{

\input tables/pipeG4O2ImpEigKrylov.tex
\eigenWaveSummaryTable
}
Table~\ref{tab:pipeze4Order2Summary} summarizes results of computing eigenpairs to second-order accuracy on a finer grid $\Gcp^{(4)}$.
The target frequency was $\omega=14$. Eighty-four eigenpairs were requested and $101$ 
eigenpairs were accurately found by the KrylovSchur
algorithm in a total of $321$ wave-solves. This corresponds to approximately 
$3.2$ wave-solves per eigenpair found.
The right graph in Figure~\ref{fig:pipeFig} shows the filter-function $\beta(\lambda;\omega)$ with the computed
eigenvalues marked.

%% file: tex/pipeGridFig.tex
{
\newcommand{\figw}{5.8cm}
\newcommand{\figh}{5.5cm}
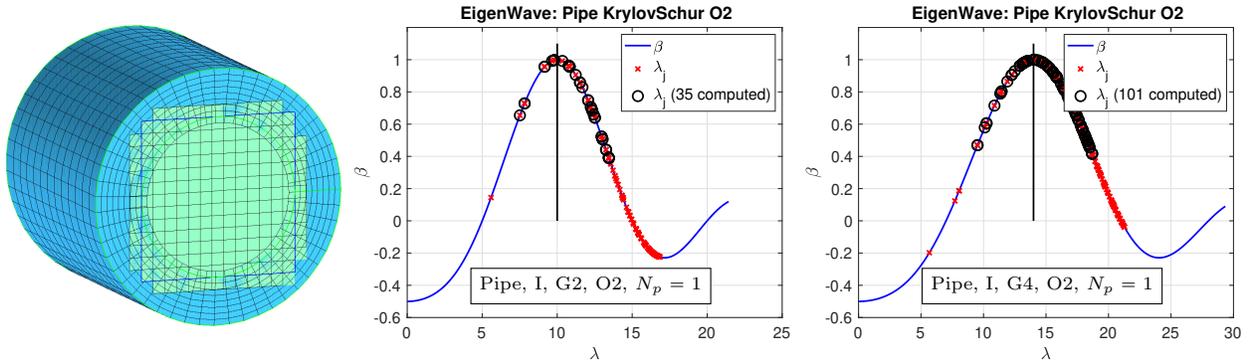
\begin{figure}[htb]
\begin{center}
\begin{tikzpicture}
  \useasboundingbox (0,.25) rectangle (17,5.25);  

  \begin{scope}[xshift=-4pt]
    \figByWidth{0.0}{0}{fig/pipeGridG2}{4.5cm}[0.025][0.015][0.065][0.065]
  \end{scope} 

  \begin{scope}[xshift=4.5cm,yshift=-12pt]
    \figByWidth{0}{0}{fig/pipeG2O2ImpEigKrylov}{\figw}[0][0][0][0]
    \draw (1.5,1) node[draw,fill=white,anchor=west,xshift=0pt,yshift=0pt] {\scriptsize Pipe, I, G2, O2, $\Np=1$};
  \end{scope} 

  \begin{scope}[xshift=10.5cm,yshift=-12pt]
    \figByWidth{0}{0}{fig/pipeG4O2ImpEigKrylov}{\figw}[0][0][0][0]
    \draw (1.5,1) node[draw,fill=white,anchor=west,xshift=0pt,yshift=0pt] {\scriptsize Pipe, I, G4, O2, $\Np=1$};
  \end{scope}    
  
\end{tikzpicture}
\end{center}
\caption{Left: overset grid $\Gcd^{(2)}$ for a cylindrical pipe. 
Middle and right: graphs of the filter function $\beta$ with the  
computed eigenvalues marked with black circles for second-order accurate computations 
on grid $\Gcp^{(2)}$ (middle) and on the finer grid $\Gcp^{(4)}$. 
 }
\label{fig:pipeFig}
\end{figure}
}

%% file: tex/pipeContoursFig.tex
{
\newcommand{\drawContour}[7]{%
\begin{scope}[#1]
\draw(0.0,0) node[anchor=south west,xshift=-4pt,yshift=+0pt] {\trimfiga{fig/#2}{\figWidtha}};
  \draw(.5,.5) node[draw,fill=white,anchor=west,xshift=2pt,yshift=1pt] {\scriptsize #3};
\begin{scope}[xshift=-.2cm,yshift=-2pt]
  \draw (\xcb,\ycb) node[anchor=south west,xshift=0.25cm,yshift=.5cm,rotate=-90] {\trimfigcb{fig/colourBarLines}{\cbWidth}{\cbHeight}};
  \draw (.8,0) node[anchor=north,xshift=+3pt,yshift=+2pt] {\scriptsize $#6$};
  \draw (4.8,0) node[anchor=north,xshift=+0pt,yshift=+2pt] {\scriptsize $#7$};
\end{scope}
\end{scope}
}
\newcommand{\cbWidth}{.2cm}
\newcommand{\cbHeight}{4cm}
\newcommand{\xcb}{.5cm}
\newcommand{\ycb}{-.2cm}
\setlength{\ycbTop}{\ycb+\cbHeight}
\setlength{\ycbMid}{\ycb+\cbHeight*\real{.5}}
\newcommand{\trimfigcb}[3]{\includegraphics[width=#2, height=#3, clip, trim=17cm 2.35cm 1.65cm 2.35cm]{#1}}
\newcommand{\figWidtha}{4.75cm}
\newcommand{\trimfiga}[2]{\trimw{#1}{#2}{.0}{.117}{.09}{.09}}
\begin{figure}[htb]
\begin{center}
\begin{tikzpicture}
   \useasboundingbox (0,.3) rectangle (15,4.5);  

   \begin{scope}[yshift=0cm]
    \drawContour{xshift= 0.cm,yshift=0.00cm}{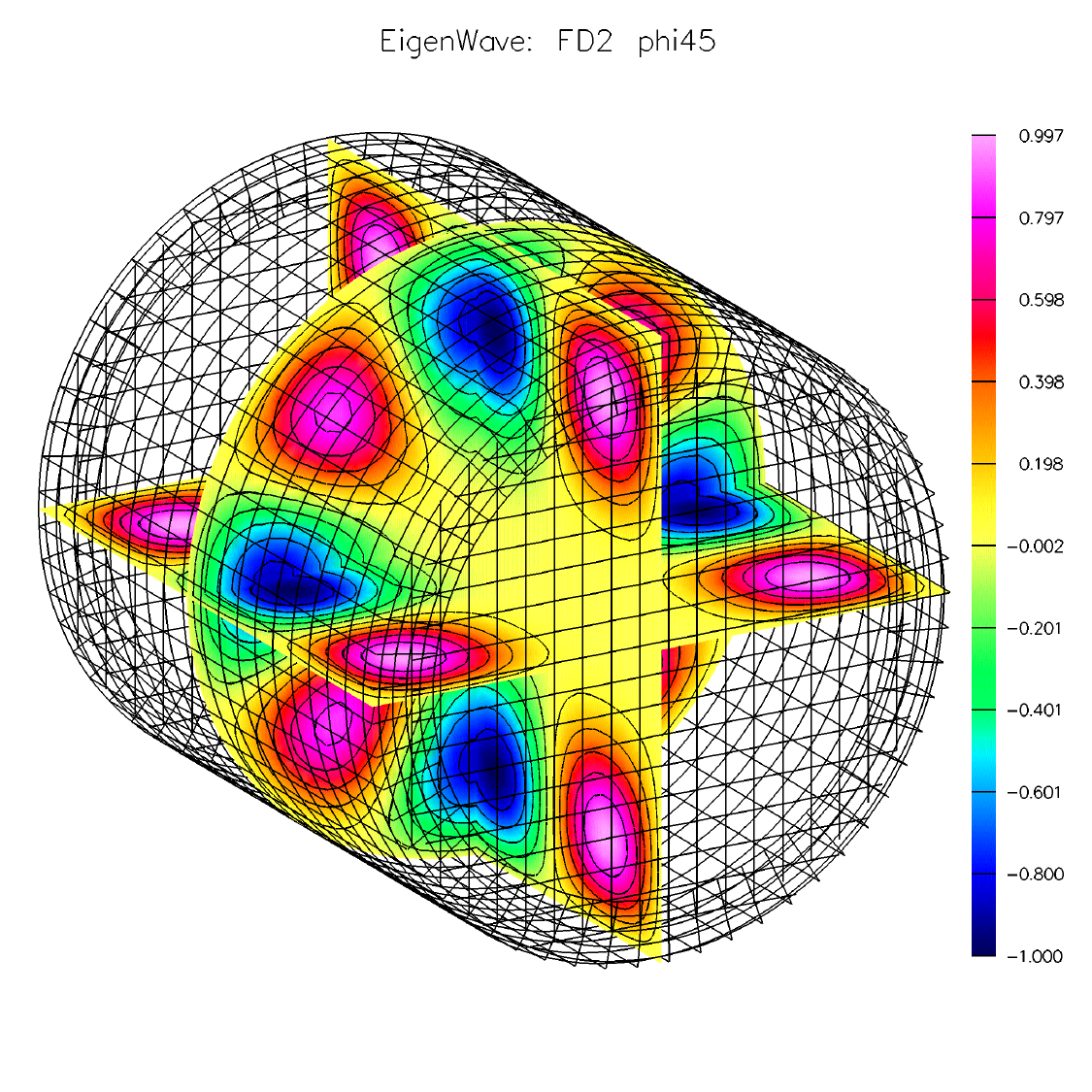}{$\phi^{(45)}$}{$v$}{$t=0.3$}{$-1.0$}{$1.0$}    
    \drawContour{xshift= 5.cm,yshift=0.00cm}{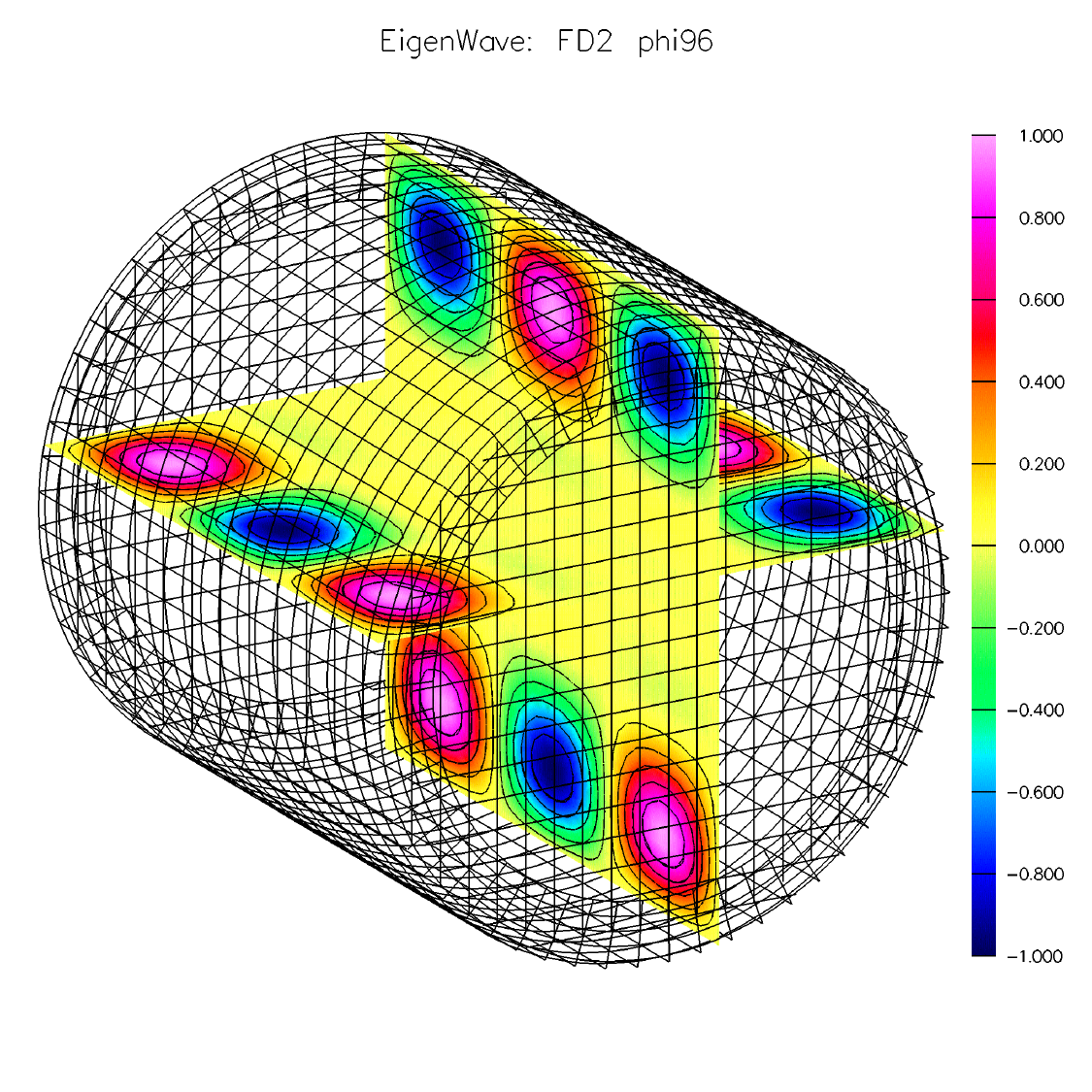}{$\phi^{(74)}$}{$v$}{$t=0.3$}{$-1.0$}{$1.0$}   
    \drawContour{xshift=10.cm,yshift=0.00cm}{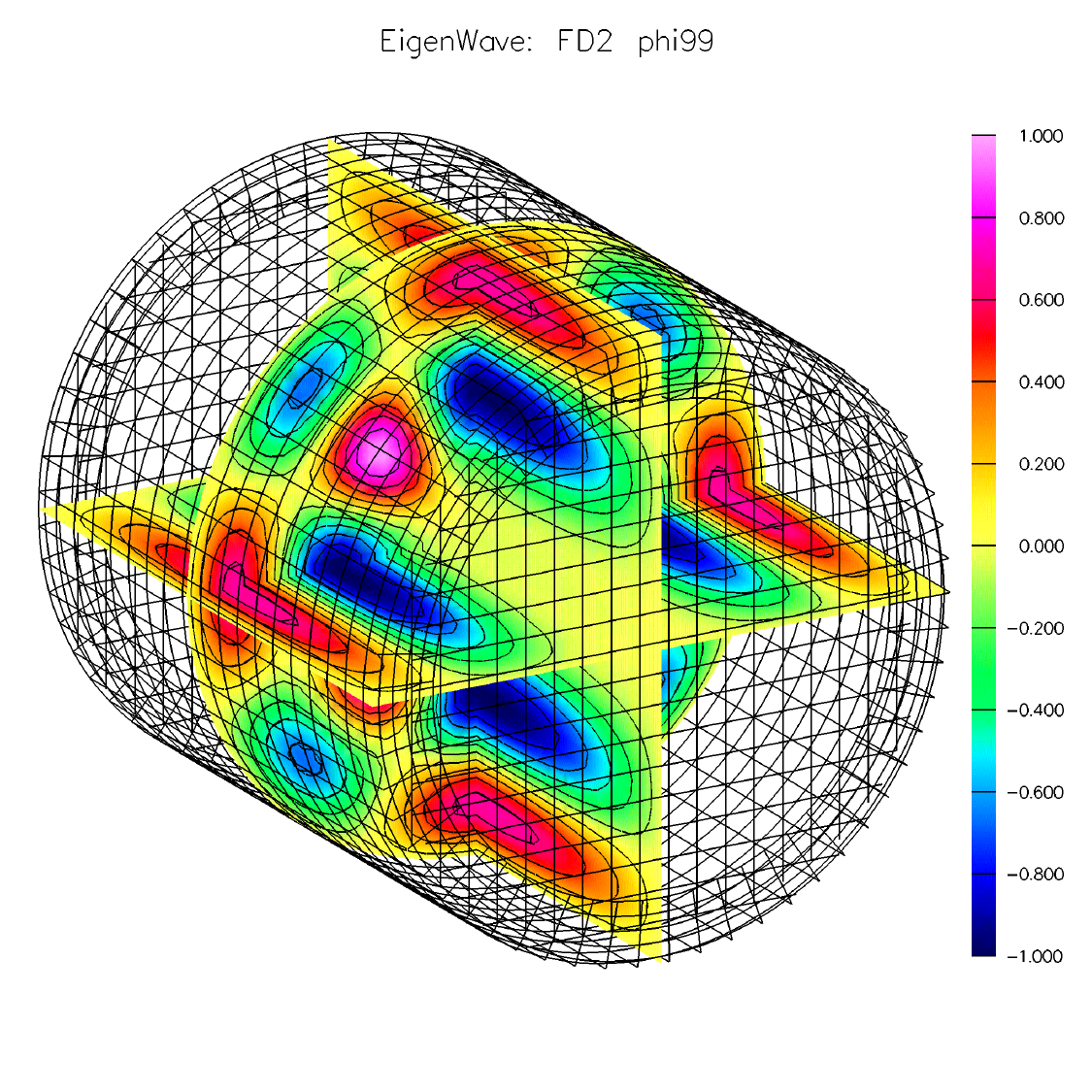}{$\phi^{(99)}$}{$v$}{$t=0.3$}{$-1.0$}{$1.0$}    
   \end{scope}

\end{tikzpicture}
\end{center}
\caption{Pipe: some computed eigenvectors. 
    }
\label{fig:pipeEigenVectors}
\end{figure}
}

%% file: tables/pipeG2O2ImpEigKrylov.tex
\newcommand{\eigenWaveLongTable}{
\begin{table}[hbt]
\begin{center}\tableFontSize
\begin{tabular}{|c|r|r|c|c|c|c|c|} \hline
 \multicolumn{8}{|c|}{EigenWave: grid=pipeze2, ts=implicit, order=2, $N_p=1$, KrylovSchur } \\ \hline 
   $j$  & \multicolumn{1}{c|}{$\lambda_j$} &  \multicolumn{1}{c|}{$\lambda_{\rm true}$}  & eig &  mult &  eig-err & evect-err & eig-res  \\ \hline
   0   &   7.885829 &   7.885829 &     1 &    1  &  4.28e-15 &  2.96e-12 &  6.39e-13 \\ 
   1   &   8.250668 &   8.250668 &     2 &    1  &  1.85e-13 &  7.78e-10 &  3.42e-13 \\ 
   2   &   8.250872 &   8.250872 &     3 &    1  &  1.15e-13 &  2.60e-09 &  2.38e-13 \\ 
   3   &   9.866775 &   9.866775 &     4 &    1  &  2.11e-14 &  4.29e-09 &  2.02e-13 \\ 
   4   &   9.866945 &   9.866945 &     5 &    1  &  3.56e-14 &  7.67e-09 &  2.70e-13 \\ 
   5   &  10.496650 &  10.496650 &     6 &    1  &  6.89e-14 &  3.84e-12 &  1.74e-13 \\ 
   6   &  10.662626 &  10.662626 &     7 &    1  &  3.42e-14 &  5.20e-11 &  1.54e-13 \\ 
   7   &  10.682982 &  10.682982 &     8 &    1  &  3.53e-14 &  7.51e-12 &  1.33e-13 \\ 
   8   &  11.395887 &  11.395887 &     9 &    1  &  1.93e-14 &  2.80e-12 &  1.02e-13 \\ 
   9   &  11.956418 &  11.956418 &    10 &    1  &  5.48e-14 &  1.22e-11 &  8.10e-14 \\ 
  10   &  11.974454 &  11.974454 &    11 &    1  &  3.71e-14 &  5.58e-11 &  2.18e-13 \\ 
  11   &  12.053323 &  12.053323 &    12 &    1  &  2.74e-14 &  4.04e-09 &  9.55e-14 \\ 
  12   &  12.053462 &  12.053462 &    13 &    1  &  7.53e-14 &  1.18e-09 &  1.10e-13 \\ 
  13   &  12.614473 &  12.614473 &    14 &    1  &  4.06e-14 &  1.96e-12 &  4.55e-14 \\ 
  14   &  13.020630 &  13.020630 &    15 &    1  &  5.06e-14 &  1.40e-10 &  9.09e-14 \\ 
  15   &  13.023784 &  13.023784 &    16 &    1  &  2.13e-14 &  4.90e-10 &  8.18e-14 \\ 
  16   &  13.250980 &  13.250980 &    17 &    1  &  4.22e-14 &  6.67e-12 &  6.38e-14 \\ 
  17   &  13.813737 &  13.813737 &    18 &    1  &  8.36e-14 &  1.06e-11 &  7.90e-14 \\ 
  18   &  13.828898 &  13.828898 &    19 &    1  &  3.15e-14 &  3.07e-11 &  1.17e-13 \\ 
  19   &  14.099027 &  14.099027 &    20 &    1  &  1.03e-14 &  9.00e-11 &  2.37e-13 \\ 
  20   &  14.101939 &  14.101939 &    21 &    1  &  8.06e-15 &  3.13e-10 &  1.49e-13 \\ 
  21   &  14.220443 &  14.220443 &    22 &    1  &  4.93e-14 &  3.89e-10 &  7.00e-14 \\ 
  22   &  14.223388 &  14.223388 &    23 &    1  &  1.10e-13 &  1.52e-11 &  5.55e-14 \\ 
  23   &  14.386254 &  14.386254 &    24 &    1  &  1.73e-14 &  2.57e-12 &  3.57e-14 \\ 
  24   &  14.507362 &  14.507362 &    25 &    2  &  4.76e-14 &  2.42e-12 &  4.93e-14 \\ 
  25   &  14.507476 &  14.507476 &    26 &    2  &  3.43e-14 &  4.26e-12 &  5.29e-14 \\ 
  26   &  15.213858 &  15.213858 &    27 &    1  &  2.58e-14 &  2.66e-10 &  2.80e-14 \\ 
  27   &  15.216611 &  15.216611 &    28 &    1  &  1.73e-14 &  6.23e-12 &  2.52e-14 \\ 
  28   &  15.307152 &  15.307152 &    29 &    1  &  5.48e-14 &  4.76e-11 &  4.52e-14 \\ 
  29   &  15.315627 &  15.315627 &    30 &    1  &  3.14e-14 &  9.78e-11 &  4.40e-14 \\ 
  30   &  15.702428 &  15.702428 &    31 &    1  &  4.34e-14 &  3.71e-10 &  6.04e-14 \\ 
  31   &  15.705038 &  15.705038 &    32 &    1  &  8.94e-14 &  5.90e-11 &  4.45e-14 \\ 
  32   &  15.993417 &  15.993417 &    33 &    1  &  5.43e-14 &  1.16e-11 &  7.00e-14 \\ 
  33   &  16.005034 &  16.005034 &    34 &    1  &  4.95e-14 &  8.32e-11 &  4.35e-14 \\ 
  34   &  16.028389 &  16.028389 &    35 &    1  &  4.21e-15 &  3.40e-11 &  4.64e-14 \\ 
 \hline 
\end{tabular}
\end{center}
\caption{grid=pipeze2, method=KrylovSchur, ts=implicit, order=2, $N_p=1$.}
\label{tab:pipeze2Order2}
\end{table}
}

\newcommand{\eigenWaveSummaryTable}{
\begin{table}[hbt]
\begin{center}\tableFontSize
\begin{tabular}{|c|c|c|c|c|c|c|c|} \hline
 \multicolumn{8}{|c|}{EigenWave: grid=pipeze2, ts=implicit, order=2, $N_p=1$, KrylovSchur } \\ \hline 
   num   &  wave       & time-steps  & wave-solves &  time-steps &   max      &  max      &  max       \\ 
   eigs  &  solves     & per period  &  per eig    &  per-eig    &   eig-err  & evect-err & eig-res    \\ 
 \hline
    35    &       124     &      10     &    3.5    &     35     &  1.85e-13 &  7.67e-09 &  6.39e-13 \\
 \hline 
\end{tabular}
\end{center}
\caption{EigenWave: grid=pipeze2, method=KrylovSchur, ts=implicit, order=2, $N_p=1$.}
\label{tab:pipeze2Order2Summary}
\end{table}
}

%% file: tables/pipeG4O2ImpEigKrylov.tex
\newcommand{\eigenWaveLongTable}{
\begin{table}[hbt]
\begin{center}\tableFontSize
\begin{tabular}{|c|r|r|c|c|c|c|c|} \hline
 \multicolumn{8}{|c|}{EigenWave: grid=pipeze4, ts=implicit, order=2, $N_p=1$, KrylovSchur } \\ \hline 
   $j$  & \multicolumn{1}{c|}{$\lambda_j$} &  \multicolumn{1}{c|}{$\lambda_{\rm true}$}  & eig &  mult &  eig-err & evect-err & eig-res  \\ \hline
   0   &   9.899117 &   9.899117 &     4 &    2  &  1.25e-13 &  1.80e-12 &  2.58e-12 \\ 
   1   &   9.899180 &   9.899180 &     5 &    2  &  1.62e-14 &  3.69e-12 &  1.83e-12 \\ 
   2   &  10.560645 &  10.560645 &     6 &    1  &  5.77e-14 &  3.54e-12 &  2.24e-12 \\ 
   3   &  10.718212 &  10.718212 &     7 &    1  &  1.19e-13 &  1.46e-10 &  1.81e-12 \\ 
   4   &  10.724761 &  10.724761 &     8 &    1  &  3.03e-14 &  2.77e-10 &  2.11e-12 \\ 
   5   &  11.454175 &  11.454175 &     9 &    1  &  1.27e-13 &  1.26e-11 &  2.15e-11 \\ 
   6   &  12.017198 &  12.017198 &    10 &    1  &  1.54e-14 &  2.30e-11 &  2.00e-12 \\ 
   7   &  12.023034 &  12.023034 &    11 &    1  &  7.98e-15 &  1.81e-11 &  3.75e-12 \\ 
   8   &  12.124743 &  12.124743 &    12 &    2  &  3.03e-14 &  1.02e-11 &  1.31e-12 \\ 
   9   &  12.124798 &  12.124798 &    13 &    2  &  1.05e-14 &  4.94e-12 &  1.01e-12 \\ 
  10   &  12.677973 &  12.677973 &    14 &    1  &  8.07e-14 &  6.00e-12 &  5.28e-12 \\ 
  11   &  13.104801 &  13.104801 &    15 &    2  &  7.92e-14 &  3.90e-11 &  2.85e-11 \\ 
  12   &  13.104874 &  13.104874 &    16 &    2  &  2.32e-14 &  3.17e-11 &  2.28e-11 \\ 
  13   &  13.406295 &  13.406295 &    17 &    1  &  1.25e-13 &  6.26e-12 &  3.72e-12 \\ 
  14   &  13.907844 &  13.907844 &    18 &    1  &  3.33e-14 &  2.08e-11 &  8.48e-13 \\ 
  15   &  13.912914 &  13.912914 &    19 &    1  &  3.03e-14 &  1.28e-11 &  6.41e-13 \\ 
  16   &  14.186922 &  14.186922 &    20 &    2  &  5.01e-16 &  4.07e-12 &  6.53e-13 \\ 
  17   &  14.186988 &  14.186988 &    21 &    2  &  3.31e-14 &  4.72e-12 &  7.72e-13 \\ 
  18   &  14.330241 &  14.330241 &    22 &    1  &  2.00e-14 &  4.73e-10 &  8.60e-13 \\ 
  19   &  14.330521 &  14.330521 &    23 &    1  &  1.00e-14 &  6.92e-10 &  7.74e-13 \\ 
  20   &  14.482613 &  14.482613 &    24 &    1  &  2.49e-14 &  3.68e-12 &  6.52e-13 \\ 
  21   &  14.669999 &  14.669999 &    25 &    2  &  1.84e-14 &  3.36e-12 &  7.44e-13 \\ 
  22   &  14.670048 &  14.670048 &    26 &    2  &  6.30e-15 &  5.24e-12 &  8.55e-13 \\ 
  23   &  15.326076 &  15.326076 &    27 &    1  &  2.41e-14 &  4.38e-10 &  5.18e-13 \\ 
  24   &  15.326338 &  15.326338 &    28 &    1  &  2.24e-14 &  2.50e-10 &  4.97e-13 \\ 
  25   &  15.437825 &  15.437825 &    29 &    1  &  2.16e-14 &  8.80e-11 &  4.31e-13 \\ 
  26   &  15.438853 &  15.438853 &    30 &    1  &  3.60e-14 &  1.62e-09 &  6.35e-13 \\ 
  27   &  15.820316 &  15.820316 &    31 &    2  &  5.66e-14 &  5.22e-12 &  6.06e-13 \\ 
  28   &  15.820385 &  15.820385 &    32 &    2  &  6.51e-15 &  5.61e-12 &  8.36e-13 \\ 
  29   &  16.174816 &  16.174816 &    33 &    1  &  7.03e-15 &  4.49e-11 &  6.26e-13 \\ 
  30   &  16.179202 &  16.179202 &    34 &    1  &  1.54e-14 &  6.24e-11 &  7.59e-13 \\ 
  31   &  16.330744 &  16.330744 &    35 &    1  &  3.02e-14 &  4.47e-12 &  5.62e-13 \\ 
  32   &  16.366408 &  16.366408 &    36 &    1  &  2.00e-14 &  9.10e-12 &  3.65e-13 \\ 
  33   &  16.367381 &  16.367381 &    37 &    1  &  1.32e-14 &  1.48e-10 &  4.77e-13 \\ 
  34   &  16.671600 &  16.671600 &    38 &    1  &  1.77e-14 &  2.02e-12 &  3.61e-13 \\ 
  35   &  16.849389 &  16.849389 &    39 &    1  &  4.11e-14 &  2.41e-09 &  5.50e-13 \\ 
  36   &  16.849621 &  16.849621 &    40 &    1  &  9.74e-14 &  5.10e-10 &  3.52e-13 \\ 
  37   &  17.032530 &  17.032530 &    41 &    1  &  6.88e-15 &  9.71e-12 &  1.01e-12 \\ 
  38   &  17.056744 &  17.056744 &    42 &    1  &  3.75e-15 &  8.20e-12 &  6.26e-13 \\ 
  39   &  17.383047 &  17.383047 &    43 &    2  &  4.44e-14 &  4.28e-11 &  1.78e-11 \\ 
  40   &  17.383092 &  17.383092 &    44 &    2  &  1.90e-14 &  6.20e-11 &  2.65e-11 \\ 
  41   &  17.503138 &  17.503138 &    45 &    1  &  1.35e-13 &  2.37e-12 &  2.10e-13 \\ 
  42   &  17.730748 &  17.730748 &    46 &    1  &  7.91e-14 &  6.49e-10 &  4.05e-13 \\ 
  43   &  17.731576 &  17.731576 &    47 &    1  &  1.30e-14 &  1.63e-10 &  2.63e-13 \\ 
  44   &  17.800886 &  17.800886 &    48 &    1  &  3.21e-14 &  2.09e-10 &  2.26e-13 \\ 
  45   &  17.801766 &  17.801766 &    49 &    1  &  2.39e-15 &  2.94e-10 &  2.94e-13 \\ 
  46   &  17.845901 &  17.845901 &    50 &    2  &  2.11e-14 &  1.72e-11 &  8.08e-12 \\ 
  47   &  17.845972 &  17.845972 &    51 &    2  &  2.27e-14 &  1.62e-11 &  5.91e-12 \\ 
  48   &  17.878473 &  17.878473 &    52 &    1  &  6.36e-15 &  1.69e-11 &  3.60e-13 \\ 
  49   &  17.901539 &  17.901539 &    53 &    1  &  1.45e-14 &  1.14e-11 &  3.17e-13 \\ 
  50   &  18.327368 &  18.327368 &    54 &    1  &  1.40e-14 &  1.39e-12 &  3.24e-13 \\ 
  51   &  18.544873 &  18.544873 &    55 &    1  &  2.49e-15 &  4.87e-11 &  2.39e-12 \\ 
  52   &  18.545664 &  18.545664 &    56 &    1  &  1.92e-15 &  1.76e-10 &  1.23e-12 \\ 
  53   &  18.670352 &  18.670352 &    57 &    1  &  4.15e-14 &  6.45e-12 &  2.30e-13 \\ 
  54   &  18.674167 &  18.674167 &    58 &    1  &  1.14e-14 &  1.19e-11 &  4.08e-13 \\ 
  55   &  18.764183 &  18.764183 &    59 &    1  &  1.72e-14 &  1.01e-11 &  2.49e-13 \\ 
  56   &  18.764383 &  18.764383 &    60 &    1  &  1.89e-16 &  2.07e-10 &  2.97e-13 \\ 
  57   &  19.102315 &  19.102315 &    61 &    1  &  4.76e-14 &  3.55e-12 &  2.02e-13 \\ 
  58   &  19.200276 &  19.200276 &    62 &    1  &  2.72e-14 &  2.05e-11 &  1.36e-11 \\ 
  59   &  19.221787 &  19.221787 &    63 &    1  &  5.73e-15 &  2.65e-12 &  2.04e-13 \\ 
  60   &  19.284389 &  19.284389 &    64 &    1  &  3.85e-14 &  2.84e-12 &  7.28e-13 \\ 
  61   &  19.618962 &  19.618962 &    65 &    1  &  1.16e-13 &  2.48e-11 &  2.10e-13 \\ 
  62   &  19.623006 &  19.623006 &    66 &    1  &  1.96e-14 &  2.42e-10 &  1.89e-13 \\ 
  63   &  19.623792 &  19.623792 &    67 &    1  &  2.28e-14 &  5.96e-10 &  1.90e-13 \\ 
  64   &  19.650418 &  19.650418 &    68 &    1  &  6.33e-15 &  1.98e-11 &  2.12e-13 \\ 
  65   &  19.652134 &  19.652134 &    69 &    1  &  5.97e-15 &  1.49e-10 &  2.03e-13 \\ 
  66   &  19.822270 &  19.822270 &    70 &    1  &  5.20e-15 &  3.49e-11 &  2.25e-13 \\ 
  67   &  19.823009 &  19.823009 &    71 &    1  &  1.70e-14 &  8.45e-11 &  1.68e-13 \\ 
  68   &  19.987588 &  19.987588 &    72 &    1  &  3.71e-14 &  8.76e-11 &  8.99e-13 \\ 
  69   &  19.991168 &  19.991168 &    73 &    1  &  2.88e-14 &  1.68e-10 &  8.98e-13 \\ 
  70   &  20.135265 &  20.135265 &    74 &    2  &  1.64e-14 &  1.45e-11 &  3.31e-13 \\ 
  71   &  20.135339 &  20.135339 &    75 &    2  &  8.82e-15 &  1.51e-11 &  4.12e-13 \\ 
  72   &  20.183297 &  20.183297 &    76 &    2  &  4.98e-14 &  1.30e-11 &  8.33e-13 \\ 
  73   &  20.183324 &  20.183324 &    77 &    2  &  8.63e-15 &  1.52e-11 &  4.12e-13 \\ 
  74   &  20.388019 &  20.388019 &    78 &    1  &  7.65e-14 &  1.04e-10 &  1.75e-13 \\ 
  75   &  20.389674 &  20.389674 &    79 &    1  &  9.23e-15 &  1.33e-10 &  2.63e-13 \\ 
  76   &  20.450533 &  20.450533 &    80 &    1  &  1.04e-15 &  4.69e-10 &  1.35e-13 \\ 
  77   &  20.452197 &  20.452197 &    81 &    1  &  2.36e-14 &  2.03e-11 &  1.41e-13 \\ 
  78   &  20.713194 &  20.713194 &    82 &    1  &  1.54e-15 &  8.01e-11 &  1.64e-13 \\ 
  79   &  20.716651 &  20.716651 &    83 &    1  &  3.24e-14 &  1.31e-10 &  2.24e-13 \\ 
  80   &  20.900680 &  20.900680 &    84 &    1  &  2.24e-14 &  2.75e-11 &  1.62e-13 \\ 
  81   &  20.920480 &  20.920480 &    85 &    1  &  1.49e-14 &  2.88e-11 &  5.97e-13 \\ 
  82   &  20.953393 &  20.953393 &    86 &    2  &  5.09e-16 &  8.28e-12 &  1.20e-13 \\ 
  83   &  20.953564 &  20.953564 &    87 &    2  &  3.05e-15 &  1.11e-11 &  2.62e-13 \\ 
  84   &  21.160258 &  21.160258 &    88 &    1  &  5.66e-14 &  3.36e-10 &  2.68e-13 \\ 
  85   &  21.161866 &  21.161866 &    89 &    1  &  3.02e-15 &  4.53e-10 &  9.29e-14 \\ 
  86   &  21.285953 &  21.285953 &    90 &    1  &  3.66e-14 &  3.52e-12 &  1.25e-13 \\ 
  87   &  21.302087 &  21.302087 &    91 &    1  &  3.34e-16 &  5.11e-11 &  1.51e-13 \\ 
  88   &  21.305325 &  21.305325 &    92 &    1  &  2.00e-14 &  9.62e-11 &  1.72e-13 \\ 
  89   &  21.473457 &  21.473457 &    93 &    1  &  7.61e-15 &  1.33e-10 &  1.52e-13 \\ 
  90   &  21.474137 &  21.474137 &    94 &    1  &  9.43e-15 &  6.98e-10 &  1.51e-13 \\ 
  91   &  21.556464 &  21.556464 &    95 &    1  &  4.78e-15 &  2.01e-10 &  1.56e-13 \\ 
  92   &  21.558020 &  21.558020 &    96 &    1  &  2.14e-15 &  2.21e-11 &  1.23e-13 \\ 
  93   &  21.681595 &  21.681595 &    97 &    1  &  8.19e-16 &  3.17e-12 &  9.82e-14 \\ 
  94   &  21.725790 &  21.725790 &    98 &    1  &  6.54e-15 &  1.10e-10 &  1.35e-13 \\ 
  95   &  21.726495 &  21.726495 &    99 &    1  &  9.81e-16 &  7.97e-11 &  1.20e-13 \\ 
  96   &  21.864235 &  21.864235 &   100 &    1  &  1.43e-14 &  2.46e-11 &  1.64e-13 \\ 
  97   &  21.867497 &  21.867497 &   101 &    1  &  8.77e-15 &  3.69e-12 &  1.71e-13 \\ 
  98   &  22.176528 &  22.176528 &   102 &    1  &  9.13e-15 &  1.14e-11 &  1.75e-13 \\ 
  99   &  22.180939 &  22.180939 &   103 &    1  &  9.29e-15 &  3.23e-11 &  1.76e-13 \\ 
 100   &  22.214770 &  22.214770 &   104 &    1  &  7.52e-15 &  6.24e-11 &  2.37e-13 \\ 
 \hline 
\end{tabular}
\end{center}
\caption{grid=pipeze4, method=KrylovSchur, ts=implicit, order=2, $N_p=1$.}
\label{tab:pipeze4Order2}
\end{table}
}

\newcommand{\eigenWaveSummaryTable}{
\begin{table}[hbt]
\begin{center}\tableFontSize
\begin{tabular}{|c|c|c|c|c|c|c|c|} \hline
 \multicolumn{8}{|c|}{EigenWave: grid=pipeze4, ts=implicit, order=2, $N_p=1$, KrylovSchur } \\ \hline 
   num   &  wave       & time-steps  & wave-solves &  time-steps &   max      &  max      &  max       \\ 
   eigs  &  solves     & per period  &  per eig    &  per-eig    &   eig-err  & evect-err & eig-res    \\ 
 \hline
    101    &       321     &      10     &    3.2    &     31     &  1.35e-13 &  2.41e-09 &  2.85e-11 \\
 \hline 
\end{tabular}
\end{center}
\caption{EigenWave: grid=pipeze4, method=KrylovSchur, ts=implicit, order=2, $N_p=1$.}
\label{tab:pipeze4Order2Summary}
\end{table}
}

%% file: tex/sphereEigenpairs.tex
\newcommand{\Gcs}{\Gc_{\rm sphere}}
\newcommand{\mPhi}{m_\phi}
\subsection{Eigenmodes of a solid sphere in three dimensions} \label{sec:sphereEigenpairs}

In this section, eigenpairs of a solid sphere are computed.
The eigenfunctions of a sphere of radius $a$ with Dirichlet boundary conditions take the form
\bse
\label{eq:sphereTrueEigenpairs}
\ba
    \phi_{\mPhi,\mr,\mTheta}(r,\phi,\theta) = r^{-1/2} J_{\mPhi+1/2}( \lambda_{\mPhi,\mr} r ) \, P_{\mPhi}^{\mTheta}(\cos(\phi)) \, 
      \begin{cases} \cos(\mTheta\theta), & \mTheta=0,1,\ldots,\mPhi, \\
                    \sin(\mTheta\theta), & \mTheta=1,2,\ldots,\mPhi,
      \end{cases} ,
\ea
where $(r,\phi,\theta)$ are the usual spherical polar coordinates, $J_{\mPhi+1/2}$ are the first kind Bessel function of fractional order, and $P_{\mPhi}^{\mTheta}$ are the associated Legendre functions. The eigenvalues are
\ba
    \lambda_{\mPhi,\mr} = {\zeta_{\mPhi,\mr}\over a}, \qquad \mPhi=0,1,2,\ldots, \quad \mr=1,2,3,\ldots, 
\ea
where $\zeta_{\mPhi,\mr}$ are the zeros of $J_{\mPhi+1/2}(\zeta)$ indexed as $\mr=1,2,3,\ldots$.
\ese
Note that eigenvalue $\lambda_{\mPhi,\mr}$ has multiplicity $2\mPhi+1$ and thus there are eigenvalues with high multiplicity which 
could cause difficulties for numerical algorithms.

The overset grid for the solid sphere, denoted by $\Gcs^{(j)}$, consist of four component grids,
each with grid spacing approximately equal to $\ds^{(j)}=1/(10 j)$.  
The sphere, of radius $a=1$, 
 is covered with three boundary-fitted patches near the surface as shown on the left in Figure~\ref{fig:sphereFig}.
There is one patch specified using spherical polar coordinates that covers much of the sphere except
near the poles. To remove the polar singularities there are two patches
that cover the north and south poles, defined by orthographic mappings.
A background Cartesian grid covers the interior of the sphere.  

\input tex/sphereGridFig

\input tex/sphereContoursFig

\begin{table}[hbt]
\begin{center}\tableFontSize
\begin{tabular}{|c|c|c|c|c|c|c|c|c|} \hline
 \multicolumn{9}{|c|}{EigenWave: sphere, ts=implicit, $\omega=   5.0$, $N_p=1$, KrylovSchur } \\ \hline 
   order & num   &  wave       & time-steps  & wave-solves &  time-steps &   max      &  max      &  max       \\ 
   & eigs  &  solves     & per period  &  per eig    &  per-eig    &   eig-err  & evect-err & eig-res    \\ 
 \hline
    2 & 29    &       123     &      10     &    4.2    &     42     &  1.91e-09 &  7.13e-05 &  2.11e-05 \\
    4 & 29    &       143     &      10     &    4.9    &     49     &  6.29e-12 &  1.68e-09 &  1.09e-08 \\
 \hline 
\end{tabular}
\end{center}
\vspace*{-1\baselineskip}
\caption{Summary of EigenWave performance for sphere grid $\Gcs^{(2)}$ using the KrylovSchur algorithm and implicit time-stepping.  The spatial order of accuracy is 2 for the top row and 4 for the bottom row, and the wave-solves use $N_p=1$ to determine the final time.}
\label{tab:spheree2Summary}
\end{table}

Table~\ref{tab:spheree2Summary} summarizes results of computing eigenpairs with target frequency $\omega=5$, using both second and fourth-order accurate discretizations on grid $\Gcs^{(2)}$. In both cases, $24$ eigenpairs are requested and $29$ 
eigenpairs are found accurately by the KrylovSchur
algorithm in a total of $123$ wave-solves. This corresponds to approximately 
$4.2$ wave-solves per computed eigenpair.
Tables~\ref{tab:spheree2Order2} and~\ref{tab:spheree2Order4} in the appendix provides more details of the eigenpairs. These tables reveal that most eigenpairs are more accurate than the worst cases reported
in the summary table. As noted above, the sphere has eigenvalues with high multiplicities, and these eigenpairs are generally found to
high accuracy.
The right graph in Figure~\ref{fig:sphereFig} shows the filter function $\beta(\lambda;\omega)$ for the second-order accurate code with the computed
eigenvalues marked. 
Figure~\ref{fig:sphereEigenVectors} shows contours on cutting planes through three of the computed eigenvectors.

%% file: tex/sphereGridFig.tex
{
\newcommand{\figw}{6.5cm}
\newcommand{\figh}{5.5cm}
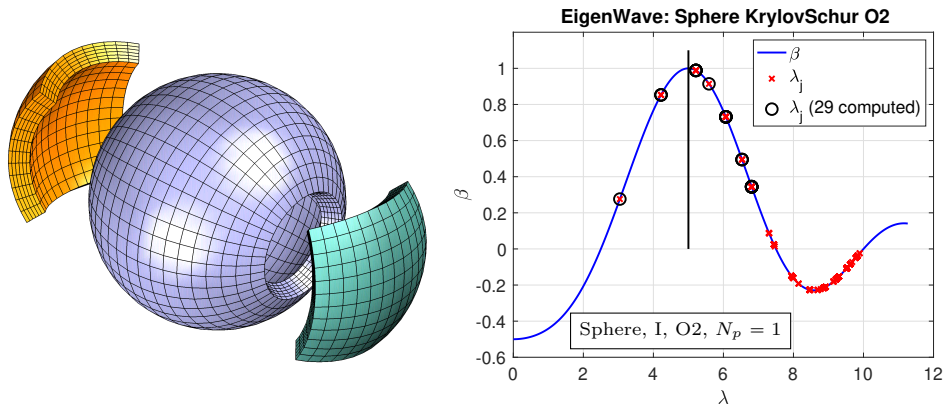
\begin{figure}[htb]
\begin{center}
\begin{tikzpicture}
  \useasboundingbox (0,.25) rectangle (12.5,5.25);  
  \figByWidth{0.0}{0}{fig/sphereGridsExploded}{5.75cm}[0.0][0.0][0.1][0.1]
   \begin{scope}[xshift=6cm,yshift=-12pt]
     \figByWidth{0}{0}{fig/sphereG2O2ImpEigKrylov}{\figw}[0][0][0][0]
     \draw (1.5,1) node[draw,fill=white,anchor=west,xshift=0pt,yshift=0pt] {\scriptsize Sphere, I, O2, $\Np=1$};
  \end{scope}  
  
\end{tikzpicture}
\end{center}
\caption{At left is exploded view of the surface patches of the overset grid for a solid sphere.
At right is a plot of the filter function $\beta$ with the  
computed eigenvalues marked with black circles for a $2nd$-order accurate computation on grid $\Gcs^{(1)}$.
The target frequency $\omega=5$ is marked as a vertical black line.
 }
\label{fig:sphereFig}
\end{figure}
}

%% file: tex/sphereContoursFig.tex
{
\newcommand{\drawContour}[7]{%
\begin{scope}[#1]
\draw(0.0,0) node[anchor=south west,xshift=-4pt,yshift=+0pt] {\trimfiga{fig/#2}{\figWidtha}};
  \draw(.5,.5) node[draw,fill=white,anchor=west,xshift=2pt,yshift=1pt] {\scriptsize #3};
\begin{scope}[xshift=-.2cm,yshift=-2pt]
  \draw (\xcb,\ycb) node[anchor=south west,xshift=0.25cm,yshift=.5cm,rotate=-90] {\trimfigcb{fig/colourBarLines}{\cbWidth}{\cbHeight}};
  \draw (.8,0) node[anchor=north,xshift=+3pt,yshift=+2pt] {\scriptsize $#6$};
  \draw (4.8,0) node[anchor=north,xshift=+0pt,yshift=+2pt] {\scriptsize $#7$};
\end{scope}
\end{scope}
}
\newcommand{\cbWidth}{.2cm}
\newcommand{\cbHeight}{4cm}
\newcommand{\xcb}{.5cm}
\newcommand{\ycb}{-.2cm}
\setlength{\ycbTop}{\ycb+\cbHeight}
\setlength{\ycbMid}{\ycb+\cbHeight*\real{.5}}
\newcommand{\trimfigcb}[3]{\includegraphics[width=#2, height=#3, clip, trim=17cm 2.35cm 1.65cm 2.35cm]{#1}}
\newcommand{\figWidtha}{4.5cm}
\newcommand{\trimfiga}[2]{\trimw{#1}{#2}{.035}{.15}{.075}{.11}}
\begin{figure}[htb]
\begin{center}
\begin{tikzpicture}
   \useasboundingbox (0,.3) rectangle (15,4.75);  

   \begin{scope}[yshift=0cm]
    \drawContour{xshift= 0.cm,yshift=0.00cm}{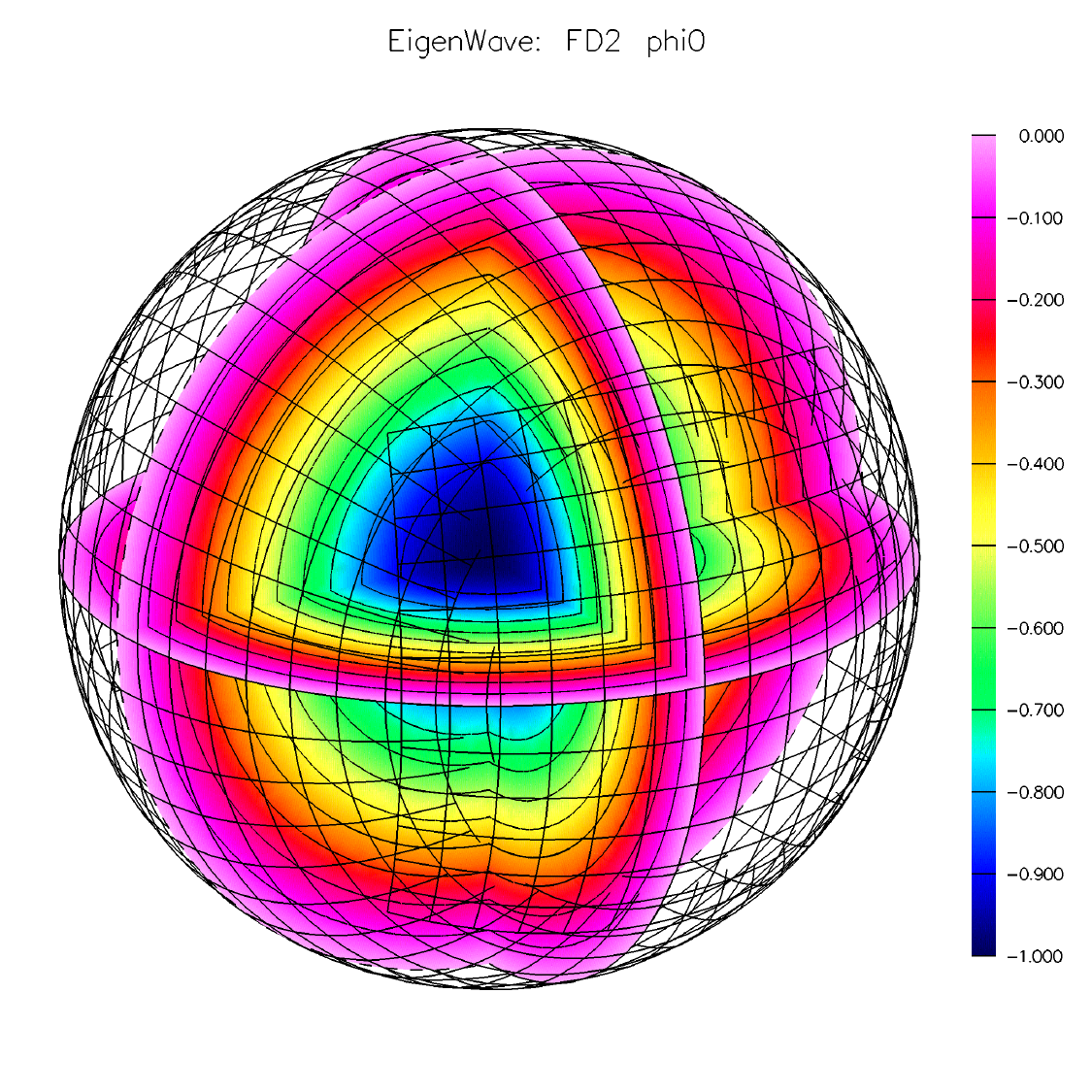}{$\phi^{(0)}$}{$v$}{$t=0.3$}{$-1.0$}{$.0$}     
    \drawContour{xshift= 5.cm,yshift=0.00cm}{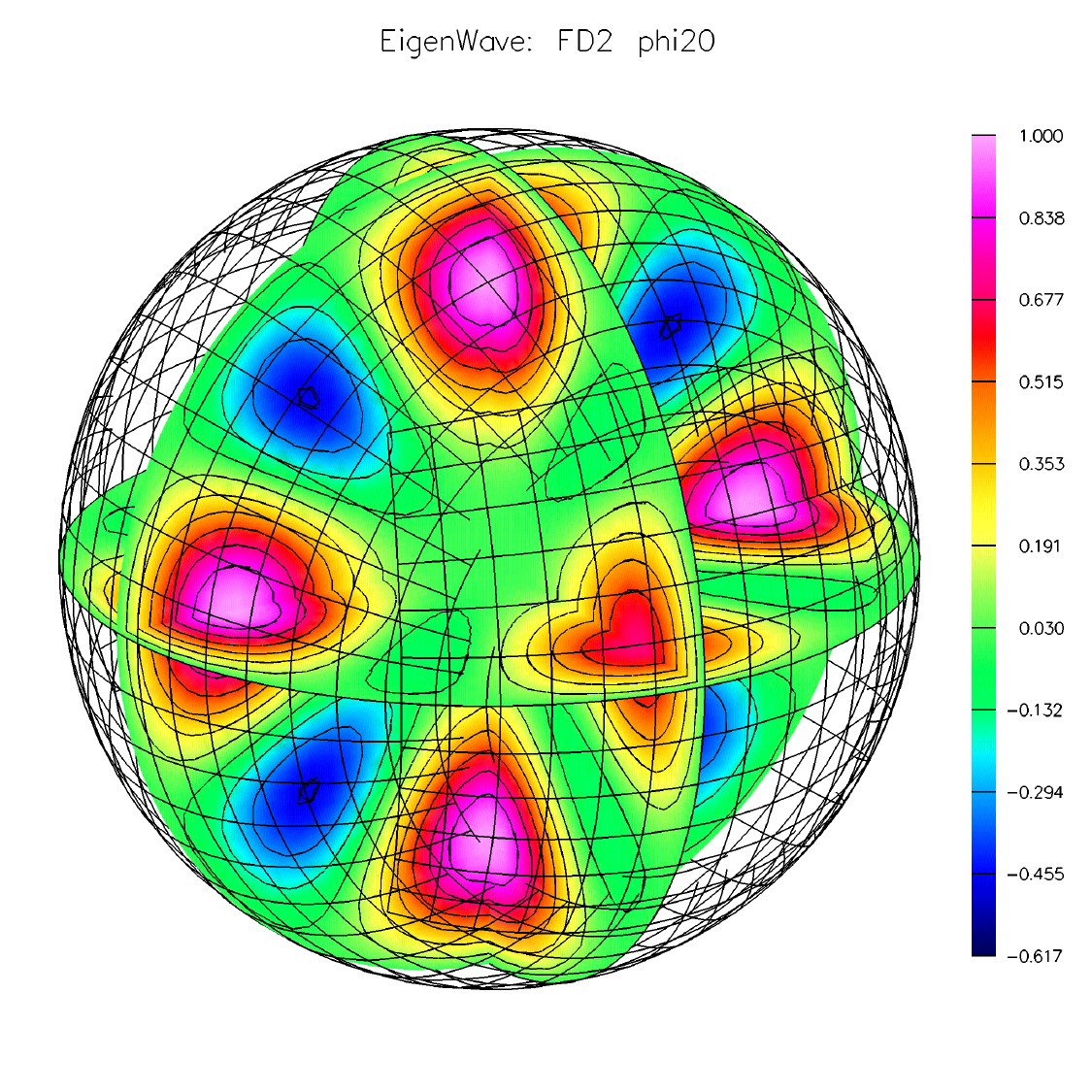}{$\phi^{(20)}$}{$v$}{$t=0.3$}{$-1.0$}{$1.0$}     
    \drawContour{xshift=10.cm,yshift=0.00cm}{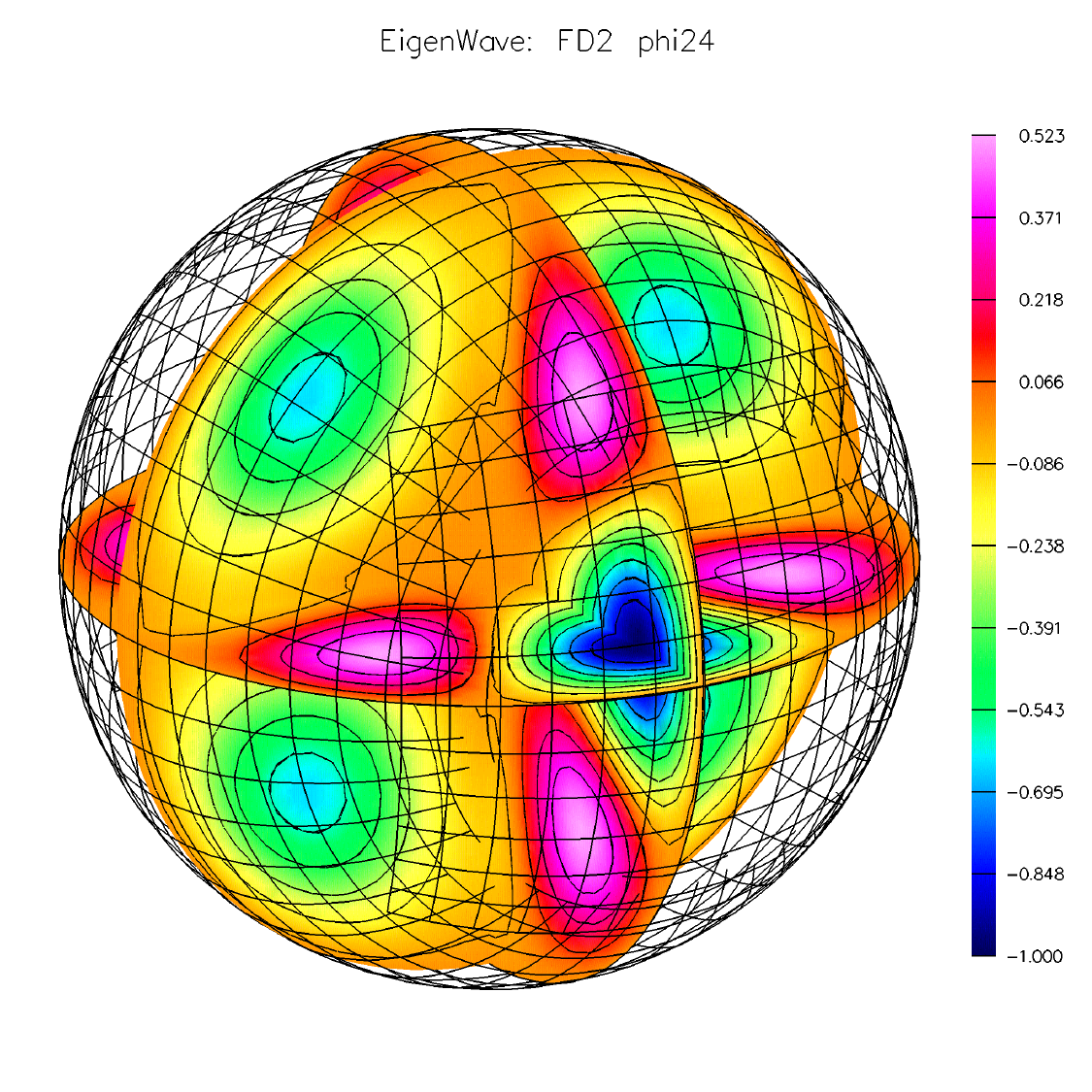}{$\phi^{(24)}$}{$v$}{$t=0.3$}{$-1.0$}{$1.0$}     
   \end{scope}

\end{tikzpicture}
\end{center}
\caption{Three computed eigenvectors on a sphere using $\Gcs^{(2)}$. 
    }
\label{fig:sphereEigenVectors}
\end{figure}
}

%% file: tex/doubleEllipsoidEigenpairs.tex
\newcommand{\Gcdes}{\Gc_{\rm des}}
\subsection{Eigenmodes of a double ellipsoid} \label{sec:doubleEllipsoidEigenpairs}

In this section eigenpairs are computed for a three-dimensional version of the Penrose unilluminable room;  
the two-dimensional case was discussed in Section~\ref{sec:doubleEllipseEigenpairs}.
The two-dimensional geometry, which consists of 
a smaller half-ellipse, with semi-axes $(a_1,b_1)=(2,1)$, and a 
larger half-ellipse, with semi-axes $(a_2,b_2)=(3,6)$,
is revolved about the axis of symmetry (i.e.~the vertical $y$-axis) to form a three-dimensional body of revolution. 
The computation is restricted to the upper half of the geometry.
Neumann boundary conditions are specified on the lower boundary which then corresponds to a symmetry plane.
Dirichlet boundary conditions are applied on all other boundaries.

\input tex/doubleEllipsoidGridFig.tex

Let $\Gcdes^{(j)}$ denote the overset grid for the double-ellipsoid geometry with typical grid spacing $\ds^{(j)}=1/(10 j)$. 
As shown in Figure~\ref{fig:doubleEllipsoidGridFig}, the overset grid is comprised of five component grids:
a large Cartesian background grid, an outer ellipsoidal shell, a green inner ellipsoidal shell, 
a light blue cylindrical grid at the top and a red patch to cover the polar singularity in the inner ellipsoid.

{
\newcommand{\figw}{6.5cm}
\newcommand{\figh}{5.5cm}

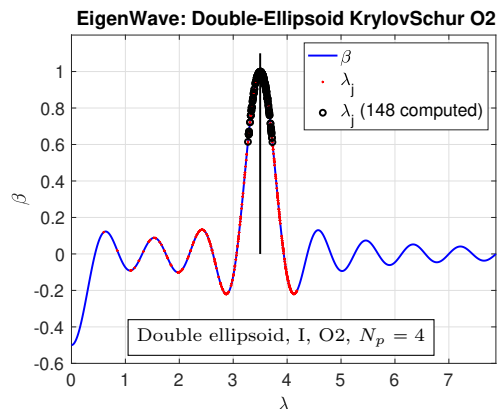
\begin{figure}[htb]
\begin{center}
\begin{tikzpicture}[scale=1]
  \useasboundingbox (0,.7) rectangle (7,5.5);  

  \begin{scope}[yshift=0cm]
     \figByWidth{0.0}{0}{fig/darkCornerRoom3dG1O2EigsKrylov}{\figw}[0][0][0][0];
     \draw (1.5,1) node[draw,fill=white,anchor=west,xshift=0pt,yshift=0pt] {\scriptsize Double ellipsoid, I, O2, $\Np=4$};
  \end{scope}

\end{tikzpicture}
\end{center}
\caption{Filter function and computed eigenvalues for the double ellipsoid.
    } 
\label{fig:doubleEllipsoidBetaAndEigs}
\end{figure}
}

\begin{table}[hbt]
\begin{center}\tableFontSize
\begin{tabular}{|c|c|c|c|c|c|c|c|} \hline
 \multicolumn{8}{|c|}{EigenWave: double ellipsoid, ts=implicit, order=2, $\omega=3.5$, $N_p=4$, KrylovSchur } \\ \hline 
   num   &  wave       & time-steps  & wave-solves &  time-steps &   max      &  max      &  max       \\ 
   eigs  &  solves     & per period  &  per eig    &  per-eig    &   eig-err  & evect-err & eig-res    \\ 
 \hline
    148    &       507     &      10     &    3.4    &     137     &  5.59e-13 &  1.75e-08 &  2.38e-09 \\
 \hline 
\end{tabular}
\end{center}
\vspace*{-1\baselineskip}
\caption{Summary of EigenWave performance for double ellipsoid grid $\Gcdes^{(1)}$ using the KrylovSchur algorithm and implicit time-stepping.  The spatial order of accuracy is 2 and the wave-solves use $N_p=4$ to determine the final time.}
\label{tab:darkCornerRoom3dGride1Order2Summary}
\end{table}


The true discrete eigenpairs for this problem are computed with SLEPc,
using GMRES (with $100$ restart vectors) to solve the implicit equations for the shifted Laplacian. 
(The problem is too large for SLEPc to solve with a direct solver on a Linux workstation with $196$ Gb of memory.)

\input tex/doubleEllipsoidContoursFig.tex

Table~\ref{tab:darkCornerRoom3dGride1Order2Summary}
summarizes results of computing eigenpairs to second-order accuracy on grid $\Gcdes^{(1)}$ (with approximately a million grid points).
Implicit time-stepping is used with $\Nits=10$ time-steps per period and $\Np=4$. 
The target frequency is $\omega=3.5$. 
A total of $128$ eigenpairs are requested and $148$
eigenpairs are found by the KrylovSchur
algorithm in a total of $507$ wave-solves. This corresponds to approximately 
$3.4$ wave-solves per eigenpair found.
Figure~\ref{fig:doubleEllipsoidBetaAndEigs} shows the filter function and the locations of the computed eigenvalues. Figure~\ref{fig:doubleEllipsoidEigenVectors} shows contours of the absolute value of some selected eigenvectors
computed on a finer grid $\Gcdes^{(2)}$ (about $6.5$ million grid points).

%% file: tex/doubleEllipsoidGridFig.tex
{
\newcommand{\figw}{6.5cm}
\newcommand{\figh}{6cm}
\begin{figure}[htb]
\begin{center}
\begin{tikzpicture}
  \useasboundingbox (0,.6) rectangle (14,.95*\figh);  

  \figByWidthb{0.0}{0}{fig/darkCornerRoom3dGridG1}{\figw}[0.][0.][0.075][0.075]

  \figByWidthb{1.1*\figw}{0.3}{fig/darkCornerRoom3dGridG1Zoom}{\figw}[0.][0.0][0.2][0.05]

\end{tikzpicture}
\end{center}
\caption{Double ellipsoid overset grid $\Gcdes^{(1)}$ 
(grid lines coarsened by a factor of $2$). 
 }
\label{fig:doubleEllipsoidGridFig}
\end{figure}
}

%% file: tex/doubleEllipsoidContoursFig.tex
{
\newcommand{\drawContour}[7]{%
\begin{scope}[#1]
\draw(0.0,0) node[anchor=south west,xshift=-4pt,yshift=+0pt] {\trimfiga{fig/#2}{\figWidtha}};
  \draw(.5,.35) node[draw,fill=white,anchor=west,xshift=2pt,yshift=2pt] {\scriptsize #3};
\begin{scope}[xshift=.1cm,yshift=-5pt]
  \draw (\xcb,\ycb) node[anchor=south west,xshift=0.25cm,yshift=.5cm,rotate=-90] {\trimfigcb{fig/colourBarLines}{\cbWidth}{\cbHeight}};
  \draw (.8,0) node[anchor=north,xshift=+3pt,yshift=+2pt] {\scriptsize $#6$};
  \draw (4.8,0) node[anchor=north,xshift=+0pt,yshift=+2pt] {\scriptsize $#7$};
\end{scope}
\end{scope}
}
\newcommand{\cbWidth}{.2cm}
\newcommand{\cbHeight}{4cm}
\newcommand{\xcb}{.5cm}
\newcommand{\ycb}{-.2cm}
\setlength{\ycbTop}{\ycb+\cbHeight}
\setlength{\ycbMid}{\ycb+\cbHeight*\real{.5}}
\newcommand{\trimfigcb}[3]{\includegraphics[width=#2, height=#3, clip, trim=17cm 2.35cm 1.65cm 2.35cm]{#1}}
\newcommand{\figWidtha}{5cm}
\newcommand{\trimfiga}[2]{\trimw{#1}{#2}{.01}{.115}{.12}{.11}}
\begin{figure}[htb]
\begin{center}
\begin{tikzpicture}
   \useasboundingbox (.25,.2) rectangle (16,4.7);  

  \begin{scope}[yshift=0.0cm]
    \drawContour{xshift=0.0cm,yshift=0.0cm}{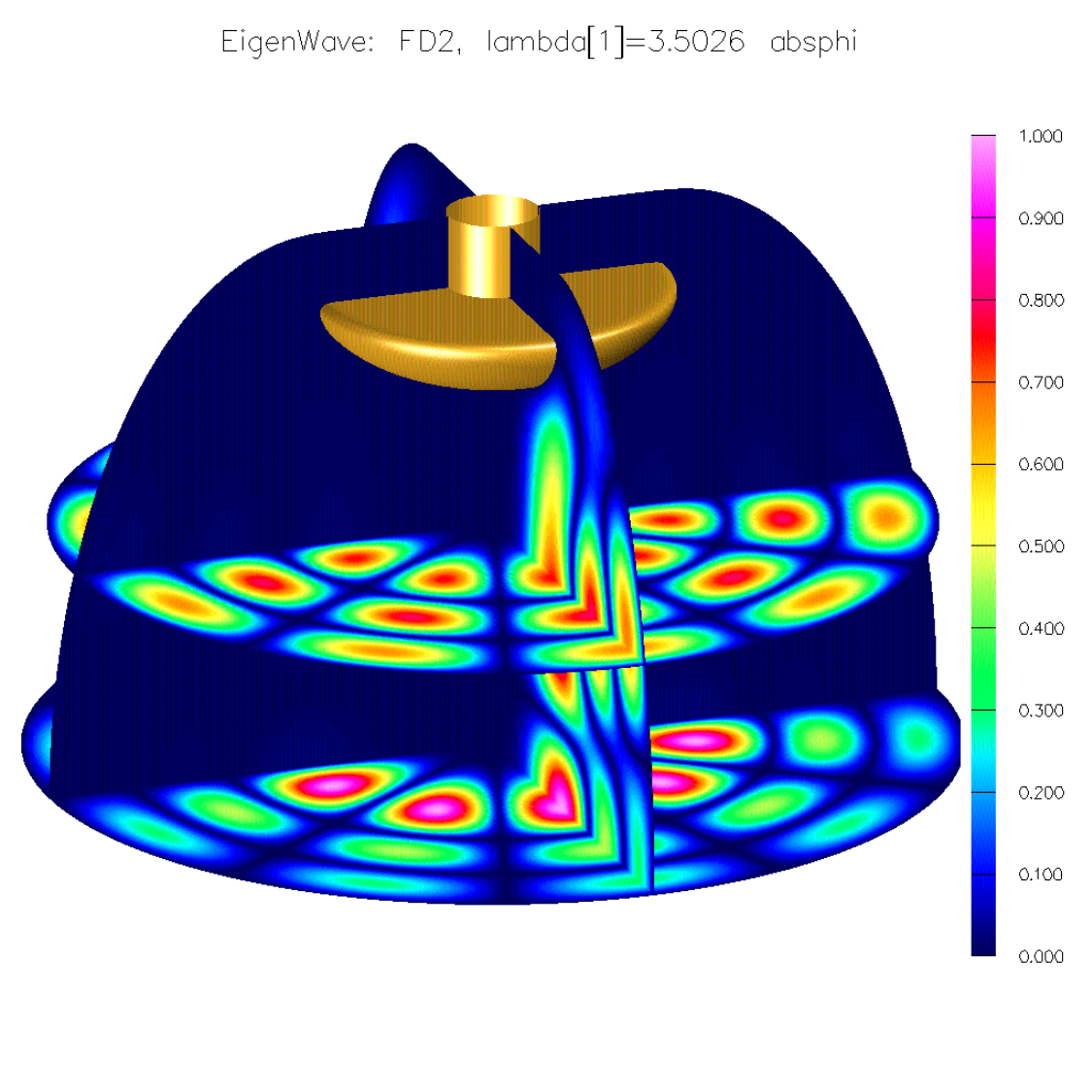}{$|\phi|$}{$v$}{$t=0.3$}{$0.0$}{$1.0$}    
    \drawContour{xshift=5.3cm,yshift=0.0cm}{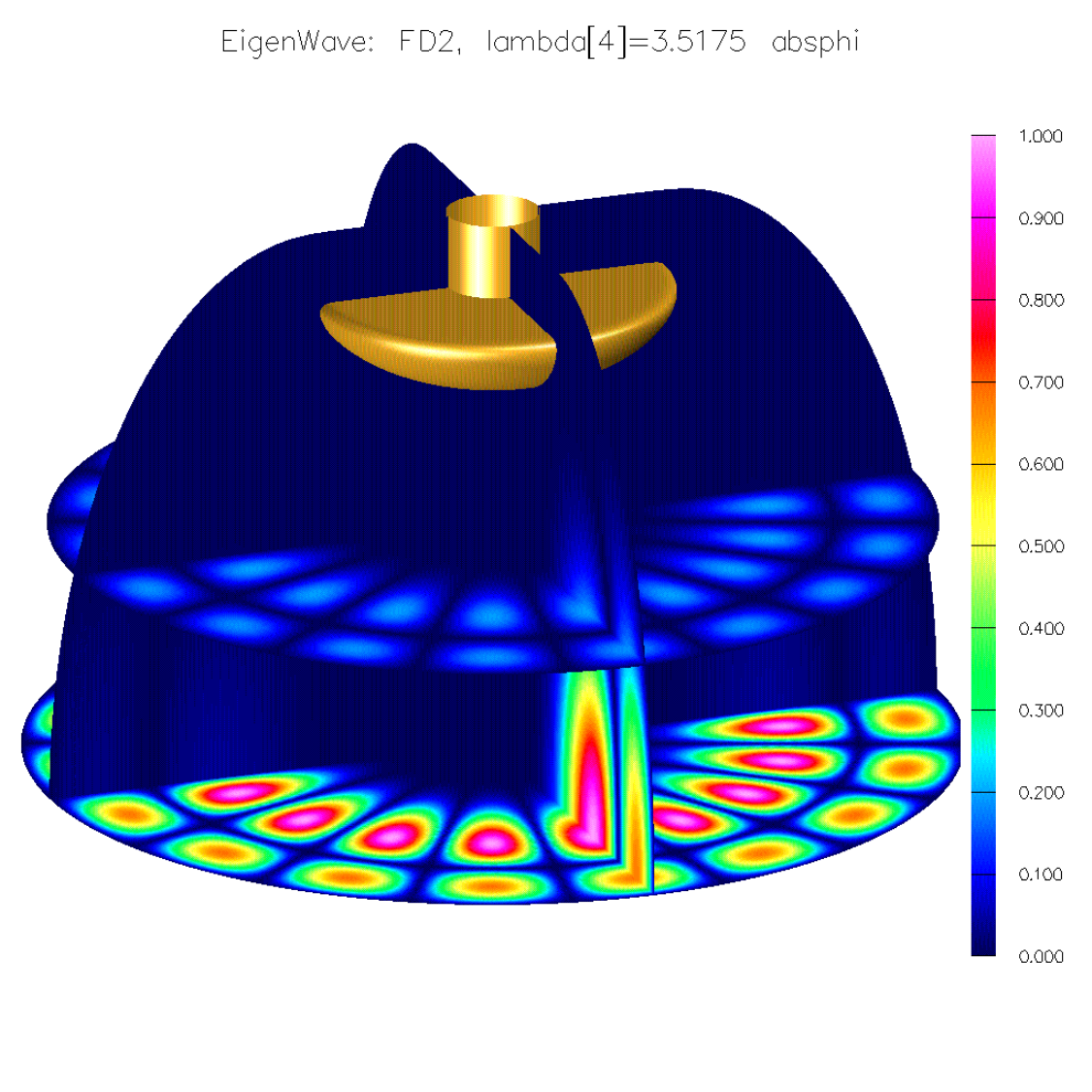}{$|\phi|$}{$v$}{$t=0.3$}{$0.0$}{$1.0$}    
    \drawContour{xshift=10.6cm,yshift=0.0cm}{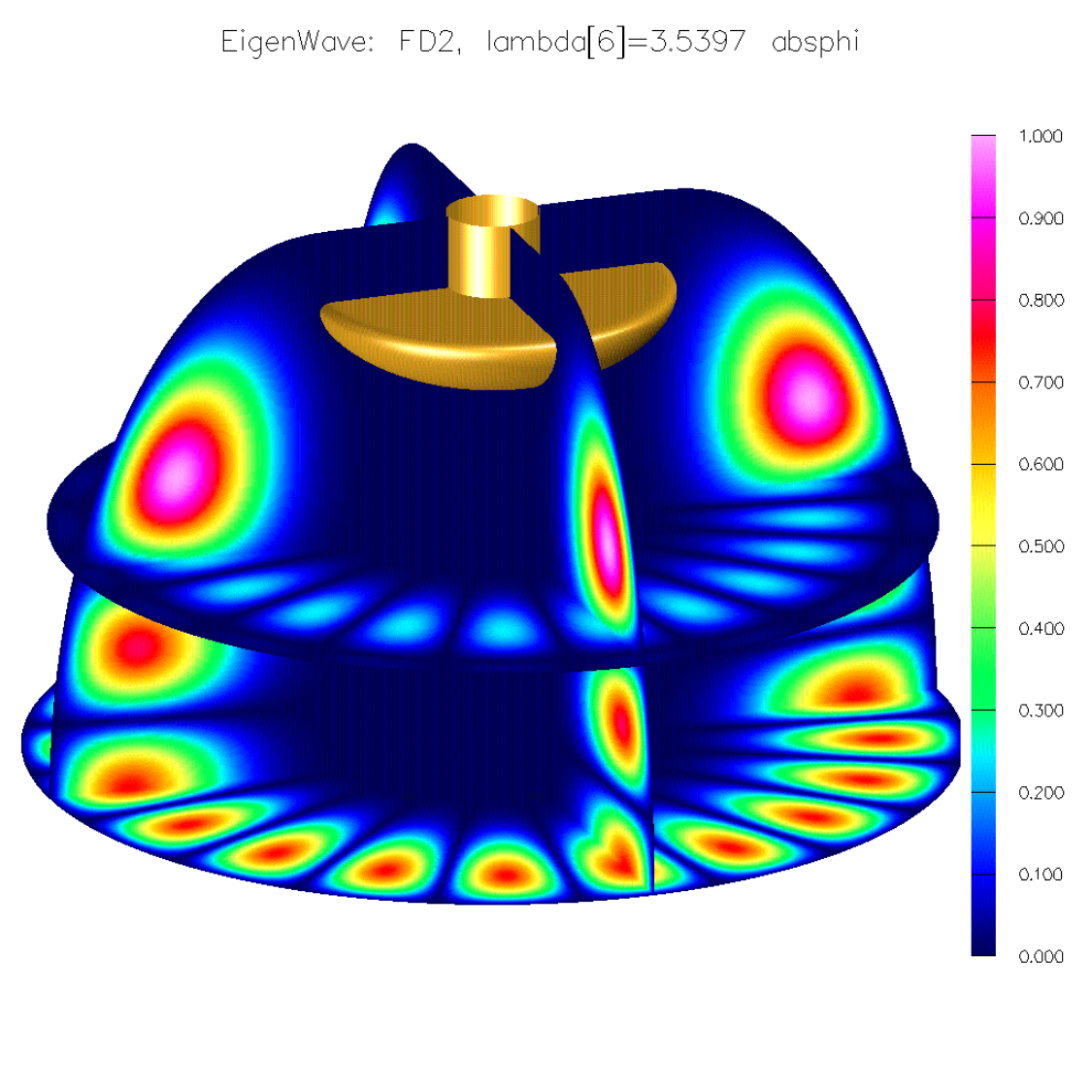}{$|\phi|$}{$v$}{$t=0.3$}{$0.0$}{$1.0$}    
  \end{scope}

\end{tikzpicture}
\end{center}
\caption{Selected eigenvectors for the double ellipsoid using $\Gcdes^{(2)}$ and the $2nd$-order accurate discretization.
    }
\label{fig:doubleEllipsoidEigenVectors}
\end{figure}
}

%% file: tex/optimalAlgorithm.tex
\section{An optimal $O(N)$ eigenvalue solver: EigenWave with implicit time-stepping and multigrid}  \label{sec:implicitTimeSteppingWithMultigrid}


In this section, results are presented for the scaling of the EigenWave algorithm 
as the mesh is refined when computing eigenpairs at a fixed target frequency.
Evidence is provided that suggests that the EigenWave algorithm using implicit time-stepping with a fixed 
number of time-steps per period, combined with a multigrid algorithm to solve the implicit time-stepping equations, is an optimal $O(N)$ algorithm.
The multigrid solver used here is the geometric multigrid solver for overset grids called Ogmg and described in~\cite{OGMG,automg,multigridWithNonstandardCoarsening2023}.
It is a matrix free solver and is thus quite memory efficient. It uses optimized red-black smoothers for Cartesian grids and zebra-line smoothers for curvilinear grids with stretched grid lines. 
It uses an adaptive cycle where the number of smooths may vary between different component grids.
It has an automatic coarsening algorithm for overset grids and performs additional smoothing near
overset grid interpolation boundaries to keep the residual smooth.

{
\newcommand{\figSize}{6.5cm}

\begin{figure}[htb]
\begin{center}
\begin{tikzpicture}[scale=1]
  \useasboundingbox (0,.7) rectangle (7,5.5);  

  \begin{scope}[yshift=0cm]
    \figByWidth{0.0}{0}{fig/betaWithCoarseFineLambda}{\figSize}[0][0][0][0]
  \end{scope}  

\end{tikzpicture}
\end{center}
\caption{The asymptotic convergence rate of the WaveHoltz filter is normally determined by the discrete eigenvalue closest to $\omega$. 
This figure shows values of $\beta$ evaluated at the coarse grid, fine grid, and continuous eigenvalues of the one-dimensional Laplacian.
As the mesh is refined the ACR approaches the value of $\beta$ at the eigenvalue $\lambda=1$.
    } 
\label{fig:betaWithCoarseAndFineLambda}
\end{figure}
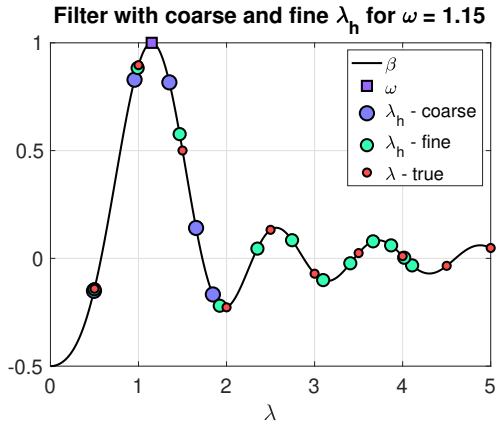
}

To understand why EigenWave might be an $O(N)$ algorithm it is useful to first discuss the WaveHoltz algorithm upon which EigenWave is based.
Evidence was provided in~\cite{overHoltzArXiv2025,overHoltzPartOne,overHoltzPartTwo} to show the optimal $O(N)$ scaling of the WaveHoltz algorithm for solving Helmholtz
problems at a fixed frequency when using implicit time-stepping.
Optimality of the WaveHoltz algorithm is based on several factors. One key factor is that the convergence rate of the WaveHoltz fixed-point iteration
is essentially independent of the mesh spacing. This is illustrated in Figure~\ref{fig:betaWithCoarseAndFineLambda} which shows the filter function~$\beta$
and the eigenvalues of the 
one-dimensional Laplacian for both coarse and fine-grid approximations. In this case, the asymptotic convergence rate (ACR)
 of the WaveHoltz fixed-point iteration is limited 
by the value of $|\beta(\lambda_h)|$ for the discrete eigenvalue near $\lambda=1$.
As the grid is refined, the eigenvalues of the discretized problem converge to the eigenvalues of the continuous problem, and as a result the rate of convergence does not depend significantly on the grid spacing. Further, the poorly resolved eigenvalues of the discrete Laplacian are large, and far away to the right along the $\lambda$-axis. These are damped rapidly during the WaveHoltz iteration.
These observations have significance with respect to the convergence of Arnoldi algorithms for computing eigenvalues since it is the eigenvalues
closest to the target frequency $\omega$ that are usually found first.
A second key factor for the $O(N)$ convergence of the WaveHoltz algorithm is that a fixed number of implicit time-steps per period can
be used independent of the mesh spacing. The linear system resulting from implicit time-stepping is then solved with an $O(N)$ multigrid
algorithm. The combination of these factors leads to the $O(N)$ algorithm.

\input tex/scalingSquareOrder2Table
\input tex/scalingSquareOrder4Table

The EigenWave algorithm, which is based on the WaveHoltz algorithm, uses some of the same components as the WaveHoltz algorithm, such as implicit time-stepping.
The EigenWave algorithm uses an Arnoldi-based eigenvalue solver, such as the KrylovSchur algorithm from SLEPSc (note that the GMRES accelerated WaveHoltz
algorithm is also based on the Arnoldi algorithm).
We first remark that the total CPU cost of the EigenWave solve is dominated by the CPU cost of the wave-solves (matrix-vector multiplies in the Arnoldi algorithm),
and the cost of each wave-solve is dominated by the cost of solving the implicit time-stepping equations.
Thus the CPU performance of the full algorithm should be dominated by the cost of solving the implicit time-stepping equations.
Since the distribution of the relevant eigenvalues (and form of the corresponding eigenvectors) being found by the Arnoldi algorithm are essentially independent of 
the mesh spacing, it is expected that the convergence of the Arnoldi algorithm should also be roughly independent of the mesh spacing (once the relevant eigenpairs are
well resolved on the grid). Thus if the cost of each Arnoldi iteration is $O(N)$, then the overall EigenWave algorithm should have an $O(N)$ cost.
We next present some computational evidence to support this conjecture.

\input tex/scalingDiskOrder2Table
\input tex/scalingDiskOrder4Table

{
\newcommand{\figSize}{6.5cm}

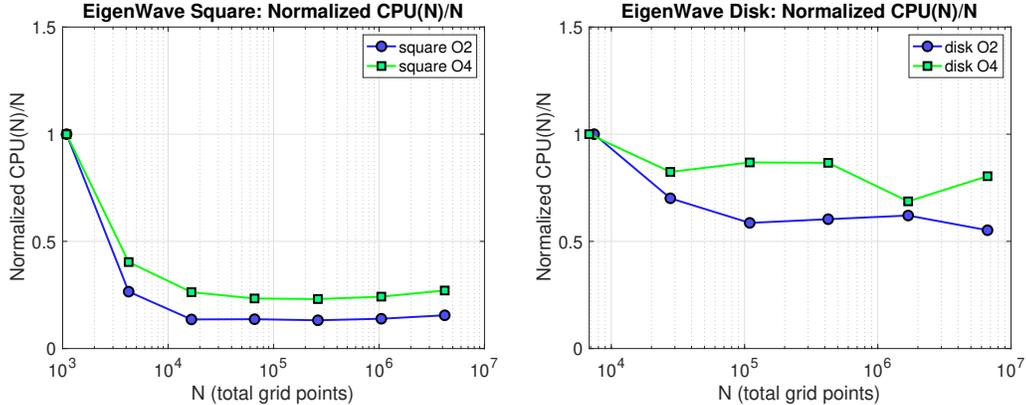
\begin{figure}[htb]
\begin{center}
\begin{tikzpicture}[scale=1]
  \useasboundingbox (0,.7) rectangle (13,5.5);  

   \begin{scope}[yshift=0cm]
    \figByWidth{0.0}{0}{fig/scalingSquareCpu}{\figSize}[0][0][0][0]
    \figByWidth{7.0}{0}{fig/scalingDiskCpu}{\figSize}[0][0][0][0]
  \end{scope}  

\end{tikzpicture}
\end{center}
\caption{Grid scaling with implicit time-stepping and multigrid. 
    Normalized values of CPU$(N)$/$N$, versus number of grid points. 
    Left: square. Right: disk.
   The results show the near optimal scaling of the EigenWave algorithm.
    }
\label{fig:squareAndDiskScaling}
\end{figure}
}

Tables~\ref{tab:scalingSquareOrder2} and~\ref{tab:scalingSquareOrder4}
show results for computing eigenpairs of the Laplacian for a two-dimensional square using EigenWave with order of accuracy equal to~two and~four.
Tables~\ref{tab:scalingDiskOrder2} and~\ref{tab:scalingDiskOrder4}
show corresponding results for a two-dimensional disk.
Sixteen eigenpairs are requested with a target frequency of $\omega=12$ (square) and $\omega=6$ (disk).
The wave equation is integrated using implicit time-stepping with $\Np=1$ and $\Nits=10$ steps per period.
The implicit time-stepping equations are solved with multigrid to a relative tolerance of~$10^{-10}$.
Results for a sequence of grids are reported with the grid spacing decreasing by a factor two from one grid to the next.
Both the number of eigenpairs found and number of wave-solves used are seen to be roughly constant. 
The average number of multigrid cycles per wave-solve increases somewhat in some cases; 
the reason for this needs further investigation. Note, however, that on an overset grid the number of smooths per component grid
is chosen dynamically and this will affect the behaviour and cost of each cycle.
Since the number of grid points increases by a factor of $4$ from one grid to the next, it is expected that
the CPU time should also increase by a factor $4$. The column titled ``CPU ratio'' gives the ratio of the CPU time for a given
grid to the CPU time for the next coarser grid. These ratios 
indicate that the CPU time does increase by roughly a factor $4$. 
To see these CPU times in an alternative fashion, Figure~\ref{fig:squareAndDiskScaling} 
shows graphs of the CPU time divided by $N$ for the sequence of grids.  The CPU time is
normalized by the CPU on the coarsest grid.
In summary, these results provide strong evidence that EigenWave, combined 
with implicit time-stepping and a multigrid solver, has near optimal $O(N)$ scaling as the grid is refined.

%% file: tex/scalingSquareOrder2Table.tex
\begin{table}[hbt]
\begin{center}\tableFontSize
\begin{tabular}{|c|c|c|c|c|c|c|} \hline
      \multicolumn{7}{|c|}{EigenWave scaling: square, order 2} \\ \hline 
              &   grid-       & eigen-  &  wave-    &  multigrid    &  CPU          &  CPU      \\
 $\ds $       &   points      & pairs   &  solves   &  cycles       &  seconds      &  ratio    \\ \hline
 $1/32   $    &   1.1e+03     &  $19$   &  $78$    &   $  6.1$      &    1.4e+00    &           \\
 $1/64   $    &   4.2e+03     &  $16$   &  $67$    &   $  5.4$      &    1.4e+00    &     1.03       \\
 $1/128  $    &   1.7e+04     &  $16$   &  $67$    &   $  5.1$      &    2.8e+00    &     2.01       \\
 $1/256  $    &   6.6e+04     &  $18$   &  $79$    &   $  5.9$      &    1.1e+01    &     4.00       \\
 $1/512  $    &   2.6e+05     &  $17$   &  $82$    &   $  6.0$      &    4.3e+01    &     3.84       \\
 $1/1024 $    &   1.1e+06     &  $17$   &  $82$    &   $  6.0$      &    1.8e+02    &     4.21       \\
 $1/2048 $    &   4.2e+06     &  $17$   &  $82$    &   $  6.0$      &    8.1e+02    &     4.46       \\
\hline
\end{tabular}
\end{center}
\vspace*{-1\baselineskip}
\caption{Scaling of EigenWave when using multigrid to solve the implicit time-stepping equations for the $2nd$-order accurate discretization.
   }
\label{tab:scalingSquareOrder2}
\end{table}

%% file: tex/scalingSquareOrder4Table.tex
\begin{table}[hbt]
\begin{center}\tableFontSize
\begin{tabular}{|c|c|c|c|c|c|c|} \hline
      \multicolumn{7}{|c|}{EigenWave scaling: square, order 4} \\ \hline 
              &   grid-       & eigen-  &  wave-    &  multigrid    &  CPU          &  CPU      \\
 $\ds $       &   points      & pairs   &  solves   &  cycles       &  seconds      &  ratio    \\ \hline
 $1/32   $    &   1.1e+03     &  $17$   &  $83$    &   $  6.0$      &    1.6e+00    &           \\
 $1/64   $    &   4.2e+03     &  $16$   &  $79$    &   $  6.0$      &    2.5e+00    &     1.57       \\
 $1/128  $    &   1.7e+04     &  $16$   &  $79$    &   $  6.0$      &    6.5e+00    &     2.57       \\
 $1/256  $    &   6.6e+04     &  $17$   &  $80$    &   $  6.0$      &    2.3e+01    &     3.52       \\
 $1/512  $    &   2.6e+05     &  $17$   &  $82$    &   $  6.0$      &    9.0e+01    &     3.94       \\
 $1/1024 $    &   1.1e+06     &  $17$   &  $82$    &   $  6.1$      &    3.8e+02    &     4.19       \\
 $1/2048 $    &   4.2e+06     &  $17$   &  $82$    &   $  6.4$      &    1.7e+03    &     4.47       \\
\hline
\end{tabular}      
\end{center}
\vspace*{-1\baselineskip}
\caption{Scaling of EigenWave when using multigrid to solve the implicit time-stepping equations for the $4th$-order accurate discretizaiton.
   }
\label{tab:scalingSquareOrder4}
\end{table}

%% file: tex/scalingDiskOrder2Table.tex
\begin{table}[hbt]
\begin{center}\tableFontSize
\begin{tabular}{|c|c|c|c|c|c|c|} \hline
      \multicolumn{7}{|c|}{EigenWave scaling: disk, order 2} \\ \hline 
              &   grid-       & eigen-  &  wave-    &  multigrid    &  CPU          &  CPU      \\
 $\ds $       &   points      & pairs   &  solves   &  cycles       &  seconds      &  ratio    \\ \hline
 $1/40   $    &   7.4e+03     &  $16$   &  $78$    &   $  6.0$      &    5.9e+00    &           \\
 $1/80   $    &   2.8e+04     &  $20$   &  $87$    &   $  6.1$      &    1.6e+01    &     2.62       \\
 $1/160  $    &   1.1e+05     &  $16$   &  $78$    &   $  6.9$      &    5.1e+01    &     3.30       \\
 $1/320  $    &   4.2e+05     &  $20$   &  $88$    &   $  7.9$      &    2.0e+02    &     3.99       \\
 $1/640  $    &   1.7e+06     &  $21$   &  $87$    &   $  8.4$      &    8.4e+02    &     4.10       \\
 $1/1280 $    &   6.7e+06     &  $21$   &  $87$    &   $  9.4$      &    2.9e+03    &     3.50       \\
\hline
\end{tabular} 
\end{center}
\vspace*{-1\baselineskip}
\caption{Scaling of EigenWave when using multigrid to solve the implicit time-stepping equations for the $2nd$-order accurate discretization.
   }
\label{tab:scalingDiskOrder2}
\end{table}

%% file: tex/scalingDiskOrder4Table.tex
\begin{table}[hbt]
\begin{center}\tableFontSize
\begin{tabular}{|c|c|c|c|c|c|c|} \hline
      \multicolumn{7}{|c|}{EigenWave scaling: disk, order 4} \\ \hline 
              &   grid-       & eigen-  &  wave-    &  multigrid    &  CPU          &  CPU      \\
 $\ds $       &   points      & pairs   &  solves   &  cycles       &  seconds      &  ratio    \\ \hline
 $1/40   $    &   6.8e+03     &  $18$   &  $86$    &   $  6.0$      &    7.8e+00    &           \\
 $1/80   $    &   2.8e+04     &  $20$   &  $86$    &   $  6.8$      &    2.6e+01    &     3.35       \\
 $1/160  $    &   1.1e+05     &  $21$   &  $88$    &   $  8.7$      &    1.1e+02    &     4.15       \\
 $1/320  $    &   4.2e+05     &  $19$   &  $90$    &   $  8.2$      &    4.2e+02    &     3.87       \\
 $1/640  $    &   1.7e+06     &  $21$   &  $88$    &   $  9.0$      &    1.3e+03    &     3.16       \\
 $1/1280 $    &   6.7e+06     &  $18$   &  $92$    &   $  9.8$      &    6.1e+03    &     4.61       \\
\hline
\end{tabular}   
\end{center}
\vspace*{-1\baselineskip}
\caption{Scaling of EigenWave when using multigrid to solve the implicit time-stepping equations for the $4th$-order accurate discretization.
   }
\label{tab:scalingDiskOrder4}
\end{table}

%% file: tex/conclusions.tex
\section{Conclusions} \label{sec:conclusions}

An algorithm, called EigenWave, has been described for 
computing discrete approximations to the eigenvalues $\lambda_j$ and eigenfunctions $\phi_j$ of elliptic boundary-value problems.
The algorithm is based on solving a related time-dependent wave equation 
whose solution is filtered in time at a target frequency $\omega$.
The EigenWave algorithm defines a new linear operator with the same eigenfunctions as the original
problem but with new eigenvalues, $\beta_j=\beta(\lambda_j,\omega)$, lying in $[-.5,1]$.
The new eigenvalues tend to be largest when the corresponding $\lambda_j$ is close to $\omega$.
Eigenvalues close to $\omega$ are the ones generally found by the algorithm.
When combined with existing high quality Arnoldi-type eigenvalue algorithms through a matrix-free interface, 
the EigenWave algorithm can be used to efficiently compute eigenpairs anywhere in the spectrum without
the need to invert an indefinite matrix, as is common with many alternative approaches.
It was shown that the wave equation can be time-stepped efficiently with an implicit scheme and amazingly only about ten time-steps per period are needed to get good results. Although implicit time-stepping requires solutions of a large matrix system, this matrix system is definite and well suited for solution by fast algorithms such as multigrid. In that case the EigenWave algorithm was shown to scale linearly with the total number of grid points $N$ leading to an optimal $O(N)$ algorithm.

Numerical results in two and three space dimensions were presented for eigenvalue problems in various geometries using finite difference approximations on overset grids.
Both second and fourth-order accurate results were obtained. The results demonstrated that EigenWave can compute multiple discrete eigenpairs to
high accuracy using only a few (e.g.~$3$--$7$) wave-solves per eigenpair.


%% file: tex/discreteAnalysis.tex
\section{Analysis of the discrete EigenWave algorithm}  \label{sec:discreteAnalysis}

An analysis of the continuous EigenWave algorithm was given previously in Section~\ref{sec:analysis}. 
Here, we perform an analysis of the algorithm for the discrete case. 
We focus on an analysis of the implicit scheme since using a large time-step $\dt$ can have potentially significant effects on the
convergence behaviour.
Comments on the analysis for the explicit scheme are made at the end of the section.

Suppose the grid function $W^n_\iv$, $n=0,1,\ldots,\NT$, $\Tbar=\NT\dt$, is given by the solution of the 
implicit scheme~\eqref{eq:implicitScheme}.  
This discrete solution at each time step can be written in terms of an expansion involving eigenvectors 
of the discrete eigenvalue problem~\eqref{eq:eigBVPdiscrete}.
Let $(\lambda_{h,j},\Phi_{j,\iv})$, $j=1,2,\ldots,\NGd$, denote the discrete eigenpairs of~\eqref{eq:eigBVPdiscrete}
where $\NGd$ is the total number of discrete eigenpairs.
The expansion for $W^n_\iv$ takes the form
\ba
   W^n_\iv = \sum_{j=1}^{\NGd} \What_j^n \, \Phi_{j,\iv}. 
\ea
Substituting this expansion into the implicit difference approximation~\eqref{eq:implicitSchemeInterior}
leads to a difference equation for the generalized Fourier coefficient $\What_j^n$
\ba
   & \Dpt\Dmt \What_j^n = - \lambda_{h,j}^2 \, \half\big( \What_j^{n+1} + \What_j^{n-1} \big) , \qquad j=1,2,\ldots,\NGd .  \label{eq:discreteModeEqn}
\ea
Note that the dependence on the spatial order of accuracy $p$ is suppressed for notational convenience.
Initial conditions~\eqref{eq:implicitSchemeIC1} and~\eqref{eq:implicitSchemeIC2} imply
\bse
   \label{eq:WhatIC}
\ba
   & \What_j^0 = \Vhat_j^{(k)}, \\
   & \What_j^1 = \What_j^{-1},
\ea 
\ese
where $\Vhat_j^{(k)}$ is the $j\sp{{\rm th}}$ coefficient in the eigenvector expansion of $V_\iv^{(k)}= v^{(k)}(\xv_\iv)$.
The general solution to the difference equation~\eqref{eq:discreteModeEqn} 
can be found by looking for homogeneous solutions of the form $\kappa^n$, for some as yet unknown constant $\kappa$.
Setting $\What_j^n=\kappa^n$ in~\eqref{eq:discreteModeEqn} leads to a quadratic equation for $\kappa$.  Taking $\What_j^n$ to be a superposition of the two roots of the quadratic, and then applying the initial conditions in~\eqref{eq:WhatIC} results in the solution 
\bse
 \label{eq:cosLambdaTilde} 
\ba
  \What_j^n  = \Vhat_j^{(k)} \, \cos(\lambdaImplicit \, t^n),
\ea
where
\ba
  \lambdaImplicit = \Lam(\lambda_{h,j},\dt),\qquad \Lam(\lambda,\dt)\eqdef \f{2}{\dt} \sin^{-1} \left( \frac{\half \lambda\dt}{\sqrt{1 + \half(\lambda\dt)^2}} \right).
  \label{eq:discreteCorr}
\ea
\ese
The function $\Lam(\lambda,\dt)$ describes the discrete time correction of the eigenvalues for the implicit time-stepping scheme, and it shows, as expected, that $\lambdaImplicit$ is a good approximation to $\lambda_{h,j}$ when 
$\lambda_{h,j}\dt$ is small. \textbf{Note that the adjusted eigenvalues $\lambdaImplicit$ are only used as part of the analysis, 
the eigenvalues computed by EigenWave are the original discrete eigenvalues $\lambda_{h,j}$.}

Following the time-continuous analysis, the discrete time filter~\eqref{eq:filterQuadrature} is applied to the expansion for $W^n_\iv$ with coefficients given by~\eqref{eq:cosLambdaTilde}.  The result is an eigenvector expansion for $\Vhat_\iv^{(k+1)}$ with coefficients given by
\ba
   \Vhat_j^{(k+1)} =\betad(\lambdaImplicit,\omega,\NT)\, \Vhat_j^{(k)} ,\qquad j=1,2,\ldots,\NGd ,
\ea
where $\betad(\lambdaImplicit,\omega,\NT)$ is the discrete WaveHoltz filter function defined by
\ba
   \betad(\lambdaImplicit,\omega,\NT) \eqdef \f{2}{\Tbar} \sum_{n=0}^{\NT} \sigma_{n} \left( \cos(\omega t^n) - \f{\alpha_d}{2} \right) \cos(\lambdaImplicit \, t^n) .
   \label{eq:betad}
\ea
We note that the discrete filter function $\beta_d$ can be written in the same form as~\eqref{eq:filterThreeSincs} using discrete analogues of the sinc functions, see~\cite{overHoltzArXiv2025,overHoltzPartOne,overHoltzPartTwo}.
As before, the discrete EigenWave operator implied by the filter in~\eqref{eq:filterQuadrature} has the same eigenvectors as the discrete operator $L_{ph}$, but with eigenvalues given by
\ba
   \beta_{d,j} \eqdef \betad(\lambdaImplicit,\omega,\NT).
\label{eq:discreteBeta}
\ea

\input tex/discreteBetaFig

Note that if $\dt$ is small (i.e.~$\NT$ is large), and assuming eigenmode $j$ is well resolved on the spatial grid, then the eigenvalues of the discrete EigenWave operator given by~\eqref{eq:discreteBeta} would be close to the values obtained using the continuous WaveHoltz filter function given in~\eqref{eq:filterThreeSincs} and shown in Figure~\ref{fig:waveHoltzBeta}.  However, with implicit time-stepping and large $\dt$ (i.e.~$\NT$ is small), the discrete values can deviate significantly from the exact curve.  Two sources of error 
arise\footnote{A third source of error arises when the implicit time-stepping equations are only solved approximately, see~\ref{sec:implicitTolerance}.}
 (still assuming well-resolved spatial modes), one from the time-stepping as given by $\Lam(\lambda_{h},\dt)$ in~\eqref{eq:discreteCorr} and the other from the approximation of the integral in the filter function.  Define 
\ba
   \betadTilde(\lambda,\omega,\NT) \eqdef \betad\bigl(\Lam(\lambda,\dt),\omega,\NT\bigr),\qquad \dt=\Tbar/\NT ,
\label{eq:correctedBeta}
\ea
as the discrete time-corrected WaveHoltz filter function.
Figure~\ref{fig:betaFunctionI} shows the behavior of $\betadTilde$ in~\eqref{eq:correctedBeta} for different values of $\NT$.  The final time for these plots is $\Tbar=2\pi/\omega$ since $\Np=1$ so that total time-steps $\NT$ is also equal to the number of implicit time-steps per period, which we denote by~$\Nits$.  
The curves of $\betadTilde$ can be compared to the exact curve of $\beta$ shown in black.
The curves for $\NT=\Nits
=4$ and~$5$ have a very broad main peak which would generally make it difficult to compute eigenvalues as there would be many values relatively close to one.
However, as $\NT=\Nits$ increases, the main peak of the time-corrected filter function narrows, and it better approximates the main peak of the exact filter function. 
Fortuitously, the filter function given by $\betadTilde$ is relatively small for large~$\lambda/\omega$.
Depending on the number of eigenvalues desired, values of $\Nits$ between $6$ and $15$ may be appropriate to use.
There does not appear to be much benefit to using values of $\Nits$ larger than~$15$.
Many of the numerical results in Section~\ref{sec:numericalResults} use $\Nits=10$, and the behavior shown in Figure~\ref{fig:betaFunctionI} explains why good results can be expected.
Note, however, that $\betadTilde(\lambda,\omega,\NT)$ in Figure~\ref{fig:betaFunctionI} 
reaches its maximum for $\lambda>\omega$. 
Thus, if eigenvalues near $\omega$ are desired then
the target frequency $\omega$ should be reduced somewhat to $\omegaTilde$ so that $\betadTilde(\lambda,\omegaTilde,\NT)=1$ for $\lambda=\omega$.
In particular $\omegaTilde$ should satisfy
\ba
  \Lam(\omega; \dtTilde) = \omegaTilde,
\ea
where $\dtTilde$ is the adjusted time-step.
That is 
\ba
	\label{eq:modifiedFrequencyEqn}
	\dfrac{2}{\dtTilde} \sin^{-1}\left(\dfrac{(\omega \dtTilde) / 2}{\sqrt{1 + (\omega \dtTilde)^2 / 2}}\right)= \omegaTilde.
\ea
Note that changing $\omega$ changes $\dt$ and $T_f$, but $\dt/T_f$, $\omega \dt$ and $\alpha_d$ ($\alpha_d$ only depends on $\omega\dt$) don't change since 
\ba
  \omega \dt = \omegaTilde \dtTilde  = \dfrac{2\pi N_p}{N_t} = \dfrac{2 \pi }{\Nits}, \text{ and } \dfrac{\dt}{T_f} =\dfrac{\dtTilde}{\tilde{T}}=\dfrac{1}{N_t}.
\ea
Solving for $\omegaTilde$ in \eqref{eq:modifiedFrequencyEqn} gives 
\ba
  \omegaTilde = \dfrac{\omega \pi}{\Nits} \sqrt{\dfrac{1-2\sin^2(\pi/\Nits)}{\sin^2(\pi/\Nits) }}. \label{eq:adjustedOmega}
\ea
Note that with this new $\omegaTilde$ 
we can still rewrite the functions $\betad(\lambda,\omegaTilde,\NT) $ and $\betadTilde(\lambda,\omegaTilde,\NT)$ as functions of $\lambda / \omega$ (rather than functions of $\lambda / \omegaTilde$ ).
The right graphs of Figure~\ref{fig:betaFunctionI} and~\ref{fig:betaFunctionII} show the corrected curves.

\input tex/discreteBetaNp2Fig

Arnoldi-type algorithms can be used to solve for multiple eigenpairs at once. 
It is observed that for good convergence of the Arnoldi algorithm, the requested number of eigenpairs 
should be some fraction of the largest eigenvalues that appear in the primary peak of the filter curve; this is discussed further 
in~\ref{sec:numberOfRequestedEigenpairs}.
If fewer eigenvalues near $\omega$ are desired (such as to avoid excessive storage requirements), 
then the number of periods, $\Np$, 
can be increased and this narrows the main peak of the filter curve.
As shown in the plots for $\Np=2$ in Figure~\ref{fig:betaFunctionII}
one should still choose about the same number of time-steps per period as with $\Np=1$ to avoid the main peak of $\betadTilde$
becoming too broad. In this case, choosing $\Nits=10$ so that $\NT = \Nits \Np = 20$, for example, would give a main peak in the curve for $\betadTilde$ that is about the same width as the
exact filter function.

\medskip
The steps in the analysis of the explicit time-stepping scheme are similar to that for implicit time-stepping. 
The main difference is that the adjusted eigenvalue in~\eqref{eq:discreteCorr} is instead
\ba
   \lambdaExplicit \eqdef \f{2}{\dt} \sin^{-1}\left(   \f{\lambda_{h,j}\dt}{2} \right). \label{eq:discreteExplicitLambdaTilde}
\ea 
Note that, for stability, the explicit time-stepping scheme must satisfy a CFL-type restriction~\cite{overHoltzArXiv2025,overHoltzPartOne,overHoltzPartTwo} with $\Nt$ increasing as the grid spacing goes to zero.
Thus, with explicit time-stepping the number of time-steps $\NT$ is usually large enough so that the discrete filter function $\beta_d$ lies very close
to the continuous one near the target frequency $\omega$.

%% file: tex/discreteBetaFig.tex
{
\newcommand{\figw}{5.25cm}
\newcommand{\figh}{4.0cm}
\begin{figure}[htb]
\begin{center}
\begin{tikzpicture}
  \useasboundingbox (0,.6) rectangle (3.*\figw,1.*\figh);  
  \begin{scope}[yshift=0*\figh]
    \figByWidth{         0}{0}{fig/betaDiscreteNt4To9}{\figw}[0.][0.][0.][0.]
    \figByWidth{1.0*\figw}{0}{fig/betaDiscreteNt10To30}{\figw}[0.][0.][0.][0.]
    \figByWidth{2.0*\figw}{0}{fig/betaDiscreteNt5To10Adjusted}{\figw}[0.][0.][0.][0.]
  \end{scope}  
\end{tikzpicture}
\end{center}
\caption{
  Discrete WaveHoltz filter $\betadTilde$ for implicit time-stepping.
  Curves are shown for different numbers of time-steps $\Nt=\Np\Nits$, with number of periods $\Np=1$.
  The continuous filter is shown in black. The vertical dashed line shows $\lambda/\omega=1$.
  Right: The maximum occurs at $\lambda=\omega$ when using the adjusted $\omegaTilde$ in~\eqref{eq:adjustedOmega}.
  }
\label{fig:betaFunctionI}
\end{figure}
}

%% file: tex/discreteBetaNp2Fig.tex
{
\newcommand{\figw}{5.25cm}
\newcommand{\figh}{4.0cm}
\begin{figure}[htb]
\begin{center}
\begin{tikzpicture}
  \useasboundingbox (0,.75) rectangle (3*\figw,1.*\figh);  
  \begin{scope}[yshift=0*\figh]
    \figByWidth{         0}{0}{fig/betaDiscreteNt8To18Np2}{\figw}[0.][0.][0.][0.]
    \figByWidth{1.00*\figw}{0}{fig/betaDiscreteNt20To60}{\figw}[0.][0.][0.][0.]
    \figByWidth{2.00*\figw}{0}{fig/betaDiscreteNt10To20Np2Adjusted}{\figw}[0.][0.][0.][0.]
  \end{scope}  
\end{tikzpicture}
\end{center}
\caption{
 Discrete WaveHoltz filter $\betadTilde$ for implicit time-stepping.
  Curves are shown for different numbers of time-steps $\Nt=\Np\Nits$, with number of periods $\Np=2$.
  The continuous filter is shown in black. The vertical dashed line shows $\lambda/\omega=1$.
  Right: The maximum occurs at $\lambda=\omega$ when using the adjusted $\omegaTilde$ in~\eqref{eq:adjustedOmega}.
   }
\label{fig:betaFunctionII}
\end{figure}
}

%% file: tex/arnoldi.tex

\section{Using Arnoldi-based algorithms to compute multiple eigenvalues}   \label{sec:Arnoldi}

EigenWave can be combined with existing Krylov-based eigenvalue solvers, e.g.~those based on the Arnoldi method,
as was discussed briefly in Section~\ref{sec:eigSoftware}.
Here we provide further details of the use of Krylov-based algorithms,
such as those available with the packages SLEPc and ARPACK, in combination with EigenWave.
It has been found that the Krylov-Schur algorithm from SLEPc and the 
IRAM algorithm from ARPACK give similarly
good results for the cases that have been computed for this paper.

Some Krylov-based eigenvalue algorithms, such as the ones noted above, do not require the matrix 
in explicit form, but can make use of a black-box routine that computes a matrix-vector product.
Such algorithms are called \textsl{matrix-free} for the current discussion.
In particular, since
the EigenWave algorithm does {\it not} require the inverse of a shifted Laplacian
(which would generally require the matrix in explicit form),
matrix-free Arnoldi algorithms using a matrix-vector product routine are very convenient.

\mni
\paragraph{EigenWave Approach} The desired eigenpairs $(\lambda_{h,j},\phi_{h,j})$ of the discretized 
BVP for $L_{ph}$ can be computed in two stages as follows. 

\mni
Stage 1. Choose a target frequency $\omega$. Using a matrix-free eigenvalue algorithm,  
compute the eigenpairs $(\beta_j,\phi_{h,j})$ of the BVP for the WaveHoltz iteration
operator $\Aw_{ph}=\Aw_{ph}(\omega)$. This operator has the same eigenvectors as $L_{ph}$ 
but has eigenvalues $\beta_j=\betadTilde(\lambda_{h,j};\omega)$.

\newcommand{\na}{{n_a}}
\mni
Stage 2. Given approximate eigenvectors, $\phi_{h,j}$, compute approximations to the true discrete
eigenvalues $\lambda_{h,j}$ of $L_{ph}$ 
using the Rayleigh quotient (although the simpler formula~\eqref{eq:simpleRayleighQuotient} can often be used in practice)
\ba
    -\lambda_{h,j}^2  = \frac{  (\phi_{h,j}, L_{ph} \phi_{h,j} )_h }{ (\phi_{h,j}, \phi_{h,j} )_h }.  \label{eq:discreteRayleighQuotient}
\ea
Here $(\cdot,\cdot)_h$ denotes a discrete approximation to the $L_2$ inner product.

\smallskip
To give a concrete outline of the approach, 
we describe the basic Arnoldi algorithm for computing eigenvalues.
Let $A\in\Real^{\na\times \na}$ denote a matrix whose eigenvalues we wish to find, 
and $\vv\in\Real^\na$ a starting vector. Here $\na$ is the total number of active grid points (see comments below).
After $m$ steps, the algorithm generates a rectangular matrix $V_m\in\Real^{\na\times m}$ with columns $\vv_i$, $i=1,2,\ldots,m$,
 and a square matrix $H_m\in\Real^{m\times m}$
which satisfies $H_m=V_m^* A V_m$. The steps in this process are shown in Algorithm~\ref{alg:arnoldi}. 
The matrix $H_m$ is an orthogonal projection of $A$ onto the space spanned by the columns of $V_m$, which span the $m$-dimensional Krylov space $\Kc_m=\text{span}\{\vv,A\vv,A^2\vv,\ldots,A^{m-1}\vv\}$ 
generated by $A$ and starting vector $\vv$.
The eigenvalues of $H_m$ are the Ritz estimates and some of these may be good approximations
to the eigenvalues of $A$. The Arnoldi algorithm tends to converge fastest to the extreme eigenvalues, largest and smallest in magnitude.
This can be seen in some of the numerical results presented in Section~\ref{sec:numericalResults}, see Figure~\ref{fig:doubleEllipseGridFig} for example, and also the plots presented later in Figure~\ref{fig:squareEigsVaryArnoldi}.

\renewcommand{\algFontSize}{\small}
\begin{algorithm}
\algFontSize 
\caption{Arnoldi algorithm.}
\begin{algorithmic}[1]
  \Function{$[\Vv_m,\Hv_m]=$Arnoldi}{$\Av$, $\vv$, $m$ }  
    \State $\vv_1=\vv/\| \vv\|_2$ \Comment Normalize starting vector, first column of $\Vv_m$
    \For{$j=1,2,\ldots,m$} 
      \State $\wv = \Av \vv$  \Comment Matrix-vector multiply
      \For{$i=1,2,\ldots,j$}        \Comment Gram-Schmidt orthogonalization
        \State $h_{ij} = \wv^*\vv_i$ \Comment Entry in the matrix $\Hv_m$
        \State $\wv = \wv - h_{ij} \vv_i$ 
      \EndFor
      \State $h_{j+1,i} = \| \wv \|_2$ 
      \State $\vv_{j+1} = \wv/h_{j+1,j}$
    \EndFor
  \EndFunction
\end{algorithmic} 
\label{alg:arnoldi}
\end{algorithm}

The basic Arnoldi algorithm has been improved in many important ways~\cite{BaiDemmel2000}.
For example, the implicitly restarted Arnoldi algorithm (IRAM) of Sorensen~\cite{Sorensen1992,BaiDemmel2000},
implemented in ARPACK, has been very successful. 
The Krylov-Schur algorithm, found in SLEPc, 
was introduced by Stewart~\cite{Stewart2002} and contains some improvements to the 
IRAM algorithm.


\renewcommand{\algFontSize}{\small}
\begin{algorithm}
\algFontSize 
\caption{Implicitly Restarted Arnoldi Method.}
\begin{algorithmic}[1]
%
    \State $\vv_1 =\vv/\| \vv \|_2;$ \Comment starting vector
    \State $\Av \Vv_\mTotal = \Vv_\mTotal\Hv_\mTotal + \fv_\mTotal \ev_\mTotal^*$ \Comment compute an $\mTotal$-step Arnoldi factorization
    \While{not converged} 
      \State Compute eigenvalues of $\Hv_\mTotal$   
       and choose  shifts $\mu_j$, $j=1,\ldots,\pExtra$
      \State $\Qv = \Iv_\mTotal$
      \For{$j=1,2,\ldots,\pExtra$} \Comment Perform $\pExtra$-steps of shifted QR algorithm
        \State $\Qv_j \Rv_j = \Hv_\mTotal - \mu_j \Iv_\mTotal$ \Comment Compute QR factorization of $\Hv_\mTotal - \mu_j \Iv_\mTotal$
        \State $\Hv_\mTotal = \Qv_j^* \Hv_\mTotal \Qv_j$
        \State $\Qv=\Qv\, \Qv_j$
      \EndFor
      \State $\beta_\kReq=\Hv_\mTotal(\kReq+1, \kReq)$; $\sigma_\kReq=\Qv(\mTotal, \kReq)$  
      \State $\fv_\kReq=\vv_{\kReq+1}\beta_\kReq +\fb_\mTotal \sigma_\kReq$
      \State $\Vv_\kReq =\Vv_\mTotal\Qv(:, 1: \kReq)$
      \State $\Hv_\kReq=\Hv_\mTotal(1:\kReq,1:\kReq)$
      \State beginning with the $\kReq$-step of Arnoldi factorization \\
     \hspace*{3.5em}   $\Av\Vv_\kReq = \Vv_\kReq\Hv_\kReq+\fv_\kReq\ev^*_\kReq$
     \State   apply $\pExtra$ additional steps of Arnoldi procedure to obtain a new $\mTotal$-step Arnoldi factorization\\
     \hspace*{3.5em} $\Av\Vv_\mTotal = \Vv_\mTotal \Hv_\mTotal + \fv_\mTotal\ev^*_\mTotal$ 
    \EndWhile   
\end{algorithmic} 
\label{alg:IRAM}
\end{algorithm}

To help better understand the numerical results presented previously in Section~\ref{sec:numericalResults}, and also the properties of the EigenWave algorithm discussed later in~\ref{sec:properties}, a brief description of the IRAM approach is now given in Algorithm~\ref{alg:IRAM}.
The number of vectors in the basic Arnoldi algorithm keeps increasing as the number of iterations increases, and
this can lead to excessive storage requirements, among other problems.
The restarted Arnoldi algorithm, however, maintains a maximum of $\Ntotal$ Arnoldi vectors, 
and so has fixed storage requirements.
When $\Nreq$ eigenvalues are requested, the IRAM algorithm maintains $\Ntotal=\Nreq + \Nextra$ Arnoldi vectors,
where the number of additional vectors, $\Nextra$, is typically chosen as $\Nextra\approx \Nreq$. 
The algorithm periodically updates the $\Nreq$ Arnoldi vectors (in a process known as restarting) so that they
converge to the desired eigenvectors. It does this by using the classical implicitly-shifted QR algorithm on the
small matrix $\Hv_m$ (of dimension $\Ntotal\times \Ntotal$) generated from the Arnoldi algorithm, i.e.~Algorithm~\ref{alg:arnoldi}.
We sort the eigenvalues of $\Hv_m$ into a ``wanted set" $\{\theta_1,\cdots, \theta_{\kReq}\}$ and an ``unwanted set" $\{\mu_1,\cdots,\mu_{\pExtra}\}$ and take the shifts to be the unwanted ones.
For EigenWave the wanted set is set of eigenvalues with largest magnitude.
The QR shifts are chosen to remove the unwanted eigenvalues. 
For more details of IRAM and the Krylov-Schur algorithms, 
see the discussion in~\cite{BaiDemmel2000} and~\cite{Stewart2002}.

\renewcommand{\algFontSize}{\small}
\begin{algorithm}
\algFontSize 
\caption{WaveHoltz matrix-vector product for matrix-free Arnoldi.}
\begin{algorithmic}[1]
  \Function{$\zv$={\mvCol matVec}}{ $\yv$ }  \Comment Compute $ \zv = \Aw_h \yv$.
    \State $\mu=0;$ \Comment Counts entries in $\yv$
    \For{$\iv \in\Active$ (active set of grid points)} \Comment Convert vector $\yv$ to grid-function $V_\iv$
      \State $V_\iv = y_\mu;$ \quad $\mu=\mu+1$;
    \EndFor
    \State $\Vv =$ \Call{applyBoundaryConditions}{$\Vv$}; \Comment Interpolate and apply BC's
    \State $\Vv =$ \Call{TakeOneWaveHoltzStep}{$\Vv$};
    \State $\mu=0;$ \Comment Counts entries in $\zv$
    \For{$\iv \in\Active$ (active set of grid points)} \Comment Convert grid-function $V_\iv$ to vector $\zv$
      \State $z_\mu = V_\iv$; \quad $\mu=\mu+1$;
    \EndFor    
 \EndFunction
\end{algorithmic} 
\label{alg:matVec}
\end{algorithm}

Algorithm~\ref{alg:matVec} outlines the form of the matrix-vector product routine
we use with Arnoldi-type methods. In what follows, \textsl{active} points, denoted by $\Active$, are those where the interior equation is applied (i.e.~these are the equations that would have a $\lambda$ in them, when posed as a eigenvalue problem). \textsl{Inactive} points are all others, such as boundary points, ghost points and interpolation points.\footnote{An overset grid consists of a collection of curvilinear component grids that cover the problem domain and overlap where they meet.  Interpolation is used to communicate solution values at grid points in regions of overlap, see~\cite{CGNS,max2006b}.}
Assume that the vector form of the active points in the solution is given by $\yv \in\Real^\na$.
First, $\yv$ is converted to a grid function $V_\iv$ with sufficient storage for both
active and inactive points, although inactive points are as yet undefined. After applying boundary conditions to $V_\iv$, which defines all inactive points, this grid function
is used as initial condition for the EigenWave algorithm. After one wave-solve step,
the active portions are then converted back to vector form, $\zv\in\Real^\na$.

\mni
\textbf{Note 1.}
An important note is that constraint equations, such as
boundary conditions and interpolation equations, are not included in the vectors used
by Arnoldi. This means the Arnoldi algorithm solves a regular eigenvalue problem 
of the form $\Aw_h \yv = \beta \yv$. Using this approach, it is very easy to eliminate constraint equations from
the eigenvalue problem.

\mni
\textbf{Note 2.} 
Our numerical results show that very accurate approximations to the discrete eigenvectors $\phi_{\iv,j}$ are often found.
In this case the discrete inner products in the Rayleigh quotient~\eqref{eq:discreteRayleighQuotient} 
can be replaced by simple unweighted sums over the set of active points, $\Active$, (or even some subset of the active points)
\ba
   -\lambda_{h,j}^2  = \frac{  \sum_{\iv\in\Active} \phi_{\iv,j} \, L_{ph}\, \phi_{\iv,j}  }{ \sum_{\iv\in\Active} (\phi_{\iv,j})^2 },  \label{eq:simpleRayleighQuotient} 
\ea
since  $L_{ph} \phi_{\iv,j}$ is a very accurate approximation to $ -\lambda_{h,j}^2 \phi_{\iv,j} $.
This is useful for overset grids since forming the weights for discrete inner products requires some work~\cite{overHoltzArXiv2025,overHoltzPartOne,overHoltzPartTwo}. 

\bigskip

\begin{table}[hbt]
\begin{center}\tableFontSize
\begin{tabular}{|c|c|c|c|c|c|c|c|} \hline
 \multicolumn{8}{|c|}{EigenWave: grid=square128, ts=explicit, order=4, $\omega=15$, $N_p=1$, KrylovSchur } \\ \hline 
   num   &  wave       & time-steps  & wave-solves &  time-steps &   max      &  max      &  max       \\ 
   eigs  &  solves     & per period  &  per eig    &  per-eig    &   eig-err  & evect-err & eig-res    \\ 
 \hline
    26    &       91     &      10     &    3.5    &     343     &  1.67e-13 &  2.38e-13 &  9.71e-12 \\
 \hline 
\end{tabular}
\end{center}
\vspace*{-1\baselineskip}
\caption{Summary of EigenWave performance for \texttt{square128} grid using the KrylovSchur algorithm and explicit time-stepping.  The spatial order of accuracy is 4 and the wave-solves use $N_p=1$ to determine the final time.}
\label{tab:square128ExpEigKrylovSummary}
\end{table}

We conclude this section by illustrating the convergence behaviour of the Krylov-Schur algorithm, which like the IRAM algorithm discussed above,  
consists of a sequence of iterations where the Krylov subspace is first expanded and then contracted.
Table~\ref{tab:square128ExpEigKrylovSummary} summarizes results of computing eigenpairs to 
fourth-order accuracy on a square domain using explicit time-stepping and the grid \texttt{square128}.
The target frequency is $\omega=15$. Twenty-four eigenpairs are requested and twenty-six eigenpairs are found by the KrylovSchur
algorithm in a total of $91$ wave-solves. This corresponds to approximately 
$3.5$ wave-solves per eigenpair found.
The top graph in Figure~\ref{fig:squareEigsOrder4KrylovSchur} shows the filter function and computed eigenvalues, while the bottom plots in the figure show the convergence behavior of the Krylov-Schur algorithm.
The bottom left graph shows the Arnoldi estimated eigenvalues $\beta_j$ at each iteration, while
the bottom right graph shows the estimated errors in $\beta_j$.
Note that the values of the converged $\beta_j$ in the lower left graph correspond to the 
heights of the circled eigenvalues in the upper graph.
The number of requested eigenvalues is $N_r=24$, while the Krylov-Schur algorithm keeps $N_e=25$ additional Arnoldi vectors in the Krylov space
for a total of $N_a=49$ vectors in the restarted Arnoldi scheme.

\input tex/squareEigsKrylovSchurConvFig.tex

From the output of the Krylov-Schur algorithm and the graphs of the convergence behavior in Figure~\ref{fig:squareEigsOrder4KrylovSchur}, we make the following observations:
\begin{enumerate}
  \item The first iteration takes place after $48$ waves-solves at which point $49$ Arnoldi vectors are known (including the initial condition). There are no converged eigenpairs at this point (the convergence tolerance is set to $10^{-14}$).
  \item The second iteration occurs after $73$ wave-solves with $16$ converged eigenpairs.
  \item The third iteration occurs after $90$ wave-solves with $26\approx N_r$ converged eigenpairs.
\end{enumerate}

%% file: tex/squareEigsKrylovSchurConvFig.tex
{
\newcommand{\figw}{6.5cm}
\newcommand{\figh}{5.5cm}

\begin{figure}[htb]
\begin{center}
\begin{tikzpicture}[scale=1]
  \useasboundingbox (0,.7) rectangle (14,11);  

  \begin{scope}[yshift=\figh*1.05]
     \figByWidth{0.0}{0}{fig/square128ExpEigKrylov}{\figw}[0][0][0][0]
     \draw (1.5,1) node[draw,fill=white,anchor=west,xshift=0pt,yshift=0pt] {\scriptsize Square128, E, O4, $\Np=1$};
  \end{scope}

   \begin{scope}[yshift=0cm]
     \figByWidth{0.0}{0}{fig/square128ExpEigKrylovEigEstimates}{\figw}[0][0][0][0]
     \figByWidth{7.0}{0}{fig/square128ExpEigKrylovEigErrEstimates}{\figw}[0][0][0][0]
  \end{scope}

\end{tikzpicture}
\end{center}
\caption{Square: computing multiple eigenpairs, order=4, explicit time stepping. 
Top left: locations of the computed eigenvalues $\lambda_j$ are marked with black circles. 
Bottom left: convergence of the (Arnoldi estimated) eigenvalues $\beta_j=\beta(\lambda_j)$ 
of the \EigenWaveb operator $\Aw_{ph}$, being computed by Krylov-Schur algorithm.
Bottom right: estimated errors in the Arnoldi estimates for $\beta_j$. 
  } 
\label{fig:squareEigsOrder4KrylovSchur}
\end{figure}
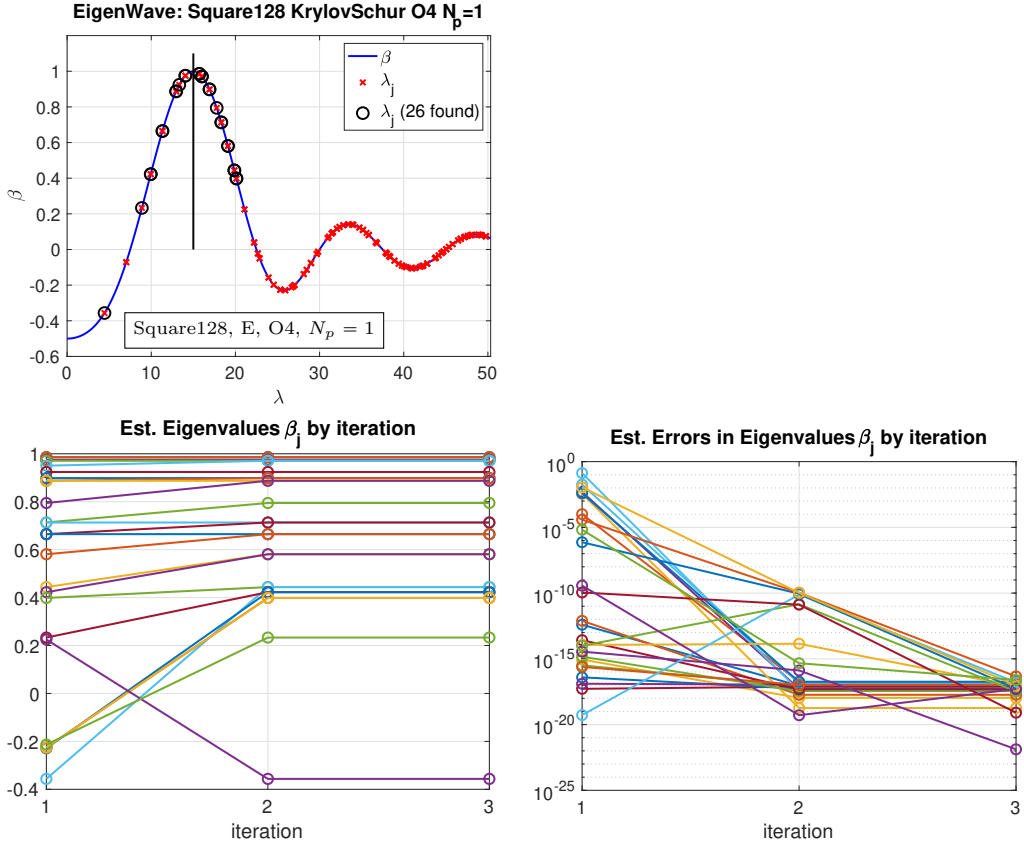
}

%% file: tex/properties.tex
\section{Properties of EigenWave and the discrete approximations} \label{sec:properties}

The discussion in this section considers various properties of the EigenWave algorithm and the discrete approximations.
These are:
\begin{description}
  \item[\ref{sec:accuracyOfEigenpairs}] shows that the computed eigenpairs converge at the expected order accuracy to the continuous eigenpairs.
  \item[\ref{sec:implicitTolerance}] shows that EigenWave behaves robustly as the iterative solution of the implicit time-stepping 
    equation is solved to different tolerances.
  \item[\ref{sec:implicitVersusExplicitComparison}] compares the CPU cost between explicit and implicit time-stepping and shows that for large enough $N$
     the implicit scheme using the $O(N)$ multigrid solver will be the fastest.
  \item[\ref{sec:numberOfRequestedEigenpairs}] shows how the computational cost depends on the number of requested eigenpairs and (if memory is available)
  then the number of wave-solves per eigenpair approaches $2$ as the number of requested eigenpairs gets large.
  \item[\ref{sec:changingImplicitStepsPerPeriod}] varies the number of implicit time-step per period, $\Nits$, and shows that 
      $\Nits\approx 10$ is a good choice to approximately minimize the computational cost.
  \item[\ref{sec:changingNp}] shows how to choose the number or filter periods, $\Np$, as a function of the number of requested eigenvalues $N_r$ so as
     to minimize the computational cost.
  \item[\ref{sec:estimatingTheConvergenceRate}] attempts to theoretically explain the convergence behaviour of the (sophisticated) Krylov-Schur algorithm found in~\ref{sec:changingNp} 
  by studying the convergence behaviour of the (simpler) simultaneous iteration algorithm.
\end{description}


\input tex/accuracyOfEigenpairs

\input tex/implicitTolerance

\input tex/implicitVersusExplicit

\input tex/convergenceProperties

\subsection{Changing the number of implicit time-steps per period} \label{sec:changingImplicitStepsPerPeriod}

The cost of solving the wave equation using implicit time-stepping depends on the number of time-steps taken. As shown in Figure~\ref{fig:betaFunctionI}, the number of implicit time-steps per period, given by $\Nits$, should be large enough so that the time-filter is effective. 
With too few implicit time-steps per period, the convergence of the EigenWave algorithm may degrade.
Here we study the effect of changing the number of implicit time-steps per period.
Figure~\ref{fig:varyTimeStepsPerPeriod} presents results for two cases.
In the first case, EigenWave is used on a square domain with the second-order accurate discretization, target frequency $\omega=10$, $\Np=1$ periods, and $12$ requested eigenvalues. 
In the second case, EigenWave is used on a disk domain using the fourth-order accurate discretization, target frequency $\omega=15$, $\Np=2$ periods, and $50$ requested eigenvalues. 
The metric used to measure the cost is 
the number of time-steps per eigenvalue.
The number of time-steps per period, $\Nits$, is varied from $5$ (the minimum needed for stability) to $15$.
As can be seen from Figure~\ref{fig:varyTimeStepsPerPeriod} the number of time-steps needed per eigenvalue 
has a minimum of about $47$ with $\Nits=11$ for the square, and $28$ with $\Nits=10$ for the disk.
The graph in the figure also shows the actual number of eigenvalues found. 
For example, $66$ eigenvalues were found ($50$ were requested) for the disk with $\Nits=10$.
These results suggest that $\Nits=10$ may be a reasonable value to choose for the number of implicit time-steps per period.

{
\newcommand{\figw}{6.5cm}
\newcommand{\figh}{5.5cm}

\begin{figure}[htb]
\begin{center}
\begin{tikzpicture}[scale=1]
  \useasboundingbox (0,.7) rectangle (7,5.5);  

  \begin{scope}[xshift=0cm]
    \figByWidth{0.0}{0}{fig/varyTimeStepsPerPeriod}{\figw}[0][0][0][0]
  \end{scope}

\end{tikzpicture}
\end{center}
\caption{Performance of the EigenWave algorithm as a function of the number of implicit time-steps per period, $\Nits$.
Case~1 (in blue) shows results for the square domain, the $2$nd order code, target frequency $\omega=10$, $\Np=1$ periods, and $12$ requested eigenvalues. 
Case 2: (in red) shows results for a disk domain, the $4$th order code, target frequency $\omega=15$, $\Np=2$ periods, and $50$ requested eigenvalues. 
The curves with shaded symbols show the total number of time-steps per eigenvalue needed.
The curves with unshaded symbols indicate the actual number of converged eigenpairs computed. 
The cost is minimized for $\Nits=11$ (case 1) and $\Nits=10$ (case 2). 
}
\label{fig:varyTimeStepsPerPeriod}
\end{figure}
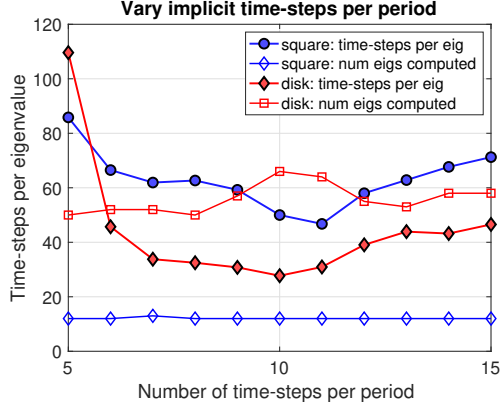
}

\input tex/varyNp.tex

\input tex/convergenceRateEstimation.tex

%% file: tex/accuracyOfEigenpairs.tex
\subsection{Accuracy of the reference discrete eigenvalues and eigenvectors}  \label{sec:accuracyOfEigenpairs}

The SLEPSc package is used to compute the ``true'' discrete eigenpairs (called \textsl{reference} eigenpairs) that are used for computing the errors in the eigenpairs computed using EigenWave. 
The discrete spatial approximations used in the SLEPSc code are the same as those used when solving the wave 
equation. As a result, it is expected that upon convergence EigenWave will give the same discrete results
as the reference values computed with SLEPc.
The reference eigenpairs are computed using the Krylov-Schur algorithm from SLEPSc~\cite{SLEPc2005}.
The procedure for doing this is described in~\cite{overHoltzArXiv2025,overHoltzPartOne} and discussed here briefly  for completeness.
The discrete approximation to the eigenvalue problem on an overset grid consists of approximations to the PDE and boundary
conditions together with interpolation equations. This is a generalized eigenvalue problem of the form $A x = \lambda B x$, since the eigenvalue does not
appear in the boundary conditions and interpolation equations. The matrix $B$ has ones on the diagonal for points where the PDE is discretized and zeros for constraint equations.
It is possible, in principle, to eliminate all constraint equations
and reduce the problem to a regular eigenvalue problem of the form $A x = \lambda x$ for a reduced matrix $A$.
For practical reasons, however, it is convenient to retain the constraint equations.
The algorithms in SLEPSc seem to work best if the matrix $B$ in the generalized form is nonsingular.
In the overset grid setting $A$ is nonsingular while $B$ is singular. To resolve this issue, the roles of $A$ and $B$ are reversed and instead we solve a related generalized eigenvalue problem $B x = (1/\lambda) A x $ for the reciprocals of the eigenvalues.
The eigenvectors returned from SLEPSc are normalized using the discrete inner product. 
For any multiple eigenvalues, an orthonormal basis for the corresponding eigen-space is found. 
It should be noted, however, that computation of the eigenmodes using SLEPSc requires the inversion of a large (often indefinite) matrix,
and generally we use a direct sparse solver to do this. This can be expensive for large problems. The EigenWave algorithm avoids the need to invert an indefinite matrix.

\input tex/diskEigenAccuracyFig.tex
\input tex/sphereEigenAccuracyFig.tex

Before proceeding to a discussion of the behavior of the EigenWave algorithm, the accuracy of these reference
eigenpairs is examined with respect to the continuous problem.  This is done by comparing the computed eigenvalues and eigenvectors with the exact continuous values
for a two-dimensional disk and a three-dimensional solid sphere.
The overset grids for the disk and solid sphere are described in Sections~\ref{sec:diskEigenpairs} and~\ref{sec:sphereEigenpairs}, respectively.
Formulas for the exact eigenvalues and eigenfunctions of the disk with Dirichlet boundary conditions 
are given in~\eqref{eq:diskTrueEigenpairs}, the corresponding formulas for the eigenpairs of the solid sphere, also with Dirichlet boundary conditions, are given in~\eqref{eq:sphereTrueEigenpairs}.
Note that both problems have eigenvalues with multiplicities larger than one (the sphere has eigenvalues with arbitrarily large multiplicities)
which could cause difficulties for numerical eigenvalue algorithms.
A grid refinement study is performed and relative max-norms errors in the eigenvalues and eigenvectors are computed.

Figures~\ref{fig:diskAccuracyOfEigenpairsOrder2} and~\ref{fig:diskAccuracyOfEigenpairsOrder4} 
show results for the eigenpairs of the disk corresponding to the smallest 10 eigenvalues to orders of accuracy two and four respectively.
Figures~\ref{fig:sphereAccuracyOfEigenpairsOrder2} and~\ref{fig:sphereAccuracyOfEigenpairsOrder4} 
show corresponding results for the solid sphere.
Note (see legends) that the results for the disk include five double eigenvalues while
results for the sphere include eigenvalues of multiplicity three and five.
In all cases the eigenvalues are seen to converge at close
to expected rates. The eigenvectors are converging at least as fast as the expected rate; 
some eigenvectors appear to be converging at a rate one order higher than expected.

%% file: tex/diskEigenAccuracyFig.tex
{
\newcommand{\figw}{5.5cm}
\newcommand{\figh}{5.5cm}

\begin{figure}[htb]
\begin{center}
\begin{tikzpicture}[scale=1]
  \useasboundingbox (0,.8) rectangle (2.2*\figw,1.05*\figh);  

  \begin{scope}[xshift=0cm]
    \figByWidth{0.0}{0}{fig/diskOrder2EigenvalueErrors}{\figw}[0][0][0][0]
  \end{scope}
  \begin{scope}[xshift=1.1*\figw]
    \figByWidth{0.0}{0}{fig/diskOrder2EigenvectorErrors}{\figw}[0][0][0][0]
  \end{scope}  

\end{tikzpicture}
\end{center}
\caption{Accuracy of computed eigenvalues and eigenvalues on a disk (order of accuracy 2) compared to the exact continuous values.
}
\label{fig:diskAccuracyOfEigenpairsOrder2}
\end{figure}
}

{
\newcommand{\figw}{5.5cm}
\newcommand{\figh}{5.5cm}

\begin{figure}[htb]
\begin{center}
\begin{tikzpicture}[scale=1]
  \useasboundingbox (0,.8) rectangle (2.2*\figw,1.05*\figh);  

  \begin{scope}[xshift=0cm]
    \figByWidth{0.0}{0}{fig/diskOrder4EigenvalueErrors}{\figw}[0][0][0][0]
  \end{scope}
  \begin{scope}[xshift=1.1*\figw]
    \figByWidth{0.0}{0}{fig/diskOrder4EigenvectorErrors}{\figw}[0][0][0][0]
  \end{scope}  

\end{tikzpicture}
\end{center}
\caption{Accuracy of computed eigenvalues and eigenvalues on a disk (order of accuracy 4) compared to the exact continuous values.
}
\label{fig:diskAccuracyOfEigenpairsOrder4}
\end{figure}
}

%% file: tex/sphereEigenAccuracyFig.tex
{
\newcommand{\figw}{5.5cm}
\newcommand{\figh}{5.5cm}

\begin{figure}[htb]
\begin{center}
\begin{tikzpicture}[scale=1]
  \useasboundingbox (0,.8) rectangle (2.2*\figw,1.05*\figh);  

  \begin{scope}[xshift=0cm]
    \figByWidth{0.0}{0}{fig/sphereOrder2EigenvalueErrors}{\figw}[0][0][0][0]
  \end{scope}
  \begin{scope}[xshift=1.1*\figw]
    \figByWidth{0.0}{0}{fig/sphereOrder2EigenvectorErrors}{\figw}[0][0][0][0]
  \end{scope}  

\end{tikzpicture}
\end{center}
\caption{Accuracy of computed eigenvalues and eigenvalues on a solid sphere (order of accuracy 2) compared to the exact continuous values.
}
\label{fig:sphereAccuracyOfEigenpairsOrder2}
\end{figure}
}

{
\newcommand{\figw}{5.5cm}
\newcommand{\figh}{5.5cm}

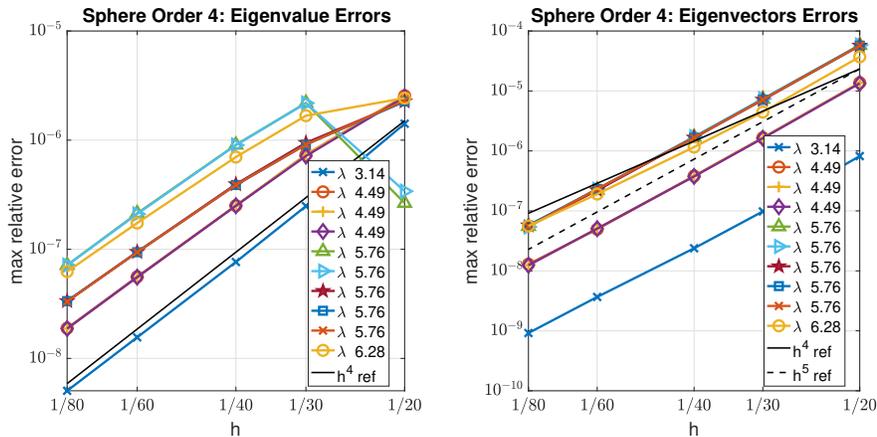
\begin{figure}[htb]
\begin{center}
\begin{tikzpicture}[scale=1]
  \useasboundingbox (0,.8) rectangle (2.2*\figw,1.05*\figh);  

  \begin{scope}[xshift=0cm]
    \figByWidth{0.0}{0}{fig/sphereOrder4EigenvalueErrors}{\figw}[0][0][0][0]
  \end{scope}
  \begin{scope}[xshift=1.1*\figw]
    \figByWidth{0.0}{0}{fig/sphereOrder4EigenvectorErrors}{\figw}[0][0][0][0]
  \end{scope}  

\end{tikzpicture}
\end{center}
\caption{Accuracy of computed eigenvalues and eigenvalues on a solid sphere (order of accuracy 4) compared to the exact continuous values.
}
\label{fig:sphereAccuracyOfEigenpairsOrder4}
\end{figure}
}

%% file: tex/implicitTolerance.tex
\newcommand{\tolMG}{\tau}
\subsection{Changing the implicit solver convergence tolerance} \label{sec:implicitTolerance}

It is of interest to know how the EigenWave algorithm behaves as a function of the convergence tolerance,~$\tolMG$, 
used to solve the implicit time-stepping equations. Three 
basic questions are how does $\tolMG$ affect
\begin{enumerate}
  \item the accuracy of the computed eigenpairs?
  \item the convergence of the Krylov-subspace eigenvalue solvers in terms of number of matrix-vector multiplies needed per eigenvalue found?
  \item the stability of the implicit time-stepping algorithm for the wave equation? In particular does implicit time-stepping remain stable when the implicit system is only approximately solved?
\end{enumerate}
Results given below indicate that, when using the multigrid solver Ogmg, the scheme is robust to changes in $\tolMG$ with the primary influence being larger values of $\tolMG$ lead to larger errors in the eigenpairs. This is not unexpected since we are computing eigenvectors of the matrix $\Aw_h$ and changing $\tolMG$ will
change the result of applying $\Aw_h$ to a vector.

\input tex/implicitToleranceFig.tex

Figure~\ref{fig:resultsVersusMultiGridTol} shows the behaviour of the EigenWave scheme, using the Krylov-Schur algorithm, as a function of the convergence tolerance
used in solving the implicit time-stepping equations with the multigrid solver.
EigenWave is used to compute eigenpairs on a square (square128, order=2) and disk ($\Gcd^{(16)}$, order=2) 
using the same EigenWave parameters as in Section~\ref{sec:implicitTimeSteppingWithMultigrid}.
The graph on the left of Figure~\ref{fig:resultsVersusMultiGridTol} shows the relative errors in the computed eigenvalues and
eigenvectors. It is seen that the accuracy of the eigenvectors is roughly scales in proportion to $\tolMG$.
The eigenvalues are more accurate than the eigenvectors but this is expected with the use of a Rayleigh quotient.
The graph on the right of Figure~\ref{fig:resultsVersusMultiGridTol} shows that the number of eigenvalues computed and 
number of matrix vector multiplies required are fairly constant as a function of~$\tolMG$.
Note that the implicit time-stepping equations are not solved very accurately. This could lead to growth in high-frequency solution modes and cause
the solution to be unstable if a large number time-steps are taken. However, 
the WaveHoltz filter is applied periodically after a relatively small number of steps and this helps damp these high-frequency modes.

Corresponding results when using bi-CG-stab to solve the implicit time-stepping equations are shown in Figure~\ref{fig:resultsVersusBiCGStabTol}.
There is some anomalous behaviour for the largest tolerance $\tau=10^{-1}$ but using such a large value for~$\tau$ is not recommended.

\mni
\textbf{Summary.} 
The EigenWave algorithm appears to be robust to the convergence tolerance used in solving the implicit time-stepping equations;
the accuracy in the eigenpairs is reduced but the number of wave-solves per eigenvalue computed is roughly constant.
Larger values of $\tolMG$ generally require fewer multigrid iterations (cycles) and this should result in a reduction in CPU time.
The accuracy of the computed eigenpairs can be independently checked by checking the residuals in the original eigenvalue problem~\eqref{eq:eigBVPdiscrete}.

%% file: tex/implicitToleranceFig.tex
{
\newcommand{\figSize}{5.5cm}
\newcommand{\figHeight}{4.8cm}

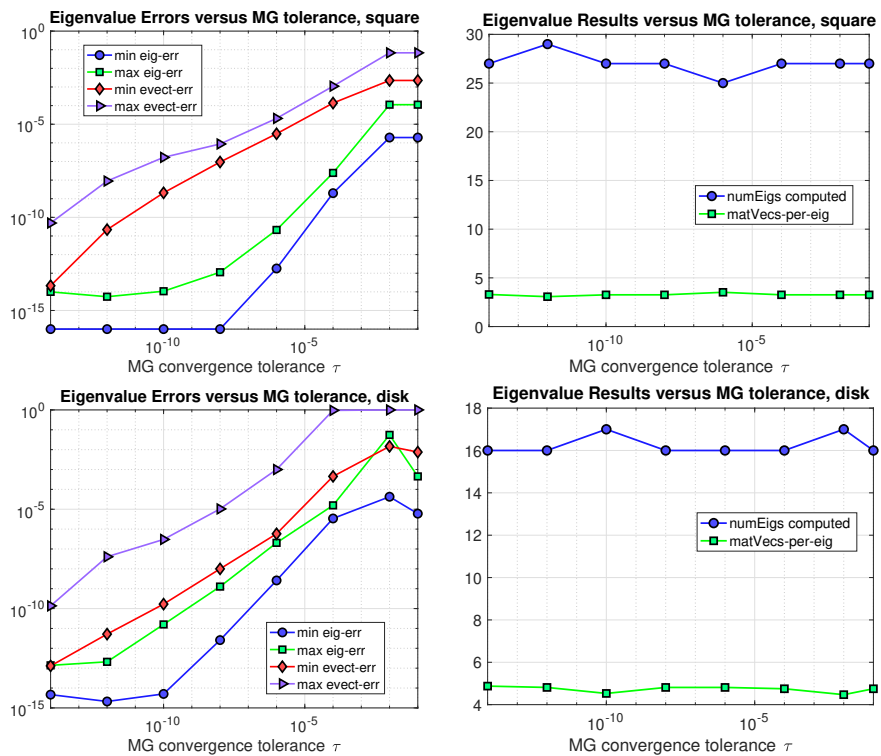
\begin{figure}[htb]
\begin{center}
\begin{tikzpicture}[scale=1]
  \useasboundingbox (0,.5) rectangle (2.1*\figSize,2.1*\figHeight);  

  \begin{scope}[yshift=1.05*\figHeight]
    \figByWidth{0.0*\figSize}{0}{fig/SquareEigenPairErrorsVersusMGTolerance}{\figSize}[0][0][0][0]
    \figByWidth{1.1*\figSize}{0}{fig/SquareEigenPairResultsVersusMGTolerance}{\figSize}[0][0][0][0]
  \end{scope}  

  \begin{scope}[yshift=0cm]
    \figByWidth{0.0*\figSize}{0}{fig/DiskEigenPairErrorsVersusMGTolerance}{\figSize}[0][0][0][0]
    \figByWidth{1.1*\figSize}{0}{fig/DiskEigenPairResultsVersusMGTolerance}{\figSize}[0][0][0][0]
  \end{scope}  

\end{tikzpicture}
\end{center}
\caption{Behaviour of the EigenWave results as a function of 
 the convergence tolerance $\tolMG$ used to solve the implicit time-stepping equations with multigrid.
 Top row: square. Bottom row: disk.
 Left: max and min relative errors in the eigenvalues and eigenvectors versus $\tolMG$.
 Right: number of computed eigenvalues and number of matrix vector multiplies (wave-solves) versus~$\tolMG$.
    }
\label{fig:resultsVersusMultiGridTol}
\end{figure}
}

{
\newcommand{\figSize}{5.5cm}
\newcommand{\figHeight}{4.8cm}

\begin{figure}[htb]
\begin{center}
\begin{tikzpicture}[scale=1]
  \useasboundingbox (0,.5) rectangle (2.1*\figSize,2.1*\figHeight);  

  \begin{scope}[yshift=1.05*\figHeight]
    \figByWidth{0.0*\figSize}{0}{fig/SquareEigenPairErrorsVersusBiCGStabTolerance}{\figSize}[0][0][0][0]
    \figByWidth{1.1*\figSize}{0}{fig/SquareEigenPairResultsVersusBiCGStabTolerance}{\figSize}[0][0][0][0]
  \end{scope}  

  \begin{scope}[yshift=0cm]
    \figByWidth{0.0*\figSize}{0}{fig/DiskEigenPairErrorsVersusBiCGStabTolerance}{\figSize}[0][0][0][0]
    \figByWidth{1.1*\figSize}{0}{fig/DiskEigenPairResultsVersusBiCGStabTolerance}{\figSize}[0][0][0][0]
  \end{scope}  

\end{tikzpicture}
\end{center}
\caption{Behaviour of the EigenWave results as a function of 
 the convergence tolerance $\tolMG$ used to solve the implicit time-stepping equations with bi-CG-stab.
 Top row: square. Bottom row: disk.
 Left: max and min relative errors in the eigenvalues and eigenvectors versus $\tolMG$.
 Right: number of computed eigenvalues and number of matrix vector multiplies (wave-solves) versus~$\tolMG$.
    }
\label{fig:resultsVersusBiCGStabTol}
\end{figure}
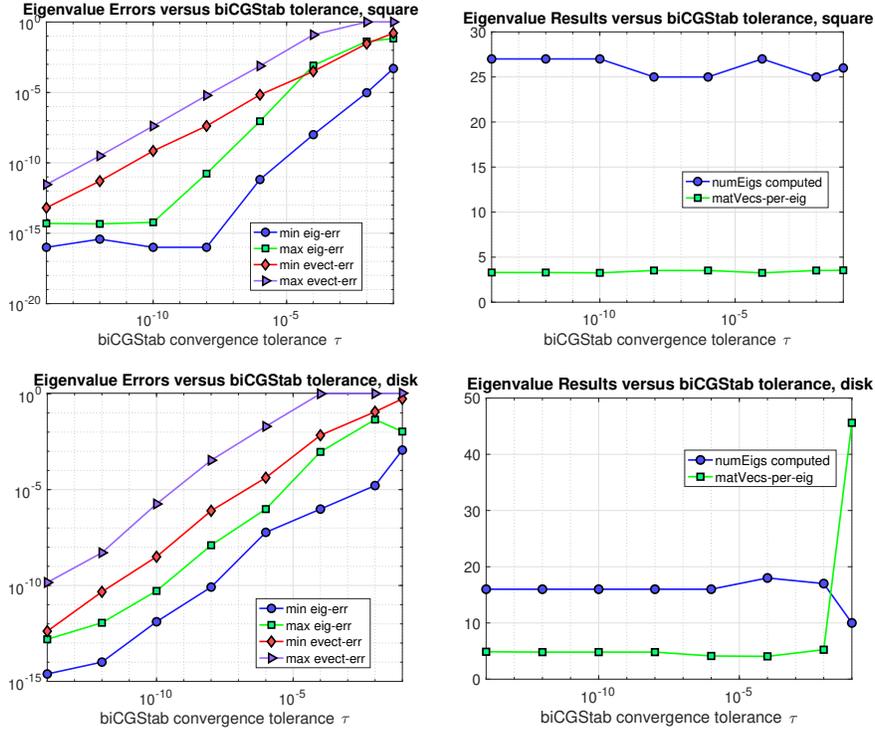
}

%% file: tex/implicitVersusExplicit.tex
\subsection{CPU time comparison for implicit versus explicit time-stepping} \label{sec:implicitVersusExplicitComparison}

In this section we compare the relative performance of using explicit versus implicit time-stepping.
When finding time-accurate solutions to the wave equation, 
explicit time-stepping methods are often preferred over implicit time-stepping methods since explicit schemes
are generally much faster per time-step. The CFL time-step restriction for explicit schemes requires the time-step $\dt$ to be proportional to 
the grid spacing. For accuracy purposes this is normally how the time-step is chosen (i.e.~the problem is not stiff as compared to solving the heat equation, say,
where implicit methods are often used).
When the problem is geometrically stiff (with local regions of small cells) a locally implicit scheme
can be advantageous~\cite{wimp2025} but an explicit scheme is still normally used over most of the domain.

The reason why implicit schemes are useful for EigenWave is because temporal accuracy is not very important.
Instead we only need to take enough time-steps so that the filter function $\beta$ in~\eqref{eq:filterFunction} is adequately resolved when discretized in time
(the time-discrete version of the filter function is given in~\eqref{eq:betad}).
When solving Helmholtz problems with WaveHoltz, 
analysis shows that at least $\Nits=5$ time-steps per-period $T=2\pi/\omega$ are required for convergence of the algorithm~\cite{overHoltzPartOne}.
Numerical experiments for EigenWave show that around $\Nits = 10$ time-steps per period is a good choice for good convergence (see Section~\ref{sec:changingImplicitStepsPerPeriod}).

Therefore, as the grids are refined, the implicit scheme can always take a constant $\Nits$ time-steps per period, 
while the number of explicit time-steps per period would increase. The break-even point when implicit methods are faster then depends on the relative
CPU costs of the implicit and explicit solvers.
The matrix $M_{ph}$ in~\eqref{eq:impMatrix} that arises from implicit time-stepping is a shifted Laplacian, but shifted in a manner to make the matrix \textbf{more definite}.
This matrix is thus well suited to fast $O(N)$ solution algorithms such as multigrid. 

\newcommand{\CFL}{K_{\ssf CFL}}
\newcommand{\Cits}{C_{\rm implicit}}
\newcommand{\Cets}{C_{\rm explicit}}
\newcommand{\PPW}{N_{\ssf PPW}}
Suppose the cost to solve the implicit system is proportional to $C_{\rm implicit} N$ while the cost for one explicit time-step is $C_{\rm explicit} N$.
The implicit method uses $N_t=\Np \Nits$ time steps per wave-solve.
If, for example, $\dt = \CFL h$ for the explicit method (which requires $\Tf/\dt$ time-steps) then the relative cost of each wave-solve is approximately 
\ba
  & \f{\text{CPU implicit}}{\text{CPU explicit}} = \f{\Cits N N_t}{\Cets N (\Tf/\dt) } 
                                          = \f{\Cits}{\Cets} \f{N_t\CFL\, \omega h}{2\pi c} = \f{\Cits}{\Cets} \f{N_t \CFL }{\PPW},
                                       \label{eq:CPUimpExp}  
\ea  
where $\PPW=2\pi c/(\omega h)$ is the number of points per wavelength.
Thus according to~\ref{eq:CPUimpExp}, for fixed target frequency $\omega$, as the mesh is refined $\PPW$ will increase and eventually
the implicit method will be faster (assuming the implicit solver has $O(N)$ scaling).

{
\newcommand{\figSize}{5cm}
\newcommand{\figHeight}{4cm}

\begin{figure}[htb]
\begin{center}
\begin{tikzpicture}[scale=1]
  \useasboundingbox (0,.6) rectangle (3.15*\figSize,1.05*\figHeight);  


  \begin{scope}[yshift=0cm]
    \figByWidth{0.0*\figSize}{0}{fig/squareOrder2CPUperWaveSolve}{\figSize}[0][0][0][0]
    \figByWidth{1.05*\figSize}{0}{fig/squareOrder4CPUperWaveSolve}{\figSize}[0][0][0][0]
    \figByWidth{2.10*\figSize}{0}{fig/nonSquareOrder4CPUperWaveSolve}{\figSize}[0][0][0][0]
  \end{scope}  

\end{tikzpicture}
\end{center}
\caption{Graphs of the CPU-time per wave-solve per-grid-point for explicit and implicit time-stepping when computing eigenpairs on the unit square
  using EigenWave with target frequency $\omega=12$.
  The implicit equations are solve using three different methods, (1) direct sparse solver, (2) multigrid, and (3) bi-CG-Stab.
   Left: Cartesian grid, second-order accuracy.
   Middle: Cartesian grid, fourth-order accuracy.   
   Right: Curvilinear grid, fourth-order accuracy.
    }
\label{fig:cpuExplicitVersusImplicit}
\end{figure}
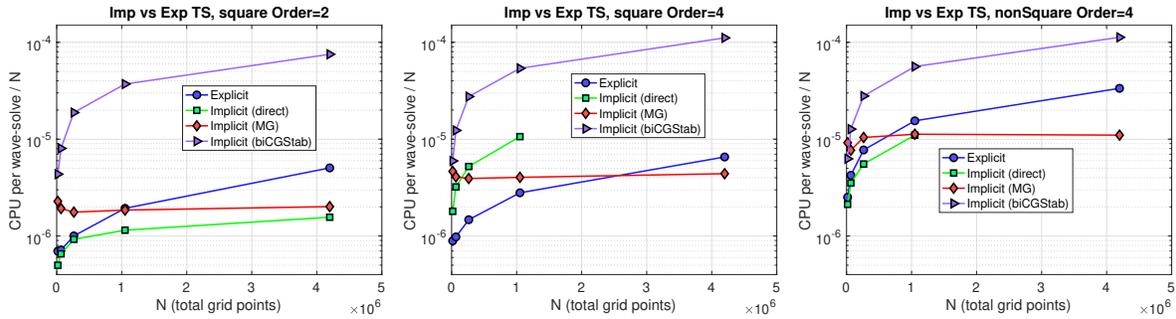
}

Figure~\ref{fig:cpuExplicitVersusImplicit} compares the CPU cost of using explicit and implicit time-stepping when using EigenWave 
to solve a problem on the unit square with target frequency $\omega=12$.
Three different implicit solvers are used, a sparse direct solver (factoring the matrix to start and then using back-substitutions for each wave-solve),
the Ogmg multigrid solver, and a bi-CG-Stab\footnote{Bi-CG-Stab is our preferred
 Krylov solver for elliptic equations on overset grids and so we this here even for a single component grid.}
Krylov solver (with ILU(3) preconditioner) from PETSc. 
The convergence tolerance for the iterative solvers was $10^{-10}$ 
(see Section~\ref{sec:implicitTolerance} for a discussion of how the iterative solvers behave as the convergence tolerance is varied).
These results show that the fastest solver depends on a number of factors such as the number of grid points $N$, the order of accuracy of the scheme, and whether the grid is Cartesian or curvilinear\footnote{The grid is the same but the code for general curvilinear grids is used instead of the code optimized for Cartesian grids.}. Implementation details of particular solvers are also important.
At second-order accuracy in two dimensions, for example, implicit time-stepping with the direct sparse solver, with reordering of equations to reduce fill-in, is the fastest over the range of $N$ considered. The direct sparse solver, which requires too much memory for the values of $N$ not shown,
 is less effective for the fourth-order accurate scheme which has a wider stencil.
Note that the red curves for multigrid in Figure~\ref{fig:cpuExplicitVersusImplicit} are nearly flat (as expected for an $O(N)$ algorithm), 
while the curves for the other methods show an increase as a function of $N$.
Thus for large enough $N$, implicit time-stepping with multigrid would ultimately be the fastest approach.

%% file: tex/convergenceProperties.tex
\subsection{Changing the number of requested eigenpairs} \label{sec:numberOfRequestedEigenpairs}

The convergence and efficiency of the EigenWave algorithm can depend strongly on the number of eigenpairs requested.
As noted in~\ref{sec:Arnoldi} the implicitly restarted Arnoldi algorithm maintains 
a basis of $\Ntotal = \Nreq + \Nextra$ Arnoldi vectors where $\Nreq$ is the number of requested vectors
and $\Nextra$ is the number of extra vectors. 
In this section we vary $\Nreq$  and determine in each case
the total number of wave-solves required for convergence of all requested eigenpairs, as well as the number of wave-solves per eigenpair.

{
\newcommand{\figw}{6.5cm}
\newcommand{\figh}{5.5cm}

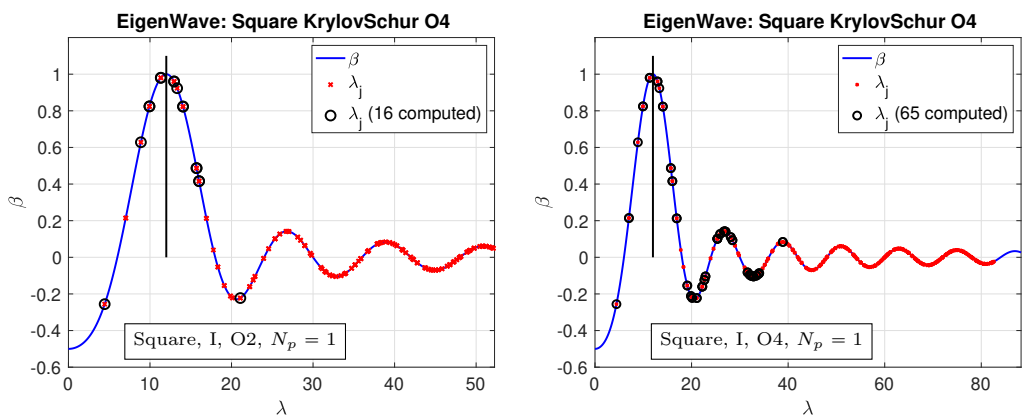
\begin{figure}[htb]
\begin{center}
\begin{tikzpicture}[scale=1]
  \useasboundingbox (0,.7) rectangle (14,5.5);  

  \begin{scope}[xshift=0cm]
    \figByWidth{0.0}{0}{fig/square128ImpEig12Krylov}{\figw}[0][0][0][0]
    \draw (1.5,1) node[draw,fill=white,anchor=west,xshift=0pt,yshift=0pt] {\scriptsize Square, I, O2, $\Np=1$};
  \end{scope}
  \begin{scope}[xshift=7cm]
    \figByWidth{0.0}{0}{fig/square128ImpEig12Num64Krylov}{\figw}[0][0][0][0]
    \draw (1.5,1) node[draw,fill=white,anchor=west,xshift=0pt,yshift=0pt] {\scriptsize Square, I, O4, $\Np=1$};
  \end{scope}  

\end{tikzpicture}
\end{center}
\caption{The problem of computing eigenpairs of a square is used to evaluate the EigenWave algorithm as the number of 
requested eigenpairs is varied. A total of $16$ eigenpairs were requested on left, while $64$ were requested at right. The plots show graphs of the filter function $\beta=\beta(\lambda;\omega)$ 
with the true discrete eigenvalues $\lambda_j$ marked
using red x's, and the computed eigenvalues marked using black circles.
}
\label{fig:squareEigsVaryArnoldi}
\end{figure}
}

Table~\ref{tab:changeNumberOfRequestedEigenpairs} shows the behavior of EigenWave when
varying the number of requested eigenpairs for the case of the unit square geometry with the \texttt{square128} grid and $\omega=12$.
For these results, the number of extra Arnoldi vectors is taken as $\Nextra=\Nreq+1$ so that
the total number of Arnoldi vectors is $\Ntotal=2\Nreq+1$.
Figure~\ref{fig:squareEigsVaryArnoldi} shows the filter function and computed eigenvalues for two representative cases with $\Nreq=16$
and $\Nreq=64$.
Table~\ref{tab:changeNumberOfRequestedEigenpairs} provides results for a range of values of~$\Nreq$, and these results suggest that the number of wave-solves per eigenpair decreases
monotonically as the number of requested eigenpairs increases, and the limit appears to be~2, at least for the range considered. Note that the IRAM algorithm uses at least $\Ntotal-1$ wave-solves, and thus in this case the number of wave-solves per eigenpair is at least $2+1/\Nreq$. Thus for the larger values of $\Nreq$ in 
Table~\ref{tab:changeNumberOfRequestedEigenpairs}, the IRAM algorithm is converging very fast and 
taking roughly one outer iteration per eigenpair.

These (somewhat limited) results suggest the following:
\begin{enumerate}
  \item The most efficient way to compute a collection of eigenpairs is to compute many of them all at once (given the available storage), rather than computing a few at a time with different values of $\omega$.
  \item Even if only a few eigenpairs are needed, it may be cheaper to compute more at once.
\end{enumerate}

\begin{table}[hbt]
\begin{center}\tableFontSize
\begin{tabular}{|c|c|c|c|c|c|c|} \hline
   Requested &  Arnoldi  & Converged &   Wave-solves &   Wave-solves  \\
       Eigs  &  Vectors  &    Eigs   &  (Mat-Vects)  &     per eig    \\ \hline
     $ 1$    &    $3$    &   $1$     &    $367$      &      $367$     \\
     $ 2$    &    $5$    &   $3$     &    $333$      &      $111$     \\
     $ 4$    &    $9$    &   $4$     &    $85$       &      $21$      \\
     $ 8$    &   $17$    &   $9$     &    $64$       &      $7.1$     \\
     $16$    &   $33$    &  $16$     &    $64$       &      $4.0$     \\
     $32$    &   $65$    &  $38$     &    $130$      &      $3.4$     \\
     $64$    &  $129$    &  $65$     &    $209$      &      $3.2$     \\
    $128$    &  $257$    & $154$     &    $411$      &      $2.7$     \\
    $256$    &  $513$    & $285$     &    $675$      &      $2.4$     \\
    $512$    &  $1025$   & $571$     &   $1325$      &      $2.3$     \\
\hline
\end{tabular}
\end{center}
\vspace*{-1\baselineskip}
\caption{Convergence of eigenpairs as a function of the number of requested eigenpairs $N_r$ for \texttt{square128}.
}
\label{tab:changeNumberOfRequestedEigenpairs}
\end{table}

In Table~\ref{tab:changeNumberOfRequestedEigenpairsII} the total number of Arnoldi vectors, $\Ntotal$, is varied.  The results here suggest that
$\Ntotal \approx 2 \Nreq$ seems to be a reasonably good choice to minimize the number of wave-solves per eigenpair.
\begin{table}[hbt]
\begin{center}\tableFontSize
\begin{tabular}{|c|c|c|c|c|c|c|} \hline
   Requested &  Arnoldi  & Converged &   Wave-solves &   Wave-solves  \\
       Eigs  &  Vectors  &    Eigs   &  (Mat-Vects)  &     per eig    \\ \hline
    $256$    &  $400$    & $285$     &    $700$      &      $2.46$     \\
    $256$    &  $475$    & $267$     &    $626$      &      $2.34$     \\
    $256$    &  $513$    & $285$     &    $675$      &      $2.37$     \\
    $256$    &  $575$    & $311$     &    $747$      &      $2.40$     \\
    $256$    &  $600$    & $315$     &    $776$      &      $2.46$     \\
\hline
\end{tabular}
\end{center}
\vspace*{-1\baselineskip}
\caption{Convergence of eigenpairs as a function of the number of Arnoldi vectors $\Ntotal$ for \texttt{square128}.
}
\label{tab:changeNumberOfRequestedEigenpairsII}
\end{table}

%% file: tex/varyNp.tex
\subsection{Changing the number of filter periods $\Np$ as a function of the number of computed eigenpairs.} \label{sec:changingNp}

In this section a study is made of how to choose the number of periods, $\Np$, over which the filter is integrated, as a function of the number of computed eigenpairs,
with the goal of minimizing the cost per eigenpair.
Taking a larger $\Np$ causes the main peak of the filter function to narrow, and leads to fewer eigenvalues 
lying near the peak. As a result, a smaller number of eigenpairs can be found more efficiently.
Reducing the number of computed eigenpairs may be desired, for example,
to avoid the storage associated with maintaining a large number of eigenvectors (the Krylov-Schur and IRAM algorithms require storage for approximately double the number of requested eigenpairs).
The disadvantage of increasing $\Np$ is that the cost per wave-solve increases by a factor of $\Np$.
To measure the performance we therefore use the number of time-steps per eigenpair found. 

\input tex/timeStepsPerEigVaryNpSummaryFig.tex

For this study we compute eigenpairs on the unit square with $128$ grid points in each direction using a target frequency of $\omega=35$.
The fourth-order accurate scheme is used with implicit time-stepping and $\Nits=10$ time-steps per period (i.e. a total of $N_t=\Np\Nits$ time-steps per wave-solve).
Figure~\ref{fig:squareEigsVaryNpSummary} shows 
the number of time-steps per eigenvalue as a function of the number of filter periods $N_p$ when requesting $6$, $12$, $24$ and $36$ eigenpairs.
It is seen that the cost per eigenvalue tends to decrease as more are requested. 
In addition, the most efficient value for $N_p$ increases as fewer eigenvalues are requested.

\input tex/timeStepsPerEigVaryNpFig.tex

Figure~\ref{fig:squareEigsVaryNp} shows more details from the computations.
The top column shows graphs of the number of time-steps per eigenvalue as a function of $\Np$ (same data as shown in Figure~\ref{fig:squareEigsVaryNpSummary}) and also indicates the actual number of eigenvalues found (Krylov-Schur type algorithms may find more converged eigenvalues than requested due to the nature of the algorithm).
The number of eigenvalues found can vary and this partially explains the up and down behaviour of some of the curves.
Below each graph in the top column, a plot is made of the corresponding  $\beta$ filter function for the value of $\Np$ that gave the best result.
These graphs show the locations of eigenvalues (red crosses) and the computed eigenvalues (black circles).
As the number of requested eigenvalues increases, $\Np$ decreases and the main peak gets wider.
A common element of the four optimal cases is that the computed eigenvalues roughly include all the eigenvalues near the main peak down
to some level where $\beta \approx 0.8$.
Further information on this behaviour is provided in~\ref{sec:estimatingTheConvergenceRate}.

%% file: tex/timeStepsPerEigVaryNpSummaryFig.tex
{
\newcommand{\figw}{6cm}
\newcommand{\figh}{5cm}

\begin{figure}[htb]
\begin{center}
\begin{tikzpicture}[scale=1]
  \useasboundingbox (0,.6) rectangle (\figw,\figh);  

  \begin{scope}[xshift=0cm]
    \figByWidth{0*\figw}{0}{fig/square128StepsPerEig}{\figw}[0][0][0][0]
  \end{scope}  

\end{tikzpicture}
\end{center}
\caption{
Number of time-steps per eigenvalue as a function of the number 
of periods $N_p$ when requesting $6$, $12$, $24$ and $36$ eigenpairs (see also Figure~\ref{fig:squareEigsVaryNp} for more details).
The cost per eigenvalue tends to decrease as more are requested. 
The most efficient value for $N_p$ increases as fewer eigenvalues are requested.
}
\label{fig:squareEigsVaryNpSummary}
\end{figure}
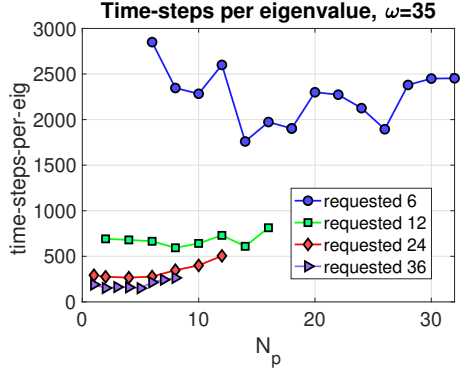
}

%% file: tex/timeStepsPerEigVaryNpFig.tex
{
\newcommand{\figw}{4cm}
\newcommand{\figh}{3.4cm}

\begin{figure}[htb]
\begin{center}
\begin{tikzpicture}[scale=1]
  \useasboundingbox (0,.5) rectangle (4*\figw,2*\figh);  

  \begin{scope}[yshift=\figh]
    \figByWidth{0*\figw}{0}{fig/square128Eig6VaryNpStepsPerEig}{\figw}[0][0][0][0]
    \figByWidth{1*\figw}{0}{fig/square128Eig12VaryNpStepsPerEig}{\figw}[0][0][0][0]
    \figByWidth{2*\figw}{0}{fig/square128Eig24VaryNpStepsPerEig}{\figw}[0][0][0][0]
    \figByWidth{3*\figw}{0}{fig/square128Eig36VaryNpStepsPerEig}{\figw}[0][0][0][0]
  \end{scope}  
  \begin{scope}[xshift=0cm]
    \figByWidth{0*\figw}{0}{fig/square128Eig6VaryNp14Krylov}{\figw}[0][0][0][0]
    \figByWidth{1*\figw}{0}{fig/square128Eig12VaryNp8Krylov}{\figw}[0][0][0][0]
    \figByWidth{2*\figw}{0}{fig/square128Eig24VaryNp4Krylov}{\figw}[0][0][0][0]
    \figByWidth{3*\figw}{0}{fig/square128Eig36VaryNp2Krylov}{\figw}[0][0][0][0]
  \end{scope}  

\end{tikzpicture}
\end{center}
\caption{
A comparison of the cost of computing different numbers of eigenvalues as a function of the number of periods $N_p$
for a problem on a square with target frequency $\omega=35$ (see also Figure~\ref{fig:squareEigsVaryNpSummary} for combined results).
The optimal values of $N_p$ depend on the number of requested eigenvalues.
Top: number of time-steps per eigenvalue found as a function of $N_p$ when requesting $6$, $12$, $24$ and $36$ eigenpairs.
The labels at each data point indicates the actual number of eigenpairs found.
Bottom: The filter functions corresponding to the best value of $N_p$ for a given requested number of eigenvalues.
For example, the bottom left filter function corresponds to the computation in the plot directly above with $N_p=14$.
}
\label{fig:squareEigsVaryNp}
\end{figure}
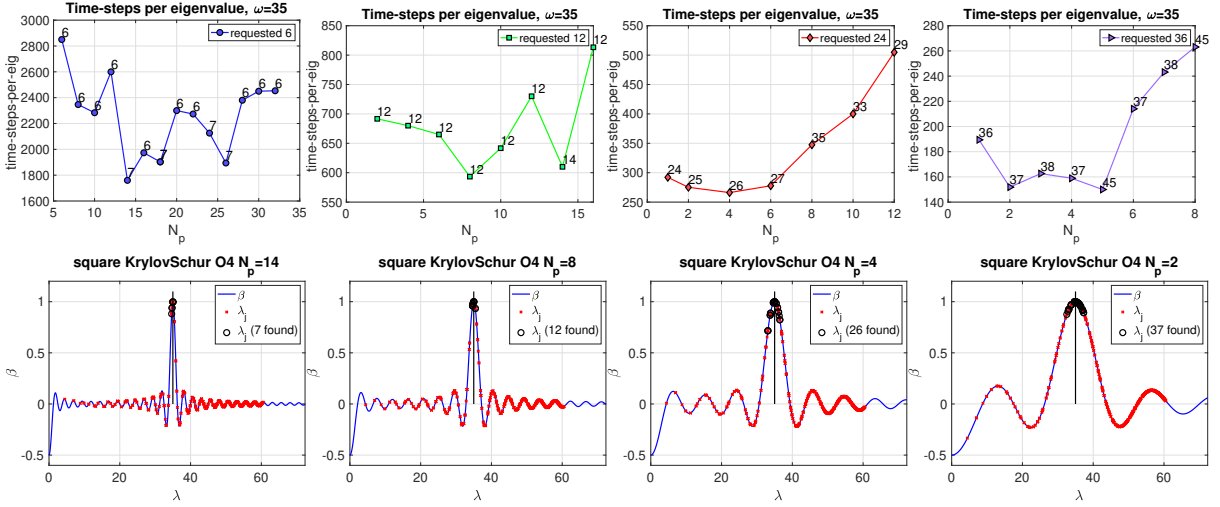
}

%% file: tex/convergenceRateEstimation.tex
\subsection{Estimating the EigenWave convergence rate based on simultaneous iteration} \label{sec:estimatingTheConvergenceRate}

In this section we provide some justification for the results observed in~\ref{sec:changingNp} for choosing the optimal value of $\Np$
in order to efficiently compute a given number of eigenpairs.
The convergence rate of EigenWave using a simple power iteration (i.e.~the fixed-point iteration)
is expected to be $\beta_2/\beta_1$ where the eigenvalues $\beta_j$ are assumed to be ordered from largest to smallest in magnitude.
The Krylov-Schur and IRAM algorithms fall into the category of eigenvalue algorithms known as Krylov subspace methods.
See~\cite{LehoucqSorensen2000}, for example, for a discussion of the expected convergence of the IRAM algorithm.
These convergence results depend, however, on a number of parameters that are dynamically found in the algorithm. 
To get some concrete estimates for the convergence rates we therefore study a simpler subspace algorithm which hopefully will provide
insight into the more complicated algorithms.

\input tex/simultaneousIterationAlgorithm.tex

A basic subspace iteration, known as simultaneous iteration (SI) (or the block power method), is given in Algorithm~\ref{alg:simultaneousIteration}.
It applies a power-like method simultaneously to a set of $N_a$ vectors, $\vv_j\in\Real^{N}$, $j=1,2,\ldots,N_a$, 
where $N$ denotes the number of grid points.
The iteration computes an $N_a$-dimensional invariant subspace whose eigenvalues
approach the $N_a$ largest eigenvalues.
Suppose, just as in the Krlov-Schur or IRAM algorithms, that $N_a$ is chosen to be larger than the number of eigenpairs we actually want.
In particular, suppose that $N_a = N_r + N_e$ where $N_r$ is the number of requested eigenpairs and $N_e$ is the number of extra eigenvalues added to the subspace.
A typical choice might be $N_e=N_r$. Thus we iterate over a subspace that has about twice the dimension of the number of eigenpairs we want. 
In this case the expected convergence rate of SI for the largest $N_r$ eigenpairs is~\cite{Gu2000,Saad2016}
\ba
  \CR \approx \f{\beta_{N_a+1}}{\beta_{N_r}}. \label{eq:subSpaceIterationCR}
\ea
Thus the CR depends on the \textsl{gap} between the eigenvalue $\beta_{N_r}$ and eigenvalue $\beta_{N_a+1}$.
The bigger this gap the faster the convergence, although the cost per iteration increases with increasing $N_e$.

\input tex/filterWidthRequestedEigsFig.tex

We would like to know how the convergence rate of EigenWave with simultaneous iteration depends on the choice of $N_r$.
To answer this question consider the diagram in Figure~\ref{fig:filterWithRequestedEigs} which shows the filter function along 
with the eigenvalues $\lambda_j$ (assumed equally spaced to be concrete).
Since the $\beta$ function only depends on the ratio $\lambda/\omega$, we take $\omega=1$ to simplify the current discussion.
We wish to find  $\beta_{N_r}$ and $\beta_{N_a+1}$ in order to use~\eqref{eq:subSpaceIterationCR} to estimate the convergence rate.
Let $\br=\beta_{N_r}$ be given,
then the number of requested eigenvalues $N_r$ is equal to the number of $\lambda_j$ where $\beta(\lambda_j)\ge \br$.
Define $\dw>0$ to satisfy $\beta(1+\dw)=\br$, and assume  $\br \in[0.7,1]$ to ensure $\dw=\dw(\br)$ will be single valued (see Figure~\ref{fig:filterWithRequestedEigs}). 
There are approximately $N_r$ eigenvalues in the interval $[1-\dw,1+\dw]$ of width $2\dw$.
If $N_a=2 N_r$ then we need to know the interval that contains twice as many eigenvalues, but this is just the 
interval $[1-2\dw,1+2\dw]$ of width $4w$ (assuming the eigenvalues are roughly evenly distributed for $\lambda$ near $\omega$). 
Therefore $\beta_{N_a+1} \approx \beta(1+2\dw)$, and the approximate convergence rate from~\eqref{eq:subSpaceIterationCR} is
\ba
  \CR(\br) \approx \f{ \beta(1+2 \dw) }{ \beta(1+\dw) }. \label{eq:subSpaceIterationEigenWaveCR}
\ea
The left graph of Figure~\ref{fig:subspaceIterationCR} shows $\dw$ as a function of $\br$ for $\Np=1,2,4,8$.
The middle graph of Figure~\ref{fig:subspaceIterationCR} plots the estimated convergence rate~\eqref{eq:subSpaceIterationEigenWaveCR}
for different values of $\Np$; note that all these curves are essentially the same.
The right graph shows the $\beta$ functions for $\Np=1,2,4,8$.

The graphs in Figure~\ref{fig:subspaceIterationCR} show that the estimated convergence rate is a strong function of $\br$ but is essentially independent of $\Np$.
Thus if the number of requested eigenpairs $N_r$ is reduced by a factor of $2$ (e.g. to reduce memory requirements) 
then by doubling $\Np$ a similar convergence rate will be obtained.

{
\newcommand{\figw}{5cm}
\newcommand{\figh}{4.4cm}

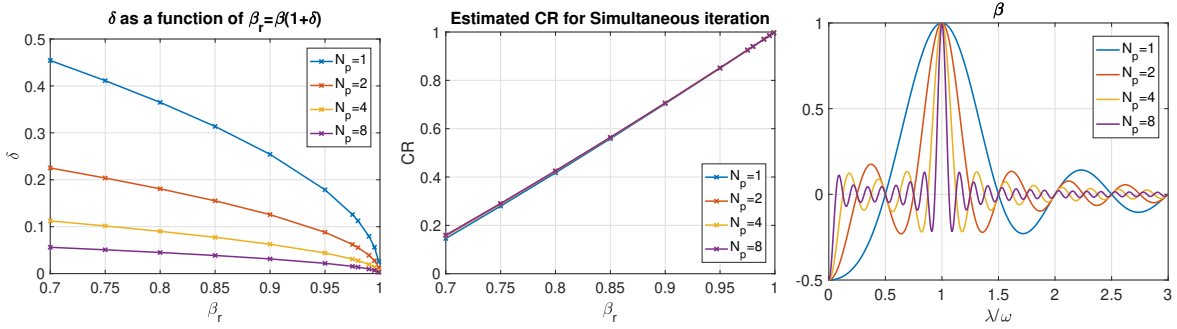
\begin{figure}[htb]
\begin{center}
\begin{tikzpicture}[scale=1]
  \useasboundingbox (0,.55) rectangle (3.1*\figw,1*\figh);  

  \begin{scope}[xshift=0cm]
    \figByWidth{0.00*\figw}{0}{fig/widthVersusBr}{\figw}[0][0][0][0]
    \figByWidth{1.05*\figw}{0}{fig/subSpaceIterationCR}{\figw}[0][0][0][0]
    \figByWidth{2.10*\figw}{0}{fig/betaFunctionsVaryNp}{\figw}[0][0][0][0]
  \end{scope}  

\end{tikzpicture}
\end{center}
\caption{
Left: $\dw$ as a function of $\br$. 
Middle: estimated convergence rate of simultaneous iteration as a function of $\br$ for different $\Np$.
Right: $\beta$ functions for $\Np=1,2,4,8$.
}
\label{fig:subspaceIterationCR}
\end{figure}
}

Consider the common situation where $\omega$ is relatively large (compared to the average spacing $\dlam$ between eigenvalues for $\lambda$ near $\omega$) 
and the requested number of eigenpairs, $\Nreq$, is not too large (this becomes more specific below).
Taking $\Np=1$ to start off, and given a desired convergence rate CR, $\br$ can be found from the middle graph
of Figure~\ref{fig:subspaceIterationCR} and $\dw$ from the left graph (e.g. $CR \approx 0.4$, $\br\approx 0.8$, $\dw\approx 0.375$). 
The corresponding value for $\Nreq^{(1)}$ (the superscript denotes $\Np=1$) 
will depend on $\dw$ and $\dlam$ by 
\ba
  \Nreq^{(1)} \approx \f{2 \dw}{\dlam}.
\ea
This value for $\Nreq^{(1)}$ may lead to excessive memory requirements and a smaller value may be desired. 
In that case, $\Np$ can be increased to reduce the memory usage, with the new number of requested eigenvalues being 
$\Nreq^{(\Np)} \approx \Nreq^{(1)}/\Np$.

{
\newcommand{\figw}{6cm}
\newcommand{\figh}{5.5cm}
\begin{figure}[htb]
\begin{center}
\begin{tikzpicture}[scale=1]
  \useasboundingbox (0,.6) rectangle (2.1*\figw,.9*\figh);  

  \begin{scope}[xshift=0cm]
    \figByWidth{0.00*\figw}{0}{fig/square256O2Freq60Nr12Np4Krylov}{\figw}[0][0][0][0]
    \figByWidth{1.05*\figw}{0}{fig/square256O2Freq60Nr12Np12Krylov}{\figw}[0.0][0][0][0]
  \end{scope}    

\end{tikzpicture}
\end{center}
\caption{
Filter function and distribution of eigenvalues for $\omega=60$ for $\Np=4$ (left) and $\Np=12$ (right).
$\Nreq=12$ eigenvalues were requested using the Krylov-Schur algorithm.
}
\label{fig:filterFunctionOmega60}
\end{figure}
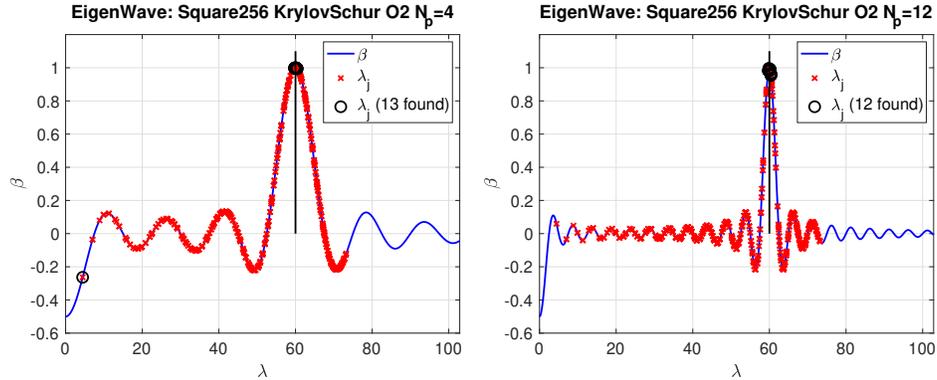
}

{
\newcommand{\figw}{7cm}
\newcommand{\figwZoom}{7cm}
\newcommand{\figh}{6.25cm}
\begin{figure}[htb]
\begin{center}
\begin{tikzpicture}[scale=1]
  \useasboundingbox (0,.6) rectangle (2.05*\figw,.95*\figh);  

  \begin{scope}[xshift=0cm]
    \figByWidth{0.00*\figw}{0}{fig/square256O2Freq60KrylovSchurVersusSIconvergeRates}{\figw}[0][0][0][0]
    \figByWidth{1.05*\figw}{0}{fig/square256O2Freq60KrylovSchurVersusSIconvergeRatesZoom}{\figwZoom}[0.0][0][0][0]
    \draw[->,very thick,black,xshift=12pt] (5.3,4) -- (8,3.75);
  \end{scope}    

\end{tikzpicture}
\end{center}
\caption{
Computed EigenWave convergence rates CR, using the Krylov-Schur algorithm, compared to the theoretical convergence rates for simultaneous iteration (SI).
Data points are coloured by $N_c$, the number of computed eigenvalues, and labeled by $\Np$, the number of filter periods.
The data is graphed versus $\br$, shown in Figure~\ref{fig:filterWithRequestedEigs}, which varies nearly linearly with the CR for the SI theory as shown
in Figure~\ref{fig:subspaceIterationCR}. The right graph is a magnified view of the left graph.
}
\label{fig:krylovSchurVersusSIconvergenceRates}
\end{figure}
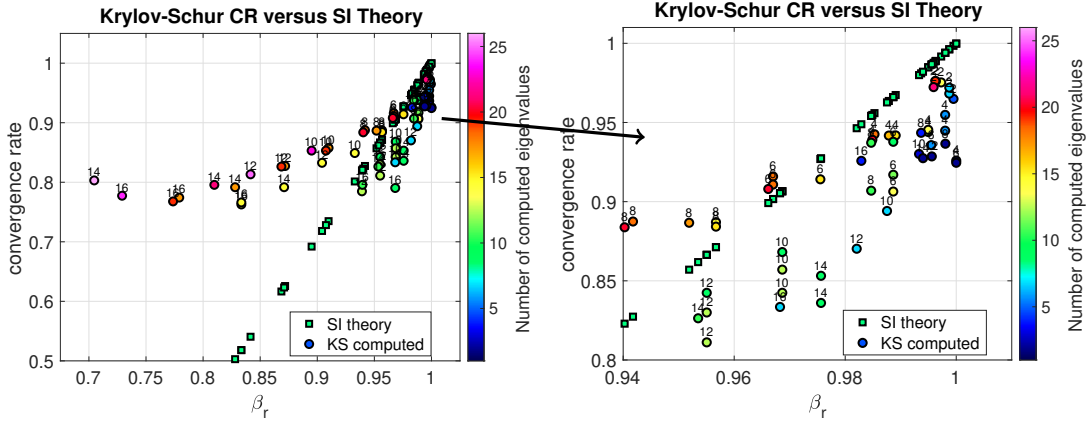
}

\medskip
To evaluate the usefulness of the SI convergence theory we consider the problem of finding eigenpairs of the Laplacian with Dirichlet boundary conditions on the unit square with $256$ grid points in each direction and a relatively high frequency $\omega=60$.
Figure~\ref{fig:filterFunctionOmega60} shows graphs of the filter function and distribution of eigenvalues for $\Np=4$ and $\Np=12$.
Computations were performed for 
$\Nreq\in\{4,6,8,10,12,14,16,18,20\}$ and $\Np\in\{2,4,6,8,10,12,14,16\}$ for a total of $9\times 8 = 72$ simulations.
Figure~\ref{fig:krylovSchurVersusSIconvergenceRates} compares actual computed
convergence rates from EigenWave to the SI theoretical results. 
The convergence rate for EigenWave is computed as 
\ba
    \CR = {\rm tol}^{1/N_{\rm mv}}
\ea
where ${\rm tol}=10^{-14}$ is the typical relative error in the computed eigenvalues 
and $N_{\rm mv}$ is the number of matrix-vector multiplies (wave-solves).
The data is plotted versus $\br$ to correspond to the middle graph of Figure~\ref{fig:subspaceIterationCR}.
Data points are coloured by $N_c$, the number of computed eigenvalues, and labeled by $\Np$, the number of filter periods.

To each computed data point (coloured disk) there is a corresponding theory point (green square) at the same value of $\beta_r$, the value of $\br$ being determined
from the computed eigenvalues using $\beta_r=\beta_{d,N_r}$ where $\beta_{d,j}$ is given by~\eqref{eq:discreteBeta}.
The results show that, despite the simplifications, the SI theory does a reasonable job at predicting the basic trends for convergence rates for $\br \ge 0.95$ 
where the Krylov-Schur rates are generally better then the SI rates.
However for smaller values of $\br$ the theory does not match computational results.
This is likely due to the fact that the SI theory ignores the oscillatory behaviour of the $\beta$ function for $\beta<0.2$ where there are many eigenvalues of similar 
size as seen in Figure~\ref{fig:filterFunctionOmega60}. These eigenvalues will slow the actual Krylov-Schur convergence behaviour and this is reflected in the results.


%% file: tex/simultaneousIterationAlgorithm.tex
\renewcommand{\algFontSize}{\small}
\begin{algorithm}
\algFontSize 
\caption{Simultaneous iteration applied to matrix $\Av$ with $N_a$ starting vectors in $\Vv^{(0)} \in\Real^{N\times N_a}$.}
\begin{algorithmic}[1]
%
    \State $\Vv \Rv = \Vv^{(0)}$       \Comment QR factorize the matrix of starting vectors $\Vv^{(0)}$.
    \For{$j=1,2,\ldots$} 
      \State $\Wv = \Av \Vv$     \Comment Multiply $\Av$ times all the columns of $\Vv$.
      \State $\Hv = \Vv^* \Wv$   \Comment Eigenvalues of $\Hv=\Vv^* A \Vv$ approximate largest eigenvalues of $\Av$.
      \If{ $\| \Wv - \Vv \Hv \|_2 \le {\rm tol}$}
          \State break from loop
      \EndIf
      \State $\Vv \Rv = \Wv$     \Comment QR factorize $\Wv$ (orthonormalization step).
    \EndFor
\end{algorithmic} 
\label{alg:simultaneousIteration}
\end{algorithm}

%% file: tex/filterWidthRequestedEigsFig.tex
{
\newcommand{\figw}{7cm}
\newcommand{\figh}{6cm}

\begin{figure}[htb]
\begin{center}
\begin{tikzpicture}[scale=1]
  \useasboundingbox (0,.5) rectangle (1*\figw,1*\figh);  

  \begin{scope}[xshift=0cm]
    \figByWidth{0*\figw}{0}{fig/filterWithRequestedEigs}{\figw}[0][0][0][0]
    \draw[black] (3,5) node[anchor=west] {\small $\beta_r=\beta(1+\dw)$};
    \draw[black] (3.3,4) node[anchor=west,yshift=-3pt] {\small $\beta(1+2 \dw)$};
    \draw[black] (3,1.75) node[xshift=-8pt] {\small $1+\dw$};
    \draw[black] (3,1.4) node[xshift=16pt] {\small $1+2 \dw$};

    \draw[black,->,thick] (2.75,1.85) -- (3,2.45); 
    \draw[black,->,thick] (3.45,1.6) -- (3.3,2.45); 
  \end{scope}  

\end{tikzpicture}
\end{center}
\caption{
Filter with $N_r$ requested and $N_e$ extra eigenvalues for $\br=.85$. 
The requested eigenvalues are those with $\beta(\lambda_j)\ge \beta_r$.
}
\label{fig:filterWithRequestedEigs}
\end{figure}
}

%% file: tex/eigenpairTables.tex
\section{Eigenpair tables} \label{sec:eigenpairTables}

The tables provided in this section give more details of the eigenpairs computed using EigenWave and their errors for many of the examples discussed previously.  
These tables are included since the worst case error for a particular example is not always representative of the typical error.
The true discrete eigenvalues $\lambda_{h,k}^{\rm true}$, $k=0,1,2,\ldots$, for a given case are sorted into increasing order by magnitude.
For a given computed eigenvalue $\lambda_{h,j}$, we find the closest true discrete eigenvalue $\lambda_{h,k}^{\rm true}$.
The tables indicate the accuracy of each eigenpair
and the multiplicity of the eigenvalue as estimated numerically to some tolerance.

Details of the eigenvalues for the square at order of accuracy $2$ are given in Table~\ref{tab:square128Order2}. 
Table~\ref{tab:sice4Order4} gives results for the disk at order $4$, while Tables~\ref{tab:spheree2Order2} and~\ref{tab:spheree2Order4}
give results for the sphere at orders $2$ and $4$, respectively.







\begin{table}[hbt]
\begin{center}\tableFontSize
\begin{tabular}{|c|r|r|c|c|c|c|c|} \hline
 \multicolumn{8}{|c|}{EigenWave: square, ts=implicit, order=2, $N_p=1$, KrylovSchur } \\ \hline 
   $j$  & \multicolumn{1}{c|}{$\lambda_{h,j}$} &  \multicolumn{1}{c|}{$\lambda_{h,k}^{\rm true}$}  & $k$ &  mult &  eig-err & evect-err & eig-res  \\ \hline
   0   &   7.024215 &   7.024215 &     1 &    2  &  3.41e-15 &  2.43e-14 &  2.60e-12 \\ 
   1   &   7.024215 &   7.024215 &     1 &    2  &  7.59e-16 &  2.45e-14 &  2.27e-12 \\ 
   2   &   8.884874 &   8.884874 &     3 &    1  &  1.80e-15 &  2.91e-14 &  1.36e-12 \\ 
   3   &   9.932544 &   9.932544 &     4 &    2  &  1.79e-16 &  2.96e-14 &  9.63e-13 \\ 
   4   &   9.932544 &   9.932544 &     5 &    2  &  3.76e-15 &  4.59e-14 &  1.01e-12 \\ 
   5   &  11.325052 &  11.325052 &     6 &    2  &  1.88e-15 &  3.91e-14 &  8.25e-13 \\ 
   6   &  11.325052 &  11.325052 &     7 &    2  &  3.14e-16 &  6.77e-14 &  7.98e-13 \\ 
   7   &  12.948204 &  12.948204 &     8 &    2  &  3.29e-15 &  1.88e-13 &  5.05e-13 \\ 
   8   &  12.948204 &  12.948204 &     9 &    2  &  5.35e-15 &  2.53e-13 &  5.60e-13 \\ 
   9   &  13.325638 &  13.325638 &    10 &    1  &  8.00e-16 &  4.89e-13 &  5.55e-13 \\ 
  10   &  14.044834 &  14.044834 &    12 &    2  &  2.15e-15 &  1.15e-13 &  5.15e-13 \\ 
  11   &  14.044834 &  14.044834 &    12 &    2  &  2.78e-15 &  1.08e-13 &  5.17e-13 \\ 
  12   &  15.702649 &  15.702649 &    14 &    2  &  3.17e-15 &  5.38e-14 &  3.57e-13 \\ 
  13   &  15.702649 &  15.702649 &    14 &    2  &  3.39e-15 &  5.10e-14 &  3.50e-13 \\ 
  14   &  16.009363 &  16.009363 &    15 &    2  &  1.33e-15 &  2.79e-14 &  3.56e-13 \\ 
  15   &  16.009363 &  16.009363 &    16 &    2  &  7.99e-15 &  4.99e-14 &  2.66e-13 \\ 
  16   &  16.908610 &  16.908610 &    17 &    2  &  4.20e-16 &  3.68e-14 &  2.39e-13 \\ 
  17   &  16.908610 &  16.908610 &    18 &    2  &  2.10e-16 &  4.41e-14 &  3.27e-13 \\ 
  18   &  17.764396 &  17.764396 &    19 &    1  &  7.60e-15 &  7.15e-14 &  2.71e-13 \\ 
  19   &  18.308930 &  18.308930 &    20 &    2  &  1.36e-15 &  2.05e-14 &  2.79e-13 \\ 
  20   &  18.308930 &  18.308930 &    21 &    2  &  9.70e-16 &  2.53e-14 &  2.54e-13 \\ 
  21   &  19.092753 &  19.092753 &    22 &    2  &  3.35e-15 &  1.99e-14 &  2.67e-13 \\ 
  22   &  19.092753 &  19.092753 &    22 &    2  &  1.86e-15 &  1.69e-14 &  2.78e-13 \\ 
  23   &  19.852824 &  19.852824 &    25 &    2  &  3.58e-16 &  8.44e-14 &  2.99e-13 \\ 
  24   &  19.852824 &  19.852824 &    25 &    2  &  1.25e-15 &  7.57e-14 &  3.44e-13 \\ 
  25   &  20.105161 &  20.105161 &    26 &    2  &  1.06e-15 &  3.74e-14 &  2.44e-13 \\ 
  26   &  20.105161 &  20.105161 &    27 &    2  &  1.77e-15 &  3.47e-14 &  2.83e-13 \\ 
 \hline 
\end{tabular}
\end{center}
\vspace*{-1\baselineskip}
\caption{Further details of the eigenpairs computed using EigenWave for a square128 grid using the KrylovSchur algorithm and implicit time-stepping (see Table~\ref{tab:square128Summary} for a summary of these results).  The spatial order of accuracy is 2 and the wave-solves use $N_p=1$ to determine the final time.
The index $k$ denotes the closest true discrete eigenvalue $\lambda_{h,k}^{\rm true}$ to the EigenWave value $\lambda_{h,j}$.
}
\label{tab:square128Order2}
\end{table}



\begin{table}[hbt]
\begin{center}\tableFontSize
\begin{tabular}{|c|r|r|c|c|c|c|c|} \hline
 \multicolumn{8}{|c|}{EigenWave: disk, ts=implicit, order=4, $\omega=  10.0$, $N_p=1$, KrylovSchur } \\ \hline 
   $j$  & \multicolumn{1}{c|}{$\lambda_j$} &  \multicolumn{1}{c|}{$\lambda_{h,k}^{\rm true}$}  & $k$ &  mult &  eig-err & evect-err & eig-res  \\ \hline
   0   &   7.015578 &   7.015578 &     8 &    2  &  3.29e-15 &  3.26e-14 &  4.06e-13 \\ 
   1   &   7.015578 &   7.015578 &     8 &    2  &  2.66e-15 &  2.30e-14 &  3.70e-13 \\ 
   2   &   7.588322 &   7.588322 &    10 &    2  &  2.93e-15 &  1.18e-13 &  4.30e-13 \\ 
   3   &   7.588322 &   7.588322 &    11 &    2  &  0.00e+00 &  8.16e-14 &  3.56e-13 \\ 
   4   &   8.417205 &   8.417205 &    12 &    2  &  2.32e-15 &  3.67e-13 &  6.17e-13 \\ 
   5   &   8.417238 &   8.417238 &    13 &    2  &  2.11e-16 &  3.35e-13 &  5.49e-13 \\ 
   6   &   8.653706 &   8.653706 &    14 &    1  &  0.00e+00 &  6.75e-13 &  9.13e-13 \\ 
   7   &   8.771438 &   8.771438 &    15 &    2  &  1.82e-15 &  2.20e-13 &  3.52e-13 \\ 
   8   &   8.771438 &   8.771438 &    15 &    2  &  2.03e-16 &  1.75e-13 &  2.97e-13 \\ 
   9   &   9.760970 &   9.760970 &    18 &    2  &  1.82e-16 &  1.62e-13 &  1.88e-13 \\ 
  10   &   9.760970 &   9.760970 &    18 &    2  &  0.00e+00 &  1.94e-13 &  2.27e-13 \\ 
  11   &   9.936020 &   9.936020 &    19 &    2  &  2.32e-15 &  1.73e-13 &  2.25e-13 \\ 
  12   &   9.936020 &   9.936020 &    20 &    2  &  3.22e-15 &  3.48e-13 &  2.33e-13 \\ 
  13   &  10.173417 &  10.173417 &    21 &    2  &  0.00e+00 &  3.23e-13 &  1.52e-13 \\ 
  14   &  10.173417 &  10.173417 &    22 &    2  &  1.75e-16 &  6.63e-13 &  2.57e-13 \\ 
  15   &  11.064600 &  11.064600 &    23 &    2  &  5.14e-15 &  4.63e-13 &  1.08e-13 \\ 
  16   &  11.064601 &  11.064601 &    24 &    2  &  3.05e-15 &  2.57e-12 &  8.43e-14 \\ 
  17   &  11.086210 &  11.086210 &    26 &    2  &  2.56e-15 &  4.16e-12 &  1.17e-13 \\ 
  18   &  11.086210 &  11.086210 &    26 &    2  &  2.08e-15 &  1.49e-12 &  9.44e-14 \\ 
  19   &  11.619654 &  11.619654 &    27 &    1  &  5.20e-15 &  1.42e-11 &  1.81e-13 \\ 
  20   &  11.619808 &  11.619808 &    28 &    1  &  2.29e-15 &  5.90e-12 &  2.43e-13 \\ 
  21   &  11.791425 &  11.791425 &    29 &    1  &  3.46e-15 &  4.03e-13 &  1.27e-13 \\ 
  22   &  12.224826 &  12.224826 &    30 &    2  &  1.45e-16 &  3.15e-13 &  1.11e-13 \\ 
  23   &  12.224827 &  12.224827 &    31 &    2  &  2.91e-16 &  3.35e-13 &  1.40e-13 \\ 
  24   &  12.338403 &  12.338403 &    32 &    2  &  1.30e-15 &  3.08e-13 &  1.36e-13 \\ 
  25   &  12.338403 &  12.338403 &    32 &    2  &  5.76e-16 &  1.70e-13 &  1.04e-13 \\ 
  26   &  13.014987 &  13.014987 &    34 &    2  &  0.00e+00 &  4.32e-14 &  9.50e-14 \\ 
  27   &  13.014987 &  13.014987 &    35 &    2  &  0.00e+00 &  1.22e-13 &  1.10e-13 \\ 
  28   &  13.323474 &  13.323474 &    36 &    2  &  3.47e-15 &  2.55e-13 &  9.08e-14 \\ 
  29   &  13.323474 &  13.323474 &    36 &    2  &  1.20e-15 &  5.05e-13 &  1.48e-13 \\ 
  30   &  13.353881 &  13.353881 &    38 &    2  &  5.72e-15 &  1.77e-12 &  5.08e-13 \\ 
  31   &  13.353881 &  13.353881 &    38 &    2  &  1.33e-16 &  1.36e-12 &  5.88e-13 \\ 
  32   &  13.588947 &  13.588947 &    40 &    2  &  5.10e-15 &  2.89e-13 &  2.01e-13 \\ 
  33   &  13.588948 &  13.588948 &    41 &    2  &  6.54e-16 &  6.94e-13 &  4.09e-13 \\ 
  34   &  14.372157 &  14.372157 &    42 &    2  &  4.45e-15 &  9.97e-14 &  7.88e-14 \\ 
  35   &  14.372157 &  14.372157 &    43 &    2  &  3.46e-15 &  1.02e-13 &  8.61e-14 \\ 
  36   &  14.474860 &  14.474860 &    44 &    2  &  1.23e-15 &  9.07e-14 &  1.00e-13 \\ 
  37   &  14.474873 &  14.474873 &    45 &    2  &  3.19e-15 &  1.59e-13 &  1.58e-13 \\ 
  38   &  14.795305 &  14.795305 &    46 &    1  &  1.92e-15 &  4.77e-12 &  5.32e-14 \\ 
  39   &  14.795795 &  14.795795 &    47 &    1  &  2.40e-15 &  1.38e-11 &  7.53e-14 \\ 
  40   &  14.820720 &  14.820720 &    48 &    2  &  5.99e-16 &  1.34e-12 &  6.35e-14 \\ 
  41   &  14.820720 &  14.820720 &    49 &    2  &  2.40e-16 &  9.65e-13 &  9.12e-14 \\ 
  42   &  14.930509 &  14.930509 &    50 &    1  &  6.78e-15 &  1.13e-13 &  5.66e-14 \\ 
  43   &  15.588923 &  15.588923 &    52 &    2  &  6.84e-16 &  1.55e-13 &  6.67e-14 \\ 
 \hline 
\end{tabular}
\end{center}
\vspace*{-1\baselineskip}
\caption{Further details of the eigenpairs computed using EigenWave for a disk with grid $\Gcd^{(4)}$ using the KrylovSchur algorithm and implicit time-stepping (see Table~\ref{tab:sice4Summary} for a summary of these results).  The spatial order of accuracy is 4 and the wave-solves use $N_p=1$ to determine the final time.
The index $k$ denotes the closest true discrete eigenvalue $\lambda_{h,k}^{\rm true}$ to the EigenWave value $\lambda_{h,j}$.
}
\label{tab:sice4Order4}
\end{table}


\begin{table}[hbt]
\begin{center}\tableFontSize
\begin{tabular}{|c|r|r|c|c|c|c|c|} \hline
 \multicolumn{8}{|c|}{EigenWave: sphere, ts=implicit, order=2, $\omega=   5.0$, $N_p=1$, KrylovSchur } \\ \hline 
   $j$  & \multicolumn{1}{c|}{$\lambda_{h,j}$} &  \multicolumn{1}{c|}{$\lambda_{h,k}^{\rm true}$}  &$k$ &  mult &  eig-err & evect-err & eig-res  \\ \hline
   0   &   3.139245 &   3.139245 &     0 &    1  &  1.32e-14 &  2.61e-13 &  6.09e-12 \\ 
   1   &   4.486704 &   4.486704 &     2 &    3  &  4.59e-11 &  2.81e-10 &  2.52e-12 \\ 
   2   &   4.486704 &   4.486704 &     2 &    3  &  4.59e-11 &  2.81e-10 &  2.52e-12 \\ 
   3   &   4.486740 &   4.486740 &     3 &    3  &  2.16e-14 &  2.19e-12 &  2.63e-12 \\ 
   4   &   5.748539 &   5.748539 &     4 &    1  &  9.12e-15 &  2.06e-10 &  9.89e-13 \\ 
   5   &   5.748753 &   5.748753 &     5 &    1  &  4.63e-16 &  1.31e-10 &  8.10e-13 \\ 
   6   &   5.751066 &   5.751066 &     6 &    3  &  4.79e-15 &  4.34e-11 &  9.37e-13 \\ 
   7   &   5.751083 &   5.751083 &     7 &    3  &  1.66e-10 &  7.87e-09 &  9.41e-13 \\ 
   8   &   5.751083 &   5.751083 &     8 &    3  &  1.66e-10 &  9.67e-09 &  9.54e-13 \\ 
   9   &   6.264269 &   6.264269 &     9 &    1  &  3.54e-15 &  3.10e-13 &  4.92e-13 \\ 
  10   &   6.963967 &   6.963967 &    11 &    2  &  1.91e-09 &  7.13e-05 &  2.11e-05 \\ 
  11   &   6.963967 &   6.963967 &    11 &    2  &  1.91e-09 &  7.13e-05 &  2.11e-05 \\ 
  12   &   6.964662 &   6.964662 &    12 &    1  &  3.32e-15 &  4.02e-11 &  2.06e-13 \\ 
  13   &   6.966035 &   6.966035 &    13 &    2  &  2.80e-12 &  1.20e-09 &  3.04e-13 \\ 
  14   &   6.966035 &   6.966035 &    14 &    2  &  2.81e-12 &  6.96e-10 &  2.85e-13 \\ 
  15   &   6.966550 &   6.966550 &    15 &    1  &  5.35e-15 &  9.19e-11 &  3.93e-13 \\ 
  16   &   6.968857 &   6.968857 &    16 &    1  &  3.31e-15 &  2.89e-11 &  4.36e-13 \\ 
  17   &   7.691124 &   7.691124 &    17 &    3  &  1.28e-11 &  1.68e-13 &  5.02e-06 \\ 
  18   &   7.691124 &   7.691124 &    17 &    3  &  1.28e-11 &  1.68e-13 &  5.02e-06 \\ 
  19   &   7.691163 &   7.691163 &    19 &    3  &  1.77e-14 &  2.00e-13 &  2.15e-13 \\ 
  20   &   8.145011 &   8.145011 &    20 &    1  &  1.42e-14 &  7.40e-12 &  4.39e-13 \\ 
  21   &   8.147207 &   8.147207 &    21 &    1  &  1.09e-15 &  2.12e-11 &  5.74e-13 \\ 
  22   &   8.147883 &   8.147883 &    22 &    2  &  2.26e-11 &  1.43e-11 &  5.15e-13 \\ 
  23   &   8.147883 &   8.147883 &    22 &    2  &  2.26e-11 &  1.41e-11 &  5.18e-13 \\ 
  24   &   8.148282 &   8.148282 &    24 &    1  &  9.16e-15 &  6.76e-12 &  3.04e-13 \\ 
  25   &   8.149261 &   8.149261 &    25 &    1  &  2.18e-15 &  1.43e-11 &  5.09e-13 \\ 
  26   &   8.152193 &   8.152193 &    26 &    1  &  3.92e-15 &  2.41e-11 &  3.86e-13 \\ 
  27   &   8.152703 &   8.152703 &    27 &    2  &  3.58e-11 &  2.50e-11 &  5.29e-13 \\ 
  28   &   8.152703 &   8.152703 &    27 &    2  &  3.58e-11 &  2.43e-11 &  5.33e-13 \\ 
 \hline 
\end{tabular}
\end{center}
\vspace*{-1\baselineskip}
\caption{Further details of the eigenpairs computed using EigenWave for a sphere with grid $\Gcs^{(2)}$ using the KrylovSchur algorithm and implicit time-stepping (see Table~\ref{tab:spheree2Summary} for a summary of these results).  The spatial order of accuracy is 2 and the wave-solves use $N_p=1$ to determine the final time.
The index $k$ denotes the closest true discrete eigenvalue $\lambda_{h,k}^{\rm true}$ to the EigenWave value $\lambda_{h,j}$.
}
\label{tab:spheree2Order2}
\end{table}

\begin{table}[hbt]
\begin{center}\tableFontSize
\begin{tabular}{|c|r|r|c|c|c|c|c|} \hline
 \multicolumn{8}{|c|}{EigenWave: sphere, ts=implicit, order=4, $\omega=   5.0$, $N_p=1$, KrylovSchur } \\ \hline 
   $j$  & \multicolumn{1}{c|}{$\lambda_{h,j}$} &  \multicolumn{1}{c|}{$\lambda_{h,k}^{\rm true}$}  & $k$ &  mult &  eig-err & evect-err & eig-res  \\ \hline
   0   &   3.141588 &   3.141588 &     0 &    1  &  7.44e-14 &  5.03e-13 &  9.87e-12 \\ 
   1   &   4.493398 &   4.493398 &     1 &    3  &  3.02e-14 &  4.24e-13 &  1.98e-10 \\ 
   2   &   4.493398 &   4.493398 &     1 &    3  &  3.02e-14 &  4.24e-13 &  1.98e-10 \\ 
   3   &   4.493398 &   4.493398 &     3 &    3  &  2.85e-14 &  4.93e-13 &  4.23e-12 \\ 
   4   &   5.763458 &   5.763458 &     4 &    5  &  4.16e-15 &  4.29e-13 &  1.48e-12 \\ 
   5   &   5.763461 &   5.763461 &     5 &    5  &  1.54e-15 &  3.62e-13 &  1.94e-12 \\ 
   6   &   5.763472 &   5.763472 &     6 &    5  &  9.05e-14 &  3.54e-13 &  5.08e-09 \\ 
   7   &   5.763472 &   5.763472 &     6 &    5  &  9.05e-14 &  3.54e-13 &  5.08e-09 \\ 
   8   &   5.763472 &   5.763472 &     8 &    5  &  5.86e-15 &  3.33e-13 &  1.97e-12 \\ 
   9   &   6.283201 &   6.283201 &     9 &    1  &  1.63e-14 &  3.04e-13 &  6.08e-13 \\ 
  10   &   6.988010 &   6.988010 &    10 &    7  &  4.80e-14 &  3.46e-13 &  1.13e-09 \\ 
  11   &   6.988010 &   6.988010 &    10 &    7  &  4.80e-14 &  3.46e-13 &  1.13e-09 \\ 
  12   &   6.988016 &   6.988016 &    12 &    7  &  1.08e-14 &  3.21e-13 &  4.52e-13 \\ 
  13   &   6.988021 &   6.988021 &    13 &    7  &  6.29e-12 &  6.37e-13 &  1.09e-08 \\ 
  14   &   6.988021 &   6.988021 &    13 &    7  &  6.29e-12 &  6.37e-13 &  1.09e-08 \\ 
  15   &   6.988032 &   6.988032 &    15 &    7  &  6.48e-15 &  6.13e-13 &  6.30e-13 \\ 
  16   &   6.988043 &   6.988043 &    16 &    7  &  3.81e-15 &  4.41e-13 &  8.72e-13 \\ 
  17   &   7.725441 &   7.725441 &    17 &    3  &  1.82e-14 &  2.45e-13 &  1.16e-10 \\ 
  18   &   7.725441 &   7.725441 &    17 &    3  &  1.82e-14 &  2.45e-13 &  1.16e-10 \\ 
  19   &   7.725452 &   7.725452 &    19 &    3  &  5.40e-15 &  2.57e-13 &  3.18e-13 \\ 
  20   &   8.182727 &   8.182727 &    20 &    6  &  7.38e-15 &  3.76e-10 &  7.11e-13 \\ 
  21   &   8.182737 &   8.182737 &    21 &    7  &  4.34e-16 &  2.10e-10 &  1.00e-12 \\ 
  22   &   8.182769 &   8.182769 &    22 &    9  &  2.76e-14 &  6.74e-10 &  4.76e-13 \\ 
  23   &   8.182789 &   8.182789 &    23 &    9  &  6.08e-15 &  6.60e-10 &  8.55e-13 \\ 
  24   &   8.182794 &   8.182794 &    24 &    9  &  7.66e-14 &  4.98e-10 &  5.05e-09 \\ 
  25   &   8.182794 &   8.182794 &    24 &    9  &  7.66e-14 &  4.98e-10 &  5.05e-09 \\ 
  26   &   8.182822 &   8.182822 &    26 &    8  &  4.99e-15 &  1.37e-09 &  1.25e-12 \\ 
  27   &   8.182837 &   8.182837 &    27 &    7  &  1.57e-12 &  1.68e-09 &  4.88e-09 \\ 
  28   &   8.182837 &   8.182837 &    27 &    7  &  1.57e-12 &  1.68e-09 &  4.88e-09 \\ 
 \hline 
\end{tabular}
\end{center}
\vspace*{-1\baselineskip}
\caption{Further details of the eigenpairs computed using EigenWave for a sphere with grid $\Gcs^{(2)}$ using the KrylovSchur algorithm and implicit time-stepping (see Table~\ref{tab:spheree2Summary} for a summary of these results).  The spatial order of accuracy is 4 and the wave-solves use $N_p=1$ to determine the final time.
The index $k$ denotes the closest true discrete eigenvalue $\lambda_{h,k}^{\rm true}$ to the EigenWave value $\lambda_{h,j}$.
}
\label{tab:spheree2Order4}
\end{table}